\documentclass[11pt, reqno]{amsart}  
\usepackage{amsmath, mathrsfs}
\usepackage{tikz}

\usepackage[all]{xy}
\SelectTips{cm}{}
\usepackage{amsthm, amsfonts}
\usepackage{geometry}
  \usepackage[overload]{empheq}
  
  \usepackage{tikz-cd}
  \usepackage{url}

  \usepackage{amssymb}
  
  \usepackage{graphicx}
  \usepackage{enumitem}
  \usepackage{color}
  \usepackage{enumitem}
  \usepackage[colorlinks=true]{hyperref}
  \hypersetup{urlcolor=cyan, citecolor=red, linkcolor=blue}
               
  \usepackage[capitalize]{cleveref}

  \usepackage{thmtools}
  \usepackage{float}

		\crefname{lemma}{Lemma}{Lemmas}
		\crefname{theorem}{Theorem}{Theorems}
		\crefname{prop}{Proposition}{Propositions}
		\crefname{cor}{Corollary}{Corollaries}
	
		\usepackage{xcolor}
\usepackage[most]{tcolorbox}
\usepackage{listings}

\definecolor{white}{rgb}{1,1,1}
\definecolor{mygreen}{rgb}{0,0.4,0}
\definecolor{light_gray}{rgb}{0.97,0.97,0.97}
\definecolor{mykey}{rgb}{0.117,0.403,0.713}

\tcbuselibrary{listings}
\newlength\inwd
\newtcblisting{pyin}[1][]{%
sharp corners,
enlarge left by=\inwd,
width=\linewidth-\inwd,
enhanced,
boxrule=0pt,
colback=light_gray,
listing only,
top=0pt,
bottom=0pt,
overlay={
  \node[
	anchor=north east,
	text width=\inwd,
	font=\footnotesize\ttfamily\color{mykey},
	inner ysep=1.5mm,
	inner xsep=0pt,
	outer sep=0pt
	] 
	at (frame.north west)
	{\refstepcounter{ipythcntr}\label{#1}In \theipythcntr:};
}
listing engine=listing,
listing options={
  aboveskip=1pt,
  belowskip=1pt,
  basicstyle=\footnotesize\ttfamily,
  language=Python,
  keywordstyle=\color{mykey},
  showstringspaces=false,
  stringstyle=\color{mygreen}
},
}

		\newtheoremstyle{1}
		{6pt} 
		{0pt} 
		{\itshape} 
		{} 
		{\bfseries} 
		{.} 
		{.5em} 
		{} 
		
		\newtheoremstyle{2}
		{6pt} 
		{0pt} 
		{} 
		{} 
		{\bfseries} 
		{.} 
		{.5em} 
		{} 
		
		\theoremstyle{1}
		\newtheorem{theorem}{Theorem}
		\newtheorem{lemma}[theorem]{Lemma}
		\newtheorem{prop}[theorem]{Proposition}
		\newtheorem{cor}[theorem]{Corollary}

    \theoremstyle{definition}
		\newtheorem{defn}[theorem]{Definition}
		\newtheorem{remark}[theorem]{Remark}
		\newtheorem{example}[theorem]{Example}
		\newtheorem{notation}[theorem]{Notation}
		\newtheorem{assumption}[theorem]{Assumption}
		\newtheorem{conjecture}[theorem]{Conjecture}

	\numberwithin{equation}{subsection}
		\numberwithin{lemma}{subsection}
		\numberwithin{theorem}{subsection}
    \numberwithin{prop}{subsection}
	\numberwithin{cor}{subsection}
		\numberwithin{defn}{subsection}
		\numberwithin{remark}{subsection}
		\numberwithin{example}{subsection}
		\numberwithin{assumption}{subsection}
    \numberwithin{conjecture}{subsection}
    \numberwithin{notation}{subsection}

        \newcommand{\F}{\mathbb{F}}
		\newcommand{\Q}{\mathbb{Q}}
		\newcommand{\Z}{\mathbb{Z}}
		\newcommand{\A}{\mathbb{A}}
		\newcommand{\Sym}{\textnormal{Sym}}

        \newcommand{\Ker}{\textnormal{Ker}}
		\newcommand{\Gal}{\textnormal{Gal}}
		\newcommand{\Hom}{\textnormal{Hom}}
		\newcommand{\gsp}{\textnormal{GSp}}
		
		\newcommand{\flag}{\mathcal{F}l}
		
		\newcommand{\KS}{\textnormal{KS}}
		\newcommand{\LL}{\mathcal{L}}
		\newcommand{\Ig}{\textnormal{Ig}}
    \newcommand{\GL}{\textnormal{GL}}

    \newcommand{\Ver}{\textnormal{Ver}}
    \newcommand{\Rep}{\textnormal{Rep}}
    \newcommand{\DiffOp}{\textnormal{DiffOp}}
    \newcommand{\crys}{\textnormal{crys}}
    \newcommand{\Std}{\textnormal{Std}}
    \newcommand{\BGG}{\textnormal{BGG}}
    \newcommand{\dR}{\textnormal{dR}}
    
    \newcommand{\et}{\textnormal{\'et}}
    \newcommand{\tor}{\textnormal{tor}}
    \newcommand{\can}{\textnormal{can}}
    \newcommand{\sub}{\textnormal{sub}}
    \newcommand{\Sh}{\textnormal{Sh}}
    \newcommand{\Shbar}{\overline{\textnormal{Sh}}}
	\newcommand{\GZip}{\mathop{\text{$G$-{\tt Zip}}}\nolimits}
	\newcommand{\GF}{\mathop{\text{$G$-{\tt ZipFlag}}}\nolimits}
	\newcommand{\Sbt}{\textnormal{Sbt}}
	\newcommand{\G}{\mathbb{G}}
	\newcommand{\SL}{\textnormal{SL}}
	\newcommand{\GSp}{\textnormal{GSp}}
	\newcommand{\Lie}{\textnormal{Lie}}
  
  \newcommand{\Fpbar}{\overline{\F}_p}
 \newcommand{\C}{\mathbb{C}}
   \newcommand{\der}{\textnormal{der}}
   \renewcommand{\H}{\textnormal{H}}
   \newcommand{\Modfg}{\textnormal{Mod}^{\textnormal{fg}}}

\begin{document}

\title{Theta operators on Hodge type Shimura varieties}
\author{Martin Ortiz}

\maketitle

\begin{abstract}
  We construct a new family of mod $p$ weight shifting differential operators on Hodge type Shimura varieties 
  at hyperspecial level. First we construct basic theta operators, labelled by positive roots, that generalize Katz's 
  theta operator for modular forms. 
  Secondly we construct theta linkage maps, these are operators between automorphic vector bundles 
  with linked weights, which can be thought of as generalizations of the classical theta cycle of Tate--Jochnowitz.
  In particular, there exist such maps within the $p$-restricted region,
  whose weight shifts are directly related to the conjectures of Herzig on the weight part
  of Serre's conjecture. We explain the relation between the two operators, and we prove some properties about them, e.g. 
  the injectivity of some of them in a generic locus of the $p$-restricted region.  As an application, we 
  produce an example of a generic entailment of Serre weights for the groups $\GL_{4,\Q_p}$ and $U(4)_{\Q_p}$,
  by combining the method of \cite{paper} with our stronger results about theta operators. 
  \end{abstract}

\tableofcontents

\clearpage

\section{Introduction}
In this article we construct a general family of theta operators on the special fibers of Hodge type Shimura varieties. Let us first 
recall the original case of mod $p$ modular forms, as it will serve us as a template to generalize. Let $p$ be a prime 
and let $X_{\Z_p}$ be the compactified modular curve at level prime to $p$. Let $\omega$ be the Hodge bundle over $X_{\Z_p}$, and $k \in \Z$. 
The first mod $p$ theta operator was defined by Katz \cite{katz-theta} for mod $p$ modular forms as a map of sheaves 
$\theta: \omega^k \to \omega^{k+p+1}$ over $X_{\Fpbar}$. In the same paper he proved some of its important properties. 
\begin{prop} \label{theta-GL2-properties}
\begin{enumerate} 
\item The map $\theta$ is Hecke equivariant up to a twist, i.e. $T_l\theta=l \theta T_l$ for $l \neq p$,
where $T_l$ are the classical 
Hecke operators. 
\item Let $H$ be the Hasse invariant. Then $\theta(H)=0$, so that $\theta$ is Hasse-linear, i.e. for any form $f \in \omega^k$, $\theta(Hf)=H\theta(f)$. 
\item Let $\theta^p$ denote the $p$-fold iterate of $\theta$, then $\theta^p=H^{p+1}\theta$.
\item Let $f \in \omega^k$ be a local section over $X_{\overline{\F}_p}$ which is not divisible by $H$, i.e. 
there does not exist a local section $g \in \omega^{k-p+1}$ such that 
$f=Hg$. Then $H \mid \theta(f)$ if and only 
if $p \mid k$. 
\item  The map $\theta: \H^0(X_{\Fpbar},\omega^k) \to \H^0(X_{\Fpbar},\omega^{k+p+1})$
is injective for $k<p-1$.
\end{enumerate}
\end{prop}

The relevance of point $(5)$ is that one can think of $\H^0(X_{\Fpbar},\omega^k)$ as spaces of mod $p$ modular forms, 
and the theta operator is producing congruences between forms of different weights. Both $(2)$ and $(4)$ are used to prove 
$(5)$. Moreover, using $(3)$ and $(4)$ Jochnowitz \cite{Jochnowitz} and Tate computed the theta cycle for modular forms of low weight. That is, 
one repeatedly applies $\theta$ to an eigenform, studying if at some point the form becomes 
divisible by the Hasse invariant. By the relation $\theta^{p}=H^{p+1}\theta$ this must eventually occur. 
Part of their result can be rephrased as follows. 
\begin{prop}
(The theta cycle) Let $1 \le k \le p$. The map $\theta^{p-k+1}: \omega^k \to \omega^{k+(p-k+1)(p+1)}$ factors through $H^{p-k+1}: \omega^{2p-k+2} \to 
\omega^{k+(p-k+1)(p+1)}$ via a map $\Theta: \omega^{k} \to \omega^{2p-k+2}$. 
\end{prop}
For us, although this is not 
standard terminology, we will say that $\Theta$ is 
the \textit{theta cycle map}. 
Moreover, in the case that $f$ is an eigenform of weight $k \le p$ which is non-ordinary, Edixhoven \cite{Edixhoven}
proved that $H \mid \Theta(f)$, thus obtaining a form of weight $p-k+3$. This partially proved the weight part of Serre's conjecture 
for modular forms, and we see it as one of the benchmarks for applications of mod $p$ theta operators. 
In this article we will generalize the two propositions above. We will denote the generalizations of $\theta$ 
as \textit{basic theta operators} and the generalizations of $\Theta$ as \textit{theta linkage maps}.
We will work with Hodge type Shimura varieties. 

These are those for which the Shimura datum $(G,X)$
embeds into the one for the Siegel Shimura variety. Let $K \subseteq G(\mathbb{A}^{\infty})$
be a neat level which is hyperspecial at $p$. That is, there exists a reductive model $\mathcal{G}/\Z_p$ such that 
$K=K^pK_p$ with $K_p=\mathcal{G}(\Z_p)$. From now on we suppress $\mathcal{G}$ from 
the notation. Given $G$, the existence of some $K$ hyperspecial at $p$ 
is equivalent to $G_{\Q_p}$ being unramified. 
By work of Kisin \cite{Kisin-integral-model} there exists 
a smooth integral model $\Sh$ for the Shimura variety attached to $(G,X,K)$, defined 
over the ring of integers $\mathcal{O}$ of a finite extension of $\Q_p$.
Let $\Shbar$ be its geometric special fiber, and $\Sh^\tor$ a choice of toroidal compactification. 

Let $\mu$ be the Hodge cocharacter and $P$ a choice of 
parabolic of $G$ corresponding to $\mu$. Let $M$ be the Levi of $P$. We fix a choice $T \subseteq B \subset P \subseteq G_{\Fpbar}$ of 
maximal torus and Borel. Let $\omega=e^*\Omega^1_{A/\Sh}$ and 
let $\omega \to \mathcal{H}:=\H^1_{\dR}(A/\Sh) \to \omega^{\vee}_{A^{\vee}}$
be the Hodge filtration. 
 Kisin \cite[Cor 2.3.9]{Kisin-integral-model} proves the existence of some Hodge tensors 
 $s_{\dR} \in \mathcal{H}^{\otimes}$ in such a way that the sheaf $G_{\dR}$ consisting of 
 trivializations of $\mathcal{H}$ respecting $s_{\dR}$ (\Cref{torsors}) becomes a $G$-torsor over $\Sh$. 
 Then the Hodge filtration defines a $P$-subtorsor $P_{\dR} \subseteq G_{\dR}$ over $\Sh$.
 For $\lambda \in X^*(T)$ let $W(\lambda) \in \Rep(M)$
be the dual Weyl module for $M$ of highest weight $\lambda$. We define $\omega(\lambda)=P_{\dR} \times^{P} W(\lambda)$
as the vector bundle 
on $\Sh^\tor$ associated to $W(\lambda)$ via $P_{\dR}$.

As in \cite{GK-stratification} we will work with the flag Shimura variety $\pi: \flag:=[B\backslash P_{\dR}] \to \Sh$ (\Cref{flag-space}), with the property that its fibers 
are isomorphic to the partial flag variety $P/B$. 
Given $\lambda \in X^*(T)$ there exists an automorphic line bundle $\LL(\lambda)$
over $\flag$ such that $\pi_* \LL(\lambda)=\omega(w_{0,M}\lambda)$, where $w_{0,M} \in W_{M}$ is the longest element. We will construct our operators on $\flag$, which 
then can be pushed forward to $\Shbar$.

\subsection{Main results} \label{section0.2}
Let $(G,X)$ be a Shimura datum of Hodge type, and $K=K^pK_p$ a neat hyperspecial level at $p$.
We assume that $(G,X)$ satisfies Deligne's axioms (SV1,2,3) \cite[\S 6]{Milne-moduli}  \cite{Deligne-Shimura}.
Assuming SV3 is harmless for our purposes. Namely, a $\Q$-simple 
factor of compact type of $G^{\text{ad}}$ does not contribute to the geometry of the Shimura variety. 
We will also use $G$ to denote the associated reductive model over $\Z_p$.
We will always consider $\flag$ over $\Fpbar$ in this section.  
We construct two kinds of theta operators, basic theta operators and theta linkage maps. Since we will repeatedly use 
this concept in this introduction, 
let us record our definition of \textit{generic weights}.
\begin{defn} \label{defn-generic}
Let $G/\Q$ be a reductive group. We fix 
a maximal torus and a Borel $T \subseteq B \subseteq G_{\overline{\Q}}$. 
\begin{enumerate}
\item Let $\epsilon \ge 0$. We say that $\lambda \in X^*(T)$ is $\epsilon$-generic with respect to a prime $p$ if 
for each $\gamma \in \Phi$, there exists an integer $n_{\gamma}$ such that 
$$
n_{\gamma}p+\epsilon< \langle \lambda, \gamma^{\vee} \rangle < (n_{\gamma}+1)p-\epsilon. 
$$
\item Fix some Shimura datum of Hodge type $(G,X)$.
Let $X=\{X_{p,K}\}_{p,K}$ be a family of statements ranging over all rational primes $p$ and levels $K \subseteq G(\mathbb{A}^{\infty})$
which are hyperspecial at $p$, concerning 
 $\Sh_{K,\mathcal{O}}$ (where $\mathcal{O}$ lives over $\Z_p$) and some weights $\lambda \in X^*(T)$. 
We say that $X$ holds for \textit{generic} $\lambda \in X^*(T)$ if there exists $\epsilon>0$ only depending on 
$G_{\Q}$ such that each statement $X_{p,K}$ holds for weights $\lambda$ which are $\epsilon$-generic with respect to $p$.
\item Similarly, fix $H$ a choice of a reductive root datum with a finite cyclic group $\Gamma$ acting on it.
For each prime $p$ let $G_{\F_p}$ be the reductive 
 group over $\F_p$ given by $(H,\Gamma)$. Choose some maximal torus $T \subset G$.
Let $X=\{X_p\}_{p}$ be a family of statements about each $G_{\F_p}$ and some weights $\lambda \in X^*(T)$. 
We say that $X$ holds for generic $\lambda$ if there exists an $\epsilon$ depending only on $(H,\Gamma)$ such that 
each $X_p$ holds for weights $\lambda$ which are $\epsilon$-generic with respect to $p$.
\end{enumerate}
\end{defn}
We refer to \Cref{rep-theory} for all the representation theoretic notation.  

\subsubsection{Basic theta operators}
The basic theta operators are the natural generalization of Katz's theta operator, as well as of all other mod 
$p$ theta operators previously defined in the literature.
For each simple root $\alpha \in \Delta$  there exists a Hasse invariant $H_{\alpha}$ on $\flag$ \cite{IK-hasses}. We note that 
$G^{\der}$ is always simply connected for Shimura varieties of PEL type A or C. 

\begin{theorem} \label{basic-theta-operaotrs}
Let $\gamma \in \Phi^{+}$ be a positive root. For each $\lambda \in X^*(T)$
there exists a basic theta operator 
$\theta_{\gamma} : \LL(\lambda) \to \LL(\lambda+\mu_{\gamma})$, which is a differential operator 
over $\flag_{\Fpbar}$, and $\mu_{\gamma}$ is an explicit weight depending on $\gamma$. They extend to toroidal 
compactifications. 
 They satisfy the following properties. 
\begin{enumerate}
\item They are Hecke equivariant away from $p$. 
\item  For each $\alpha \in \Delta$, $\theta_{\gamma}(H_{\alpha})=0$. 
\item 
For $\gamma_1,\gamma_2 \in \Phi^{+}$ we have the commutation relations 
$[\theta_{\gamma_1},\theta_{\gamma_2}]=\theta_{\gamma_1+\gamma_2}$ up to multiplying the right-hand side 
by an explicit collection of Hasse invariants and a small explicit constant. If $\gamma_1+\gamma_2 \notin \Phi^{+}$, 
the right-hand side is zero by convention. 
\item Assume that $G^{\textnormal{der}}$ is simply connected. Let $\{G_i\}$ be the almost-simple factors of 
$G^{\der}_{\Fpbar}$, and fix such an $i$.
The Frobenius in $\textnormal{Gal}(\Fpbar/\F_p)$
 maps $G_i$ isomorphically to $G_j$ for some $j$. 
Suppose that 
there exists $\alpha \in \Delta_{G_i}$ not contained in 
$\Lie(M)$, and let $\delta$ be the longest root of $G_j$. Then if $G_i$ has rank at least $2$
$$
\theta^p_{\alpha}=H^p_{\alpha}\theta_{\delta}.
$$
\item Let $\alpha \in \Delta$, and assume that $H_{\alpha}$ has a simple zero 
(this holds 
if $G^{\textnormal{der}}$ is simply connected).
Let $f \in \LL(\lambda)$ be a local section, then
$H_{\alpha} \mid \theta_{\alpha}(f)$ if and only if $p \mid \langle \lambda,\alpha^{\vee} \rangle$ or $H_{\alpha} 
\mid f$. 
\item Assume that $G^{\der}$ is simply connected, and that \Cref{assumption-vanishing}(2) holds. 
Then the map
$$
\theta_{\gamma}: \H^0(\Shbar^\tor,\omega(\lambda)) \to \H^0(\Shbar^\tor,\omega(\lambda+w_{0,M}\mu_{\gamma}))
$$
is injective for generic $\lambda \in X_1(T)$.
\end{enumerate}
\end{theorem}
Here \Cref{assumption-vanishing}(2) concerns the vanishing of $\H^0(\Shbar^\tor,\omega(\lambda))$ 
for weights $\lambda$
which are generically away from the dominant cone.
It is implied in many cases by the cone conjecture of Goldring--Koskivirta 
\cite{general-cone-conjecture1}, which has been proved in some low rank examples, as well as in the Siegel case 
\cite{cone-conjecture-Siegel}. This theorem is the amalgamation of 
\Cref{theta-kills-Hasse}, \Cref{pth-power-relation}, \Cref{restriction-simple-theta}, and \Cref{basic-theta-injective}
in that order, 
where we give more precise statements. See \Cref{table-lie} and \Cref{Hasse-dfn} for how to compute 
$\mu_{\gamma}$.
 We would like to draw parallels between 
the case of the modular curve in \Cref{theta-GL2-properties} and this theorem. Point $(3)$ above
does not appear in the case of the modular curve since it is only relevant  
for groups of higher rank. 
 We also note that some of the operators become linear 
after pushing them forward to $\Shbar$, these are precisely $\theta_{\gamma}$
for $\gamma \in \Lie(U_{B})\backslash
\Lie(U_{P})$.
\begin{remark}
To prove points $(5)$ and $(6)$ we in fact prove a more precise result about the restriction of basic theta operators 
to codimension $1$ strata, e.g. \Cref{general-restriction-strata}. This could be seen as a first step towards the 
study of differential operators on EO strata. 
\end{remark}

We illustrate the weight shifts of the basic theta operators and Hasse invariants in the case of $G=\GSp_4$.
Weights are given by a tuple $(k,l,w) \in \Z^3$ with $w=k+l \bmod{2}$, and then for $k\ge l$ 
we have $\omega(k,l,w)=\Sym^{k-l}\omega \otimes \det^l \omega \otimes \delta^{(w-k-l)/2}$, where 
$\delta:=\omega(0,0,2)$. It is a non-canonically trivial line bundle which twists the Hecke action, and 
$\delta^{p-1}$ is canonically trivial over $\Shbar$. 
The simple roots are $\alpha=(1,-1,0)$ and $\beta=(0,2,0)$, and the basic theta operators are
\begin{itemize}
\item $\pi_*\theta_{\beta}: \omega(\lambda) \to \omega(\lambda+(p+1,p-1,0))$
\item $\pi_*\theta_{\alpha}: \omega(\lambda) \to \omega(\lambda+(p-1,0,p-1))$
\item $\pi_*\theta_{\alpha+\beta}: \omega(\lambda) \to \omega(\lambda+(2p,p-1,p-1))$
\item $\pi_*\theta_{2\alpha+\beta}: \omega(\lambda) \to \omega(\lambda+(2p,0,2p-2))$.
\end{itemize}
Note that only the weight shift of $\pi_*\theta_{\alpha}$ does not introduce a twist by $\delta$.
The arrows in \Cref{figure1} illustrate these weight shifts, where 
we are suppressing the central character. The large parallelogram
roughly corresponds to the 
$p$-restricted weights. It will be a phenomenon that the operators only produce congruences between 
$p$-restricted weights for non-generic weights, i.e. those very close to the boundary of the $p$-restricted region. 
Here $\pi_* \theta_{\beta}$ was already defined by \cite{Yamauchi-1} and \cite{Eischen-Mantovan-1}, but the others 
are new. In the setup of \cite{Eischen-Mantovan-1}, $\pi_*\theta_{\alpha}$ appearing in point $(4)$ agrees 
with their theta operator. 
\begin{figure}
  \centering
  \begin{tikzpicture}
\node [color=black] (1) at (0,0) (snlabel) {$\bullet$};
\draw (0,0) -- (4,0);
\draw[dashed] (4,0) -- (8,0);
\draw (0,0) -- (4,4);
\draw (2,2) --(4,0);
\draw (4,0) --(4,4);
\draw (4,4) --(8,4);
\draw (4,0) -- (8,4);
\draw (4,4) --(6,2);
\draw[dashed] (0,0) --(0,4);
\node (2) at (8.6,0) (snlabel) {$k$};
\node (3) at (0,4.2) (snlabel) {$l$};
\node[color=black] (4) at (-0.3,-0.3) (snlabel) {$0$};
\node (5) at (4,-0.5) (snalabel) {$(p,0)$};
\node (6) at (4,4.3) (label) {$(p,p)$};
\node (7) at (8.2,4.3) (label) {$(2p,p)$};
\node[red] (H1) at (3.4,3.9) {$H_{\beta}$};
\draw[dashed,thick, red,->] (0,0) -- (3.9,3.9);
\draw[dashed,thick,red,->] (0,0) --(4,-0.2);
\node[red] (H2) at (2,-0.4) {$H_{\alpha}$};
\draw[dashed,blue,->] (0,0) --(4.3,3.9);
\node[blue] at (2.9,2.2) {$\theta_{\beta}$};
\draw[dashed,blue,thick,->] (0,0) --(3.9,0);
\node[blue] at (2,0.3) {$\theta_{\alpha}$};
\draw[dashed,blue,->] (0,0) --(8,3.9);
\node[blue] at (8.6,3.8) {$\theta_{2\alpha+\beta}$};
\draw[blue,dashed,->] (0,0) -- (8,0);
\node[blue] at (7.4,-0.3) {$\theta_{\alpha+\beta}$};
\draw[pink,thick,dashed,->] (0,0) --(0,1.5);
\node (-beta) at (-0.3,1.7) {$\beta$};
\draw[pink,thick,dashed,->] (0,0) --(0.9,-0.9);
\node (alpha) at (0.9,-1.2) {$\alpha$};
\end{tikzpicture}
  \caption{}
  \label{figure1}
  \end{figure}

\subsubsection{The theta linkage maps}
In order to get maps between $p$-restricted weights we can use \textit{theta linkage maps}. On one hand, these generalize 
the theta cycle map $\Theta: \omega^k \to \omega^{2p-k+2}$
from the modular curve case, and on the other hand they pair nicely with the modular representation theory appearing 
in the setup of mod $p$ algebraic automorphic forms. 
Let $\rho$ be the half sum of positive roots for $G_{\Fpbar}$. 
We can divide weight space $X^*(T) \otimes \mathbb{R}$ into 
$\rho$-shifted alcoves, defined as the connected components of weight space after dividing it by the 
hyperplanes $H_{\gamma,n}:\{\lambda: np=\langle \lambda+\rho,\gamma^{\vee} \rangle \}$ for all $n \in \Z$
and $\gamma \in \Phi^{+}$. The affine Weyl group $W_{\text{aff}}=p\mathbb{Z}\Phi^{+} \rtimes W$ 
acts transitively on the alcoves, and 
it is generated by reflections $s_{\gamma,n}: \lambda \mapsto \lambda-
(\langle \lambda+\rho,\gamma^{\vee} \rangle-pn)\gamma$. We say that $\lambda \uparrow \mu$ ($\mu$ is linked 
to $\lambda$) if there is a sequence $\lambda=\lambda_0 \le \lambda_1 \le \ldots \le \lambda_k=\mu$ such that 
$\lambda_{i+1}=s_{\gamma_i,n_i} \cdot \lambda_i$ for some $\gamma_i$ and $n_i$. 
We say that $\lambda \uparrow_{\gamma} \mu$ if $\mu=s_{\gamma,n} \cdot \lambda$ such that $n$ is the smallest integer
for which $\lambda \le s_{\gamma,n} \cdot \lambda$. In this case $\mu$ is a reflection across 
the closest $\gamma$-hyperplane to $\lambda$, such that it moves in the positive root direction. 

\begin{theorem}[\Cref{Verma-linkage}, \Cref{linkage-map-composition-basic}]
Assume that $G^{\text{der}}$ is simply connected. 
Let $\lambda,\mu \in X^*(T)$ such that $\lambda \uparrow_{\gamma} \mu$. Then there exists a differential 
 operator 
$\theta^{\lambda \uparrow \mu}: \LL(-\mu) \to \LL(-\lambda)$ on $\flag_{\Fpbar}$. It is Hecke 
equivariant away from $p$, and after post-composing $\theta^{\lambda \uparrow \mu}$ with 
some number of Hasse invariants it is identified with a combination of  basic theta operators. 

\end{theorem}
We denote these operators by \textit{theta linkage maps}. For general $\lambda \uparrow \mu$, 
by composing the maps above we can construct a map $\LL(-\mu) \to \LL(-\lambda)$, but 
we don't know if this is independent of the chain of reflections going from $\lambda$ to $\mu$.
Even in the setup of the theorem we cannot prove in general that these operators are non-zero. 
The last sentence in the statement means 
that there is an element $f \in U\mathfrak{u}_{B}$, which can be thought of as a 
non-commutative polynomial on the positive roots $\{x_{\gamma}\}$, 
such that when evaluating it on the basic theta operators $\{\theta_{\gamma}\}$ we get 
$$
f([\theta_{\gamma}]_{\gamma \in \Phi^{+}})=\prod_{\alpha \in \Delta} H^{n_{\alpha}}_{\alpha} \theta^{\lambda \uparrow \mu}
$$
for some 
$n_{\alpha} \ge 0$. This clearly generalizes the behaviour of the theta cycle map $\Theta=\theta^{k-1}/H^{k-1}: 
\omega^k \to \omega^{2p-k+2}$
 for the modular curve,
but in general multiple 
basic theta operators are required in this combination. Only in the case that $\lambda \uparrow_{\gamma} \mu$ and 
$\gamma$ is a simple root, 
then $\theta^{\lambda \uparrow \mu}$ will be a power of a single basic theta operator.  We call 
these \textit{simple theta linkage maps}. \Cref{figure2} illustrates the $3$ linkage maps within 
the $p$-restricted region for $\GSp_4$. Because of the twist $\pi_* \LL(\lambda)=\omega(w_{0,M}\lambda)$ we get 
maps $\omega(\lambda_i+\eta) \to \omega(\lambda_{i+1}+\eta)$, where $\eta=(3,3)$. 
The one going from $\lambda_1$ to $\lambda_2$ 
is the pushforward of a simple theta linkage map ($w_{0,M}$ corresponds to the 
reflection across the $x=y$ axis up to a shift by $\rho$).
\begin{figure}
\centering
\begin{tikzpicture}
    \node [color=black] (1) at (0,0) (snlabel) {$\bullet$};
    \node[below left] at (0,0) {$-\rho$};
    \draw (0,0) -- (4,0);
    \draw[dashed] (4,0) -- (8,0);
    \draw (0,0) -- (4,4);
    \draw (2,2) --(4,0);
    \draw (4,0) --(4,4);
    \draw (4,4) --(8,4);
    \draw (4,0) -- (8,4);
    \draw (4,4) --(6,2);
    \draw[dashed] (0,0) --(0,4);
    \draw[pink,thick,dashed,->] (0,0) --(0,1.5);
    \node (-beta) at (-0.3,1.7) {$\beta$};
    \draw[pink,thick,dashed,->] (0,0) --(0.9,-0.9);
    \node (alpha) at (0.9,-1.2) {$\alpha$};
    \node (lambda0) at (3,0.3) {$\scriptstyle \lambda_0$};
    \node (wlambda0) at (3,0.5) {$\cdot$};
   \node (wlambda1) at (3.5,1) {$\cdot$};
\node (wlambda2) at (4.5,1) {$\cdot$};
\node (wlambda3) at (7,3.5) {$\cdot$};
\node (lambda1) at (3.65,0.75) {$\scriptstyle \lambda_1$};
\node (lambda2) at (4.3,0.8) {$\scriptstyle \lambda_2$};
\node (lambda3) at (6.5,3.5) {$\scriptstyle \lambda_3$};
\draw[green,dotted,thick,->] (3,0.5) --(3.5,1);
\draw[green,dotted,thick,->] (3.5,1) --(4.5,1);
\draw[green,dotted,thick,->] (4.5,1) --(7,3.5);
    \end{tikzpicture}
\caption{}
\label{figure2}
\end{figure}

In general theta linkage maps are more difficult to study than 
basic theta operators. For instance, they are of higher degree as differential operators,
and they don't necessarily commute with Hasse invariants 
in any reasonable way, so it does not make sense to restrict them to codimension $1$ strata. 
In the case of simple theta linkage maps we can say a bit more, thanks to 
point $(4)$ in \Cref{basic-theta-operaotrs}. We will use a relabelling of the theta linkage maps which is more 
suitable for applications. 
\begin{theorem}[\Cref{injective-linkage}] \label{injective-linkage-intro}
  Assume that $G^{\der}$ is simply connected and that \Cref{assumption-vanishing}(2) holds. 
  Let $G_i$ be an almost-simple factor of $G^{\der}$
  and suppose that there exists 
  $\alpha \in \Delta_{G_i}$ which is not in $\Lie(M)$. Let $\lambda \in X_1(T)$ and $\mu \in X^*(T)$
  such that $\lambda \uparrow_{w_{0,M}\alpha} \mu$. 
  Then 
  $$
  \theta_{\lambda \uparrow \mu}:=\theta^{- w_{0,M}\mu-\eta \uparrow  -w_{0,M}\lambda-\eta} : \H^0(\Shbar^{\tor},\omega(\lambda+\eta)) \to \H^0(\Shbar^{\tor},\omega(\mu+\eta))
  $$
  is injective for generic $\lambda \in X_1(T)$. Here $\eta$ is the weight of the canonical bundle of $\Shbar$. 
\end{theorem}

In the setup from above, such an $\alpha$ is unique.
For the other simple roots $\beta \in \Delta_{G_i}$, we have that $\theta^p_{\beta}=0$ and 
$\pi_* \theta^{- w_{0,M}\mu-\eta \uparrow  -w_{0,M}\lambda-\eta}$ is a linear map between automorphic vector bundles
coming from a map of $P$-representations,
so those are not as interesting.

\subsubsection{(Generic) entailments}
As an application of all of the previous machinery we prove some examples of entailments of Serre weights, generalizing \cite{paper}. We remark that these
could not have been proved using only the techniques of \cite{paper}. 
Let $\mathfrak{m} \subset \mathbb{T}$ be a mod $p$ Hecke eigensystem appearing in the \'etale cohomology 
of $\Sh_{\overline{\Q}_p}$. We define its set of modular Serre weights to be
$$
W(\overline{\rho}_{\mathfrak{m}}):=\{F(\lambda): \H^{\bullet}(\Sh_{\overline{\Q}_p},\underline{F}(\lambda))_{\mathfrak{m}} \neq 0 \}.
$$
An entailment of Serre weights occurs when 
there exist $p$-restricted weights $\lambda \neq \lambda' \in X_1(T)$ such that for all $\mathfrak{m}$ we have
$F(\lambda) \in W(\overline{\rho}_{\mathfrak{m}}) \implies F(\lambda') \in W(\overline{\rho}_{\mathfrak{m}})$. We will  
relax the definition by allowing ourselves to restrict to sufficiently generic $\mathfrak{m}$. 
We will consider generic non-Eisenstein $\mathfrak{m}$ as defined below. 
\begin{defn} \label{generic-non-Eisenstein}
Let $\mathfrak{m} \subseteq \mathbb{T}$ be a maximal ideal. Assume that all the almost simple factors of $G^{\der}_{\overline{\Q}_p}$ are of type $A$ or $C_2$
and that we are in the setup of \Cref{generalities-SW}. 
We say that $\mathfrak{m}$ is non-Eisenstein if its associated mod $p$ global Galois representation 
$\overline{r}_{\mathfrak{m}}: \Gal_{F} \to \check{G}(\Fpbar)$ is irreducible. In the case that $\check{G}=\GSp_4$ we say that 
$\overline{r}_{\mathfrak{m}}$ is irreducible if it is so as a $\GL_4$-valued representation.
 We say that 
$\mathfrak{m}$ is generic if there exists an auxiliary prime $l \neq p$ such that the restriction of $\overline{r}_{\mathfrak{m}}$
to the decomposition groups of primes above $l$ is generic in the sense of \cite[Def 1.1]{Hamann-Lee}.
\end{defn}
If for some $\lambda \in X_1(T)$ 
there exists a set $\{\lambda_i\} \subseteq X_1(T)$ not containing $\lambda$ such that 
for all generic non-Eisenstein $\mathfrak{m}$ we have that
$F(\lambda) \in W(\overline{\rho}_{\mathfrak{m}}) \implies F(\lambda_i) \in W(\overline{\rho}_{\mathfrak{m}})$ 
for some $i$ (which may depend on $\mathfrak{m}$), we say 
that it is a \textit{weak entailment}. Finally, we will say a (weak) entailment is generic if
the implication $F(\lambda) \in W(\overline{\rho}_{\mathfrak{m}}) \implies F(\lambda_i) \in W(\overline{\rho}_{\mathfrak{m}})$ 
holds for \textit{generic} $\lambda \in X_1(T)$ as in \Cref{defn-generic}. 
We prove some examples of generic weak entailments using theta linkage maps. 

In \cite{paper} we proved a generic weak entailment for $\GSp_{4,\Q_p}$.
In the next theorem, let $\Shbar$ be a compact unitary Shimura variety such that $G_{\overline{\Q}}=\GL_4 
\times \mathbb{G}_m$.
When working with unitary Shimura varieties we will often identify the representation theory of $G_{\Fpbar}$ with 
the one of $\GL_n$, suppressing the twist coming from the $\mathbb{G}_m$.
In this case there are $6$ $p$-restricted alcoves $C_i$ for $0 \le i \le 5$. Given $\lambda_0 \in C_0$, 
let $\lambda_i \in C_i$ be the corresponding affine Weyl translates. 
\begin{theorem} \label{entailment-intro} [\Cref{entailment-GL4}]
  Let $\Sh/\mathcal{O}$ be a compact unitary Shimura variety of signature $(3,1)$ 
  with $p$ split in the quadratic imaginary field.
  Let $\mathfrak{m} \subseteq \mathbb{T}$ be a generic non-Eisenstein maximal ideal.
  Then the statement
  $$
	F(\lambda_0) \in W(\overline{\rho}_{\mathfrak{m}}) \implies F(\mu) \in W(\overline{\rho}_{\mathfrak{m}})
  $$
  for some $F(\mu) \in \textnormal{JH}[V(\lambda_5)_{\overline{\F}_p}]$ holds for $\lambda_0 \in C_0$ generic.
  Moreover,
	$F(\mu) \neq F(\lambda_0)$.
  The same result holds for $p$ inert if \Cref{assumption-vanishing}(1) holds.
	\end{theorem}

  In the case that $p$ is inert $G_{\Q_p}=GU(4)$ with respect to $\Q_{p^2}/\Q_p$
  is a non-split unitary group, which have been much less studied 
  than split groups. A result like this might suggest that the phenomenon of entailments 
  is more related to the representation theory of $G_{\Fpbar}$ as opposed to the one of $G(\F_p)$.
   For $GU(4)_{\Q_p}$ \Cref{assumption-vanishing}(1) is not known, but it is nevertheless 
  within reach with the current techniques \footnote{In fact \Cref{assumption-vanishing}(1) is implied 
  by the long-announced work of \cite{ampleness}.}. 
  
  \begin{remark}
  Assuming the weight part of Serre's conjecture as in \cite[Conj 3.2.7]{Gee-Herzig-Savitt} these entailments are 
  related to certain Breuil--Mezard cycles containing more than one irreducible component 
  of the reduced Emerton--Gee stack. Then \Cref{entailment-intro} would conjecturally imply that for generic $\lambda_0 \in C_0$ we have a containment 
  of Breuil--Mezard cycles
   $$
  \mathcal{Z}_{F(\lambda_0)} \subseteq \mathcal{Z}_{V(\lambda_5)}.
  $$
  Based on \cite{Le-Hung-Feng}, in work in progress \cite{Le-Hung-Lin} prove that if the 
  Breuil--Mezard conjecture is 
  true, then for generic  $\lambda_0 \in C_0$ (see \Cref{defn-generic}(3)) one must have  
  $\mathcal{Z}_{F(\lambda_0)} \subseteq \mathcal{Z}_{F(\lambda_5)}$ for $\GL_{4,\Q_p}$. 
  Moreover, their results would imply that this is the only generic entailment for $\GL_{4,\Q_p}$,
  the generic (weak) entailment for $\GSp_{4,\Q_p}$ proved in \cite{paper} is the only one, 
  and no generic entailment exists for $\GL_3/\Q_p$. 
  Thus, we are able to detect the generic (weak) entailments for the two unramified examples of lowest rank. 
  \end{remark}

  We can also produce non-generic (weak) entailments,
   using the linkage maps, basic theta operators, 
  and Hasse invariants. Here a non-generic entailment is defined as an entailment that holds 
  for a class of weights 
  which is not generic. Informally, they only hold for a particular class of weights close to the boundary 
  of the alcoves.
   Moreover, we can construct examples of interesting congruences between forms with non-regular 
  weights and forms with regular weights. We give examples of these for $\GL_3$ in
  \Cref{non-generic-entailment-GL3}.

  \subsubsection{The de Rham realization functor}
  As an outcome of the way we construct the theta linkage maps, in \Cref{section7} we will also construct
exact functors
\begin{equation} \label{realization-functor}
\Psi: D^{b}(\mathcal{O}_{P,\Fpbar}) \xrightarrow{f} D^{b}((\Shbar/\Fpbar)_{\crys}) \xrightarrow{g} D^{b}(C_{\Shbar/\Fpbar}).
\end{equation}
Here $\mathcal{O}_{P,\Fpbar}$ is a category of $(U\mathfrak{g},P)_{\Fpbar}$-modules modelled after the classical complex category $\mathcal{O}$, 
$(\Shbar/\Fpbar)_{\crys}$ is the crystalline topos over $\Fpbar$, and $C_{\Shbar/\Fpbar}$ is the category 
of $\mathcal{O}_{\Shbar}$-modules with $\Fpbar$-linear maps. The image of $f$ lands in the subcategory 
of complexes of crystals on 
$(\Shbar/\Fpbar)_{\crys}$.
 The map $g$ in \eqref{realization-functor} is the pushforward 
from the crystalline topos 
to the Zariski topos. Moreover, $\Psi$ extends to toroidal compactifications. 
For $W \in \Rep(P)$ let $\Ver_{P}(W)=U\mathfrak{g} \otimes_{U\mathfrak{p}} W$ be a parabolic Verma module.
The functor satisfies that for  V a $G$-representation $\Psi([V])=\mathcal{V}^{\vee} \otimes \Omega^{\bullet}_{\Shbar}$ is the de Rham complex of the associated 
vector bundle with connection $\mathcal{V}:=G^{\dR} \times^{G} V$. Thus, $\Psi$ computes de Rham cohomology. 
 It also satisfies that $\Psi(\Ver_{P}(W)[0])=\mathcal{W}^{\vee}[0]$ where $\mathcal{W}:=P_{\dR} \times^{P} W$ is the
 vector bundle associated to $W$. Moreover,
  maps between parabolic Verma modules are sent by $\Psi$ to 
differential operators (these are often theta linkage maps, see \Cref{intro-Verma}). 
We will use this functor in a companion paper \cite{de-Rham-paper} to compute de Rham cohomology of some low dimensional Shimura varieties 
in terms of generalized mod $p$ BGG decompositions.
\subsection{The methods} \label{section0.3}
\subsubsection{Verma modules and differential operators} \label{intro-Verma}
The construction of all the operators becomes a simple representation/group theoretic problem
 once we establish 
the following relation between Verma modules and differential operators on the (flag) Shimura variety. 
For a smooth scheme $X/S$ such that $p^nS=0$ let $P_{X/S}$ be the sheaf of PD principal parts. It is the structure sheaf 
of the PD envelope of $\Delta: X \to X \times_{S} X$. It is a natural $\mathcal{O}_{X}$-bimodule. 
Let $P^n_{X/S}$ be the $(n-1)$th formal neighbourhood of the diagonal. It is defined for any $S$.
Then (crystalline) differential operators between 
two $\mathcal{O}_{X}$-modules $V_{1},V_2$
are defined as $\DiffOp_{X/S}(V_1,V_2):= \cup \Hom_{\mathcal{O}_{X}}(P^n_{X/S} \otimes V_1,V_2)$, where $P^n_{X/S} \otimes V_1$
carries its left $\mathcal{O}_{X}$ structure. Write $D_{X}:=\DiffOp(\mathcal{O}_X,\mathcal{O}_X)$.
There is a map $\DiffOp_{X/S}(V_1,V_2) \to \Hom_{\mathcal{O}_S}(V_1,V_2)$ which generally fails to be injective. 
We then prove the following. Let $\Sh/\mathcal{O}$ with residue field $k$, and let $R \in \{\mathcal{O},k\}$.
Let $F_{B} : \text{Mod}_{R}(B) \to \textnormal{QCoh}(\mathcal{O}_{\flag_R})$ be the map defined by the canonical $B$-torsor 
over $\flag_R$, we will also use the notation $\mathcal{V}:=F_{B}(V)$.
For $V \in \Rep_{R}(B)$ let $\Ver_{B}(V):=U\mathfrak{g} \otimes_{U\mathfrak{b}} V$.
\begin{theorem} \label{Verma-diff-intro}
  [\Cref{canonical-iso-Pm}, \Cref{creator-diff}]
\begin{enumerate}
\item There exist canonical
 isomorphisms $e_{V}:F_{B}(\Ver_{B}(V)) \cong D_{\flag/R} \otimes \mathcal{V}$ compatible with the filtrations on both sides. 
\item They satisfy a number of compatibilities, for instance the map 
$\Ver_{B}(\Ver_{B}(1)) \to \Ver_{B}(1)$ corresponds to the composition map $D_{\flag}  \otimes D_{\flag} \to D_{\flag}$.
\item There is a natural map 
$$
\Hom_{(U\mathfrak{g},B)}(\Ver_{B}V_1,\Ver_{B}V_2) \to \DiffOp_{\flag_R}(\mathcal{V}^{\vee}_2,\mathcal{V}^{\vee}_1)
$$
producing Hecke equivariant maps away from $p$, 
which is compatible with composition on both sides. 
\end{enumerate}
The analogue result works for $P$ and $\Sh/R$, and both extend to toroidal compactifications. 
\end{theorem}

Part $(3)$ follows formally from parts $(1)$ and $(2)$. Over the complex numbers, the functor in $(3)$ was known 
from \cite{faltings-chai}, and used in the integral setting by \cite{Mokrane-Tilouine}, \cite{Lan-Polo}.
Our treatment emphasizes the use of crystalline differential operators, and we give a more detailed description of 
these isomorphisms, 
which is needed to get a finer control on the operators produced in this way. Over the complex numbers one
can prove the theorem above using 
 the Borel embedding from the universal cover of $\Sh_{\C}$
to $G/P$. Integrally, we lack an analogous global period map, so we will use the pointwise Grothendieck--Messing period map instead. 
\begin{remark}
There is also a notion of non-crystalline differential operators $\DiffOp_{G}$ 
\cite[\href{https://stacks.math.columbia.edu/tag/09CH}{Tag 09CH}]{stacks-project}, and a "PD"-Verma 
module using the algebra of distributions $\Ver^{PD}_{B}(V)$.  To illustrate the difference between the two,
for $X=\mathbb{A}^1$, we have that
$D_{G,X}$ is generated by $\frac{1}{n!}(\frac{d}{dx})^n$,
while $D_{X}$ is generated by $(\frac{d}{dx})^n$.
 However, the functor from point $(3)$ 
does not extend to a functor 
$\Hom_{(U(G),B)}(\Ver^{PD}_{B}V_1,\Ver^{PD}_{B}V_2) \to \DiffOp_{G}(\mathcal{V}^{\vee}_2,\mathcal{V}^{\vee}_1)$
in any reasonable sense, see \Cref{PD-diff-ops-dont-work}. 
Similarly, there is a notion of baby Verma modules in characteristic $p$, and 
"Frobenius differentials" on $\flag_{\Fpbar}$, but the above functor does not relate the two. See \Cref{Baby-Vermas-dont-work}.
This justifies that crystalline differential operators are the right kind of operators for our purpose.
\end{remark}
We prove $(1)$ and $(2)$ in \Cref{section2.2}, and in \Cref{section2.1} we introduce some of necessary
 ingredients.  
First, we construct the maps $e_{V}$ on PD formal completions of $\Fpbar$-points, 
where it is given by the differential 
of the Grothendieck--Messing period map, which is \'etale in an appropriate sense.
We use that on the flag variety $G/B$ there is an analogous relation 
between Verma modules and $G$-equivariant differential operators which is much easier to prove i.e. 
we use the equivalence between $B$-representations and $G$-equivariant sheaves on $G/B$.
All the compatibilities will be transferred from the flag variety in this way.
Then we construct the isomorphism $e_{V}$. 
First, we can reduce to constructing it at the level of the degree at most $1$ part $\Ver^{\le 1}(V)$ using 
the compatibility of $\Ver_{P}(\Ver_{B}) \twoheadrightarrow
\Ver_{B}$ with composition on $D_{\flag}$. We construct $e^{\le 1}_{V}$ for $V \in \Rep(G)$ using the Gauss--Manin connection 
attached to $\mathcal{V}^{\vee}$ and the Kodaira--Spencer isomorphism. Namely, for $V \in \Rep(G)$ the associated 
bundle on $G/B$ is equipped with a connection, which defines an isomorphism of $B$-representations
$$
\Ver^{\le 1}_{B}(V)\cong V \otimes \Ver^{\le 1}_{B}(1).
$$
The dual of the stratification $P^{1} \otimes \mathcal{V}^{\vee} \cong \mathcal{V}^{\vee} \otimes P^{1}$
associated to the Gauss--Manin connection on $\mathcal{V}^{\vee}$ and the Kodaira--Spencer isomorphism 
defines $e^{\le 1}_{V}$ after applying $F_{B}$ to the isomorphism above.
For general $V \in \Rep(B)$ we construct it by devissage. A priori this construction 
is not well-defined, but since it agrees with the construction on PD formal completions, it is well-defined and all the required properties hold.

\subsubsection{Basic theta operators} \label{subsubsection-basic-theta}
We construct the basic theta operators in \Cref{section3}.
Consider the map $\flag \to \GF \to \Sbt_{G}=[B \times B \backslash G]$ as in \cite[\S 3.1]{IK-hasses}. All our operators will come from 
pullback from this map. We remark that by construction of the map, all automorphic vector bundles arise 
as pullbacks along this map.
Let $U \subset \flag$ be the pullback of the open Bruhat cell $\Sbt_{1}:=[B \times B \backslash Bw_0B]$. 
Over $\Sbt_1$, 
there is a $T$-reduction of the $B$-torsor $[1 \times B \backslash G] \to\Sbt_{G}$, which we can pull back 
to a $T$-reduction over $U$ of the canonical $B$-torsor over $\flag$. 
Therefore, we obtain a functor $F_{T}: \Rep(T) \to \text{Coh}(U)$ which is compatible with $F_{B}: \Rep(B) \to \text{Coh}(\flag)$
in the obvious sense. 
Thus, applying $F_T$ and using \Cref{Verma-diff-intro}(1,2) yields a functor 
\begin{equation} \label{U-diff}
\Hom_{T}(\lambda,\Ver_{B}(\mu)) \to \cup_{n} \Hom_{\mathcal{O}_{U}}(P^n \otimes \LL(-\mu),\LL(-\lambda))=\DiffOp_{U}(\LL(-\mu),\LL(-\lambda)).
\end{equation}
The left-hand side is in bijection 
with elements of $U\mathfrak{u}_{B}$ of weights $\mu-\lambda$. In this way, for each $\gamma \in \Phi^{+}$ 
we get operators $\tilde{\theta}_{\gamma}: \LL(\lambda) \to \LL(\lambda+\gamma)$ over 
$U$ corresponding to $x_{-\gamma} \in \mathfrak{u}_{B}$.
\begin{remark}
Even though it is framed in a group-theoretic language, the $\tilde{\theta}_{\gamma}$ can also be defined 
more classically in terms of the Gauss--Manin connection, see \Cref{easy-lemma-theta}(1).
\end{remark}
 To extend the $\tilde{\theta}_{\gamma}$ to $\flag$, we must determine the appropriate Hasse 
invariants by which to post-compose them. Crucially, via \Cref{Verma-diff-intro} the functor above \eqref{U-diff}
is pulled-back 
from the functor 
$$
\Hom_{T}(\lambda,\Ver^{\le 1}_{B}(\mu)) \to \Hom_{\Sbt_{1}}(\underline{\Ver^{\le 1,\vee}_{B}(\mu)},\underline{-\lambda}),
$$
where $\Sbt_{1} \subseteq \Sbt_{G}$ is the open cell, and the underline denotes the functor 
sending representations of $B \times B$
to vector bundles on $\Sbt$ (here both representations factor through the first projection). 
The Hasse invariants $H_{\alpha}$ are also defined via pullback from $\Sbt$, so 
to extend $\tilde{\theta}_{\gamma}$ we need to figure out for which powers $n_{\alpha,\gamma} \ge 0$ the map
$$
\underline{\Ver^{\le 1,\vee}_{B}(\lambda+\gamma)} \to \underline{-\lambda} 
\xrightarrow{\prod H^{n_{\alpha,\gamma}}_{\alpha}} \underline{(-\mu_1,-\mu_2)}
$$
extends to $\Sbt$. Then the pullback of this composition to $\flag$ defines $\theta_{\gamma}$ via
 \Cref{Verma-diff-intro}.
 This becomes a simple problem in Schubert calculus, which we can compute explicitly, following 
 \cite{IK-hasses}.
For $(5)$ and $(6)$ of \Cref{basic-theta-operaotrs} we prove some relations between the restriction of basic 
theta operators to codimension $1$ strata. These restrictions are well-defined when the associated Hasse invariant has a simple zero
(which always holds if $G^{\text{der}}$ is simply connected), 
by point $(2)$. 
By our construction, we can 
just check these relations on $\Sbt$, i.e. by checking that two maps on closed codimension $1$ stratum $\overline{\Sbt}_{\alpha}$
agree. 
 We illustrate $(5)$ and $(6)$ for $G=GU(2,1)$, with simple roots 
$\alpha_1,\alpha_2$. In \Cref{restriction-simple-theta} and \Cref{general-restriction-strata} we prove that 
the restriction to the closed codimension $1$ stratum $\overline{D}_{\alpha_i} \subseteq \flag$ 
of $\theta_{\alpha_i}$ on $\LL(\lambda)$
is
$\theta_{\alpha_i,\overline{D}_{\alpha_i}}=\langle \lambda, \alpha^{\vee}_i \rangle B_{\alpha_i}$,
where $B_{\alpha_i}$ is a section on $\overline{D}_{\alpha_i}$ which is non-vanishing on the interior.
This proves $(5)$.
Furthermore, $\theta_{\alpha_1+\alpha_2,\overline{D}_{\alpha_2}}=
B_{\alpha_2}\theta_{\alpha_1,\overline{D}_{\alpha_2}}$. Now, let $f \in \H^0(\Shbar,\omega(\lambda))$ for a generic 
$\lambda \in X_1(T)$. Suppose that $\theta_{\alpha_1+\alpha_2}(f)=0$,
since the weight is small we may assume it is not divisible by a Hasse invariant, by \Cref{assumption-vanishing}(2).
From the above 
\begin{align*}
\theta_{\alpha_1+\alpha_2, \overline{D}_{\alpha_2}}(f)=0 &\implies H_{\alpha_2} \mid \theta_{\alpha_1}(f)
\implies \theta_{\alpha_1}(f)=0 \\
 &\implies \theta_{\alpha_1,\overline{D}_{\alpha_1}}(f)=0 \implies H_{\alpha_1} \mid f
\implies f=0,
\end{align*}
where the second implication also follows from \Cref{assumption-vanishing}(2). 
In general, point $(6)$ is proved via 
an inductive process 
like this one. 

Points $(2)$, $(3)$, $(4)$ of \Cref{basic-theta-operaotrs} can be proved over $U$. 
For this, we give an explicit description of the isomorphism $e_{V}$ of \Cref{Verma-diff-intro}
over $U$, which by its 
construction reduces to an explicit description of the Gauss--Manin connection on some local basis.
This is treated in \Cref{section2.3}. 
The crucial fact that we need (\Cref{special-embedding}), 
is that there is a $G$-representation $V$ such that after some compatible choice of 
Borels the induced full Hodge and conjugate flag on $\mathcal{V}$ are in general position. That is, 
there exists some choice of Borels of $\GL(V)$ such that the embedding $\Sbt_{G} \hookrightarrow \Sbt_{\GL(V)}$ 
maps the open cell to the open cell. In fact, one can take $V$ to be the Hodge embedding. 
Then we can describe $\nabla$ on $\mathcal{V}_{U}$ with respect 
to an appropriate section of the $T$-torsor over $U$ which we call \textit{adapted elements} 
(\Cref{adapted-element}),
in terms of the Kodaira--Spencer isomorphism. 
The definition of an adapted element uses the property that the full conjugate filtration of $\mathcal{V}$, 
which by our assumption can be expressed in terms of the Hodge filtration, 
is flat for $\nabla$ with trivial $p$-curvature on its graded pieces. As a toy example,
for the modular curve, an adapted element would be given by
a local basis $\{e_1,e_2\}$ of $\mathcal{H}$ with $e_2 \in \Ker(V)$ such that $\nabla(e_2)=0$, and
 $e_1 \in \omega$ such 
that $\nabla(e_1)=e_2 \otimes \KS(e^2_1)$. This basis is well-defined up to $p$th powers. 
In such a basis we get a simple expression for $e_{V}$ in \Cref{canonical-iso-U},
which proves point $(2)$. Point $(3)$ 
follows from compatibility of \eqref{U-diff} with composition in an appropriate sense. Roughly, 
if a map of Verma modules 
is given by an element $f(x_{-\gamma})$ in $U\mathfrak{u}^{-}_{B}$, then its corresponding theta operator on $U$ will 
be $f(\tilde{\theta}_{\gamma})$.
To prove this compatibility with composition we use the explicit description of $e_{V}$ when fixing 
an adapted element. 
Furthermore, $(4)$ reduces to a computation involving the $p$-curvature map
on $\mathcal{V}$, where the $p$th iterate of a derivation 
shows up. In this computation we use the fact that $\mathcal{V}$ has the structure of a \textit{de Rham F-gauge}.
That is, the Kodaira--Spencer map on $\mathcal{V}$ agrees with 
the map induced by the $p$-curvature and the Cartier isomorphism, see \Cref{conjugate-p-curvature}.

\subsubsection{The linkage maps}
We construct the theta linkage maps via \Cref{Verma-diff-intro}, 
by constructing maps of Verma modules in characteristic $p$.
The key observation is that 
$$
\Hom_{(U\mathfrak{g},B)}(\Ver_{B}(\lambda)_{\Fpbar},\Ver_{B}(\mu)_{\Fpbar})
$$ 
is invariant under $p$-translation of the weights,
 so that we can reduce to the case where
$\lambda \uparrow s_{\gamma} \cdot \lambda=\mu$, in which case there exists a map in characteristic $0$, by the classical 
theory. See \Cref{Verma-linkage}. Then we simply define theta linkage maps as the differential operators coming from maps 
of Verma modules, and similarly for parabolic theta linkage maps and parabolic Verma modules. 
To prove that theta linkage maps are related to combinations of basic theta operators, we prove that when restricted to $U$ the linkage maps 
are a combination of $\tilde{\theta}_{\gamma}$. This follows from
the compatibility with composition of \eqref{U-diff}, as explained above. 
 To prove the generic injectivity of
\Cref{injective-linkage-intro}, by the compatibility of \eqref{U-diff} with composition,
the simple theta linkage map is given by repeated application of a 
basic theta operator. Then we use point $(4)$ in \Cref{basic-theta-operaotrs} to reduce the injectivity of
the linkage map to the injectivity of  
a basic theta operator. 
\subsubsection{The entailments}
To prove the entailments we need $3$ facts. 
\begin{enumerate}
\item The injectivity of a certain theta linkage map in \Cref{injective-linkage-intro}, or the injectivity 
of any operator we are using.
\item The fact that for generic non-Eisenstein $\mathfrak{m}$ and generic $\lambda \in X_1(T)$ we have the relation 
$$
\H^0(\Shbar^\tor,\omega(\lambda+\eta))_{\mathfrak{m}} \neq 0 \iff \H^{d}_{\et}(\Sh_{\overline{\Q}_p},V(\lambda)_{\Fpbar})_{\mathfrak{m}} \neq 0,
$$
where $\eta$ is the weight of the canonical bundle of $\Shbar$. This is proved in \Cref{coherent-to-betti-GL4}
for instance.
 The forward direction follows 
from the vanishing results for $\H^1$ coherent cohomology of \cite{alexandre}, \cite{Deding-unitary} (where 
we use that $\mathfrak{m}$ is non-Eisenstein to identify canonical and subcanonical coherent cohomology)
to lift coherent cohomology to characteristic $0$, and then use 
the BGG decomposition over $\C_p$. The other direction is more delicate, first one 
uses the vanishing results of \'etale cohomology of \cite{Hamann-Lee} \cite{caraiani-scholze-compact} 
to lift \'etale cohomology to characteristic $0$.
Then one has to check that for non-Eisenstein eigensystems, if an element of the archimedean L-packet 
contributes to the middle degree \'etale cohomology in characteristic $0$,
then so does the holomorphic component of the L-packet, which 
contributes to coherent cohomology in degree $0$. In the $\GSp_4/\Q$ case and the Harris--Taylor case 
one can prove this by
using the fact that, under the assumptions on $\mathfrak{m}$ 
the Galois representation $\overline{r}_{\mathfrak{m}}$ is contained in middle degree \'etale cohomology. 
\item 
Finally, we need the decomposition of dual Weyl modules $V(\lambda)_{\Fpbar}$ for
 $\lambda \in X_1(T)$ into Serre weights.
 This is a difficult 
 problem, so even for $G=\GL_5$, there are no results written in the literature \footnote{Although Lusztig's conjecture 
 \cite[\S II.7.20]{Janzten-book}
 describes the constituents of $V(\lambda)_{\Fpbar}$ in terms of Kazhdan--Lusztig polynomials 
 for $p$ sufficiently large, which can in principle be computed.}, but 
in our examples these decompositions have been explicitly figured out.
\end{enumerate}
 For simplicity let us illustrate the generic weak entailment for $\GSp_4$. The one for $\GL_4$ in
 \Cref{entailment-GL4} is proved in the same way. Let $\mathfrak{m}$ be generic non-Eisenstein, and $\lambda_0 \in C_0$. 
 Let $\lambda_i \in C_i$ be its affine Weyl reflections in the other $p$-restricted alcoves.
  We want to prove that the statement '$F(\lambda_0) \in W(\overline{\rho}_{\mathfrak{m}})$ implies
  $F(\lambda_{i}) \in W(\overline{\rho}_{\mathfrak{m}})$
 for either $i=1$ or $i=2$' holds for $\lambda_0 \in C_0$ generic. 
For point $(1)$ we know that the theta linkage map
 $\H^0(\Shbar,\omega(\lambda_1+\eta))_{\mathfrak{m}} \xrightarrow{\theta_{\lambda_1 \uparrow \lambda_2}}
  \H^0(\Shbar,\omega(\lambda_2+\eta))_{\mathfrak{m}}$
is injective, and in point $(3)$ we have following relations in 
the Grothendieck group of $\Rep_{\Fpbar}G(\F_p)$: $[V(\lambda_1)_{\Fpbar}]=F(\lambda_0) + F(\lambda_1)$, 
$[V(\lambda_2)_{\Fpbar}]=F(\lambda_1)+ F(\lambda_2)$. 
Then $F(\lambda_0) \in W(\overline{\rho}_{\mathfrak{m}})$ implies 
\begin{align*}
& \H^3_{\et}(F(\lambda_0))_{\mathfrak{m}} \neq 0 \implies \H^3_{\et}(V(\lambda_1)_{\Fpbar})_{\mathfrak{m}} \neq 0
\implies \H^0(\Shbar,\omega(\lambda_1+\eta))_{\mathfrak{m}} \neq 0 \\
& \implies \H^0(\Shbar,\omega(\lambda_2+\eta))_{\mathfrak{m}} \neq 0
 \implies \H^3_{\et}(V(\lambda_2)_{\Fpbar})_{\mathfrak{m}} \neq 0,
\end{align*}
which means $F(\lambda_i) \in W(\overline{\rho}_{\mathfrak{m}})$ for $i=1$ or $i=2$. The first arrow is implied by 
$[V(\lambda_1)_{\Fpbar}]=F(\lambda_0) + F(\lambda_1)$ and the concentration of \'etale cohomology in middle degree. 
The second is point $(2)$ above, the third follows by point $(1)$, and the last by point $(2)$ again. 
Finally, the decomposition $[V(\lambda_2)_{\Fpbar}]=F(\lambda_1)+ F(\lambda_2)$ yields the final result. 
 Thus, using the more general injectivity result \Cref{injective-linkage-intro}, we could 
prove more general entailments if we understood the constituents of certain linked dual Weyl modules, 
but this gets very hard as the rank of $G$ grows.

\subsection{Organization}
In \Cref{section1} we explain all of the geometric tools that we will need to construct the operators
in the setting of Hodge type Shimura varieties. In particular, we prove a general version of 
the Kodaira--Spencer isomorphism, and  
its compatibility with a Grothendieck--Messing period map. Although not difficult, this result does not appear explicitly in the literature.
 In \Cref{section2} we prove the relation between Verma modules and 
differential operators on the Shimura variety, and we give an explicit description 
of this on the open stratum of $\flag_{\Fpbar}$. In \Cref{section3} we construct the basic theta operators, 
and we prove a number of properties about them. In \Cref{section4} we construct the theta linkage maps, and 
we prove the injectivity of a particular simple theta linkage map. In \Cref{section5} we prove that 
duality of differential operators is compatible with duality of (parabolic) theta linkage maps, defined by 
Serre duality on the flag variety. In \Cref{section6} we prove some generic and non-generic entailments. 
In \Cref{section7} we construct the de Rham realization functor, and we reprove the lowest alcove integral BGG
decomposition of \cite{Lan-Polo}.

This article supersedes most of the author's first version of \cite{paper}, which 
concerns the case of $\GSp_4$. The main novelties, apart from dealing with general 
Hodge type Shimura varieties, are as follows. \Cref{section2.1} and \Cref{section2.2} are mostly the same, but we also consider objects
coming from category $\mathcal{O}$. In \Cref{section2.3} we use a more general method, not relying on the existence 
of an ordinary locus. 
In \Cref{section3} we give a new, more conceptual construction of 
the basic theta operators, 
a proof of \Cref{pth-power-relation} without using Serre--Tate coordinates, 
and we prove relations between restrictions of basic theta operators to codimension $1$ 
strata. \Cref{section5} is new and so is the relation between Borel and parabolic theta linkage maps in \Cref{section4}.
\Cref{section6} is new, although the method to obtain the entailments is the same. 
In \Cref{section7} we introduce the de Rham realization functor.

\subsection{Acknowledgments}
The contents of this article are part of the author's PhD thesis. For that I thank my PhD advisor George Boxer. This was 
supported by the Engineering and Physical Sciences Research Council [EP/S021590/1]. 
The EPSRC Centre for Doctoral Training in Geometry and Number Theory 
(The London School of Geometry and Number Theory), University College London, King's College London and 
Imperial College London. Part of its writing was done while the author was a member of the Max Planck Institute for Mathematics 
in Bonn, that I thank for their hospitality.

\section{The Shimura variety} \label{section1}

\subsection{Representation theoretic notation} \label{rep-theory}
Let $G/\F_p$ be a reductive group, such 
that it becomes split over some finite extension $k/\F_p$. Fix a choice $T \subset B \subset P \subset G$
of maximal torus, Borel subgroup and a parabolic subgroup, with $T,B$ defined over $\F_p$, and $P$ over $k$.
Let $P=MU_{P}$ be the Levi decomposition of 
$P$. Let $B^{-}$ be the opposite Borel, uniquely determined by $B \cap B^{-}=T$.
\begin{itemize}
\item Let $X^{*}(T)=\Hom_{\Fpbar}(T,\G_m)$ and $X_{*}(T)=\Hom_{\Fpbar}(\G_m,T)$ be the group of characters and cocharacters respectively.
Let $W=N_{G}(T)/Z(T)$ be the Weyl group of $G/k$.
The Galois group $\Gamma=\Gal(k/\F_p)$ acts on them.
Let $\sigma \in \Gamma$ be the Frobenius $x \mapsto x^{p}$.
\item Let $\Phi \subset X^*(T)$ be the set of roots and $\Phi^{\vee} \subset X_{*}(T)$ the set of coroots.
\item Let $\mathfrak{g}$ be the Lie algebra of $G$.
For $\alpha \in \Phi$ choose $x_{\alpha}\in \mathfrak{g}$ so that 
$\{x_{\alpha}: \alpha \in \Phi\}$
is a Chevalley basis for $G_{\Fpbar}$. Concretely, $h_{\alpha}:=[x_{\alpha},x_{-\alpha}]$ acts 
as the scalar $\langle \lambda, \alpha^{\vee} \rangle$ on a vector of weight $\lambda$. 
Also for $\alpha,\beta \in \Phi$,  $[x_{\alpha},x_{\beta}]=(c+1)x_{\alpha+\beta}$ up to a sign 
when $\alpha+\beta \in \Phi$, 
and $c$ is the greatest positive integer such that $\alpha-c\beta \in \Phi$. 
\item For each $\alpha \in \Phi$  fix a root morphism 
$u_{\alpha}: \G_{a} \to U_{\alpha} \subset G$ over $\Fpbar$, determined by 
the property that $du_{\alpha}(1)=x_{\alpha}$.
They satisfy 
$$
t u_{\alpha}(x) t^{-1}=u_{\alpha}(\alpha(t)x)
$$
for all $x \in \G_a$, $t \in T$. 
\item 
Let $\Phi^{+} \subset \Phi$ be the set of positive roots, defined by the condition that $\alpha \in \Phi^{+}$
if $U_{\alpha} \subset B$, and let 
$\Delta \subset \Phi^{+}$ be a choice of simple roots. 
\item For $\alpha \in \Phi$ let 
$\phi_{\alpha}: \SL_2 \to G_k$ be the map corresponding to $u_{\alpha}$. It satisfies 
$$
\phi_{\alpha}\begin{pmatrix}
1 & t \\
0 & 1
\end{pmatrix}=u_{\alpha}(t),  \;\;\;
\phi_{\alpha}\begin{pmatrix}
	1 & 0 \\
	t & 1
	\end{pmatrix}=u_{-\alpha}(t).
$$
\item For $\alpha \in \Phi$ let $\alpha^{\vee} \in \Phi^{\vee}$
be the corresponding coroot. The tuple $(X^{*},\Phi, X_{*}(T),\Phi^{\vee})$ together with their natural pairings and Galois actions 
is the root system of $G/\F_p$.
\item For $\alpha \in \Phi$ let $s_{\alpha} \in W$ be the corresponding reflection. It acts on $X^*(T)$ as 
$\lambda \mapsto \lambda-\langle \lambda, \alpha^{\vee} \rangle \alpha$. Let $w_0 \in W$ be the longest element 
according to the Bruhat order. 
\item The parabolic $P$ corresponds to a subset $I \subset \Delta$, determined by 
the Lie algebra $\mathfrak{p}$ containing $x_{-\alpha}$ for $\alpha \in \Delta$ 
precisely when $\alpha \in I$. We will use the notation $P=P_I$. Let $W_{I}=W_{M}$ be the subgroup of $W$ 
generated by $s_{\alpha}$ for $\alpha \in I$, or equivalently, the Weyl group of $M$. Let
 $W^{I}=W^{M}$ be the 
set of minimal length representatives of $W_I \backslash W$.
\item (Dot action) Let $\rho \in X^*(T)$ be a choice of half sum of positive roots for $G_{\Fpbar}$. We define the 
 dot action of $W$ on $X^*(T)$ by $w \cdot \lambda:=w(\lambda+\rho)-\rho$.
 \item (Affine Weyl group) Let $W_{\text{aff}}=p\mathbb{\Z}\Phi^{+} \rtimes W$ be the affine Weyl group.
 The dot action of $W$ on $X^*(T)$ extends to an action of $W_{\text{aff}}$ where $p\mathbb{\Z}\Phi^{+}$ acts by 
 translation. Then $W_{\text{aff}}$ is generated by reflections $s_{\gamma,n}$ for $\gamma \in \Phi$ and $n \in \Z$,
 which are characterized by acting as 
 $s_{\gamma,n} \cdot \lambda=\lambda+(pn-\langle \lambda+ \rho,\gamma^{\vee} \rangle) \gamma$.
 \item We divide $X^{*}(T)-\rho$ into $\rho$-shifted alcoves, defined as the 
 interior of the regions defined by the hyperplanes $H_{n,\gamma}=\{\langle \lambda+\rho, \gamma^{\vee} \rangle=np\}$
 for $\gamma \in \Phi$ and $n \in \Z$.
 \item  We say that 
 $\lambda \in X^*(T)$ is $p$-small if $|\langle \lambda+ \rho, \gamma^{\vee} \rangle| <p$ for 
 all $\gamma \in \Phi$.
 The affine Weyl group acts simply transitively on $\rho$-shifted alcoves 
  and the stabilizer of the region of $p$-small weights is the finite Weyl group $W$.
\item  Let $X_{+}(T) \subset X^{*}(T)$ be the set of dominant weights $\lambda$, satisfying 
$\langle \lambda, \alpha^{\vee} \rangle \ge 0$ for every $\alpha \in \Delta$. Similarly, let $X_{M,+}(T)$
be the set of $M$-dominant weights, consisting of the $\lambda$ satisfying 
$\langle \lambda, \alpha^{\vee} \rangle \ge 0$ for all $\alpha \in \Delta_{M}$. 
\item For a linear algebraic group $H$ over a ring $R$, let $\Rep_{R}(H)$ be the category of algebraic representations 
of $H$ 
which are finite projective $R$-modules. It is equivalent to the category of 
modules for the algebra of distributions 
$U(H)$ which are also
Let $\Modfg_{R}(H)$ be the category of algebraic representations of $H$ which are finitely generated as $R$-modules. 
\end{itemize} 

Since it will be a key part of this article we highlight our definition of generic weights.

\begin{defn}[\Cref{defn-generic}] \label{defn-generic1}
Let $\lambda \in X^*(T)$ and $\epsilon \ge 0$. We say that $\lambda$ 
is $\epsilon$-generic with respect to $p$ if 
for each $\gamma \in \Phi$, there exists an integer $n_{\gamma}$ such that 
$$
n_{\gamma}p+\epsilon< \langle \lambda, \gamma^{\vee} \rangle < (n_{\gamma}+1)p-\epsilon. 
$$
\end{defn}

Similarly, for their relevance in this article, we state the definition of dual Weyl modules, irreducible algebraic 
representations, and Serre weights. 
\begin{defn} \label{dual-Weyl-module}
Let $G/\Z_p$ be a reductive group scheme, and $T \subseteq B$ a choice of maximal torus and Borel defined over $\Z_p$.
 Let $P \subseteq G$ be a parabolic subgroup defined over some $W(k)$, for $k/\F_p$ finite. 
Let $M$ be the Levi of $P$. 
\begin{enumerate}
\item For $\lambda \in X_{M,+}(T)$ let 
$$
V_{P}(\lambda):=\text{Ind}^{P}_{B} w_{0,M} \lambda=
\H^0(P/B,P \times w_{0,M}\lambda/B),
$$
where $B$ acts as $b(p,v)=(pb,w_{0,M}\lambda(b)v)$. Then $V_{P}(\lambda) \in \Rep_{W(k)} G$. It 
is naturally inflated from the induced representation $\text{Ind}^{M}_{B^{-} \cap M} \lambda$ of $M$. We say that  
$V_{P}(\lambda)$ is the \textit{dual Weyl module} for $M$ of highest weight $\lambda$.
When $P=G$ we will simply write $V(\lambda)=V_{G}(\lambda)$.
In the setup of \Cref{integral-model}, when $P$ is the Hodge parabolic we will use
the notation $W(\lambda):=V_{P}(\lambda)$.
\item Let $L_{M}(\lambda) \subset V_{P}(\lambda)_{\Fpbar}$ be the socle as an $M_{\Fpbar}$-representation.
We will use the notation $L(\lambda):=L_{G}(\lambda)$.
The
$L(\lambda)$ for $\lambda \in X_{+}(T)$ parametrize all irreducible representations of $G_{\Fpbar}$ \cite[Cor II.2.7]{Janzten-book}.
\item The set of $p$-restricted weights $X_1(T) \subseteq X_{+}(T)$ consists of those $\lambda$ such that 
$0 \le \langle \lambda, \alpha^{\vee} \rangle<p$ for all $\alpha \in \Delta$. Given $\lambda \in X_1(T)$
let $F(\lambda)\in \Rep_{\Fpbar}G(\F_p)$ be the restriction of $L(\lambda)$ to $G(\F_p)$. 
Then $\{F(\lambda): \lambda \in X_1(T)\}$ contains all the \textit{Serre weights}, i.e. irreducible representations in 
$\Rep_{\Fpbar}G(\F_p)$. Moreover, $F(\lambda)=F(\lambda')$ only if $\lambda-\lambda' \in 
(\sigma-\text{id})X^0(T):=(\sigma-\text{id})\{ \lambda \in X^*(T): 
\langle \lambda, \alpha^{\vee} \rangle=0 \text{ for all } \alpha \in \Delta\}$ 
\cite[Lem 9.2.4]{Gee-Herzig-Savitt}. 
\end{enumerate}
\end{defn}

\subsection{The integral model} \label{integral-model}
We recall some of the properties of the integral model of our Shimura variety. 
Consider the Siegel Shimura datum
$(\GSp_{2g},S^{\pm})$. We take $\GSp_{2g}/\Q$ to be the base change to $\Q$
of the reductive group over $\Z$
defined as the group of symplectic similitudes of the free abelian group $\Lambda$ of rank $2g$ together with a 
symplectic form $\langle \cdot,\cdot \rangle$. We choose a basis $\{e_1,\ldots,e_{2g}\}$ for $\Lambda$
so that $\langle \cdot, \cdot \rangle$ 
 is given by the matrix 
$
J=\begin{pmatrix}
0 & S \\
-S & 0 
\end{pmatrix}
$
with $S$ the matrix with $1$s on the antidiagonal and $0$s elsewhere. Denote by 
$$
0 \to L \to \Lambda \to L^{\vee} \to 0 
$$
the exact sequence where $L=\langle e_1,e_2,\ldots, e_{g} \rangle$ is maximal isotropic for $\langle \cdot,\cdot \rangle$. 
The Siegel double space $S^{\pm}$
is defined as the set of maps $h:\mathbb{S} \to \GSp_{2g,\mathbb{R}}$ such that 
$h$ induces a Hodge structure of weights $(-1,0)$, $(0,-1)$ on $\Lambda_{\mathbb{R}}$, and 
the pairing $(x,y) \mapsto \langle x, h(i)y \rangle$ is positive or negative definite on $\Lambda_{\mathbb{R}}$.
Let $(G,X)$ be a Shimura datum of Hodge type, that is, we fix an embedding of Shimura data 
$i: (G,X) \hookrightarrow (\GSp_{2g},S^{\pm})$ for some $g \ge 1$. For us, $(G,X)$ satisfies axioms 
(SV1,2,3) in Deligne's definition \cite[\S 6]{Milne-moduli} \cite{Deligne-Shimura}, with SV3 saying that $G^{\text{ad}}$ has no $\Q$-simple factor such 
that $h$ is trivial on it. 
Let $E_0/\Q$ be the reflex field of $(G,X)$.
Fix a prime $p$, and a prime $\mathfrak{p}$ of $E_0$ above $p$.
Let $E$ be the completion  of $E_0$ at the prime $\mathfrak{p}$, with 
ring of integers $\mathcal{O}$ and residue field $k$. 
Assume that there is a reductive model $G_{\Z_p}$ whose generic fiber is $G_{\Q_p}$, and fix such a choice.
Let $K \subseteq G(\A^{\infty})$
be a compact open neat subgroup of the form $K=K_pK^p$ with $K_p=G_{\Z_p}(\Z_p)$, i.e. it is a level hyperspecial 
at $p$. By the existence of $G_{\Z_p}$ we have that $p$ is unramified in $E_0$, since $G_{\Z_p}$ will be 
split over $W(k')$ for some $k'/\F_p$ finite. Therefore, $\mathcal{O}=W(k)$.
After possibly replacing $\Lambda$ by another lattice in $\Lambda \otimes \Q$, there exist tensors $s_{\alpha}
\in \Lambda_{\Z_p}^{\otimes}$ such that $G_{\Z_p} \subset \GL(\Lambda)$ is the group of transformations leaving $s_{\alpha}$
invariant \cite[Lem 2.3.1]{Kisin-integral-model}. For a module $M$ we will use $M^{\otimes}$ to denote the direct sum of all the combinations of tensor powers, duals,
 symmetric powers, and wedge powers of $M$. In particular, the symplectic pairing restricts to a pairing
 $\Lambda \times \Lambda \to \Z_p$. 

 Let $K'_p \subseteq \GSp(\Lambda)(\Q_p)$ be the stabilizer of $\Lambda_{\Z_p}$. Then we can choose a level 
 $K'=K^{'p}K^{'}_{p} \subseteq \GSp(\Lambda)(\A^{\infty})$ such that there is a closed embedding of Shimura varieties 
 \cite[Lem 2.1.2]{Kisin-integral-model}
 $$
\Sh_{K}(G,X) \hookrightarrow \Sh_{K'}(\GSp,S^{\pm})_{E}.
 $$
Let $\Sh_{\GSp,K'}/\Z_p$ 
be the canonical integral model for the Shimura variety for $(\GSp_{2g},S^{\pm})$. Then Kisin \cite{Kisin-integral-model},
Vasiu \cite{Vasiu}, and \cite{models-p2} (for $p=2$)
construct a canonical smooth integral model $\Sh_{K}/\mathcal{O}$ of $\Sh_{K}(G,X)$ together with a closed embedding 
$i: \Sh_{K} \hookrightarrow \Sh_{\GSp,K',\mathcal{O}}$. 
It is defined as the flat closure of $\Sh_{K,E} \to \Sh_{K',\GSp,\mathcal{O}}$. In \cite{Kisin-integral-model}
one also takes the normalization, but \cite{xu2021normalizationintegralmodelsshimura} proved 
that this was not necessary. 
Let $\Shbar:=\Sh_{\Fpbar}$ be the geometric special fiber. With the integral model defined, we can state our definition of generic weights for a statement involving the 
Shimura variety. 
\begin{defn}(\Cref{defn-generic}) \label{generic-Sh}
Fix a Shimura datum $(G,X)$ of Hodge type. We fix 
a maximal torus and a Borel $T \subseteq B \subseteq G_{\overline{\Q}}$. 
Let $X=\{X_{p,K}\}_{p,K}$ be a family of statements ranging over all rational primes $p$ and levels $K \subseteq G(\mathbb{A}^{\infty})$
which are hyperspecial at $p$, concerning 
 $\Sh_{K,\mathcal{O}}$ (where $\mathcal{O}$ lives over $\Z_p$) and some weights $\lambda \in X^*(T)$. 
We say that $X$ holds for \textit{generic} $\lambda \in X^*(T)$ if there exists $\epsilon \ge 0$, depending only on 
$G_{\Q}$, such that each statement $X_{p,K}$ holds for weights $\lambda$ which are $\epsilon$-generic with respect to $p$.
\end{defn}

Associated to $h$ there is a conjugacy class of cocharacters 
$\mu_{0} : \G_m \to G$ defined over $E$. Since $G_{\F_p}$ is quasi-split we can choose $B \subseteq G_{\Z_p}$
such that both its generic and special fibers are Borel subgroups. 
Choose a maximal torus $T \subseteq B_{\F_p}$, it lifts to a torus $T$ over $\Z_p$. Then 
there is an associated cocharacter $\tilde{\mu}: \mathbb{G}_m 
\to G$ defined over $\mathcal{O}$ as in \cite[2.2.3]{Zhang-EO} whose base change to $E$ lies in
 the conjugacy class of $\mu_{0}$. There exists a unique parabolic $P \subseteq G_{\mathcal{O}}$ 
 such that $B \subseteq P$ and such that it is 
 the parabolic associated to a
 conjugate of $\tilde{\mu} \in X_{*}(T)$. Namely, $P=P_{\mu}:=\{g \in G: \lim_{t \to 0} \mu(t)g\mu(t^{-1}) \;\;\text{exists} \}$,
 where $\mu \in X_*(T)$ is the dominant conjugate of $\tilde{\mu}$. 
We will fix such a tuple $(T,B,P,G)$ with all of them defined over $\Z_p$ except $P$, which is defined over $\mathcal{O}$.

\subsubsection{Hodge embeddings}
We will need some structural results on what the possible Hodge embeddings are, since that will restrict what $G_{\Q}$
can be. The most important result in this subsection is \Cref{special-embedding} ahead, 
which will be used repeatedly. 
\begin{defn} \label{simple-reductive}
We say that a 
connected semisimple group over a field $k$ is almost-simple if every proper normal algebraic subgroup is finite. 
Given a semisimple group $G$ over $k$, let $\{G_i: i \in I\}$ be the set of almost-simple normal subgroups of $G_{\overline{k}}$.
 Then there is an isogeny 
$\prod G_i \to G_{\overline{k}}$. Given $i \in I$ we define $G^i$ as the cokernel of $\prod_{j \neq i} G_j \to G_{\overline{k}}$, 
it is isogenous to $G_i$. 
 We will refer to $G_i$ as the almost-simple factors, and to 
$G^i$ as the almost-simple quotients. 
Then $\Gal(\overline{k}/k)$ permutes the $G_i$, and thus also the $G^i$. 
\end{defn}

\begin{lemma} \label{lift-rep}
Let $G$ be a reductive group over a field $k$, and $G^{\textnormal{der}} \to \GL(V)$ a representation over $\overline{k}$.
Then 
we can extend $V$ to a representation of $G_{\overline{k}}$.
\begin{proof}
We have that $G_{\overline{k}}=\mu \backslash(T \times G^{\text{der}})$ where $T$ is a split torus and $\mu=Z(G) \cap G^{\der} \subseteq T$ is
a finite
subgroup scheme. Thus, to extend $V$ we just need to extend the character by which $\mu$ acts on $V$ to $T$. 
We have that $\mu$ decomposes into a product
$\mu=\prod \mu_{n_i} \subseteq \prod \mathbb{G}_m=T$, with $\mu_{n_i}$ the multiplicative group of order $n_i$.
It is easy to see that a character of $\mu_{n_i}$ extends to a character of $\mathbb{G}_m$.
\end{proof}
\end{lemma}

Given $G$ a reductive group, we say that a cocharacter $\mu$ is minuscule if 
$\lvert \langle \mu, \gamma \rangle \rvert \le 1$
for every root $\gamma \in \Phi$. 
For a minuscule cocharacter $\mu$ of $G$, we say that $\mu$ is non-trivial on an almost-simple factor $G_i$ of 
$G^{\text{der}}$ if there exists a simple root $\alpha$ of $G_i$ such that $\langle \mu, \alpha \rangle=\pm 1$.
If moreover $\mu$ is dominant, i.e. $\langle \mu, \alpha \rangle \ge 0$ for all simple roots $\alpha$,
then there exists a unique $\alpha \in \Delta_{G_i}$ such that $\langle \mu, \alpha \rangle=1$, and
we say that it is the \textit{special node} for $\mu$ and $G_i$.
When talking about the special node for a conjugacy class of cocharacters we will 
do it with respect to a dominant representative $\mu$ in the conjugacy class. 
We will use the notation of \Cref{table-lie} for the roots of the different semisimple Lie algebras. 
We recall that we are assuming that our Shimura variety satisfies the axiom SV3.
\begin{prop}
  \label{Hodge-embeddings}
Consider the Hodge embedding $G_{\C} \hookrightarrow \GL(\Lambda)_{\C}$. 
The image of $G$ contains the scalars, and the associated cocharacter on $\GL(\Lambda)$ 
is conjugate to one of the standard minuscule cocharacters. Also, $G^{\textnormal{der}}$ does not contain 
an exceptional group as an almost-simple factor.  
Further, $G^{\der}_{\C}$ is a quotient by a finite group of a product $\prod H_i$, which we now describe.    
The action of $H_i$ on $\Lambda$ is a direct sum of minuscule representations $\Lambda_i$.
Moreover, $H_i$ is one of the following. 
\begin{itemize}
\item (Type $A$) $H_i$ is simply connected of type $A_{n}$ for $n \ge 1$. 
\begin{enumerate}
\item The special node of $\mu$ for $H_i$ can be any simple root. 
Assume that the special node for $\mu$ and $H_i$ is neither $\alpha_1$ nor $\alpha_n$. Then
 $\Lambda_i$ can only be the standard representation $\SL_{n+1} \hookrightarrow \GL_{n+1}$.
\item If the special node for $\mu$ is one of the above, then $\Lambda_i$ can be any minuscule representation, corresponding 
to the wedge powers of the standard representation.
\end{enumerate}
\item (Type $B$)  $H_i$ is simply connected of type $B_n$ for $n \ge 2$.
Then the special node of $\mu$ for $H_i$ must be $\alpha_1$.
There is only one option for $\Lambda_i$, corresponding to the spin representation. 
\item (Type $C$) $H_i$ is simply connected of type $C_n$ for $n \ge 2$. The special node 
for $\mu$ must be $\alpha_n$.
 There is only one minuscule representation $\Lambda_i$, which is the
 standard representation $\textnormal{Sp}_{2n} \hookrightarrow \GL_{2n}$.
\item (Type $D$) $H_i$ is of type $D_n$ for $n \ge 4$.
There are two options. Either $H_i=\textnormal{SO}_{2n}$ is of type $D_n$ for $n \ge 5$, 
with special node $\alpha_{n-1}$ or $\alpha_{n}$, in which case $\Lambda_i$ can only be 
the standard representation. Or $H_i$ is simply connected of type $D_n$, in which case 
$\Lambda_i$ can be either of the two half-spin representations. 
For $n \ge 5$ this implies 
that the special node is $\alpha_1$, and for $n=4$ either one of $\{\alpha_1,\alpha_{n-1},\alpha_{n}\}$
can be the special node. 

\end{itemize}
\begin{proof}
Let $H$ be the simply connected cover of $G^{\der}$. Then each almost-simple factor $H^{\text{sc}}_i$
of $H$ admits the complexification of a
\textit{symplectic representation} with finite kernel as in \cite[Def 10.2]{Milne-moduli} via the Hodge embedding.
Fix some $i$ and assume first that $\mu$ is non-trivial on it, and that $H^{\text{sc}}_i$ is contained in the 
complexification of a factor of non-compact
type of $H_{\mathbb{R}}$.
Then under these conditions, in \cite[Summary 10.8]{Milne-moduli} and the classification right before, they
describe the possibilities for what the special node for $H^{\text{sc}}_i$ and $\Lambda_i$ 
can be, as in the statement. 
To deduce the same properties about the action of $H^{\text{sc}}_i$ on $\Lambda$ for all $i$, 
we use the $\Gal_{\Q}$ action on the almost-simple factors,
and their subrepresentations of $\Lambda$. 
Since we are assuming that our Shimura variety satisfies axioms SV1,2,3, we have that 
each $\Q$-simple factor of $G^{\der}$ has a non-compact factor over $\mathbb{R}$.
Similarly, by SV3 each $\Q$-simple factor of $H$ contains 
an almost-simple factor
such that $\mu$ is non-trivial on it, so we deduce that each $H^{\text{sc}}_i$ acts 
on $\Lambda$ as in the statement. Therefore, if $\mu$ is non-trivial on $H_i$ its special node must be as stated.
Only in type $D_n$ with $n \ge 5$ and special node $\alpha_{n-1}$ or $\alpha_{n}$ it must happen that $\Lambda_i$ 
factors through $\textnormal{SO}_{2n}$. 
 Therefore, $G^{\der}$ is the quotient of $\prod H_i$ by a finite group. 
\end{proof}

\end{prop}

\begin{defn} \label{filtration-functoriality}
Let $G$ be a reductive group over a field, and $Q \subseteq G$ a parabolic subgroup. It is defined by a cocharacter 
$\mu \in X_*(G)$ in the sense that 
$Q=\{g \in G: \lim_{t \to 0} \mu(t)g\mu(t^{-1}) \;\;\text{exists} \}$. For a representation $\rho: G \to \GL(V)$ 
we say that the parabolic $Q_{V,\mu}$ associated to $(Q,\mu)$ is the 
one which is induced by $\rho \circ \mu$. 
In the case that $Q=P$ is the Hodge parabolic, there is a unique $\mu$ in the conjugacy class of the Hodge 
cocharacter lying in $X_{*}(T)$ such that $P$ is associated to $\mu$, 
so we will omit $\mu$ from the notation. Similarly, over $\Fpbar$, for $Q=P^{-(p)}$ we take its cocharacter to be 
$-\sigma \mu$.
\end{defn}

In order to reduce some results to the case of $\GL_n$ we will need the following important lemma. 
We will use it 
repeatedly in \Cref{section2.3}, and in the subsequent results that depend on it. 
Recall that $G$ is defined over $\mathbb{F}_p$ and $P$ over $k$.
Let $P^{(p)}=P \times_{k,\sigma} k$ denote the base change along the $p$th power Frobenius $\sigma$.
We consider it sitting inside $G$ via the base change map $P^{(p)} \to P$. For $z \in W$, let $B^{z}=\dot{z}B\dot{z}^{-1}$
for some representative $\dot{z} \in N(T)$ of $z$. 
 Before proceeding with the statement, let us give some motivation. 
 In \Cref{section2.3} we will see that there exists a map $\flag \to [B \times B^z \backslash G]$. The lemma below induces 
 a map $[B \times B^z \backslash G] \hookrightarrow [B_0 \times B^{z'}_{0} \backslash \GL(V)]$ sending the open cell 
 to the open cell, so that we can relate the Hasse invariants on both sides. Moreover, we will need $B_0 \subseteq P_{V}$
 to make sense of a Kodaira--Spencer isomorphism for $\mathcal{V}_{\flag}$.
\begin{lemma} \label{special-embedding}
Fix a Borel pair $(T,B)$ of $G_{\overline{\F}_p}$
such that $B \subseteq P$. 
Let $z \in W$ be such that $B^{z}$ over $\Fpbar$ is the Borel generated by $U_{P^{-(p)}}$ and $(B \cap M)^{(p)}$
\cite[Lem 2.3.4]{GK-stratification}.  
Let $V \in \Rep_{\Fpbar}(G)$ be an embedding of the form $V=\oplus V_j$ such that each weight space of $V_j$ is $1$-dimensional. There exists a Borel pair $(T_0,B_0)$ of 
 $\GL(V)_{\overline{\F}_p}$ and $z' \in W_{\GL(V)}$ 
such that under the map $G \hookrightarrow \GL(V)$ we have $B \subseteq B_0 \subseteq P_{V}$, 
$B^{z} \subseteq B^{z'}_0$ 
 and $Bw_0 z^{-1}B^{z} \subseteq B_{0} w_{0,\GL} z^{'-1} B^{z'}_0$. 
 \end{lemma}

\begin{proof}
   Let $\mu$ be a dominant representative of the Hodge 
cocharacter. Choose some basis $\{v_i: i \in I\}$ of $V$ such that each $v_i$ 
lies in a single weight space of $V$. This defines a maximal torus $T_0 \subseteq \GL(V)$ 
such that $T \subseteq T_0$. Define $B_0$ 
by an ordering $\{v_1,v_2,\ldots,v_n\}$ of the basis above (with $v_1$ generating the rank $1$ object) via the following procedure. 
First choose $v_1$ to be a highest weight vector such that $\langle \mu, w(v_1) \rangle$ is 
maximal among all the basis vectors, where $w(v_i)$ denotes the weight of $v_i$.
This is possible since $\mu$ is dominant. 
Let $i \ge 1$ and suppose that $\{v_1,\ldots,v_i\}$ have been chosen with the following properties: 
$\langle \mu, w(v_1) \rangle \ge \langle \mu, w(v_2) \rangle \ge \ldots  \langle \mu, w(v_i) \rangle$,
and no pair with $(1 \le j < k \le i)$
exists such that $w(v_k)>w(v_j)$. 
To choose $v_{i+1}$, let $S \subset I$ consist of the  $j \in I \setminus \{1,2,\ldots,i\}$ such 
that $w(v_j)$ is maximal under the $\ge$ order. Choose $v_{i+1}$ to be $v_j$ for some $j \in S$ such that 
 $\langle \mu, w(v_{j})\rangle$ is maximal among the elements of $S$. 
 Then we cannot have $w(v_{i+1}) > w(v_k)$ for some $1\le k < i+1$ by the way 
 $v_k$ was chosen, and $\langle \mu, w(v_{i})\rangle \ge \langle \mu, w(v_{i+1})\rangle$ by the way $v_i$ was chosen and since $\mu$ is dominant. 
Therefore, we have constructed $\{v_1,v_2,\ldots,v_n\}$ such that if $w(v_i)\ge w(v_j)$ then $j \ge i$, and 
$\langle \mu, w(v_i) \rangle \ge \langle \mu, w(v_{i+1}) \rangle$. The second condition is equivalent to 
$B_0 \subseteq P_{V}$. The first one implies that $B \subseteq B_0$. This can be checked at the level of Lie algebras.  
For $\alpha \in \Phi^{+}$
we have that $x_{\alpha}v_i \in \langle v_1, v_2, \ldots, v_{i}\rangle$ by the way the set $\{v_j\}$ is defined. 
Since the $x_{\alpha}$ for $\alpha \in \Phi^{+}$ generate $\Lie(B)$,
this means that $\Lie(B) \subseteq \Lie(B_0)$ as desired. The choice of $(T,T_0)$ induces a map 
$W \to W_{\GL(V)}$. This is because the weight spaces of each $V_j \subseteq V$ are $1$-dimensional, so that $W$ permutes the 
$1$-dimensional subspaces given by the basis $\{v_i\}$. 
Let $z'=zw_{0}w_{0,\GL} \in W_{\GL(V)}$ under this map.
Proving $B^{z} \subseteq B^{z'}_{0}$ is equivalent to proving the inclusion of opposite Borels $B^{-} \subseteq B^{-}_{0}$, 
which follows in the same way as $B \subseteq B_0$ does. Then since $w_0z^{-1}$ maps to $w_{0,\GL}z^{'-1} \in W_{\GL(V)}$
we obtain that $Bw_{0}z^{-1}B^{z} \subseteq B_0w_{0,\GL}z^{'-1}B^{z'}_0$.
\end{proof}

We note that an irreducible minuscule representation $V$ satisfies that all of its weight spaces are $1$-dimensional, 
so that \Cref{special-embedding} applies to the Hodge embedding.
By essentially the same proof we get the following result integrally. Let $\check{\Z}_p=W(\Fpbar)$.
 \begin{lemma}
Let $\rho: G \hookrightarrow \GL(\Lambda)$ be the Hodge embedding over $\Z_p$.
 There exists tuple $T_0 \subseteq B_0 \subseteq P_0 \subseteq \GL(\Lambda)$
 defined over $\check{\Z}_p$
 of a maximal torus, Borel and $P_0=P_{\rho \circ \mu}$ the parabolic associated to $\rho \circ \mu$
  such that 
 $G \cap (T_0,B_0,P_0,\GL(\Lambda))=(T,B,P,G)_{\check{\Z}_p}$ under the Hodge embedding.
 \begin{proof}
 We construct $T_0, B_0$ over $\Fpbar$ as in \Cref{special-embedding} and then we lift them 
 to $\check{\Z}_p$.
 \end{proof}
 \end{lemma}

 When working over $\check{\Z}_p$ we will also fix 
a tuple
$(T_0,B_0,P_0,\GL(\Lambda))$ with the property that $G \cap (T_0,B_0,P_0,\GL(\Lambda))=(T,B,P,G)_{\check{\Z}_p}$.

\subsubsection{Automorphic vector bundles}
Recall that we have fixed an embedding of 
the tuples $(T,B,P,G) \hookrightarrow (T_0,B_0,P_0,\GL(\Lambda)=\GL_n)$, with all the members of the first 
tuple defined over $\Z_p$
except $P$, and the second tuple defined over $\check{\Z}_p$. We will invoke \Cref{special-embedding} to choose another embedding and tuple
over $\Fpbar$ whenever it is needed.  
Let $0 \to L \to \Lambda \to L^{\vee} \to 0$ be the filtration on $\Lambda$ induced by $P_0$. 
Similarly, $B_0$ defines a full descending flag $F^{\bullet}_{\Lambda} \subseteq \Lambda_{\check{\Z}_p}$ compatible 
with the Hodge filtration.
Let $\{v_i\}$ be a basis of $\Lambda_{\check{\Z}_p}$ such that $v_{i}$ generates $F^{n-i}_{\Lambda}/F^{n-i+1}_{\Lambda}$. 
 Let $h: A \to \Sh$ be the pullback of the universal abelian variety 
over $\Sh_{\GSp}$, and let $\mathcal{H}:=\H^1_{\dR}(A/\Sh)$. It is a vector bundle equipped with the Gauss--Manin 
connection $\nabla$. In general, the Gauss--Manin connection is defined for a triple $X/Y/T$ of schemes with 
$Y/T$ smooth, 
$\nabla: \H^{i}_{\dR}(X/Y) \to \H^{i}_{\dR}(X/Y) \otimes \Omega^1_{Y/T}$ \cite{Katz-Oda}.
 It comes by pullback from the same object $\mathcal{H}_{\GSp}$ over $\Sh_{\GSp}$. From the degeneration 
of the Hodge de Rham spectral sequence it is equipped with a Hodge filtration 
$$
0 \to \omega \to \mathcal{H} \to \omega^{\vee}_{A^{\vee}} \to 0
$$
where $\omega=e^* \Omega^1_{A/\Sh}$, with $e: \Sh \to A$ the zero section. Kisin
\cite[Cor 2.3.9]{Kisin-integral-model} proves that there exist some natural tensors 
$s_{\alpha, \dR} \in \mathcal{H}^{\otimes}$ such that on every $\mathcal{O}$-point of $\Sh$ 
they define the group $G_{\mathcal{O}}$, and on $\mathbb{C}$-points 
they are naturally identified with the tensors induced by $s_{\alpha}$ under the complex de-Rham-Betti comparison 
theorem \cite[\S 2.2]{Kisin-integral-model}. They are moreover horizontal for $\nabla$. 

\begin{defn} \label{torsors}
Let 
$$
G_{\dR}=\text{Isom}_{\Sh}((\mathcal{H},s_{\alpha,\dR}),(\Lambda, s_{\alpha}) \otimes \mathcal{O}_{\Sh}).
$$
It is an \'etale (left) $G$-torsor over $\Sh$, with $G$ acting on $\Lambda$ \cite[Thm 2.4.1]{Zhang-EO}.
 Similarly, let 
$$
P_{\dR}=\text{Isom}_{\Sh}((\omega \subseteq \mathcal{H},s_{\alpha,\dR}),(L \subseteq \Lambda, s_{\alpha}) \otimes \mathcal{O}_{\Sh}).
$$
It is an \'etale $P$-torsor. 
\end{defn}

As an auxiliary space we can define the flag Shimura variety.
\begin{defn} \label{flag-space}
Let $\flag:=P_{\dR}/B$, it comes with a map $\pi: \flag \to \Sh$, whose fibers are isomorphic to $P/B$.
Let $B_{\dR}$ be ${P}_{\dR} \to \flag$ considered as a $B$-torsor over $\flag$. 
\end{defn}

There is a $P$-equivariant embedding $B_{\dR} \hookrightarrow P_{\dR} \times_{\Sh} \flag$ over $\flag$, and in fact 
$B_{\dR} \times^{B} P \cong P_{\dR} \times_{\Sh} \flag$ as $P$-torsors over $\flag$.
 Consider the embedding of pairs $(G,B) \to (\GL(\Lambda),B_0)$
 over $\check{\Z}_p$.
Then the pushout $B_{\dR} \times^{B} B_0 \subseteq G_{\dR} \times^{G} \GL(\Lambda)$ defines a descending filtration $F^{i}_{\mathcal{H}}
\subset \mathcal{H}$, such 
that $\LL_{i}:=\text{gr}^{i} F^{\bullet}_{\mathcal{H}}$ are line bundles. It is easy to see that if 
$F^{\bullet}_{\Lambda} \subseteq \Lambda_{\check{\Z}_p}$ is the flag induced by $B_0$, then 
since $B_{\check{\Z}_p}=G \cap B_0$
$$
B_{\dR, \check{\Z}_p}=\text{Isom}_{\flag}(
	(F^{\bullet}_{\mathcal{H}},s_{\alpha,\dR}),(F^{\bullet}_{\Lambda}, s_{\alpha}) \otimes \mathcal{O}_{\flag}).
$$

\begin{defn}(Automorphic vector bundles) \label{automorphic-vector-bundles}
The torsors above define exact functors 
$$
F_P: \Rep_{\mathcal{O}}(P) \to \text{Coh}(\Sh), \; \; \; F_B: \Rep_{\mathcal{O}}(B) \to \text{Coh}(\flag),
\;\; F_G: \Rep_{\mathcal{O}}(G) \to \text{Coh}(\Sh)
$$
by e.g. $F_P(V)=P_{\dR} \times^{P} V$, which by definition locally consists of pairs 
$(\phi, v)$ with $\phi \in P_{\dR}, v \in V \otimes \mathcal{O}$
up to the equivalence relation $(\phi,v)=(p\phi,pv)$ for $p \in P$. These are compatible with the base change map $\mathcal{O} \to k$
on both sides,
and the functors naturally extend to 
representations which are a filtered colimit of finitely generated representations. 
 Similarly, we can define by pushout $F_{\GL(\Lambda)}$, $F_{P_0}$ and $F_{B_0}$. When the domain 
is clear we will also use the notation $\mathcal{V}:=F_{B}(V)$. We define 
$\LL(\lambda):=F_{B}(\lambda)$ as a line bundle over $\flag$,
 and $\omega(\lambda):=F_P(W(\lambda))$, where $W(\lambda)$ is the dual Weyl module for $M$ as in
\Cref{dual-Weyl-module}. 
\end{defn}

For $V \in \Rep(P)$, when considered as a $B$-representation we have that $\pi^*F_{P}(V)=F_{B}(V)$, since 
$B_{\dR} \times^{B} P \cong P_{\dR} \times_{\Sh} \flag$. 

\begin{lemma}
We have $\pi_* \LL(\lambda)=\omega(w_{0,M}\lambda)$.
\begin{proof}
This follows as in \cite[Prop 3.20]{alexandre} using that $\text{ind}^{P}_{B} \lambda=\text{ind}^{M}_{B \cap M} \lambda=
\text{ind}^{M}_{B^{-} \cap M} w_{0,M} \lambda=W(w_{0,M}\lambda)$, since $(B \cap M)^{w_{0,M}}
:=w_{0,M}(B\cap M)w^{-1}_{0,M}=B^{-} \cap M$.
\end{proof}
\end{lemma}

The next proposition shows that automorphic vector bundles coming from $G$-representations 
also carry a flat connection. The proof follows 
\cite[\S 5.2.1]{Mokrane-Tilouine}, but we record it here since we will use its construction often. 
\begin{prop} \label{GM-definition}
Let $R \in \{\mathcal{O},\mathcal{O}/p^n\}$, and let $V \in \Modfg_{R}(P)$ with a compatible action of $\mathfrak{g}$, i.e.
the derivative of $P$ (well-defined since $V$ is finitely generated) and the action of $\mathfrak{p}=\Lie(P)$ agree.
Then $F_{P}(V)$ is equipped with a natural flat connection over $\Sh_{R}$. 
The same property holds replacing $P$ by $B$. 
\begin{proof}
Let $\phi \in P_{\dR}$ be a local section. By transporting the Gauss--Manin connection along $\phi$ we get a connection 
$\tilde{\nabla}_{\phi}: \Lambda \otimes \mathcal{O}_{\Sh} \to \Lambda \otimes \Omega^1_{\Sh}$. The restriction 
to $\Lambda$ is  $\mathcal{O}$-linear, which can be seen as an element
 $\xi_{\phi} \in \mathfrak{g} \otimes \Omega^1_{\Sh}$, since 
the tensors $s_{\alpha,\dR}$ are parallel for $\nabla$ and $\phi$ maps them to $s_{\alpha}$. 
The image of $\xi_\phi$ under $\mathfrak{g} \otimes_{R} \Omega^1_{\Sh} \xrightarrow{\rho_V \otimes \text{id}}
  \text{End}(V) \otimes_{R} \Omega^1_{\Sh}$, where $\rho_V$ is the $\mathfrak{g}$ action on $V$, can be extended
  to a connection on $V \otimes \mathcal{O}_{\Sh}$ in a unique way, which we denote by $\nabla_{\phi}$.
  Define a connection $\nabla: F_P(V) \to F_P(V) \otimes \Omega^1_{\Sh}=P_{\dR} \times^P (V \otimes \Omega^1_{\Sh})$,
  where 
  in the second equality $P$ acts trivially on $\Omega^1_{\Sh}$, by
  \[
  \nabla(\phi,v)=(\phi, \nabla_\phi(v)).
  \]
  Checking the compatibility with the $P$ action reduces to the identity $p^{-1} \circ
   \nabla_{p \phi} \circ p=\nabla_{\phi}$
  for all $p \in P$. 
  From the definition of $\nabla_\phi$
  it follows from the identity $\xi_{p\phi}=p \xi_{\phi} p^{-1}$, which can be easily checked.
   It is clear that 
  $\nabla$ is a flat connection 
  since each $\nabla_{\phi}$ is.
\end{proof}
\end{prop}

\begin{remark} \label{connection-g-action}
Consider the following setup. Let $P \subseteq G$ be a parabolic inside a reductive group, $E/X$ is a vector bundle 
with a flat connection $\nabla$ over a scheme $X$, $I_{G}$ is a $G$-torsor on $X$ consisting of trivializations of $E$
respecting some tensors of $E$ that are horizontal for $\nabla$, and $I_{P} \subseteq I_{G}$ a sub $P$-torsor.
Then for $V$ a $(U\mathfrak{g},P)$-module as above, the proof above shows that
$I_{P} \times^{P} V$ is equipped with a flat connection. 
\end{remark}

\subsection{Examples}
We list a few examples of Hodge type Shimura varieties which can serve to illustrate our general setup. We will be use some of them in \Cref{section6}.
\begin{example}(PEL Shimura varieties)
The Shimura varieties that we will use for applications are of PEL type, which 
have the advantage of having an integral model which has a moduli interpretation. The data to define 
a PEL Shimura variety (over a number field) is a tuple $(D,*,\Lambda,\langle ,\rangle,h)$ where 
$(D,*)$ is a $\Q$-algebra with a positive involution $*$, $\Lambda$ is a finitely generated $D$-module 
equipped with a non-degenerate $\Q$-valued alternating form $\langle, \rangle$ such that 
$\langle bv,w \rangle=\langle v, b^*w\rangle$ for all $b \in D$, $v,w \in \Lambda$. This induces an involution 
on the algebra $C=\text{End}_{D}(\Lambda)$. Finally, $h: \mathbb{C} \to C_{\mathbb{R}}$ is a $*$-homomorphism 
such that the symmetric bilinear form $\langle \cdot, h(i) \cdot \rangle$ on $\Lambda_{\mathbb{R}}$ is positive-definite.
 Then the group $G$ associated to the data is defined by $G(R):=\{ (x,\nu) \in C_{R} \times R^{\times} : xx^*=\nu\}$. 
 We say that $\nu: G \to \mathbb{G}_m$ is the similitude character. One can then define a 
 Shimura datum $(G,X)$
 and a moduli problem of PEL type defined over the reflex field. Given a rational PEL data 
 as above a choice of integral PEL data consists
 of a tuple $(\mathcal{O}_{D},*,\Lambda_{\Z},\langle,\rangle)$ where $\mathcal{O}_{D}$ is an order on $D$ stable 
 under $*$ and $\Lambda_{\Z}$ is a lattice in $\Lambda$ stable by $\mathcal{O}_{D}$, and such that $\langle, \rangle$
 takes integer values on $\Lambda_{\Z}$. Then an integral PEL datum defines an integral
 PEL moduli problem, and hence an integral model for the Shimura variety. 
 We will use this to specify the integral model of a PEL
 Shimura variety, where we will assume that $p$ is a good prime for the integral PEL data, as in \cite[Def 2.1.6]{boxer-thesis}.
We refer to \cite{Kottwitz-PEL} for more details.
\end{example}

\begin{example}(Compact unitary Shimura varieties) \label{unitary-compact}
For applications it can be convenient to consider compact unitary varieties. Let $E_0/\Q$ be quadratic imaginary and $n,a,b$ as before, 
and let $B/E_0$ a division algebra of dimension $n^2$ over $E_0$ such 
that $E_0$ is its center, admits a positive involution $*$ of the second kind 
($z^*=z^{c}$ for all $z \in E_0$), and $B$ splits over the primes above $p$. Let $\mathfrak{p}$ be such a prime,
and $E:=E_{0\mathfrak{p}}$. 
Then $B_{\mathfrak{p}}\cong M_{n}(E)$, so we let $\mathcal{O}_{B}/\mathcal{O}_{E_0}$ the maximal order stable under 
$*$ such that
$\mathcal{O}_{B,\mathfrak{p}}=M_{n}(\mathcal{O}_E)$.
Let $\epsilon \in \mathcal{O}_{B}$ be the idempotent corresponding to the matrix with a $1$ in the entry $(1,1)$
and zeros elsewhere. Finally, let $\Lambda=B/E_0$ as a module over itself together with a non-degenerate 
$*$-hermitian alternating form $\langle, \rangle$ such that $\Lambda_{\mathbb{R}}$ has signature $(a,b)$.
We choose as a lattice $\Lambda_{\Z}=\mathcal{O}_{B}$. 
Let $\Sh$ be the PEL Shimura variety associated to the integral data 
$(\mathcal{O}_{B},*,\Lambda_{\Z},\langle, \rangle)$.
The integral model is defined over $\mathcal{O}_{E}$ if $a \neq b$ and $\Z_p$ otherwise. 
It is
proper smooth.  
Applying the idempotent $\epsilon$ to the $\sigma$-isotypic component of 
the Hodge filtration of $A/\Sh$ (in the case that $p$ is split in $E$ one can also consider $\epsilon$ 
for the prime corresponding to $\sigma$) we get 
$$
0 \to \omega_{\epsilon} \to \mathcal{H}_{\epsilon} \to \omega^{\vee}_{A^\vee,\epsilon} \to 0,
$$
where the first term has rank $a$ and the last $b$. We can identify $\omega^{\vee}_{A^\vee,\epsilon}$ 
with $\omega^{\vee}_{\epsilon}$. This defines $P_{\dR}$, we will use 
the notation $\omega(\lambda)=\omega_{\epsilon}(\lambda_a) \otimes \omega^{\vee}_{\epsilon}(\lambda_b)$.
We have THAT $G_{\Q}=GU(a,b)$ with respect to $E/\Q$. If $p$ splits in $E$ we have $G_{\Q_p}=\GL_n \times \mathbb{G}_m$ and if $p$ is 
 inert $G_{\Q_p}=GU(n)$ with respect to $\Q_{p^2}/\Q_p$. 
 In \Cref{Hodge-embeddings} the relevant Hodge embedding is
   $\GL_n \times \mathbb{G}_m \hookrightarrow 
 \GL_{2n}$, sending $(A,v)$ to the block diagonal matrix given by $(A,vA^{-1,t})$. 
For more details see \cite[\S 2.1]{Mantovan-unitary}. When the signature is $(n-1,1)$ these are also referred to as 
Harris--Taylor Shimura varieties, \cite{Harris-Taylor}.
\end{example}

\begin{example}(Spin Shimura varieties) \label{Spin-Shimura-varieties}
To get examples of Shimura varieties where $G$ is of type $B$ or $D$ we introduce the following example. 
It is also an example of a Hodge type Shimura variety which is not of PEL type. 
For more details see \cite{Spin}.
Let $n \ge 1$, and $p>2$. Let $V$ be a finite dimensional $\Q$-vector space of dimension $n+2$
with a bilinear form $B$ whose quadratic form $Q$ is non-degenerate of signature $(n,2)$ over $\mathbb{R}$.
Let $L \subseteq V \otimes \Q_p$ be a $\Z_p$ lattice with such that $B$ is a perfect form on it, i.e. 
it induces an isomorphism $L \cong L^{\vee}$. 
Let $C$ be its Clifford algebra, defined to be the universal algebra over $\Z_p$ with a map 
$L \to C$ that satisfies 
$v^2=Q(v)$ for $v \in C$. By construction it has a $\Z/2$-grading $C=C^{+}\oplus C^{-}$.
Then $G=\text{GSpin}(L)=\{ x \in (C^{+})^{\times}: xLx^{-1}=L \}$.  The action by conjugation of $C$ on $L$ 
defines the standard representation $G \hookrightarrow \GL(L)$.
 Under these conditions 
one can define a symplectic form $\psi$ on 
$C$ such that it admits a Shimura datum of Hodge type $(G,X) \hookrightarrow (\GSp(C,\psi),S^{\pm})$.
The reflex field is $\Q$, $\Shbar$ has dimension $n$ and $G_{\Z_p}=\text{GSpin}_{n+2}$, 
which for $n$ even it is of type $D$ and for $n$ odd it is of type $B$. In the language of \Cref{Hodge-embeddings}
 the special node corresponding to $\mu$ is $\alpha_1$, which means that the induced Hodge filtration on the standard 
 representation is a $2$-step filtration of ranks $1,n+1,n+2$. The Hodge embedding is the spin representation, 
 which for type $B$ it is irreducible, and for type $D$ it is the sum of the two half-spin representations. 
\end{example}

\subsection{Hecke operators away from $p$} \label{Hecke-section}
Fix a level $K=K^pK_p \subset G(\mathbb{A}^{\infty}_{\Q})$, and let $\mathbb{A}^{\infty,p}$
be the finite adeles away from $p$. Consider $
\mathcal{H}^{p}=C^{\infty}_{c}(G(\mathbb{A}^{\infty,p})//K^p,\Z_p)
$, the abstract Hecke algebra away from $p$
defined as a convolution algebra of locally constant, compactly supported, 
bi-$K^p$ invariant functions on $G(\mathbb{A}^{\infty,p})$ with coefficients in $\Z_p$.
We define the action of $\mathcal{H}^{p}$ on coherent cohomology via Hecke correspondences.
Let $g \in G(\mathbb{A}^{\infty,p})$, and $K_g=K \cap gKg^{-1}$. There are two finite \'etale maps  
$p_1,p_2: \text{Sh}_{K_g} \to \text{Sh}_{K}$ over $\mathcal{O}$ defined as follows. On the generic fiber 
$p_1$ is defined by the inclusion $K_{g} \subseteq K$, and $p_2$ is defined 
by composing $\text{Sh}_{K_g} \to \text{Sh}_{gKg^{-1}}$ defined as in $p_1$ with the isomorphism 
$[g]: \text{Sh}_{gKg^{-1}} \to \text{Sh}_{K}$, defined on $\text{Sh}_{gKg^{-1}}(\mathbb{C})=G(\Q)\backslash 
X \times G(\A^{\infty})/gKg^{-1}$ by $(x,h) \mapsto (x,hg)$.
By \cite[Lem 2.1.2]{Kisin-integral-model} we can also find a level $K'$ of $\GSp_{2g}(\A^{\infty})$ hyperspecial at $p$
such that it fits in the diagram 
$$
\begin{tikzcd}
\Sh_{K_{g}} \arrow[d,"p_i"]  \arrow[r,hook] & \Sh_{K'_{g'},\GSp} \arrow[d,"p_i"] \\
\Sh_{K} \arrow[r,hook] & \Sh_{K',\GSp} 
\end{tikzcd}
$$
over the generic fiber, where $g' \in \GSp(\mathbb{A}^{\infty,p})$ is the image of $g$. 
We can define both projections integrally on $\Sh_{\GSp}$ as follows. 
The points of the Siegel Shimura variety are tuples $(A,\lambda,\eta_{K})$ of an abelian variety 
$A$ with prime-to-$p$ polarization $\lambda$ and level structure $\eta_{K}$ as in 
\cite[Def 1.3.8]{Lan-thesis}. Then $p_1$ is defined by 
$(A,\lambda,\eta_{K_g}) \to (A,\lambda,\eta_{K})$ where $\eta_K$ is obtained by taking a $K$-orbit of any 
$\eta$ in $\eta_{K_g}$, and $[g]$ is defined by sending $A$ to the unique 
$A'$ prime to $p$ quasi-isogenous to $A$ via $f: A \to A'$ satisfying the conditions in
\cite[Prop 1.4.3.4]{Lan-thesis}. By functoriality of flat closure  we get (unique) maps 
$p_i: \Sh_{K_g} \to \Sh_{K}$ over $\mathcal{O}$ which can be used to define the square above integrally. 
Moreover, these projections are finite \'etale. 
We can check that this correspondence only depends on the coset $K^pgK^p$ 
(by checking it both in the Siegel case and on the generic fiber), so in particular
one can check that 
it induces 
an action of $\mathcal{H}^{p}$ by correspondences, and it readily extends to $\flag_{K}$.

We claim that for $Q \in \{G,P\}$, we have natural 
isomorphisms of torsors $T_{g}: p^*_1 Q_{\dR} \cong p^*_2 Q_{\dR}$ on $\Sh$ or $\flag$. Over 
$\Sh_{\GSp}$ this is induced by the map $f: A \to A'$ used to define $[g']$. By definition of $Q_{\dR}$
we just need to check that the isomorphism
$\H^1_{\dR}(A/\Sh) \cong \H^1_{\dR}(A'/\Sh)$ respects the de Rham tensors, which can 
be checked on formal completions of points, by the construction of the tensors \cite[1.5.4]{Kisin-integral-model}. 
This induces an action in cohomology: for $V \in \Rep_{\mathcal{O}}(P)$
$$
T_{g}: \H^i(\text{Sh}_K, F_{P}(V)) \xrightarrow{p^*_2} \H^i(\text{Sh}_{K_g}, p^*_2 F_{P}(V))
\xrightarrow{T_{g}} \H^i(\text{Sh}_{K_g}, p^*_1 F_{P}(V)) \xrightarrow{\text{Tr } p_1 } \H^i(\text{Sh}_K, F_{P}(V)),
$$
and similarly over $\flag$. 
We say that a collection of maps $F_{Q}(V)_{K} \to F_{Q}(W)_{K}$
of automorphic vector bundles, ranging across $K=K^pK_p$ neat with $K_p$ hyperspecial 
is Hecke equivariant away from $p$ if the associated maps on cohomology are equivariant for $\mathcal{H}^{p}$.
From the definition of $T_{g}$ it is enough to check that the maps are compatible with base change and 
prime to $p$ quasi-isogenies $f: A \to A'$, in the 
sense that they are compatible under replacing the torsors $Q_{\dR}$ with the ones defined by $\H^1_{\dR}(A'/\text{Sh}_K)$. Similarly, 
there is an action of $\mathcal{H}^{p}$ on de Rham cohomology of $\Sh$ with coefficients in $\mathcal{V}$ for 
$V \in \Rep(G)$. This is because the isomorphisms $T_{g}: p^*_1 G_{\dR} \cong p^*_2 G_{\dR}$ are compatible with the way the Gauss--Manin connection is defined on 
$\mathcal{V}$, since the projections $p_{i}$ are \'etale.

Given a level $K$ let $S$ be a finite set of places containing $p$ 
such that $G(\Z_{\ell})$ is contained in $K$ (for some reductive model over $\Z_{\ell}$) for $\ell \notin S$.
 Let $\mathbb{T}/\Z_p$ be the spherical Hecke algebra 
away from $S$, defined as a restricted tensor product of local Hecke algebras 
$\bigotimes^{'}_{\ell \notin S} C^{\infty}_{c}(G(\Q_l)//G(\Z_{\ell}),\Z_p)$.
It is a commutative subalgebra of $\mathcal{H}^{p}$, so it also acts on coherent cohomology. 
Although it depends on the level $K$ we will suppress it from the notation.

\subsection{Toroidal compactifications}
Many Shimura varieties are not proper, so we introduce the machinery of toroidal compactifications, with 
the aim of extending all our operators to them. Here $\Sh/\mathcal{O}$ is always considered with respect to $K_p$ 
hyperspecial. 

\begin{theorem} 
  Let $\Sigma$ be a complete smooth admissible rational polyhedral cone decomposition 
  with respect to $(G,X,K)$ as in \cite[\S 2.1.23]{compactification}. 
  By \cite{compactification} there exists an associated toroidal compactification
   $\Sh^{\tor, \Sigma}$ over $\mathcal{O}$. It is a smooth 
  proper scheme satisfying the following properties. 
  \begin{enumerate} \label{toroidal}
  \item The boundary $D\coloneqq \Sh^{\tor,\Sigma}-\Sh$ with its reduced structure 
  is a Cartier divisor with simple normal crossings. 
  \item The universal abelian scheme on $\Sh$
  extends to a semi-abelian scheme $\pi^{\tor}: A^{\tor} \to \Sh^{\tor}$.
  The prime to $p$ polarization extends to a prime to $p$ isogeny 
  $\lambda: A^{\tor} \to A^{\tor,\vee}$. 
  Define $\omega^{\tor}=e^* \Omega^1_{A^{\tor}/\Sh^{\tor}}$, it extends $\omega/\Sh$.
  In characteristic $p$ the Frobenius and Verschiebung maps extend to $\omega^{\tor}$.
  \item There is a canonical extension $\mathcal{H}^{\tor}$ of $\mathcal{H}$ to a vector bundle over $
  \Sh^{\tor}$ in such a way 
  the de Rham tensors extend to $\Sh^{\tor}$, and we can define 
  an extension of the $G$ torsor $G_{\dR}$ to $\Sh^{\tor}$ \cite[Prop 4.3.7]{compactification}. 
  Moreover, it fits in the exact sequence
  $$
  0 \to \omega^{\tor} \to \mathcal{H}^{\tor} \to \omega^{\tor, \vee}_{A^{\vee}} \to 0 
  $$
  which can be used to define an extension of the $P$ torsor $P_{\dR}$ \cite[Prop 4.3.9]{compactification}.
  \item We extend $\flag \to \Sh$ to $\pi^{\tor}: \flag^{\tor} \to \Sh^{\tor}$ by defining it to be 
  $P_{\dR}/B$.
  Let $D_{\flag}=\pi^{\tor,-1}(D)$ be its boundary divisor. The functors defining automorphic vector bundles 
extend to $F_{B}^{\can}: \Rep_{R}(B) \to \textnormal{Coh}(\flag^{\tor}_R)$ using $B_{\dR}$ over $\flag^{\tor}$.
We define subcanonical extensions by $F^{\sub}_{B}(V)=F^{\can}_{B}(V)(-D_{\flag})
$, and similarly for $\Rep(P)$. This defines extensions $\LL^{?}(\lambda)$, $\omega^{?}(\lambda)$ for
$? \in \{\can, \sub\}$,
and we still have $\pi^{\tor}_{*} \LL^{?}(\lambda)=\omega^{?}(w_{0,M}\lambda)$. 
From now on we will largely drop the indices for canonical 
extensions. 
\item For some appropriate hyperspecial level the Hodge embedding extends to a closed embedding 
$\Sh^{\tor} \hookrightarrow \Sh^{\tor}_{\GSp}$.
 \item There exists a log connection 
 $\nabla^{\tor}: \mathcal{H}^{\tor} \to \mathcal{H}^{\tor} \otimes \Omega^1_{\Sh^{\tor}/\Z_p}(\log D)$ extending 
 the Gauss--Manin connection, and the de Rham tensors in $(3)$ are parallel under it. This
 induces log connections on both $F^{\can}_G(V)$ and $F^{\sub}_G(V)$ for $V \in \Rep(G)$. 
 \item The coherent cohomology groups $\H^i(\Sh^{\tor,\Sigma},\omega^{?}(\lambda))$ for $? \in \{\can, \sub\}$
 are independent of the cone decomposition $\Sigma$. 
 \item In the PEL setting when the codimension of $\Sh$ in its minimal compactification is at least $2$ we have 
 $\H^0(\Sh^{\tor}_{R},\omega(\lambda))=\H^0(\Sh_{R},\omega(\lambda))$ for $R \in \{\Z_p,\F_p\}$, by 
 \cite{Higher-Koecher}. This is satisfied when $\Sh$ decomposes up to finite morphisms 
 as a product of factors which are 
 either proper or of dimension at least $2$ over $\mathcal{O}$.
 \item The Hecke operators away from $p$
 extend to operators 
 $$
 T_{g}: \H^i(\Sh^{\tor},\omega^{?}(\lambda)) \to \H^i(\Sh^{\tor},\omega^{?}(\lambda))
 $$
 for $? \in \{\sub,\can\}$.
  \end{enumerate}
  \end{theorem}

\subsection{Deformation theory and the Kodaira--Spencer isomorphism} \label{KS-section}
Here we explain a version of the Kodaira--Spencer isomorphism that works for the connection attached 
to any $V \in \Rep(G)$, and 
we upgrade the isomorphism to $\flag$. This is proved using the deformation theory of the integral model,
as studied in \cite{Kisin-integral-model}.
Let $x \in \Sh(\overline{\F}_p)$. Let $\mathbb{X}$ be the pullback of the $p$-divisible group $A[p^\infty]$
to $\mathcal{O}^{\wedge}_{\Sh,x}$, and $\mathbb{X}_0$ the reduction over $\Fpbar=\kappa(x)$. Let $W=W(k(x))$.
Let $\mathbb{D}(\mathbb{X}_0)(W)=(M^1_{0} \subseteq M_0,\phi_{M_0})$ be the filtered Dieudonne module of $\mathbb{X}_0$.
Then Faltings and Kisin \cite[Prop 2.3.5]{Kisin-integral-model}
 prove that $\mathcal{O}^{\wedge}_{\Sh,x} \cong R_{G}$, where the latter
is the ring of functions of the completion at the identity of the opposite unipotent $U^{-}_{P,W}$.
This is a quotient of $R$,
the versal deformation ring of $\mathbb{X}_0$. Then $M:=\mathbb{D}(\mathbb{X})(R_{G})=M_0 \otimes R_{G}$
with its filtration $M^1$ being induced by $M_0$, and the Frobenius is given by the universal element in $U^{-}_{P}(R_{G})$. 
We can also identify $M$ with the pullback $\mathcal{H}_{R_{G}}$.
By virtue of being a crystal it is also equipped with a flat topologically nilpotent connection $\nabla_{R_G}$
with the property that the Frobenius $\phi^*M \to M$  is horizontal for it. 

\begin{lemma} \label{GM-compatibility-Dieudonne}
The connection $\nabla_{R_G}$ is the pullback connection from the Gauss--Manin connection on $\mathcal{H}$. 
\end{lemma}

Let $(G,P) \hookrightarrow (\GL_n,P_0)$ be induced by the Hodge embedding. Let $U^{-}_{\GL}$ be the opposite unipotent 
to $P_0$. We can identify $\mathfrak{gl}_n/\mathfrak{p}_0$ with $\Hom_{W}(L,L^{\vee})$.

\begin{prop}(Kodaira--Spencer on formal completions) \label{KS-local}
There exists a commutative diagram 
$$
\begin{tikzcd}
T_{R_G} \arrow[dotted,rr,"\KS_{R_G}"] \arrow[d, hook] & & F_P(\mathfrak{g}/\mathfrak{p}) \arrow[d,hook] \\
T_R \otimes R_{G} \arrow[rr,"\KS_{\GL,R} \otimes R_G"] & & F_{P_0}(\mathfrak{gl}_n/\mathfrak{p}_0) \otimes R_{G}
\end{tikzcd}
$$
here $\KS_{\GL,R}$ is defined by $D \mapsto \nabla_{R,D} \in \Hom(M^1, M/M^{1})$. The two horizontal maps 
are isomorphisms.
\begin{proof}
Since on both formal completions the Hodge filtration is induced from the one on $M$ we have identifications 
$F_{P_0}(\mathfrak{gl}_n/\mathfrak{p}_0)=\text{Lie}(U^{-}_{\GL}) \otimes R$, and $
F_{P}(\mathfrak{g}/\mathfrak{p})=\text{Lie}(U^{-}_{P}) \otimes R_{G}$. The fact that $\KS_{\GL,R}$ is an isomorphism 
is a classical fact, since $R$ is the versal deformation ring of $\mathbb{X}_{0}$. Moreover, 
\cite[\S 1.5]{Kisin-integral-model}
shows that $\nabla_{R}$ restricted to $M_0$ is given by an element $\omega$ of $\Lie(U^{-}_{\GL}) \otimes \Omega^1_R$.  
Since $\nabla_{R_G}$ is the pullback connection, and the tensors $s_{\alpha,\dR}$ are parallel under it, this shows 
that the restriction of $\nabla_{R_G}$ to $M_0$ is given by an element 
in $\Lie(U^{-}_{P}) \otimes \Omega^1_{R_G}$, which is also the projection of $\omega \otimes 1$ via 
 $\Omega^1_R \otimes R_{G}\to \Omega^1_{R_G}$.
 This shows that $\KS_{R_G}$ exists, and it fits in the diagram.
It is automatically injective, and the cokernels of the vertical inclusions are locally free, since $R \to R_{G}$
is induced from $U^{-}_{P} \hookrightarrow U^{-}_{\GL}$ and both of them are smooth. 
After choosing compatible splittings we see that $\KS_{R_G}$ must be an isomorphism, since $\KS_{\GL,R}$
is given by a block diagonal matrix. 
\end{proof}
\end{prop}

We will also need a more general version of the previous proposition. 
\begin{prop} \label{general-KS-local}
Let $V \in \Rep_{W}(G)$ such that $\Lie(G^{\der}) \to \Lie(\SL(V))$ is an embedding with a free quotient. Then $P_{V}$ from \Cref{filtration-functoriality} defines a descending filtration 
$F^{\bullet}_{V}$ of $V$. We identify $\mathfrak{gl}(V)/\mathfrak{p}_{V}$ with 
$\oplus_{k} \Hom(F^k_{V},V/F^k_{V})$. Define $\KS_{R_G,V}: T_{R_G} \to F_{P_{V}}(\mathfrak{gl}(V)/\mathfrak{p}_{V})$ 
by $D \mapsto \nabla_{D}$. It fits in the diagram 
$$
\begin{tikzcd}
& F_{P}(\mathfrak{g}/\mathfrak{p}) \arrow[d,hook] \\
T_{R_G} \arrow[r,swap,"\KS_{R_G,V}"] \arrow[ru,"\KS_{R_G}"]
 & F_{P_{V}}(\mathfrak{gl}(V)/\mathfrak{p}_{V}),
\end{tikzcd}
$$
so in particular the diagonal is an isomorphism independent of $V$.
\begin{proof}
We follow \cite[Prop 4.16]{Spin}. Note that the vertical embedding only depends on the Lie algebra of $G^{\der}$.
First we prove that the diagonal arrow exists. Since the cokernel of vertical 
arrow is locally free we can prove the factorization after tensoring with $R_{G} \to W$. We can identify 
$T_{R_{G}} \otimes W$ with $\Lie(U^{-}_{P}) \otimes W$, then it is enough to show that 
$\KS_{R_G,V} \otimes W$ corresponds to the natural map
$\Lie(U^{-}_{P}) \to \mathfrak{gl}(V)/\mathfrak{p}_{V}$. 
For $V=\Lambda$ this can be extracted from the proof of \Cref{KS-local}.
For a general $V$ the construction of the Gauss--Manin 
connection on $V \otimes R_{G}$ in \Cref{connection-g-action} goes by applying
$\nabla_{\Lambda} \in \mathfrak{g} \otimes \Omega^1_{R_G}$ to $V$, 
so we can reduce to the case $V=\Lambda$. 
We prove that the diagonal map is independent of $V$, so that it must be $\KS_{R_G}$. For $V_1,V_2 \in \Rep_{W}(G)$
consider $V=V_1\oplus V_2$, and we can see that both $\KS_{R_{G},V_i}$ factor through 
$F_{P_V}(\mathfrak{gl}(V)/\mathfrak{p}_{V})$. 
\end{proof}
\end{prop}

The same proposition as above applies to $V \in \Rep_{\Fpbar}(G)$ and the corresponding special fibers. 
From these local results we get a description of $T_{\Sh}$. 
\begin{prop}(Kodaira--Spencer isomorphism)
  Let $V \in \Rep_{\mathcal{O}}(G)$ such that the differential $\Lie(G^\der) \to \Lie(\SL(V))$ is an embedding 
  with a free quotient. 
	Define $\KS_{\GL(V)}$ to be the linear map $T_{\Sh} \to 
	F_{P_0}(\mathfrak{gl}(V)/\mathfrak{p}_V)$
	given by $D \mapsto (x \in F^k_{\mathcal{V}}\mapsto \nabla_{D}(x) \bmod F^k_{\mathcal{V}})$. It factors through 
	a Hecke equivariant away from $p$ isomorphism 
	$$
	\KS: T_{\Sh} \cong F_{P}(\mathfrak{g}/\mathfrak{p}),
	$$ 
  which is moreover independent of the choice of $V$. We can always take 
  $V$ to be the Hodge embedding. Using the log Gauss--Manin connection over $\Sh^{\tor}$ we can 
  define $\KS_{\GL(V)}: T_{\Sh^\tor}(-\log D) \to F_{P_0}(\mathfrak{gl}(V)/\mathfrak{p}_V)$. It factors through 
  a Hecke equivariant away from $p$ isomorphism 
  $$
	\KS: T_{\Sh^\tor}(-\log D) \cong F_{P}(\mathfrak{g}/\mathfrak{p}).
	$$ 
	\begin{proof}
	It is enough to check it on formal completions of $\overline{\F}_p$-points, where it follows from \Cref{KS-local},
	 \Cref{GM-compatibility-Dieudonne} and \Cref{general-KS-local}. We prove the Hecke equivariance.
  The map $\KS_{\GL(V)}$
   is Hecke equivariant since the Hodge filtration and the Gauss--Manin connection are compatible with base change 
   and prime-to-$p$ isogenies. Since $\KS$ is fully determined by $\KS_{\GL(V)}$ it is also Hecke equivariant. 
   For the statement over $\Sh^{\tor}$, the existence of $\KS$ can be checked on a dense open subset, since 
   the cokernel of $F_{P}(\mathfrak{g}/\mathfrak{p}) \to F_{P_V}(\mathfrak{gl}(V)/\mathfrak{p}_V)$ is a vector bundle, 
   and each connected component of $\Sh$ is integral. The fact that it is an isomorphism holds in the Siegel case by 
   \cite{faltings-chai}, and since we can choose the toroidal compactifications compatible with the Hodge embedding
   $i: \Sh^{\tor} \hookrightarrow \Sh^{\tor}_{\GSp}$
   we get that $\KS$ is injective, hence an isomorphism, since 
   the cokernel of $T_{\Sh^{\tor}}(-\log D) \hookrightarrow i^*T_{\Sh^{\tor}_{\GSp}}(-\log D')$ is locally free. It is Hecke equivariant 
   since the log Gauss--Manin connection is. 
  \end{proof}
	\end{prop}

	We can upgrade this isomorphism to one for $T_{\flag}$. 
  For $V \in \Rep_{\check{\Z}_p}(G)$ choose a Borel $B \subseteq B_0 \subseteq P_{V} \subseteq \GL(V)$. 
  That is, it is given by a full flag $F^{\bullet}_{B,V}$ refining the one induced by $P_V$.
	We identify $\mathfrak{gl}(V)/\mathfrak{b}_0$ with a system of compatible elements $\phi_i
	\in \Hom(F^i_{B,V}, V/F^i_{B,V})$. 
	We first note that the Kodaira--Spencer map is also an isomorphism on $G/B$. Let 
  $F_{G/B}: \Rep(B) \to \text{Coh}(G/B)$ be the functor defined by the canonical $B$-torsor over 
  $G/B$.

	\begin{lemma} \label{KS-flag}
  Let $(B,G)$ be a Borel and a reductive group. 
	The map $\KS: T_{G/B} \to F_{G/B}(\mathfrak{g}/\mathfrak{b})$ defined by $D \mapsto (\nabla_{D} \in 
	\Hom(F^i_{\Lambda} \otimes \mathcal{O}_{G/B}, \Lambda \otimes \mathcal{O}_{G/B}/
  F^i_{\Lambda} \otimes \mathcal{O}_{G/B})$ is an isomorphism. 
	\end{lemma}

	We can prove a local Kodaira--Spencer isomorphism for $T_{\flag}$.
  Let $y \in \flag(\overline{\F}_p)$ mapping to $x \in \Sh(\overline{\F}_p)$. Then 
	$\flag^{\wedge}_y \to \Sh^{\wedge}_x$ is identified with 
	$U^{-\wedge}_{B,\text{id}} \to U^{-\wedge}_{P,\text{id}}$. On functions write it as 
	$R_{G} \to R_{B}$. 
  Let $(B_0,P_{\Lambda},\GL(\Lambda))$ be a choice of Borel contained in $P_0=P_{\Lambda}$.
  Write $R \to R_{B_0}$, 
	where $R_{B_0}=P_{0,\dR,R}/B_0$.
	\begin{prop} \label{KS-flag-local}
  \begin{enumerate}
 \item There is a commutative diagram 
	$$
    \begin{tikzcd}
	T_{R_B} \arrow[r,"\KS_{R_B}"] \arrow[d] & F_{B}(\mathfrak{g}/\mathfrak{b}) \arrow[d,hook] \\
	T_{R_{B_0} \otimes R_{B}} \arrow[r,"\KS_{R_{B_0}} \otimes R_{B}"] & 
	F_{B_0}(\mathfrak{gl}_n/\mathfrak{b}_0) \otimes R_{B}
	\end{tikzcd}
	$$
	where both horizontal maps are isomorphisms, and $\KS_{R_{B_0}}$ is defined by 
	$D \mapsto \nabla_{R_{B_0},D}$
	seen as a compatible family of maps $\Hom_{R_{B_0}}(F^i_{\Lambda} \otimes R_{B_0}, M_{R_{B_0}}/
  F^i_{\Lambda} \otimes R_{B_0})$.
   \item  Moreover, 
	$\KS_{R_{B}}$ maps the exact sequence $0 \to T_{R_{B}/R_{G}} \to T_{R_{B}} \to \pi^* T_{R_{G}} \to 0$
  to $0 \to F_{B}(\mathfrak{p}/\mathfrak{b}) \to F_{B}(\mathfrak{g}/\mathfrak{b})
	\to F_B(\mathfrak{g}/\mathfrak{p}) \to 0$, 
	where the quotient map is $\KS_{R_{G}}$.
  \item More generally, let $V \in \Rep_{W}(G)$ as in \Cref{general-KS-local}. Choose a Borel $B_0 \subseteq \GL(V)$
  such that $B \subseteq B_0 \subseteq P_{V} \subseteq \GL(V)$. Then the natural map 
  $T_{R_{B}} \to F_{B_0}(\mathfrak{gl}(V)/\mathfrak{b}_0) \otimes R_{B}$ factors through 
  $F_{B}(\mathfrak{g}/\mathfrak{b})$ 
  via $\KS_{R_B}$.
\end{enumerate}
	\begin{proof}
	The map $\KS_{R_{B_0}}$ is well-defined by the property of being a connection. After fixing 
	a trivialization of $M_0$, we have that $R_{B_0}=R \times P_0/B_0$, and $R_{B}=R_{G} \times P/B$.
	We claim that $\KS_{R_{B_0}}$ maps the exact sequence 
	$T_{R_{B_0}/R} \to T_{R_{B_0}} \to \pi^* T_{R}$ to the one induced by $\mathfrak{p}_{0}/\mathfrak{b}_0
	\to \mathfrak{gl}_n/\mathfrak{b}_0 \to \mathfrak{gl}_n/\mathfrak{p}_{0}$. This is because for $D 
	\in T_{R_{B_0}/R}$, 
  $\nabla_{D}(\omega) \subseteq \omega$
   by functoriality of $\nabla$ along $R \to R_{B_0}$, so $D$ maps to zero 
   in $F_{P_{0}}(\mathfrak{gl}_0/\mathfrak{p_0})$. The two extremes 
	of the sequences are mapped isomorphically by \Cref{KS-local} 
	and \Cref{KS-flag} with respect to $P/B=M/(M\cap B)$, so that $\KS_{R_{B_0}}$ is an isomorphism. 
	As in the proof of \Cref{KS-local}, the restriction of $\nabla_{R_{B_0}}$ to $M_0$
	is given by an element of $\text{Lie}(U^{-}_{B_0}) \otimes \Omega^1_{R_{B_0}}$.
	By preservation of the tensors, $\nabla_{R_{B}}$
	is then given by an element of $\text{Lie}(U^{-}_{B}) \otimes \Omega^1_{R_B}$,
	 which proves that $\KS_{R_B}$
	is well-defined, and by smoothness it is an isomorphism. Part $(2)$ 
	follows from the statement for $\GL(\Lambda)$. Finally, for $(3)$, the map 
  $T_{R_{B}} \to F_{B_0}(\mathfrak{gl}(V)/\mathfrak{b}_0) \otimes R_{B}$ also satisfies the analogous 
  property to $(2)$, by functoriality of the connection, and that $B_0 \subseteq P_V$. Both graded pieces of 
  the map are isomorphisms, by \Cref{general-KS-local} and \Cref{KS-flag}. Independence of $V$ follows 
  as in  the proof of \Cref{general-KS-local}. 

	\end{proof}
	\end{prop}
	
	Globally we get the following result describing $T_{\flag}$.
	
	\begin{prop} \label{basic-iso}
    Let $V \in \Rep_{\check{\Z}_p}(G)$ as in \Cref{general-KS-local}. Choose a Borel $B_0 \subseteq \GL(V)$
    such that $B \subseteq B_0 \subseteq P_{V} \subseteq \GL(V)$.
		There exists a commutative diagram 
		 $$
        \begin{tikzcd}
			T_{\flag/\mathcal{O}} \arrow[r,"e^1_1"] \arrow[rd,"e"]
			& F_{B}(\mathfrak{g}/\mathfrak{p}) \arrow[d,hook]\\
			& F_{B_0}(\mathfrak{gl}(V)/\mathfrak{b}_0)
		\end{tikzcd},
		 $$
		 where $e$ is defined by sending $D \in T_{\flag}$ to the compatible family induced by $\nabla_{D}$
     on $F^{\bullet}_{\mathcal{V}}$. 
		 Then $e^1_1$ is a Hecke equivariant away from $p$ isomorphism, it is independent of $V$, and it
		fits in the diagram of exact sequences
		$$
		\begin{tikzcd}
		0 \arrow[r] &  F_{B}(\mathfrak{p}/\mathfrak{b}) \arrow[d] \arrow[r] & F_B(\mathfrak{g}/\mathfrak{b}) \arrow[d,"e^1_1"]
		 \arrow[r] & F_{B}(\mathfrak{g}/\mathfrak{p}) \arrow[d,"\pi^*\KS"] \arrow[r] &  0 \\
		0 \arrow[r] &  T_{\flag/\Sh} \arrow[r] & T_{\flag/\check{\Z}_p} \arrow[r] &
		 \pi^*T_{\Sh/\check{\Z}_p} \arrow[r] & 0.
		\end{tikzcd}
		$$
    Moreover, $e^1_1$ extends to a
    Hecke equivariant away from $p$ isomorphism $e^1_1: T_{\flag^\tor}(-\log D) \cong F_{B}(\mathfrak{g}/\mathfrak{b})$.
		\begin{proof}
		We can check all the required properties about $e^1_1$ on formal completions of points,
		where it follows from \Cref{KS-flag-local}.
    The extension to $\flag^{\tor}$ follows immediately from the extension of the Kodaira--Spencer isomorphism.

		\end{proof}
		\end{prop}

	\subsection{The Grothendieck--Messing period map}
We reinterpret the isomorphisms of \Cref{basic-iso} on divided power formal 
completions of points in terms of the Grothendieck--Messing period map. This will be crucially used 
in \Cref{section2}.
 First recall an important property of crystals.  

\begin{lemma}\cite[\href{https://stacks.math.columbia.edu/tag/07J6}{Tag 07J6}]{stacks-project} \label{lemma-crystal}
Let $k/\F_p$ be a perfect field, $R \to k$ a PD thickening, and $W=W(k)$. Let $\mathbb{X}_0$
be a $p$-divisible group over $k$. Then $M=\mathbb{D}(\mathbb{X}_0)(R)$ comes equipped with a crystalline connection 
$\nabla_{M}: M \to M \otimes \Omega^1_{R/W,\delta}$. There is a map of 
PD pairs $(W,k) \to (R,k)$ induced by the natural section of $R/p \to k$. By the crystal property 
it induces an isomorphism
$M \cong \mathbb{D}(\mathbb{X}_0)(W) \otimes R$. Then any element $m \in \mathbb{D}(\mathbb{X}_0)(W)$ 
satisfies $\nabla_{M}(m \otimes 1)=0$.
\end{lemma}

Let $y \in \flag(\overline{\F}_p)$ 
and $x \in \Sh(\overline{\F}_p)$ its projection. Let $\mathbb{X}_0$ be the $p$-divisible over $\kappa(y)$ coming from the universal abelian variety, 
and $M_0=\mathbb{D}(\mathbb{X}_0)(W)$.
 Choose a trivialization 
$\phi$ of $(F^{\bullet}_0 \subseteq M_0,s_{\alpha,\dR,y})$ over $W$ such 
that over $\Fpbar$ it identifies $F^{\bullet}_0$ with $\infty=[B] \in G/B$. Let $\Sh^{\sharp}_{x}$ be the divided power
formal completion
at $x$, and 
let $\flag^{\sharp}_y$ be the divided power 
formal completion of $\flag^{\wedge}_y$ along the ideal defined by $x$. Observe that in the notation of \Cref{KS-section} $\flag^{\wedge}_y$
is the pullback of $\Sh^{\sharp}_x \to \text{Spf}(R^{\sharp})$ along $R^{\sharp}_{B_0}:=P_{0,\dR,R^{\sharp}}/B_0$. Here $R^{\sharp}$ is the PD completion of $R$ along $x$. 
Then $R^{\sharp}_{B_0}$ parametrizes full flags $M^{\bullet} \subseteq M_{R^{\sharp}}$ containing $M^1$.
Let 
$I_{P}=[G \to G/P], I_{B}=[G \to G/B]$  be the tautological $P$ and $B$-torsors over the respective 
flag varieties. They induce functors $F_{G/P}: \Rep(P) \to \text{Coh}(G/P)$ and  
$F_{G/B}: \Rep(B) \to \text{Coh}(G/B)$. After pulling back to the formal completion $G/P^{\wedge}_{\infty}$ 
we make the identification
$I_{P}=[G^{\wedge}_{P} \to G/P^{\wedge}_{\infty}]$.

\begin{prop} \label{GM-maps-P}
Fix a trivialization $\phi$ as above. There exist  Grothendieck--Messing period maps (depending on $\phi$)
$$
\begin{tikzcd}
\Sh^{\sharp}_{x} \arrow[r, "\pi_P"] \arrow[hook,d] & G/P^{\wedge}_{\infty} \arrow[hook,d] & 
\flag^{\sharp}_{y} \arrow[r, "\pi_{B}"] \arrow[hook,d] & G/B^{\wedge}_{\infty} \arrow[hook,d] \\
\textnormal{Spf} R^{\sharp} \arrow[r,"\pi_{P_0}"] & {\GL_n/P_0}^{\wedge}_{\infty} & 
\textnormal{Spf} R^{\sharp}_{B_0} \arrow[r,"\pi_{B_0}"] & {\GL_n/B_0}^{\wedge}_{\infty}
\end{tikzcd}
$$
with
$$
\begin{tikzcd}
	\flag^{\sharp}_{y} \arrow[r, "\pi_{B}"] \arrow[hook,d] & G/B^{\wedge}_{\infty} \arrow[hook,d] \\
	\Sh^{\sharp}_{x} \arrow[r, "\pi_P"]  & G/P^{\wedge}_{\infty}  
	\end{tikzcd}
$$
being cartesian. They satisfy the following properties. 
\begin{enumerate}
\item There are natural isomorphisms $\pi_{P}^* I_P \cong P_{\dR}$ and  
$\pi_{B}^* I_B \cong B_{\dR}$.
\item The connection $\nabla_{\Sh^{\sharp}_{x}}$ on $\mathcal{H}$ 
is the pullback connection along $\pi_{P}$ of the 
trivial connection on $F_{G/P}(\Lambda)=\Lambda \otimes \mathcal{O}_{G/P}$.
The analogous statement holds for $B$.
\end{enumerate}
\begin{proof}
The definition of $\pi_{P_0}$ is as follows. It is enough to define it on points $(R \to k)$
of $\Sh^{\sharp}_{x}$
that are PD thickenings. By the crystal property 
\begin{equation} \label{crystal}
\mathcal{H}_R \cong \mathbb{D}(\mathbb{X}_0)(R) \cong M_0 \otimes R \cong \Lambda \otimes R,
\end{equation}
where the last isomorphism is given by $\phi$. Therefore, 
it induces a trivialization of the torsor $\GL_{\dR}$ over $R$. 
The map $\pi_{P_0}$ sends $R$ to $P_{0,\dR,R} \subseteq \GL_{\dR,R} \cong \GL_{R}$, which defines 
a point of $G/P^{\wedge}_{0,\infty}(R)$ by the choice of $\phi$, and since 
the second isomorphism respects the Hodge filtration over $k$.  
Similarly, the second isomorphism in \eqref{crystal} is compatible with the tensors $s_{\alpha, \dR}$ and $s_{\alpha}$
by construction of $s_{\alpha,\dR}$ \cite[1.5.4]{Kisin-integral-model},
so that it induces a trivialization of $G_{\dR,R}$, and $\pi$ is defined by 
$P_{\dR,R} \subseteq G_{\dR,R} \cong G_{R}$. They clearly fit in the commutative diagram. 
For $(1)$, as a space over $\Sh^{\sharp}_x$ we have 
$$
\pi_{P}^* I_P(R \to k)=\{ y \in \Sh^{\sharp}_x(R) , g \in G^{\wedge}_{P}(R) : \pi_P(y)=[gP \subseteq G] \}
$$
with $P$ acting on $g$ on the right. 
We define a map $\pi_{P}^* I_P \cong P_{\dR}$ by sending $(y,g)$ to the element of $P_{\dR,R}$ corresponding to $g$ 
under $gP \cong P_{\dR,R}$. It is clearly $P$-equivariant. Conversely, an element of $P_{\dR,R}$ defines 
an element of $G(R)$ via $P_{\dR,R} \subseteq G_{R}$. It factors through $G^{\wedge}_{P}(R)$ since $G_{\dR,R} \cong G_{R}$
respects the Hodge filtration over $k$. 
For $(2)$, by definition $\pi^*F_{G/P}(\Lambda)(R)=\mathcal{H}_{R} \cong M_0 \otimes R$,
and the pullback connection is defined by $\nabla(m \otimes 1)=0$ for $m \in M_0$. 
By \Cref{lemma-crystal} this agrees with the Gauss--Manin connection $\nabla_{R}$ on $\mathcal{H}_{R}$. 
Doing this for every $R$
ensures the equality over $\Sh^{\sharp}_x$.
The definition of $\pi_{B_0}$ is the same but using $B_{0,\dR}$ instead of $P_{0,\dR}$, and by the same reasoning 
$\pi_{B}$ is well-defined and fits in the diagram. The cartesianness of the diagram involving $\pi_{P}$
and $\pi_{B}$ follows since $[G/B^{\wedge}_{\infty} \to G/P^{\wedge}_{\infty}]=I_{P}/B$, 
and properties $(1), (2)$ follow in a similar way. 

\end{proof}
\end{prop}

Let $j: \Sh^{\sharp}_x \hookrightarrow \Sh$ and $i: \flag^{\sharp}_y \hookrightarrow \flag$. 
If $(R,I,\gamma_i)/W$ is a ring with divided powers, define $\Omega^1_{(R,I)/W}$
as the quotient of $\Omega^1_{R/W}$ by the relations $d \gamma_i(x)=\gamma_{i-1}(x)dx$. 
Define $T_{(R,I)/W}$ as the dual of $\Omega^1_{(R,I)/W}$.
Then we define $T_{\Sh^{\sharp}_x}$ and $T_{\flag^{\sharp}_y}$
as the direct limit of $T_{(\mathcal{O}^{\sharp}_{\Sh,x}/\mathfrak{m}^n_x, \mathfrak{m}_x)/W}$
and $T_{(\mathcal{O}^{\sharp}_{\flag,y}/\mathfrak{m}^n_y, \mathfrak{m}_x)/W}$. Using smoothness 
and \'etale local coordinates 
we see that $T_{\flag^{\sharp}_y} \cong i^* T_{\flag}$ and $T_{\Sh^{\sharp}_x} \cong j^* T_{\Sh}$.

\begin{prop} \label{GM-KS}
Fix a trivialization $\phi$ of $M_0$ as before. 
The differential of the Grothendieck--Messing map fits in the commutative diagram 
$$
\begin{tikzcd}
T_{\Sh^{\sharp}_{x}} \arrow[r,"d\pi_P"] \arrow[rd,swap,"j^*\KS"] & 
\pi_{P}^* T_{G/P^{\wedge}_{\infty}} \arrow[d,"\cong"] \\
& j^*F_P(\mathfrak{g}/\mathfrak{p}),
\end{tikzcd}
$$
where the vertical isomorphism is induced from the isomorphism $\pi_{P}^* I_P \cong P_{\dR}$
and $F_{G/P}(\mathfrak{g}/\mathfrak{p}) \cong T_{G/P}$. In particular, 
$d\pi_{P}$ is an isomorphism. Similarly, $d\pi_{B}$ is an isomorphism fitting in the commutative diagram  
$$
\begin{tikzcd}
	T_{\flag^{\sharp}_{y}} \arrow[r,"d\pi_B"] \arrow[rd,swap,"i^*e^1_1"] & 
	\pi_{B}^* T_{G/B^{\wedge}_{\infty}} \arrow[d,"\cong"] \\
	& i^*F_B(\mathfrak{g}/\mathfrak{b}),
	\end{tikzcd}
$$
with $e^1_1$ from \Cref{basic-iso}.
\begin{proof}
For the first diagram, since $j^*\KS=\KS_{R^{\sharp}_{G}}$ and $\pi_{P}$ are induced from $\KS_{\GL(\Lambda)}$
and $\pi_{B_0}$ it is enough to prove the analogous proposition for $(B_0,P_{0},\GL(\Lambda))$
and $R^{\sharp}$. 
Let $\{e_i : i=1,\ldots, n\}$ be the basis of $M_0$ induced by $\phi$, with 
$\{e_1, \ldots, e_m\}$ corresponding to the Hodge filtration $M^1_0$. 
Then $M^1 \subset M_{R^{\sharp}}$ is given by a basis 
$\{\tilde{e}_i\}=\{e_1 + \tau_{1,m+1}e_{m+1} + \ldots + \tau_{1,n}e_n, 
\ldots, e_m + \tau_{1,m+1}e_{m+1} + \ldots + \tau_{1,n}e_n\}$, where $\tau_{i,j} \in R^{\sharp}$. 
Similarly, let $x_{ij}$ be the coordinates on $\GL(\Lambda)/P_0$
 parametrizing the one-step filtration. Then $\pi_{P_0}$
sends $x_{ij}$ to $\tau_{ij}$. Identifying 
$F_{P_0}(\mathfrak{gl}_n/\mathfrak{p}_0)$ with $\Hom_{R^{\sharp}}(M^1, M/M^1)$
the vertical isomorphism sends $\partial x_{ij}$ to the map $\xi$ defined
 by $\tilde{e}_i \mapsto \tilde{e}_j$
and every other basis element to $0$. This is because under the isomorphism $\pi^* I_P \cong P_{\dR}$ 
in the proof of \Cref{GM-maps-P}
$(y,\pi^{*}_{0}g)$ with 
$$
g=
\begin{pmatrix}
\text{id}_{m} & 0 \\
(x_{ij}) & \text{id}_{n-m}
\end{pmatrix}
$$
($\pi^{*}_{P_0}g$ denotes the specialization of $g$ along $\pi_{P_0}$, i.e. replacing $x_{ij}$ with $\tau_{ij}$)
is sent to the trivialization of $M_{R^{\sharp}}$ given by the basis 
$\{\tilde{e}_1, \ldots, \tilde{e}_{m}, e_{m+1}, \ldots, e_{n}\}$.
Since $\nabla(e_i)=0$ by \Cref{lemma-crystal} we see 
that $\KS(\partial \tau_{ij})=\xi$. The commutativity follows since $d\pi_{P_0}(\partial \tau_{ij})=\partial x_{ij}$.

For the second diagram, it is also enough to prove it for $B_0$. Using the convention that 
$e_i \in M_0$ generates the graded $F^{n-i}_{\Lambda}/F^{n-i+1}_{\Lambda}$ of $\Lambda$ under $\phi$, let 
 $F^{n-i}_{M}$ be generated by $\{f_1, f_2, \ldots, f_i\}$ 
 where $f_i=e_i+y_{i,i+1}e_{i+1} + \ldots + y_{i,n}e_n$, and $y_{ij} \in R^{\sharp}_{B_0}$.
 Let $z_{ij}$ be the analogous coordinates in $\GL_n/B_0$. Then $\pi_{B_0}$ maps $z_{ij}$ to $y_{ij}$.
 A similar analysis to before shows that
 the vertical isomorphism sends $\partial z_{ij}$ to the tuple of elements in $\Hom(F^k_{M}, M/F^k_{M})$
 given by $f_i \mapsto \overline{f_j}$, and every other basis element to $0$. Then since $\nabla(e_i)=0$,
 $\KS_{R_{B_0}}(\partial y_{ij})$ is given by the compatible system sending $f_i \to \overline{f}_j$
 and the rest of the basis elements to $0$. This proves the commutativity. 

\end{proof}
\end{prop}
We note that the same result holds over the special fiber of $\flag^{\sharp}_{x}$, even though it 
is a highly non-reduced ring. 

\subsection{Vanishing of coherent cohomology}
For applications it will be crucial to have good vanishing results for mod $p$ coherent cohomology. 
Most of these results are not available in full generality, so we will state a general expected theorem as an 
assumption, and then specify the available results. This expected theorem is closely related to Goldring--Koskivirta's 
cone conjecture, which we now state. We refer to \Cref{section2.3} for the definition of $\GZip$.

\begin{conjecture}[\cite{general-cone-conjecture0}]
Let $C_{\Shbar}$ be the cone of weights $\lambda \in X^*(T)$ such that there exists an integer $n \ge 1$ satisfying 
$\H^0(\Shbar^\tor,\omega(n\lambda)) \neq 0$. Let $C_{\text{Zip}}$ be the cone of weights $\lambda \in X^*(T)$ such that there exists an integer $n \ge 1$ satisfying 
$\H^0(\GZip,\underline{W(n\lambda)}) \neq 0$. Here $\underline{W(\lambda)}$ is 
the vector bundle on $\GZip$ such that its pullback to $\Sh$ is $\omega(\lambda)$. 
Then $C_{\Shbar}=C_{\text{Zip}}$. 
\end{conjecture}

The surjectivity of the Zip map $\Shbar^\tor \to \GZip$ implies the inclusion $C_{\text{Zip}} \subseteq C_{\Shbar}$, 
and the conjecture states that this is an equality. Since $C_{\text{Zip}}$ can be studied with simpler group theoretic techniques 
one can get estimates on $C_{\Shbar}$ assuming this conjecture. 

\begin{assumption} \label{assumption-vanishing}
Let $(G,X,K)$ be a Shimura datum of Hodge type and $\Shbar$ its associated Shimura variety.
Then for $\lambda \in X^*(T)$.
\begin{enumerate}
\item There exists some integer $n \ge 0$ depending only on $G_{\Q}$ such that for any $i \ge 1$
$$
\H^i(\Shbar^\tor,\omega(\lambda)^{\sub})=0
$$
for all $\lambda \in X_{1}(T)$ satisfying $\langle \lambda-2w_{0,M}\rho, \alpha^{\vee} \rangle > n$ 
for all $\alpha \in \Delta$.
\item There exists an integer $m \ge 0$ depending only on $G_{\Q}$ such that the following is satisfied. 
For all $\alpha \in \Delta$ and $\lambda \in X^*(T)$ such that 
$\langle \lambda,\beta^{\vee} \rangle \le 3p+m$ for all $\beta \in \Delta \setminus \{\alpha\}$,
then 
 $\H^0(\Shbar^{\tor},\omega(\lambda)) \neq 0$ implies that  
 $\langle \lambda, \alpha^{\vee} \rangle>-m$. 
\end{enumerate}
\end{assumption}

The first point has to do with the ampleness of the line bundles $\LL(\lambda)$ over $\flag^{\tor}$, 
by applying a type of Kodaira vanishing result. In fact, in all cases known the vanishing holds for $\lambda$ in a larger 
region than the one stated. More precisely, it holds for $\lambda+2w_{0,M}\rho$
where $\lambda$ lies in a certain cone (the ample cone) contained in the dominant cone.
The second point is related to the cone conjecture in an obvious way. The constant $3$ appearing in its statement is rather arbitrary, 
but it is one that makes our proof of \Cref{basic-theta-injective} work. 
In fact in \cite[\S 3]{GK-GM} they prove that in many cases the cone conjecture
implies \Cref{assumption-vanishing}(2), but we remark that the cone conjecture is a much finer result. 
\begin{theorem} \label{known-vanishing}
Let $\Shbar$ of Hodge type, then 
\begin{enumerate}
\item Let $p> \textnormal{dim}(\flag)$. 
Part $(1)$ in \Cref{assumption-vanishing} is proved in the Siegel case by \cite{alexandre}
and in the case of compact unitary Shimura varieties over a CM field $F/F^{+}$ such that 
the primes of $F^{+}$ above $p$ split over $F$ in \cite{Deding-unitary}. 
\item Part $(2)$ in \Cref{assumption-vanishing} is proved for the Siegel case in \cite{vanishing-Siegel},
 for Hilbert-Blumenthal 
Shimura varieties  \cite{cone-conjecture-Hilbert}, and
for $G_{\F_p}=\GL_{4,\F_p}$, $G_{\F_p}=\GL_{3,\F_p}$, $G_{\F_p}=GU(3)_{\F_p}$,
and $G_{\F_p}=GU(4)_{\F_p}$ such that $\mu$ has signature $(3,1)$ in \cite{general-cone-conjecture0} and 
\cite{general-cone-conjecture1}.
\item Assume that $\Shbar$ is of PEL type, that $p > \textnormal{dim}(\Shbar)$, and 
that $G^{\der}_{\overline{\Q}}$ is almost-simple.
For generic weights $\lambda \in C_0$ in the lowest 
alcove and $w \in W^{M}$
we have that 
$\H^{\bullet}(\Shbar^{\tor},\omega^{\sub}(w \cdot \lambda)^{\vee})$ vanishes above  
degree $d-l(w)$, and $\H^{\bullet}(\Shbar^{\tor},\omega^{\can}(w \cdot \lambda)^{\vee})$
vanishes below degree $d-l(w)$ \cite[Prop 8.12]{Lan-Suh-non-compact}, where $d$ is the dimension 
of $\Sh$. 
\end{enumerate}
\end{theorem}

In \cite{Deding-unitary} they
 moreover obtain sharp bounds on the ample cone of automorphic line bundles, which makes the genericity conditions
 on the weights 
 explicit, i.e. one can determine what $n$ is in \Cref{assumption-vanishing}(1). Similarly, in the results of 
 \Cref{known-vanishing}(2) one is able to explicitly compute a bound for $m$ in \Cref{assumption-vanishing}(2).
  We also note that the PEL assumption in $(3)$ is not essential, and one could reprove their result in the Hodge case. 
 In \cite{Lan-Suh-non-compact} they also prove a more general result without the restriction on 
 $G^{\der}_{\overline{\Q}}$ being almost-simple, where the weights considered are in a smaller region than the 
 lowest alcove.

\section{Verma modules and differential operators} \label{section2}
Let $(G,X)$ be a Shimura datum of Hodge type, and $\Sh$ its associated integral model.
The goal of this section is to relate 
Verma modules to differential operators on the integral model of the (flag) Shimura variety. 
First we define the kind of differential operators and Verma modules that we will use. 
\subsection{Sheaves of differential operators and Lie algebras} \label{section2.1}

In this section we consider $G$ a split reductive group scheme over a ring $R$, which will be either $\mathcal{O}:=W(k)$ for some finite $k/\F_p$, or $\mathcal{O}/p^n$
for some $n \ge 1$. Let $(T,B)$ be a Borel pair of $G$.
We start by reviewing 
various definitions of differential operators, and the relation to Verma modules in the case of flag varieties. 
\begin{defn}
\begin{itemize}
\item (Universal enveloping algebra)
Let $U \mathfrak{g}$ be the $R$-algebra $\bigoplus \mathfrak{g}^{\otimes n}/\langle x \otimes y -y \otimes x -[x,y] 
\rangle$ with $x,y,$ running along $\mathfrak{g}$. 
\item (Restricted universal enveloping algebra) Let $R=\mathcal{O}/p^n$, and write $q=p^n$.
 Over $R$ the Lie algebra has a $q$-operation $x \mapsto x^{[q]}$
corresponding to considering $x$ as a left invariant derivation of $G$ and composing it $q$ times. For each $m\ge 1$,
the
image of the map 
$$
\Sym^{\bullet}\mathfrak{g}^{(q^m)} \hookrightarrow U\mathfrak{g}_{R}
$$
sending $X$ to $X^{q^m}-X^{[q^m]}$ lands in the center of $U \mathfrak{g}_{R}$ 
\footnote{To see why, apply \cite[Lem 1.3.1]{BMR} to $G$ and take its fiber at the identity.}. Denote it by $Z^{\text{Fr}^m}$. 
Define $U^m \mathfrak{g}=U\mathfrak{g} \otimes_{Z^{\text{Fr}^m}} R$, where the character $Z^{\text{Fr}^m} \to R$
is induced by augmentation map $\Sym^{\bullet}\mathfrak{g}^{(q^m)} \to R$. 
\end{itemize}
The first one is a filtered algebra, for $m \ge 0$ we denote $U^{\le m} \mathfrak{g}$  the pieces of 
degree at most $m$.
\end{defn}

\begin{remark} One can also consider the algebra of distributions $U(G)=\Hom_{R}(\mathcal{O}^{\wedge}_{G,e},R)$. 
There is a map of algebras $U \mathfrak{g} \to U(G)$ induced by the
isomorphism $\mathcal{O}_G \otimes \mathcal{O}_{G}/I^2 \to \mathcal{O}_G/\mathfrak{m}_{e}^2$
sending $x \otimes 1-1\otimes x$ to $x$. For $R=\mathcal{O}$ this is an injection that becomes an isomorphism 
after inverting $p$, but it is not 
injective over $k$. In fact, the surjection $U \mathfrak{g}_{\Fpbar} \to U^0 \mathfrak{g}$ can be identified 
with the image of $U \mathfrak{g} \to U(G)$. 
One can see this explicitly: 
let $\{H_i : i=1,\ldots k\}$ be some generators of the Cartan algebra for the maximal torus 
$\mathfrak{h}$, 
and extend it to a basis
$\{H_i, X_j\}$ of $\mathfrak{g}$. Then $U(G)_{\Z_p}$ is generated by elements of the form 
${ H_i \choose k }=\frac{H_i(H_i-1)\ldots(H_i-k+1)}{k!}$ and $\frac{X^n_i}{n!}$. Since $X^{[p]}_i=0$
and
$H^{[p]}_i=H_i$ we see that $U \mathfrak{g}_{\Fpbar} \to U(G)_{\Fpbar}$ factors through $U^0 \mathfrak{g}$.
\end{remark}

We want to define an integral/mod $p$ analogue of category $\mathcal{O}$
over the complex numbers. Let us first recall its definition for motivation. 

\begin{defn}(Complex category $\mathcal{O}$)
Let $\mathfrak{g}$ be the Lie algebra of a complex reductive group, $\mathfrak{h}$ a Cartan subalgebra, $\Phi^{+}$
a choice of positive roots
and $\mathfrak{u}=\oplus_{\alpha \in \Phi^{+}} \mathfrak{g}_{\alpha}$.
Then the category $\mathcal{O}$ with respect to $(\mathfrak{g},\mathfrak{h},\mathfrak{u})$
consists of $U\mathfrak{g}$-modules $M$ such that 
\begin{enumerate}
\item $M$ is finitely generated as a $U\mathfrak{g}$-module.
\item The action of $\mathfrak{h}$ on $M$ is semisimple.
\item For each $v \in M$, the $U\mathfrak{u}$-submodule generated by $v$ is finite dimensional.
\end{enumerate}
\end{defn}

Fix a choice $T \subseteq B \subseteq Q \subseteq G$ of maximal torus, Borel, and some parabolic
$Q$, with its opposite $Q^{-}$. Let $\Modfg_{R}(Q)$ be the category of algebraic representations of $Q_{R}$ which are 
finitely generated as $R$-modules. 

\begin{defn}(Verma modules and category $\mathcal{O}$ integrally) \label{verma-defn}
A $(U \mathfrak{g},Q)$-module over $R$ is a $R$-module $M$ equipped with an action of $U\mathfrak{g}$
and an algebraic action of $Q$, i.e. $M$ is a filtered union of elements in $\Modfg_{R}(Q)$.
Moreover, $M$ 
 satisfies the following properties.
\begin{enumerate}
\item $M$ is a finitely generated $U\mathfrak{g}$-module. 
\item Since the action of $Q$ is algebraic, the derivative of $Q$ on $M$ is well-defined.
It agrees with the restriction of the
$U\mathfrak{g}$ action to $U\mathfrak{q}$.
\end{enumerate}
We denote this category of $(U \mathfrak{g},Q)$-modules by $\mathcal{O}_{Q,R}$. Note that 
the next two properties are automatically satisfied by $M \in \mathcal{O}_{Q,R}$, which are analogues of 
the ones of complex category $\mathcal{O}$.
\begin{itemize}
\item [3.] The action of $T \subseteq Q$ on $M$ is semisimple.
\item [4.] Let $v \in M$, then the $Q$-submodule generated by $v$ is in 
$\Modfg_{R}(Q)$.
\end{itemize}
Similarly, in characteristic $p$ a $(U^0 \mathfrak{g},Q)$-module is an element of $\Rep_{\Fpbar}(Q)$
with an action of $U\mathfrak{g}$ factoring through $U^{0}\mathfrak{g}$, such that the derivative of the $Q$-action 
factors through $U^0 \mathfrak{q}$. Let $V \in \Modfg_{R}(Q)$.
\begin{itemize}
\item Let $\Ver_{Q}(V)\coloneqq U \mathfrak{g} \otimes_{U \mathfrak{q}} V
\in \mathcal{O}_{Q,R}$, with 
$U \mathfrak{g}$ acting on the left and $Q$ by the adjoint action on the left, and by its action on $V$ on the right. 
\item (Baby Verma modules) For $R=\mathcal{O}/p^n$ and $m \ge 1$ let 
$\Ver^m_{Q,R}(V)\coloneqq U^m \mathfrak{g} \otimes_{U^m \mathfrak{q}} V$ as
a $(U^m\mathfrak{g},Q)_{R}$-module.
\item For $V \in \Modfg_{R}(B)$ we can define the variants 
$\Ver_{Q/B}(V)\coloneqq U\mathfrak{q} \otimes_{U \mathfrak{b}} V$
 and $\Ver^{m}_{Q/B,R}(V)\coloneqq U^m\mathfrak{q} \otimes_{U^m \mathfrak{b}} V$.
\end{itemize}
The first one has a filtration induced by the one on $U \mathfrak{g}$, i.e. 
$\Ver_{Q}(V)^{\le m}$ is generated by simple tensors $x \otimes v$ with $x \in U^{\le m} \mathfrak{g}$. 
These filtered pieces are preserved by the action of $Q$.
\end{defn}

We have $\text{gr}^{\bullet} \Ver_{Q}(V) \cong \Sym^{\bullet} \mathfrak{g}/\mathfrak{q} \otimes V$
as $Q$-representations. In general the graded pieces don't split, 
except in the case 
$$
\Ver^{\le 1}_{Q}(1) \cong  1 \oplus \mathfrak{g}/\mathfrak{q}.
$$
It is clear that $\mathcal{O}_{Q,R}$ is an abelian category, since $U\mathfrak{g}$ is a Noetherian ring, by the PBW theorem. 
Our definition of $\mathcal{O}_{Q,R}$ mirrors the definition of category $\mathcal{O}$ over 
the complex numbers, with two main changes. One is that we work with a general parabolic instead of a Borel. 
The second is that when working integrally one cannot just work with representations of Lie algebras, 
since in characteristic $p$ the Lie algebra does not distinguish between weights
which are $p$-translations of each other.
When $Q=B$ this category also appears in \cite{modular-O} under the name modular category $\mathcal{O}$, 
and some of its basic properties are also discussed in \cite[\S 2]{Quan-O},
 which appeared when the contents of this paper were already written. 
 We define the admissible dual of an element of $\mathcal{O}_{Q,R}$.
 
\begin{defn}(Duality) \label{dual-O}
Let $V \in \mathcal{O}_{Q,R}$. We define its dual 
$V^{\vee}$ as the linear functionals $\phi: V \to R$ with the property that $\phi$ factors through a surjection $V \twoheadrightarrow W$
with $W$ a finitely generated $R$-module. 
We define a $U\mathfrak{g}$ action on $V^{\vee}$ given by $x f(v):=f(-xv)$ for $x \in \mathfrak{g}$ and $v \in V$.
\end{defn}

In complex category $\mathcal{O}$ it is customary to define the $\mathfrak{g}$ action on $V^{\vee}$ by 
$xf(v)=f(-\tau(x)v)$, where $\tau$ is the anti-involution on $\mathfrak{g}$ acting as $\text{-id}$ on $\mathfrak{h}$ and as the identity
on the root spaces. We do not do this, and as a result  
the weight spaces are related by the relation $(V^{\vee})_{\lambda}=V_{-\lambda}$, which will be key for us. 
Importantly, duality behaves differently mod $p^n$ compared to characteristic zero. For that we need to define 
a variant of $\mathcal{O}_{Q,R}$. 

\begin{defn} \label{category-o-2}
Let $R=\mathcal{O}/p^n$.
Let $\mathcal{O}^{+}_{Q,R}$ be the category of $R$-modules $M$ with an action of $U\mathfrak{g}$ and
an algebraic action of $Q$, i.e. $M$ is a filtered union of elements in $\Modfg_{R}(Q)$ satisfying 
property $(2)$ from \Cref{verma-defn}. Moreover, the weights of $M$ must be bounded below 
in the $\le$ order. 
\end{defn}
A Verma module $\Ver_{Q}W$ for $W \in \Modfg_{R}(Q)$ belongs to $\mathcal{O}_{Q,R}$, and we will see that its 
dual belongs to $\mathcal{O}^{+}_{Q,R}$. 

\begin{lemma}
Let $R=\mathcal{O}/p^n$ for some $n \ge 1$, and $V \in \mathcal{O}_{Q,R}$. Then the $U\mathfrak{g}$
action on $V^{\vee}$ integrates to an algebraic action of $Q$, in a way that $V^{\vee} \in \mathcal{O}^{+}_{Q,R}$.
\begin{proof}
Let $V=\langle v_1, v_2, \ldots, v_n \rangle_{U\mathfrak{g}}$, and let $W=\langle v_i \rangle_{Q}$. 
It is an element of $\Modfg_{R}(Q)$ since $V$ is a filtered union of elements of $\Modfg_{R}(Q)$. 
There is a natural map 
$\Ver_Q(W) \to V$. Its image is 
$U \mathfrak{q}^{-}W=\langle v_i \rangle_{U\mathfrak{q}^{-}Q} \supseteq \langle v_i \rangle_{U\mathfrak{q}^{-}U\mathfrak{q}}=V$,
so it is surjective. 
Thus, induces an injection $V^{\vee} \hookrightarrow \Ver_Q(W)^{\vee}$ which 
is $U\mathfrak{g}$-equivariant. Therefore, to prove that $V^{\vee} \in \mathcal{O}^{+}_{Q,R}$ we may assume that $V=\Ver_{Q}(W)$.
For each $m\ge 1$ we have surjections $\Ver_{Q}(W) \twoheadrightarrow \Ver^{m}_{Q,R}(W)$ as $Q$-modules. From the PBW theorem 
we see that only the zero element in $\Ver_{Q}(W)$ vanishes in all of the $\Ver^{m}_{Q,R}(W)$,
since $U^{m}\mathfrak{g}$ is generated 
as a vector space by $\{\prod x^{n_{\gamma}}_{\gamma} : 0 \le n_{\gamma}<q^m\}$. For any given 
$\phi \in \Ver_{Q}(W)^{\vee}$ one can then find some $m$ such that  $\phi$ factors through $\Ver^{m}_{Q,R}(W)$, since $\phi$ factors thorough a finitely generated quotient of 
$\Ver_{Q}(W)$.
Therefore, 
$\Ver_{Q}(W)^{\vee}=\cup_{m} \Ver^{m}_{Q,R}(W)^{\vee}$ is a union of elements in 
$\Modfg_{R}(Q)$. We can then see that the properties of \Cref{category-o-2} are satisfied. 

\end{proof}
\end{lemma}

Note that the proof implies that if $V \in \mathcal{O}_{Q,R}$ then its weights are bounded above in the $\le$ order. 
Thus, dual Verma modules are not in $\mathcal{O}_{Q,R}$ since their weights are not bounded above. 
Alternatively, one could prove directly that over $R=\mathcal{O}/p^n$ dual Verma modules
are not finitely generated as $U\mathfrak{g}$ modules. For instance, for $Q=B$, $R=\Fpbar$, and $\gamma \in \Phi^{+}$
the elements $x^p_{\gamma} \in Z(U\mathfrak{g})_{\Fpbar}$ act trivially on $\Ver_{B}(1)^{\vee}$, so that
the $U\mathfrak{g}$ submodule generated by any vector is finite dimensional. 

\begin{remark}
Let $V \in \mathcal{O}_{Q,\mathcal{O}/p^n}$. Even though the action of $\Lie(Q^{-})$ is locally finite on $V^{\vee}$ it does not upgrade 
to an algebraic action of $Q^{-}$. For an example take $G=\GL_2$, $R=k$, $Q=B$, $V=\Ver_{B}(1)$,
and let $x_{+},x_{-} \in \mathfrak{g}$
be the positive and negative root. 
Then for $n \ge p$, $x^p_{-}(x^n_{-} \otimes 1)^{\vee}=(x^{n-p}_{-} \otimes 1)^{\vee} \neq 0$, but $x^{p}_{-}$ is zero in $U(B^{-})_{k}$. This 
means that one cannot integrate the $\mathfrak{g}$ action into an action of $B^{-}$. 
If $V \in \mathcal{O}_{Q,\Q_p}$ then $V^{\vee} \in \mathcal{O}_{Q^{-},\Q_p}$ since the action of $\Lie(Q^{-})$ is locally finite, 
but $V^{\vee} \notin \mathcal{O}^{+}_{Q,\Q_p}$. In the example of $G=\GL_2$ and $V=\Ver_{B}(1)$ 
we have that $x_{+}(x^n_{-} \otimes 1)^{\vee}=
\frac{1}{(n+1)(-n)}(x^{n+1}_{-}\otimes 1)^{\vee}$, so that the action of $\Lie(B)$ is not locally finite. 
\end{remark}

Note that we can still define the admissible dual of an element of $\mathcal{O}^{+}_{Q,R}$ word by word 
from \Cref{dual-O}. 

\begin{lemma} \label{lemma-category-O}
Let $R=\mathcal{O}/p^n$. Then
\begin{enumerate}
\item The functor $(-)^{\vee}: \mathcal{O}_{Q,R} \to \mathcal{O}^{+}_{Q,R} $ is exact. 
\item Let $V \in \mathcal{O}_{Q,R}$. Then $(V^{\vee})^{\vee}=V$.
\item Let $V,W \in \mathcal{O}_{Q,R}$. Then $(V\otimes W)^{\vee}=V^{\vee} \otimes W^{\vee}$.
\end{enumerate}
\begin{proof}
For $(1)$, since $V^{\vee}$ is a submodule of $\Hom_{R}(V,R)$ we just have to prove that given
an embedding $V_1 \hookrightarrow V_2$ the map $V^{\vee}_2 \to V^{\vee}_1$ is surjective. Let $\phi_1 \in V^{\vee}_1$
such that it factors through $V_1 \twoheadrightarrow W_1$ with $W_1$ finitely generated.
 We can construct a finitely generated quotient $V_2 \twoheadrightarrow W_2$ such that $W_1 \hookrightarrow W_2$.
  Thus, we have reduced to proving exactness duality on 
$\Modfg(R)$. This follows since each finitely generated $R$ module is a direct sum of free and torsion modules, by the structure theorem for modules over a PID. 
For $(2)$ the $(U\mathfrak{g},Q)$-equivariant injection $V \hookrightarrow (V^{\vee})^{\vee}$ is an isomorphism.
 Write $V \twoheadrightarrow V_n$ for an exhaustive cofiltration
with $V_n \in \Modfg_{R}(Q)$, we claim that $(V^{\vee})^{\vee} \twoheadrightarrow (V^{\vee}_n)^{\vee}$ is also exhaustive. 
Then one uses the observation that finitely generated $R$-modules are reflexive.
 We can prove the claim as follows.  
Let $\phi \in (V^{\vee})^{\vee}$, then it factors through 
some projection $V^{\vee}\twoheadrightarrow W$ with $W$ finitely generated. Let $W$ be generated 
by some 
$\phi_1, \ldots \phi_k \in V^{\vee}$, and for each $1 \le i \le k$ let $\phi_i$ factor through $V \twoheadrightarrow W_{i}$ with
$W_i \in \Modfg_{R}(Q)$. 
Consider $V \twoheadrightarrow \tilde{V}$ with $\tilde{V} \in \Modfg_{R}(Q)$ so that all $\phi_i$ factor through it.
Then $\phi$ is not in the kernel of $(V^{\vee})^{\vee} \to (\tilde{V}^{\vee})^{\vee}$. Observe that $(V^{\vee})^{\vee}$ is countably generated, 
so that we can inductively construct an exhaustive cofiltration 
$V \twoheadrightarrow V_n$ by repeating the process along the generators of $(V^{\vee})^{\vee}$. The above argument shows that 
$(V^{\vee})^{\vee} \twoheadrightarrow (V^{\vee}_n)^{\vee}$ is exhaustive. 
For $(3)$ there is a natural map $e: W^{\vee} \otimes V^{\vee} \to (W \otimes V)^{\vee}$.
By writing $V=\cup_n V_n$, $W=\cup_n W_n$ with $V_n,W_n \in \Modfg_{R}(Q)$, the map $e$ factors through an isomorphism
$W^{\vee}_n \otimes V^{\vee}_n \cong (W_n \otimes V_n)^{\vee}$. This implies that $e$ is injective. Given 
$\phi \in (W \otimes V)^{\vee}$, let it factor through $W \otimes V \twoheadrightarrow U$ for some $U \in \Modfg_{R}(Q)$.
By finding some simple tensors in $W \otimes V$ that generate $U$ 
we can construct projections $W \twoheadrightarrow U_1$, $V \twoheadrightarrow U_2$ and $\phi_1 \in W^{\vee}$ $\phi_2 \in V^{\vee}$ factoring through these 
such that $\phi_1 \otimes \phi_2$ maps to $\phi$.
Thus, $e$ is an isomorphism.
\end{proof}
\end{lemma}

Now we define the sheaves of differential operators that we will use.  

\begin{defn}(Sheaves of differential operators) \label{diff-ops-defn}
Let $Y/S$ be a smooth map of schemes. Suppose first that $Y/S$ is separated. 
Let $P_{Y}$ be the divided power envelope of the diagonal $\Delta: Y \to Y \times_{S} Y$, i.e. it is the divided power 
envelope of $\mathcal{O}_{Y}$ with respect to $I$, the ideal sheaf of $\Delta$. If $p^nS=0$ for some $n \ge 1$ 
then $P_{Y}$ is well-defined without any assumption on $Y/S$ \cite[Rem 3.31]{Berthelot-Ogus}. 
We consider $P_{Y}$ as a $\mathcal{O}_{Y}$-bimodule via the two projections. 
Let $J$ be the kernel of $P_Y \to \mathcal{O}_Y$. Denote by $J^{[n]}$
  the ideal of $P_Y$ generated by elements $\prod \gamma_{n_i}(x_i)$ for $x_i \in I$ and 
  $\sum n_i \ge n$.
\begin{itemize}
  \vspace{-0.5em}
\item (Crystalline differential operators) 
For an integer $m \ge 0$, let $P^{m}_Y=P_{Y}/J^{[m+1]}$ as a $\mathcal{O}_Y$-bimodule via each 
projection. It is well-defined for any $Y/S$. 
The sheaf of crystalline differential operators 
of degree at most $m$ is $D^{m}_Y=\Hom_{\mathcal{O}_Y}(P^{m}_Y, \mathcal{O}_Y)$, 
where the $\Hom$ is taken as left $\mathcal{O}_Y$-modules.
Then the sheaf of crystalline differential operators is $D_Y=\cup_{m \ge 0} D^{m}_Y$.
More explicitly, $D_Y$ is the sheaf of rings generated by $\mathcal{O}_Y$ and $T_{Y/S}$, under 
the relations $AB-BA=[A,B]$, $Af=fA+A(f)$ for all $A,B \in T_{Y/S}$ (seen as derivations),
and $f \in \mathcal{O}_Y$. 
\item(Log crystalline differential operators) Let $D \hookrightarrow Y$ be a relative Cartier divisor with normal crossings.
Define $D^{\log}_{(Y,D)}$ as the sheaf generated by $T_{Y/S}(-\log D)$ and $\mathcal{O}_Y$
together with the relations described above. Define $P^{\log}_{(Y,D)}$ as in \cite[\S 4.2]{Mokrane-Tilouine}, using the 
fiber product in the category of log schemes. 
\end{itemize}
\end{defn}

\begin{remark}
 One can also consider Grothendieck's sheaf of differential operators 
$\tilde{D}_Y$, defined by taking the formal completion of $Y \times Y$ along the diagonal. 
To illustrate the difference between the two, for $Y=\mathbb{A}^1$ we have that $\tilde{D}_{Y}$ is generated by the divided powers 
$\frac{1}{n!}\frac{d^n}{dx}$, while $D_Y$ is generated by the $\frac{d^n}{dx}$. In this way we see that $D_{Y}$ 
does not act faithfully on $\mathcal{O}_{Y}$ in characteristic $p$. 
We won't make use of $\tilde{D}_Y$, except to remark that crystalline differential operators are the ones that 
make our formalism work, as opposed to $\tilde{D}$.
\end{remark}

Then $P_Y$ is a cofiltered $\mathcal{O}_Y$ bi-module, and 
one has $\text{gr}^{\bullet} P_Y=\Sym^{\bullet} \Omega^1_{Y}$. It also carries a natural coalgebra structure,
given by the maps $\delta: P_Y \to P_Y \otimes_{\mathcal{O}_Y} P_Y$
of $\mathcal{O}_Y$-bimodules, defined on local coordinates by 
$$
a \otimes_R b \mapsto a \otimes_R 1 \otimes_{\mathcal{O}_Y} 1 \otimes_R b.
$$
For $x \in 
\mathcal{O}_Y$ let $\xi_x=x \otimes 1-1\otimes x \in I$. If $Y/S/(\Z/p^n)$ is smooth then 
$P_Y$ is locally freely generated 
by divided powers of elements 
like this, using \'etale coordinates \cite[Prop 3.32]{Berthelot-Ogus}. In these coordinates $\delta$ is given by  
$\delta: \prod_i \frac{\xi^n_{x_i}}{n!} \mapsto \prod_i (\frac{\xi^n_{x_i}}{n!} \otimes 1+1\otimes 
\frac{\xi^n_{x_i}}{n!})$, and one sees
 that 
$P^{n+m}_Y$ lands in $P^n_Y \otimes P^m_Y$.

\begin{defn} Let $Y/S$ be a scheme, and 
 $\mathcal{F}_{1,2}$ sheaves of $\mathcal{O}_{Y}$-modules. 
Define the sheaf of (crystalline) differential operators of degree at most $m$ from $\mathcal{F}_1$ to $\mathcal{F}_2$ as
$$
\DiffOp^{m}(\mathcal{F}_1,\mathcal{F}_2)\coloneqq \Hom_{\mathcal{O}_Y}(P^{m}_{Y/S} \otimes_{\mathcal{O}_Y} \mathcal{F}_1,
\mathcal{F}_2)
$$
where in the tensor product $P^{m}_Y$ is a $\mathcal{O}_Y$-module on the right,
and the $\mathcal{O}_Y$-module structure is given by the left $\mathcal{O}_Y$ action on $P^m_Y$.
Then $\DiffOp(\mathcal{F}_1,\mathcal{F}_2):=\cup_m \DiffOp^{m}(\mathcal{F}_1,\mathcal{F}_2)$. 
There is a map 
\begin{equation} \label{diffops}
\DiffOp(\mathcal{F}_1,\mathcal{F}_2)
\rightarrow \Hom_{\mathcal{O}_S}(\mathcal{F}_1,\mathcal{F}_2),
\end{equation}
defined by precomposing with the (not $\mathcal{O}$-linear) map $\mathcal{F}_1 \to P^m_Y \otimes \mathcal{F}_1$ given by 
$v \mapsto 1 \otimes 1 \otimes v$.
Given $\phi \in \DiffOp^{m}(
\mathcal{F}_1,\mathcal{F}_2
), \psi \in \DiffOp^{n}(
  \mathcal{F}_2,\mathcal{F}_3
  )
$
its composition is defined by 
$$
\psi \circ \phi\coloneqq P^{n+m}_Y \otimes_{\mathcal{O}_Y} \mathcal{F}_1 \xrightarrow{\epsilon \otimes \text{id}} P^n \otimes_{\mathcal{O}_Y} P^m
\otimes \mathcal{F}_1
\xrightarrow{\text{id} \otimes \phi} P^n \otimes_{\mathcal{O}_Y} \mathcal{F}_2 \xrightarrow{\psi} \mathcal{F}_3,
$$
This makes \eqref{diffops} compatible with composition on the right-hand side.
When $p^nS=0$ we will occasionally need the notion of a HPD differential operator between $\mathcal{F}_1$ and $\mathcal{F}_2$: this 
is an element of $\Hom_{\mathcal{O}_Y}(P_Y \otimes \mathcal{F}_1,\mathcal{F}_2)$. Thus, every HPD differential operator 
is a differential operator, but not conversely, see \Cref{stratification-defn}(4) ahead. 
We define log differential operators on $(Y,D)$ as 
$\DiffOp^{\log}(\mathcal{F}_1,\mathcal{F}_2)=\cup_{n} \Hom_{\mathcal{O}_Y}(P^{\log,n}_{(Y,D)} \otimes \mathcal{F}_1,\mathcal{F}_2)$,
and similarly for HPD log differential operators.
\end{defn}

\begin{notation}
For a left $\mathcal{O}$-module $\mathcal{F}$ and $\mathcal{G}$ a bimodule
we will always use $\mathcal{G} \otimes \mathcal{F}$ for the tensor product with the right 
$\mathcal{O}$ structure on $\mathcal{G}$, and $\mathcal{F} \otimes \mathcal{G}$ where the tensor product is taken 
with respect to the left $\mathcal{O}$ structure on $\mathcal{G}$.
\end{notation}

Differential operators are tightly related to vector bundles with connections, and more generally to stratifications. 

\begin{defn}(Stratifications) \label{stratification-defn}
For a quasi-coherent sheaf $E$ on a smooth scheme $Y/S$ we define the following. 
\begin{enumerate} 
\item A connection on $E$ 
can be seen as a $\mathcal{O}$-linear isomorphism $P^{1}_Y \otimes E \cong E \otimes P^{1}_Y$ 
such that
the associated projection map $E \to E$ is the identity, and the connection is flat 
if the isomorphism satisfies some cocycle condition described in the point below
\cite[\S 2]{Berthelot-Ogus}. 
\item 
An extension to a collection of compatible isomorphisms of left $\mathcal{O}_Y$-modules
 $(\epsilon_n: P^n_Y \otimes E \cong E \otimes P^n_Y)_{n \ge 1}$ satisfying a cocycle condition with respect to 
 $Y \times_S Y \times_S Y$
is said to be a PD stratification on $E$. We also have the notion 
of a stratification when using the non-PD version $\tilde{P}^n$ as opposed to $P^n$.
The cocycle condition is the commutativity
of the following diagram 
\[
\begin{tikzcd} \label{cocycle-stratification}
P^n_Y \otimes P^m_Y \otimes E \arrow[d,"\text{id}\otimes \epsilon_n"]
\arrow[r,"\delta^{*} \epsilon_{n+m}"] & E \otimes P^m_Y \otimes P^n_Y ,\\
P^m_Y \otimes E \otimes P^n_Y \arrow[ru,swap,"\epsilon_m \otimes \text{id}"] & 
\end{tikzcd}
\]
where $\delta: P^{n+m} \to P^n \otimes P^m$ is the coalgebra structure map \cite[2.10]{Berthelot-Ogus}.
Here by $\delta^{*} \epsilon_{n+m}$ we mean the extension of scalars of 
$\epsilon_{n+m}: P^{n+m} \otimes E \cong E \otimes P^{n+m}$ along the map of algebras 
$\delta$. The left-hand side of $\epsilon_{n+m}$ is seen as a left $P^{n+m}$-module, and the right-hand side as a 
right $P^{n+m}$-module. 
\item Suppose that $p^mS=0$ for some $m \ge 1$. An extension to an isomorphism $P_Y \otimes E \cong E \otimes P_{Y}$ of $P_Y$-modules
satisfying the analogue of the diagram above \eqref{cocycle-stratification}
 is said 
to be a HPD stratification \footnote{HPD stands for hyperstratification.}.
Given $E$ any quasi-coherent sheaf there is a canonical HPD stratification 
on $P_Y \otimes E$  
$$
\epsilon_{P_Y \otimes E}: P_Y \otimes P_Y \otimes E \cong P_Y \otimes E \otimes P_Y
$$
given by $(a\otimes b) \otimes (c\otimes d) \otimes m \mapsto (ac \otimes b) \otimes m \otimes (1\otimes d)$.
\item The data of a flat connection on $E$ is equivalent to the data of a PD stratification. 
The data of a flat quasi-nilpotent connection is equivalent to the data of a HPD stratification. A connection $(\nabla,E)$
over $Y/S$ with $p^nS=0$
is quasi-nilpotent if for any local section $s \in E$ there exist \'etale local coordinates $x_1, x_2, \ldots, x_n$ and a
positive integer $k$ such that $\nabla(\partial/\partial x_i)^k s=0$ for all $i=1,\ldots,n$
\cite[Thm 4.12]{Berthelot-Ogus}. In particular, the condition of being quasi-nilpotent is independent of the choice of coordinates. 
 All these results can be extended to the log setting of $(Y,D)$
as in \cite[\S 4.2]{Mokrane-Tilouine}. 
\end{enumerate}
\end{defn}

Crucially the Gauss--Manin connection on $\Shbar$ is quasi-nilpotent. 
\begin{lemma} \cite{Katz-p-curvature} \label{GM-quasi-nilpotent}
Let $X/Y/S$ schemes with $Y/S$ smooth, $X/Y$ proper, and $S$ is annihilated by a power of $p$.
Then the Gauss--Manin connection on $\H^i_{\dR}(X/Y)$ is
flat and quasi-nilpotent. 
\end{lemma}

In characteristic $p$ crystalline differential operators interact nicely with the relative Frobenius. 
\begin{prop}(Frobenius differentials) \label{Frobenius-differentials}
  Let $Y/S$ be a smooth map of schemes over $\F_p$ and $E$ a $\mathcal{O}_Y$-module. Let $F: Y \to Y^{(p)}$
  be the relative Frobenius. Then for all $m \ge 1$
  there is a natural map 
  \begin{equation} \label{aaad}
  F^{*}F_{*} E \to P^{m}_{Y/S} \otimes E
\end{equation}
  which is injective for sufficiently large $m$. For 
  $\mathcal{V}_{1,2}$ sheaves of $\mathcal{O}_Y$-modules
  define the sheaf of Frobenius differential operators as
  $$
  D^{[p]}(\mathcal{V}_1,\mathcal{V}_2)\coloneqq \textnormal{Hom}_{\mathcal{O}_Y}(F^* F_* \mathcal{V}_1,\mathcal{V}_2).
  $$
  The composition of two maps $f: F^* F_* \mathcal{V}_1 \to \mathcal{V}_2$ and $g: F^* F_* \mathcal{V}_2 \to \mathcal{V}_3$
  is defined as $F^*F_* \mathcal{V}_1 \xrightarrow{\textnormal{id} \circ s \circ \textnormal{id}} F^*F_*F^*F_* \mathcal{V}_1 
  \xrightarrow{F^*F_* f} F^*F_* \mathcal{V}_2 \xrightarrow{g} \mathcal{V}_3$, where $s: 1 \to F_*F^*$ is the unit of the adjunction.
  Then \eqref{aaad} induces a surjection $\DiffOp(\mathcal{V}_1,\mathcal{V}_2) \to D^{[p]}(\mathcal{V}_1,\mathcal{V}_2)$
  compatible with composition.
  \begin{proof}
  Define a map $E \to  P^{m}_{Y/S} \otimes E$ is given by $v \mapsto 1 \otimes 1 \otimes v$, 
  which is Frobenius linear since $(f \otimes 1 -1\otimes f)^p=0$ on $P_{Y/S}$ for all 
  $f \in \mathcal{O}_Y$. By adjunction it defines the map \eqref{aaad}.
  The rest of the statements can easily be checked on \'etale local coordinates. 
  \end{proof}
  \end{prop}

  The above also says that crystalline differential 
  operators are always Frobenius linear in characteristic $p$, and that the image of \eqref{diffops}
  is precisely $D^{[p]}$.

  \begin{example}
  Let $X=\mathbb{A}^1_{\Fpbar}$. Then $D_{X/\Fpbar}=\DiffOp_{X}(\mathcal{O}_{X},\mathcal{O}_{X})$
  is generated by the elements $(\frac{d}{dx})^n$ for $n \ge 0$. On the other hand 
  $D^{[p]}(\mathcal{O}_{X},\mathcal{O}_{X})$ is generated by $(\frac{d}{dx})^n$ for $0 \le n \le p-1$. 
  We see that for $n \ge p$, $(\frac{d}{dx})^n$ acts trivially as an operator in 
  $\Hom_{\Fpbar}(\mathcal{O}_{X},\mathcal{O}_{X})$, so that the image of $D_{X} \to \Hom_{\Fpbar}(\mathcal{O}_{X},\mathcal{O}_{X})$
  is indeed 
  $D^{[p]}$.
  \end{example}

  \subsubsection{Differential operators on flag varieties}
  Let $Q \subset G$ be a reductive group and a parabolic subgroup, defined over $R$, as in the beginning 
  of the section. 
  Differential operators on the flag variety $G/Q$ are intimately related to Verma modules.

  \begin{lemma} \label{vb-on-flag}
  Let $R \in \{\mathcal{O},\mathcal{O}/p^n\}$.
  There is an equivalence of categories between $\Modfg_{R}(Q)$
  and $G$-equivariant coherent $\mathcal{O}_{G/Q}$ modules. The functor in one direction 
  $F_{G/Q}: \Modfg_{R}(Q) \to \textnormal{Coh}_{G}(\mathcal{O}_{G/Q})$ is given by 
  $V \mapsto (G \times V)/Q$ with $Q$ acting as 
  $h \cdot (g,v)=(hg,hv)$. The functor in the other direction is taking the fiber 
  at $\infty=[Q] \in G/Q$. Both are exact and tensor functors, in particular $F_{G/Q}$ commutes with duality.
 
  \end{lemma}

  This equivalence readily extends to intermediate flag varieties 
  like $P/B$, and to the situation where the representations are a countable 
  union of finite free pieces. We will use the notation $\mathcal{V}:=F_{G/Q}(V)$ whenever it is clear we are 
  working on $G/Q$. 

  \begin{prop} \label{verma-flag}
  Let $R=\mathcal{O}$ or $R=\mathcal{O}/p^n$. 
  Let $V \in \Rep_{R}(Q)$. Under the equivalence of \Cref{vb-on-flag} we have the following  
  identifications.
  \begin{enumerate}
  \item There exists a canonical isomorphism
  $$
  F_{G/Q}(\Ver_Q(V))=D_{G/Q} \otimes F_{G/Q}(V)
  $$
  compatible with the filtrations on both sides.
  If $R=\mathcal{O}/p^n$
  $$
  F_{G/Q}(\Ver^{\vee}_{Q}(V))=P_{G/Q} \otimes \mathcal{V}^{\vee}.
  $$
  The second isomorphism extends to $V \in \mathcal{O}_{Q,R}$, where now 
  $\mathcal{V}^{\vee}:=F_{G/Q}(V^\vee)$.
  \item If $V \in \Rep_{k}(G)$, then $F_{G/Q}(\Ver^0_{Q,k}(V))^{\vee}=F^* F_* \mathcal{V}^{\vee}$, where $F: G/Q \to (G/Q)^{(p)}$
  is the Frobenius. The surjection 
  $$
  D_{G/Q} \otimes \mathcal{V}=\DiffOp_{G/Q}(\mathcal{V}^\vee,\mathcal{O}) \twoheadrightarrow 
  D^{[p]}_{G/Q}(\mathcal{V}^\vee,\mathcal{O})
  $$
  from \Cref{Frobenius-differentials} corresponds under the equivalence of \Cref{vb-on-flag}
  to the surjection 
  \begin{equation} \label{Ver to Ver^0}
  \Ver_Q(V) \twoheadrightarrow \Ver^0_{Q,k}(V)
  \end{equation}
  of $Q$-modules induced by $U \mathfrak{g} \to U^0 \mathfrak{g}$.
  \item Let $V,W \in \Rep_{k}(Q)$, and $f: \Ver_{Q}(V^{\vee}) \to \Ver_{Q}(W^{\vee})$ a map of
  $(U\mathfrak{g},Q)$-modules.
  It descends to a map $\overline{f}: \Ver^0_Q(V^{\vee}) \to \Ver^0_Q(W^{\vee})$. By $(1)$
  the dual of $f$ induces a map 
  $\phi(f) \in \DiffOp(\mathcal{W},\mathcal{V})$ and the dual of $\overline{f}$ induces a map 
  $g: F^*F_* \mathcal{W} \to F^* F_* \mathcal{V}$ by part $(2)$. Then $g=F^* F_* \phi(f)$ where 
  $\phi(f): \mathcal{W} \to \mathcal{V}$ is the map of sheaves associated to $\phi(f)$. In particular 
  $$
  \Ker F^*F_* \phi(f)=F_{G/Q}(\textnormal{coker} \overline{f})^{\vee}.
  $$
\end{enumerate}
\begin{proof}
Let $\pi: G \to G/Q$. For the first part of $(1)$ we have to prove that $(D_{G/Q} \otimes \mathcal{V})_{\kappa(\infty)}\cong \Ver_{Q}(V)$
as $(U\mathfrak{g},Q)$ modules. 
There is a $G$-equivariant map (hence also $\mathfrak{g}$-equivariant)
$\pi_{*}(D_{G} \otimes_{R} V) \to D_{G/Q} \otimes \mathcal{V}$ given as a map of bundles 
by $(g,D \otimes v) \mapsto (\overline{g}, d\pi(D) \otimes v)$. 
Taking its fiber at $\infty$ we get a map 
$\pi_{*}(D_{G} \otimes_{R} V)_{\kappa(\infty)} \to (D_{G/Q} \otimes \mathcal{V})_{\kappa(\infty)}$ 
which is $(U\mathfrak{g},Q)$-equivariant. There is also a $(U\mathfrak{g},Q)$-equivariant map 
$U\mathfrak{g} \otimes_{R} V \to \pi_{*}(D_{G} \otimes_{R} V)_{\kappa(\infty)}$ by identifying $U\mathfrak{g}$ as 
$G$-equivariant differential operators in $D_{G}$. Thus, we get a $(U\mathfrak{g},Q)$-equivariant 
map $\phi: U\mathfrak{g} \otimes_{R} V \to (D_{G/Q} \otimes \mathcal{V})_{\kappa(\infty)}$.
We claim that it factors through $\Ver_{Q}(V)$, and then it becomes an isomorphism. Since $\phi$ is compatible 
with the composition maps $U\mathfrak{g} \otimes U\mathfrak{g} \otimes V \to U\mathfrak{g} \otimes V$ and 
$D_{G/Q} \otimes D_{G/Q} \otimes \mathcal{V} \to D_{G/Q} \otimes \mathcal{V}$, it is enough to prove the factorization of 
$\phi$ 
for the degree at most $1$ part. The latter follows from the fact that $T_{G/Q,\kappa(\infty)}\cong \mathfrak{g}/\mathfrak{q}$, 
and keeping track of the $Q$-action. Then $\Ver_{Q}(V) \to (D_{G/Q} \otimes \mathcal{V})_{\kappa(\infty)}$ is an isomorphism 
since on the associated graded module it is given by the symmetric powers of the degree at most $1$ map $\phi^{\le 1}$.
We now prove the second part of $(1)$. We compute the fiber of $P_{G/Q,R}$ at $\infty$ on the open Bruhat cell 
containing $\infty$. Namely, let $M$ be the Levi of $Q$, and $w \in W^{M}$ the longest element, then
$\infty \in w^{-1}BwQ/Q$. We have that $w^{-1}BwQ/Q \cong \mathbb{A}^d$ where the map is given  
 by composition of the $u_{\alpha} : \mathbb{G}_a \to G$ for $\alpha \in \Phi$ not contained in $\Lie(Q)$. 
 Therefore, $P_{G/Q,R,\kappa(\infty)}$ as an $R$-module is the PD polynomial algebra on those
  $\xi_{\alpha}:=x_{\alpha}\otimes 1-1\otimes x_{\alpha}$. 
 We define the map $(P_{G/Q,R} \otimes \mathcal{V})_{\kappa(\infty)} \to \Ver_{Q}(V)^{\vee}$  as follows. 
 By choosing appropriate projections we may assume that $V$ is free, let $\{v_i\}$ be a basis for $V$. The map sends
 $\prod \frac{\xi^{n_{\alpha}}_{\alpha}}{n_{\alpha}!} \otimes v_i$ to the dual of 
 $(\prod x^{n_{\alpha}}_{\alpha} \otimes v_i) \in \Ver_{Q}(V)$
 with respect to the PBW basis induced by $\{x_{\alpha}\}$ and $\{v_i\}$. One can check that this descends to maps 
 $(P^{n}_{G/Q,R} \otimes \mathcal{V})_{\kappa(\infty)} \to \Ver_{Q}(V)^{\le n, \vee}$, which coincide with the dual 
 of the maps defined in the first part. Therefore the map defined is an isomorphism, and it is $(U\mathfrak{g},Q)$-equivariant. 
 
 Part $(2)$ was originally proved in \cite{original-baby-verma} in the case of $Q=B$ and $V$ a character,
we sketch a proof in general. We have that the fiber at $\infty$ of 
 $(F^*F_* \mathcal{V}^\vee)^{\vee}=\underline{\textnormal{Hom}}_{G/Q}(F^*F_*\mathcal{V}^{\vee},\mathcal{O})=
 \underline{\textnormal{Hom}}_{(G/Q)^{(p)}}(F_*\mathcal{V}^{\vee},F_*\mathcal{O})
 $ 
 is $\Hom_{G_1Q/Q}(\mathcal{V}^{\vee},\mathcal{O})$. Here $G_1=G[F] \subseteq G_{\Fpbar}$
  is the Frobenius kernel of $G$. 
 The latter in turn is identified with $\text{Ind}^{G_1Q}_{Q} V^{\vee}$. 
 Moreover, the
 the category of $(U^{0} \mathfrak{g},Q)_{k}$-modules is equivalent to the category of 
$G_1Q$-representations. 
 Under this equivalence, we have that $\text{Ind}^{G_1Q}_{Q} V^{\vee}\cong \Ver^0_{Q,k}(V)$ by \cite[Lem 9.2]{Janzten-book}, 
 since the baby Verma is an example of a coinduced representation.
 Since $F_{G/Q}$ is compatible with duality we obtain the desired isomorphism.  The second part of $(2)$
 can be proved noting that in the \'etale local coordinates of the proof of $(1)$, the fiber 
 $F^*F_* \mathcal{V}^{\vee}$ at $\infty$ is generated as a vector space by $\prod x^{n_{\alpha}}_{\alpha} \otimes v$ for 
 $\alpha \in \Phi$ not contained in $\Lie(Q)$, $0 \le n_{\alpha}<p$ and $v \in V^{\vee}$. Then the isomorphism 
 $F_{G/Q}(\Ver^0_{Q,k}(V))^{\vee}=F^* F_* \mathcal{V}^{\vee}$ is given by the natural map given these coordinates. 
Part $(3)$ is equivalent to 
proving that the following square commutes
 $$
 \begin{tikzcd}
 \Hom_{Q}(V^{\vee},\Ver_{Q}(W^{\vee})) \arrow[r,"\sim"] \arrow[d] & \DiffOp_{G/Q}(\mathcal{W},\mathcal{V})
 \arrow[d]\\
  \Hom_{Q}(V^{\vee},\Ver^{0}_{Q,k}(W^{\vee})) \arrow[r,"\sim"] & 
  D^{[p]}_{G/Q}(\mathcal{W},\mathcal{V}).
 \end{tikzcd}
 $$
 For this, one needs to check that the dual of $\Ver_{Q}(W^{\vee}) \to \Ver^{0}_{Q,k}(W^{\vee})$ corresponds under point $(1)$
 and $(2)$ to the map
 $F^*F_* \mathcal{W} \to P_{G/Q} \otimes \mathcal{W}$ defined in \Cref{Frobenius-differentials}.
 This is part of point $(2)$.
\end{proof}
  \end{prop}
  With the above we can upgrade \Cref{vb-on-flag} to the category of $(U\mathfrak{g},Q)$-modules. 
  We work over $R=\mathcal{O}/p^n$.
  \begin{prop} \label{flag-g-modules}
  Let $R=\mathcal{O}/p^n$.
	\begin{enumerate}
	\item  Let $V \in \mathcal{O}_{Q,R}$. Then $F_{G/Q}(V^{\vee})$ is a $G$-equivariant quasi-coherent sheaf on $(G/Q)_{R}$
	which is moreover equipped with a $G$-equivariant HPD stratification. Conversely, if 
  $\mathcal{W}$ is a $G$-equivariant quasi-coherent sheaf (which is a countable union of coherent sheaves)
  equipped with a $G$-equivariant HPD stratification, 
  then $\mathcal{W}=F_{G/Q}(W)$ for $W$ a $(U\mathfrak{g},Q)$-module satisfying property $2$ in \Cref{verma-defn}.
  Moreover, whenever they are both defined, these two procedures are inverse to each other. 
	\item If $V \in \Rep_{R}(G)$ then $\mathcal{V} \cong V \otimes \mathcal{O}_{G/Q}$
	 $G$-equivariantly and the associated 
	HPD stratification by part $(1)$ is the trivial one. 
	\item (Tensor identity) For $V \in \mathcal{O}_{Q,R}$, there exists an isomorphism of 
	$(U\mathfrak{g},Q)$-modules $\phi_{V}: \Ver_{Q}(V) \cong V \otimes \Ver_{Q}(1)$ such that
  under \Cref{verma-flag}(1) its dual induces 
	the HPD stratification on $\mathcal{V}^{\vee}:=F_{G/Q}(V^{\vee})$ given by part $(1)$.
	 \item For $W$ a filtered colimit of elements of $\Modfg_{R}(Q)$, the isomorphism 
	$F_{G/Q}\Ver_{Q}(W)^{\vee} \cong P_{G/Q} \otimes \mathcal{W}^{\vee}$
	from \Cref{verma-flag}(1) is an isomorphism of $G$-equivariant 
	sheaves with HPD stratifications, where we are using the stratification defined by $(1)$
	on the left-hand side, and the stratification of 
	\Cref{stratification-defn}(3) on the right-hand side. 
	\end{enumerate}
  
   \begin{proof}
   We can identify the $Q$ torsor $G \mapsto G/Q$ with trivializations of 
   $\Lambda \otimes \mathcal{O}_{G/Q}$ respecting the Hodge tensors $s_{\alpha}$ and the flag corresponding to 
   $Q$. Therefore, by \Cref{connection-g-action}
	if $V \in \mathcal{O}_{Q,R}$, then $F_{G/Q}(V^{\vee})$ is 
   equipped with a $G$-equivariant flat connection, using the trivial connection on
   $\Lambda \otimes \mathcal{O}_{G/Q}$. We check that it is quasi-nilpotent, so that it extends to a HPD 
   stratification, which will be automatically $G$-equivariant. 
  Given a local section $\phi \in Q$, together with the trivial connection 
  it induces the element $\xi_{\phi} \in \mathfrak{g} \otimes \Omega^1_{G/Q}$. 
  There exist some \'etale local sections $\{x_i\}$ for $G/Q$ such that each $x_{\phi,i}:=\partial_{x_i} \circ \xi_{\phi} \in \mathfrak{g}$
is nilpotent, since the trivial connection on $\Lambda \otimes \mathcal{O}_{G/Q}$ is
 quasi-nilpotent. 
That is, $x_{\phi,i}$ is in the span of $\{x_{\gamma} : \gamma \in \Phi\}$. If $x_{\gamma} \in \Lie(Q)$ then it acts 
locally nilpotently on $V^{\vee}$ since the action of $Q$ is algebraic. If $x_{\gamma} \in \Lie(B^{-})$, let
$\psi \in V^{\vee}$ factoring through
$\pi: V \twoheadrightarrow V_1 $ with $V_1 \in \Modfg_{R}(Q)$. Then, for $n$ large enough 
$\pi(x^n_{\gamma}V)=0$, since the weights of $V$ are bounded above. Therefore 
$x^n_{\gamma}\psi=0$, so that $x_{\phi,i}$ acts locally nilpotently on $V^{\vee}$. Since 
$\nabla_{\partial_{x_i}}(\phi,\psi)=
(\phi,x_{\phi,i}\psi)$
this implies that 
the connection is quasi-nilpotent. Now let $\mathcal{W}$ as in the second part of $(1)$. We can write it as 
$\mathcal{W}=F_{G/Q}(W)$ for $W$ a $Q$-module. The HPD stratification on $\mathcal{W}$ induces a map 
$\mathcal{W}\to \mathcal{W} \otimes P_{G/Q} \cong P_{G/Q} \otimes \mathcal{W}$.
 Under the identification \Cref{verma-flag}(1) its dual defines a map 
$\Ver_{Q}(W) \to W$, which defines the $U\mathfrak{g}$ action on $W$.

    For $(2)$ the isomorphism is given by $(g,v) \mapsto (\overline{g},gv)$. 
   The fact that the associated connection on $\Lambda \otimes \mathcal{O}_{G/Q}$ is the trivial one
   follows from \Cref{connection-g-action}. 
   For $(3)$ the map $\phi_{V}$ is defined by $x_1 x_2 \ldots x_n \otimes v \mapsto (x_1 \otimes 1 + 1 \otimes x_1)(x_2 \otimes 1 + 1 \otimes x_2)
  \ldots (x_n \otimes 1 + 1 \otimes x_n)v \otimes 1$
  for $x_i \in \mathfrak{g}$, and $1 \otimes v \mapsto  v \otimes 1$. This is an isomorphism by 
  \cite[Prop 1.7]{MR414645},
  and it can easily checked to be $\mathfrak{g}$-equivariant. In fact $\phi_{V}$ is an isomorphism even for $R=\Z$. 
  Since $\Z$ is torsion-free we see that it is also $U(Q)$-equivariant over $\Z$. 
  Then using the PBW theorem we deduce that $\phi_{V}$ is $U(Q)$-equivariant for any $R$, hence 
  it is $Q$-equivariant. To prove that its dual induces the HPD stratification of point $(1)$, we only need to check it at the 
  degree at most $1$ part. There it follows from the way the connection is defined on $F_{G/Q}(V^{\vee})$, and the way its associated 
  degree $1$ stratification is built from the mentioned connection. 
   For $(4)$, by the second part of $(1)$ we just need to check that the $U\mathfrak{g}$-action 
   on the fiber of $P_{G/Q} \otimes \mathcal{W}^{\vee}$ induced by the HPD stratification $\epsilon_{P_{G/Q} \otimes \mathcal{W}^{\vee}}$ 
   from \Cref{stratification-defn}(3) agrees with the $U\mathfrak{g}$-action on $\Ver_{Q}(V)^{\vee}$.
   In the notation of the proof of \Cref{verma-flag}(1), on an open Bruhat cell the map 
   $P_{G/Q} \otimes \mathcal{W}^{\vee} \to P_{G/Q} \otimes P_{G/Q} \otimes \mathcal{W}^{\vee}$ is given by 
   $\prod \frac{\xi^{n_{\alpha}}_{\alpha}}{n_\alpha!} \otimes v \mapsto \prod \frac{\xi^{n_{\alpha}}_{\alpha}}{n_\alpha!} \otimes 1 \otimes  v$. 
   Under the (dual) isomorphism in the proof of \Cref{verma-flag}(1) the dual of this map agrees with the natural map 
   $\Ver_{Q}\Ver_{Q}(W) \to \Ver_{Q}(W)$. This means that the fiber at $\infty$ of $P_{G/Q} \otimes \mathcal{W}^{\vee}$
   is isomorphic to $\Ver_{Q}(W)^{\vee}$ as elements of $\mathcal{O}^{+}_{Q,R}$.

   \end{proof}
  \end{prop}

  \begin{remark} \label{PD-tensor-identity}
  For $V \in \Rep(Q)$ we can define $\Ver^{\text{PD}}_{Q}(V)=U(G) \otimes_{U(Q)} V$. 
  Then for $V \in \Rep(G)$ we also get a tensor identity of the form 
  $\Ver^{\text{PD}}_{Q}(V) \cong V \otimes \Ver^{\text{PD}}_{Q}(1)$ since the trivial connection 
  on $\mathcal{V}$ also extends to a (non-PD) stratification.
  \end{remark}

\subsection{Differential operators on $\flag$ and Verma modules} \label{section2.2}
\begin{lemma}  \label{nilpotent-connection}
Let $R=\mathcal{O}/p^n$
and $V \in \mathcal{O}_{P,R}$. Then $F_{P}(V^{\vee})$ is equipped with a flat quasi-nilpotent connection
on $\Sh_{R}$, 
which agrees with the one of \Cref{GM-definition} whenever $V \in \Rep_{R}(G)$. It extends to a log quasi-nilpotent 
connection on $\Sh^{\tor}_{R}$.
\begin{proof}
To define a flat connection we mimic the construction of \Cref{GM-definition}.  We prove that it is quasi-nilpotent on $\Sh_{R}$, 
the proof on the toroidal compactification is similar. It follows the proof of \Cref{flag-g-modules}(1).
Let $\psi \in V^{\vee}$. 
Given a local section $\phi \in P_{\dR}$, recall $\xi_{\phi} \in \mathfrak{g} \otimes \Omega^1_{\Sh_{R}}$ from the proof of 
\Cref{GM-definition}. There exist some \'etale local sections $\{x_i\}$ such that $x_{\phi,i}:=\partial_{x_i} \circ \xi_{\phi} \in \mathfrak{g}$
is nilpotent, since the Gauss--Manin connection on $\mathcal{H}$ is quasi-nilpotent \Cref{GM-quasi-nilpotent}. 
That is, $x_{\phi,i}$ is in the span of $\{x_{\gamma} : \gamma \in \Phi\}$. If $x_{\gamma} \in \Lie(Q)$ then it acts 
locally nilpotently on $V^{\vee}$ since the action of $Q$ is locally finite. If $x_{\gamma} \in \Lie(B^{-})$, let 
let
$\psi \in V^{\vee}$ factoring through
$\pi: V \twoheadrightarrow V_1 $ with $V_1 \in \Modfg_{R}(Q)$. Then, for $n$ large enough 
$\pi(x^n_{\gamma}V)=0$, since the weights of $V$ are bounded above. Therefore 
$x^n_{\gamma}\psi=0$, so that $x_{\phi,i}$ acts locally nilpotently on $V^{\vee}$. Since 
$\nabla_{\partial_{x_i}}(\phi,\psi)=
(\phi,x_{\phi,i}\psi)$
this implies that 
the connection is quasi-nilpotent. 
\end{proof}
\end{lemma}
The goal of this subsection 
is to prove the following theorem, which will immediately allow us to construct differential operators out of 
maps of Verma modules. 

\begin{theorem} \label{canonical-iso-Pm}
Let $R \in \{ \mathcal{O}, \mathcal{O}/p^n\}$. 
For any $V \in \Modfg_{R}(B)$ and $W \in \Modfg_{R}(P)$ there exist canonical isomorphisms
$$
e_V: F_B\Ver_{B}(V) \cong D_{\flag_{R}} \otimes \mathcal{V}, \;\;\; 
$$
$$
e_{W}: F_P\Ver_{P}(W) \cong D_{\Sh_{R}} \otimes \mathcal{W}, \;\;\; 
$$
$$
e^{P/B}_V: F_B \Ver_{P/B}(V) \cong D_{\flag_R/\Sh_R} \otimes \mathcal{V}
$$
which are compatible with the filtrations on both sides, i.e. they induce isomorphisms 
$e^{\le n}_{V}: F_{B}\Ver^{\le n}_{B}(V)\cong D^{\le n}_{\flag_R} \otimes \mathcal{V}$.
For $R=\mathcal{O}/p^n$, let $V \in \mathcal{O}_{B,R}$ or $V \in \Modfg_{R}(B)$,
and $W \in \mathcal{O}_{P,R}$ or $W \in \Modfg_{R}(P)$.
Then we have dual isomorphisms 
$$
e^{\vee}_{V} : F_{B} \Ver_{B}^{\vee}(V) \cong P_{\flag_{R}} \otimes \mathcal{V}^{\vee} \;\;\; 
e^{\vee}_{W} : F_{P} \Ver_{P}^{\vee}(W) \cong P_{\Sh_{R}} \otimes \mathcal{W}^{\vee}.
$$ 
Moreover, for $V \in \Rep_{\Fpbar}(B)$ there exist isomorphisms 
$$
e^0_V: F_B \Ver^0_{P/B,\Fpbar}(V)^{\vee} \cong F^{*}F_{*} \mathcal{V}^{\vee},
$$
where $F: \flag \to \flag^{(p)}$ is the relative Frobenius with respect to $\flag/\Shbar$. They are 
all Hecke equivariant away from $p$, and they satisfy the following properties. 
\begin{enumerate}
	\item The map $\Ver_{P/B}(V) \to \Ver_{B}(V)$ is identified with the map 
	$D_{\flag/\Sh}\otimes \mathcal{V} \to D_{\flag} \otimes \mathcal{V}$ induced by 
	$T_{\flag/\Sh} \to T_{\flag/\mathcal{O}}$.
	\item $e^{0}_V$ is uniquely determined by the commutative square 
	$$
	\begin{tikzcd}
	F_B(\Ver_{P/B}(V)) \arrow[d,"e^{P/B}_V"] \arrow[r,"F_B(\pi)"] &  F_B(\Ver^0_{P/B}(V)) \arrow[d,"e^0_V"] \\
	D_{\flag/\Shbar} \otimes \mathcal{V} \arrow[r,"\Pi"] & (F^*F_* \mathcal{V}^\vee)^{\vee}
	\end{tikzcd}
	$$
	where $\pi: \Ver_{P/B}(V) \to \Ver^0_{P/B}(V)$ is the map \eqref{Ver to Ver^0}, and 
	$\Pi$ is defined in \Cref{Frobenius-differentials}.
	\item If $\Ver_{Q}(V) \to \Ver_{Q}(W)$ is induced
	by a map $f: V \to W$ of $Q$-representations, then the associated map
	$D \otimes \mathcal{V} \to D \otimes \mathcal{W}$ is 
	$\textnormal{id} \otimes F_Q(f)$. 
	\item The map $\Ver_{Q}(\Ver_{Q}(V)) \to \Ver_{Q}(V)$ given by $x \otimes (y \otimes v) \mapsto xy \otimes v$ is sent 
	to $D \otimes (D \otimes \mathcal{V}) \to D \otimes \mathcal{V}$ induced by composition of differential 
	operators on $\Shbar$ or $\flag$. 
	\item On graded pieces of $\Ver_{B}(V)$ and $D \otimes \mathcal{V}$, $e_V$ induces the map
	$F_{B}(\textnormal{gr}^n \Ver_{B}(V)) \cong F_{B}(\Sym^n \mathfrak{g}/\mathfrak{b} \otimes V) \xrightarrow{\Sym^{\bullet} e^1_1 \otimes \text{id}}
	\Sym^n T_{\flag} \otimes \mathcal{V} \cong 
	\textnormal{gr}^n D_{\flag} \otimes \mathcal{V}$, where $e^1_1$ is defined in \Cref{basic-iso}.
	\item Let $R=\mathcal{O}/p^n$. For $V \in \mathcal{O}_{Q,R}$, consider the isomorphism 
	$\phi_V: \Ver_{Q}(V) \cong V \otimes \Ver_{Q}(1)$ from \Cref{flag-g-modules}(3).
	Then $e^{\vee}_{V}$ fits into the following 
	commutative diagram 
	$$
	\begin{tikzcd}[
	  ar symbol/.style = {draw=none,"#1" description,sloped},
	  isomorphic/.style = {ar symbol={\cong}},
	  equals/.style = {ar symbol={=}},
	  ]
	F_{Q}(\Ver_{Q}(V)^{\vee}) \arrow[r,"F(\phi^{\vee}_V)"] \arrow[d,"e^{\vee}_V"] & 
  F_{Q}(V^{\vee} \otimes \Ver_{Q}(1)^{\vee}) \arrow[r,"\textnormal{id} \otimes e^{\vee}_{1}"] & 
	\mathcal{V}^{\vee} \otimes P \\
	P \otimes \mathcal{V}^{\vee} \arrow[rru,"\nabla"] & & 
	\end{tikzcd}
	$$
	where the diagonal map is the isomorphism corresponding to the HPD stratification associated to 
  the Gauss--Manin connection on 
	$\mathcal{V}^{\vee}$ by \Cref{nilpotent-connection}. The same holds for $V \in \Modfg_{R}(G)$
  and $e^{\le n,\vee}_{V}$ using the PD stratification.
	\item Let $R=\mathcal{O}/p^n$. For $V \in \mathcal{O}_{Q,R}$,
  the dual of $\phi_{\Ver_{Q}(V)}: \Ver_Q \Ver_Q(V) \cong \Ver_{Q}(V) \otimes \Ver_{Q}(1)$ induces a 
  HPD stratification on 
	$P \otimes \mathcal{V}^{\vee}$ under $e^{\vee}_{V}$. It coincides with $\epsilon_{P \otimes \mathcal{V}^{\vee}}$ of 
	\Cref{stratification-defn}(3).
	\item All the isomorphisms extend to (sub)canonical extensions on toroidal compactifications by replacing $D$
	with $D^{\log}$ and $P$ with $P^{\log}$. In $(2)$ one just considers $D_{\flag^\tor/\Shbar^\tor}$. They satisfy all the previous properties. 
	\end{enumerate}
\end{theorem}

The proof of \Cref{canonical-iso-Pm} will take the rest of the subsection. Here
we sketch the strategy of the proof. First we construct local isomorphisms on divided power
formal completions
of points by using the Grothendieck--Messing period map to a flag variety, which allows to transport 
\Cref{verma-flag} to the Shimura variety. Then we construct candidate maps $e_V$, and we show that they agree
with the local ones on divided power formal completions, so that we can check all their properties locally. 
We could only manage 
to construct $e_V$ in an indirect way. First we reduce to the case of  $\Ver^{\le 1}(V)$
by an inductive process relating the maps
$\Ver^{\le 1}(\Ver^{\le n}(V)) \to \Ver^{\le n+1}(V)$ to the algebra structure 
on $D$. If $V$ is a $G$-representation we can naturally construct 
$e^{\le 1}_V$ via property $(6)$. For an arbitrary $V \in \Rep_{R}(Q)$ we prove that it is a subquotient
of the restriction to $Q$ of a $G$-representation, so that we can construct $e_V$ by exactness of $F$ and $\Ver_{Q}(-)$.
A priori the inductive process might not be well-defined, but the compatibility with the local isomorphisms 
ensures that it is. To construct $e^{\vee}_{V}$ we use that the PD stratification given by the 
Gauss--Manin connection extends to a HPD stratification, so that we can lift all the isomorphisms 
of vector bundles $e^{\le n, \vee}_{V}$ to the desired isomorphism. 

\begin{remark} \label{PD-diff-ops-dont-work}
Let $\tilde{D}$ be Grothendieck's sheaf of differential operators, and define the variant 
$\Ver_P^{PD}(V)=U(G)\otimes_{U(P)} V$. 
It is not true that $
F_P(\Ver_P^{PD}(V)) \cong \tilde{D}_{\Sh} \otimes F_P(V)
$ for all $V$ in a functorial way. If we had such isomorphisms extending $e_V$
along the embeddings $\Ver_P(V) \hookrightarrow \Ver_P^{PD}(V)$ and 
$D_{\Sh} \hookrightarrow \tilde{D}_{\Sh}$ for both $V$ the trivial and 
the standard representation $\Lambda$, using the diagram of part 6) together with \Cref{PD-tensor-identity}
would imply that the Gauss--Manin connection on $\mathcal{H}$ over $\Sh_{\mathcal{O}}$ extends to a stratification 
$(\tilde{P}^n_{\Sh} \otimes H \cong H \otimes \tilde{P}^n_{\Sh})_n$.
Take for example $\Shbar$ the Siegel Shimura variety. 
Over a lift of the ordinary locus $W \subset \Sh_{\Z_p}$ there is a lift 
of the Frobenius isogeny $F: A \to A'$, and by functoriality the existing PD stratification 
respects this morphism. 
Since the cokernel of $P^n_{\Sh} \to \tilde{P}^{n}_{\Sh}$ is $p$-torsion, and all the sheaves are locally free,
we see that the potential stratification on $W$ would have to respect $F$ too. Therefore, it would 
also be compatible with Frobenius on $\Shbar^{\text{ord}}$, and hence in all of $\Shbar$ since 
being compatible with Frobenius is 
an open condition. Such an extension of $\nabla$ to a stratification on $\Shbar$ compatible with the Frobenius
is well-known not to exist, by considering the rank of powers of Verschiebung, 
 see \cite[Ex 2.18]{Berthelot-Ogus}.
\end{remark}

As an immediate corollary from \Cref{canonical-iso-Pm} we get a functor from 
maps of Verma modules to differential operators on the (flag) Shimura variety. 

\begin{theorem} \label{creator-diff}
	Let $R \in \{\mathcal{O},k\}$,  $V_{1},V_{2} \in \Rep_{R}(P)$, and $W_{1},W_{2} \in \Rep_{R}(B)$.
  There exist functorial 
	embeddings
	$$
	\Phi_{P}: \Hom_{(U \mathfrak{g},P)_R}(\Ver_{P}(V_1), \Ver_{P}(V_2)) \hookrightarrow \DiffOp_{\Sh_R}(\mathcal{V}^{\vee}_2,\mathcal{V}^{\vee}_1),
	$$
	$$
	\Phi_{B}: \Hom_{(U\mathfrak{g},B)_R}(\Ver_{B}(W_1), \Ver_{B}(W_2)) \hookrightarrow \DiffOp_{\flag_R}(\mathcal{W}^{\vee}_2,\mathcal{W}^{\vee}_1),
	$$
	$$
	\Phi_{P/B}: \Hom_{(U\mathfrak{p},B)_R}(\Ver_{P/B}(W_1), \Ver_{P/B}(W_2))
	 \hookrightarrow \DiffOp_{\flag/\Sh_R}(\mathcal{W}^{\vee}_2,\mathcal{W}^{\vee}_1).
	$$
	If $W_i \in \Rep_{\Fpbar}(B)$ there also exist functorial embeddings
	$$
	\Phi^{\textnormal{Fr}}_{P/B}: \Hom_{(U^0 \mathfrak{p},B)_{\Fpbar}}(\Ver^0_{P/B,\Fpbar}(V_1), \Ver^0_{P/B,\Fpbar}(V_2)) \hookrightarrow
	D^{[p]}_{\flag/\Shbar}(\mathcal{V}^{\vee}_2,\mathcal{V}^{\vee}_1).
	$$
	All the functors extend to functors $\Phi^{\can,\sub}$ on (sub)canonical extensions on
	toroidal compactifications,
	by using 
	log crystalline differential operators. 
	All the differential operators produced in this way are Hecke equivariant away from $p$.
	Moreover, they satisfy the following properties.
	\begin{enumerate}
	\item For $f \in \Hom_{(U \mathfrak{p},B)}(\Ver_{P/B}(V_1), \Ver_{P/B}(V_2))$, let $\overline{f}$
	 be the induced map on baby Verma modules. Then $\Phi_{P/B}(f): \mathcal{V}^{\vee}_2 \to \mathcal{V}^{\vee}_1$
	is Frobenius linear with respect to $\Shbar$, and it satisfies 
	$$
	F^*F_* \Phi_{P/B}(f)=F_{B}(\overline{f}^{\vee})
	$$
	under the canonical isomorphism $e^0_{P/B}$, where $F: \flag \to \flag^{(p)}$ is the relative Frobenius 
	with respect to $\Shbar$.
	\item 
	For $\Phi_{B}, \Phi_P, \Phi_{P/B}$, and $\Phi^{\textnormal{Fr}}_{P/B}$
	composition as maps of $(U\mathfrak{g},Q)$-modules on the left matches with composition 
	of differential operators on the right.
  \end{enumerate}
  
	\begin{proof}
	Note that $\Hom_{(U \mathfrak{g},Q)}(\Ver_{Q}(V_1), \Ver_{Q}(V_2))=\Hom_{Q}(V_1,\Ver_{Q}(V_2))$
	and similarly for the $\Ver^0_B$. Moreover, since the $V_i$ are finitely generated any map must factor through some 
  $\Ver^{\le n}_{Q}(V_2)$. 
   Applying 
	$F_P$ or $F_B$  respectively, and then dualizing gives the desired differential operators as in 
	\Cref{diff-ops-defn}, by \Cref{canonical-iso-Pm}. The maps are embeddings since $F_{P}$ is
  exact and compatible with duality on $\Rep_{R}(P)$, 
  and $\Hom_{Q}(V_1,\Ver_{Q}(V_2))=\cup_{n} \Hom_{Q}(V_1,\Ver^{\le n}_{Q}(V_2))$. 
  The functors extend 
	to toroidal compactifications since \Cref{canonical-iso-Pm}
	does, and they are Hecke equivariant away from $p$. Properties $(1)$ and $(2)$ can be checked on the interior,
	being closed conditions. 
  
	Part $(1)$ follows directly from
	\Cref{canonical-iso-Pm}(2) and the way that $\Pi$ is defined in \Cref{Frobenius-differentials}.
   For part $(2)$, given two $Q$-equivariant maps
	 $f_i: V_i \to \Ver_Q(V_{i+1})$
	its composition 
	is given by
	$f_2 \circ f_1: V_1 \xrightarrow{f_1} \Ver_Q(V_2) \xrightarrow{\Ver(f_2)} \Ver_Q(\Ver_Q(V_3)) \xrightarrow{j}
	\Ver_Q(V_3)$, where $j$ is induced from multiplication on $U\mathfrak{g}$.
	By part $(3)$ and $(4)$ of \Cref{canonical-iso-Pm}, after dualizing,
  applying $F_Q$ and the canonical isomorphisms, 
	the previous map is identified with 
  $P^{n+m} \otimes \mathcal{V}^{\vee}_3 \xrightarrow{\delta} P^{n} \otimes P^{m} \otimes \mathcal{V}^{\vee}_3  
  \xrightarrow{\text{id} \otimes \Phi_{Q}(f_2)} P^n \otimes \mathcal{V}^{\vee}_2 
  \xrightarrow{\Phi_{Q}(f_1)} \mathcal{V}^{\vee}_1$,
  where $\delta$ is the coalgebra map on $P$. 
  This is precisely how composition is defined for elements of 
	$\DiffOp(-,-)$.
	\end{proof}
	\end{theorem}

	\begin{remark} \label{restriction-creator-diff}
	Suppose that there exists some open embedding $U \hookrightarrow \flag$ and a subgroup $Q \subseteq B$ together 
	with a $Q$-reduction $Q_{\dR} \subseteq B_{\dR}$ over $U$,
	i.e. $B_{\dR}$ is the pushout of $Q_{\dR}$ along $Q \to B$. This defines a functor 
  $F_{Q}: \Rep(Q) \to \text{Coh}(U)$, and 
	for $V \in \Rep(B)$ by restricting \Cref{creator-diff} we have
	$F_Q \Ver_{B}(V) \cong F_{B} \Ver_{B}V \cong D_U \otimes \mathcal{V}$. Therefore we still get a functor 
	$$
    \Phi_Q: \Hom_{Q}(V_1, \Ver_{B}(V_2)) \to \DiffOp_U(\mathcal{V}^{\vee}_2, \mathcal{V}^{\vee}_1),
	$$
	where $V_{i} \in \Rep(Q)$. This situation will occur on the special fiber of $\flag$, 
  and by this procedure we will construct 
	the basic theta operators. 
	\end{remark}

\subsubsection{The isomorphism of divided formal completions of points}
Here we construct the isomorphisms of \Cref{canonical-iso-Pm} on PD formal completions of points, using 
the Grothendieck--Messing period map.
Fix $y \in \flag(\overline{\F}_p)$ mapping to 
$x \in \Sh(\overline{\F}_p)$. 

\begin{defn} \label{PD-differentials}
	Let $(R,I,\gamma_i)$ be a ring with divided powers over $\Z_p$.
	\begin{itemize}
	\item Define the sheaf of PD differential forms
	$\Omega^1_{(R,I)}$ as the quotient of the sheaf of differentials $\Omega^1_{R/\Z_p}$
	by imposing the relations $d \gamma_n(x)=\gamma_{n-1}(x)dx$ for all $x \in I$ and $n \ge 1$.
	\item Let $I_{\Delta} \subset D_{I_{\Delta}}(R \otimes R)$ 
	be the ideal of 
	the diagonal inside its divided power envelope, which we also denote as $P_{(R,I)}$.
	Then the sheaf of divided power differentials can be identified with 
	$I_{\Delta}/(I^{[2]}_{\Delta}+K^{[2]})$, where $K$
	is generated by elements of the form $\gamma_{m}(x) \otimes 1 - 1\otimes \gamma_{m}(x)$ for 
	$x \in I$, $m\ge 1$
	\cite[\href{https://stacks.math.columbia.edu/tag/07HT}{Tag 07HT}]{stacks-project}.
	\item  Define $P^{n}_{(R,I)}\coloneqq P_{(R,I)}/(I^{[n+1]}_{\Delta}+K^{[n+1]})$
	with the natural bimodule structure and coalgebra structure,
	and let $D^{n}_{(R,I)}$ be its dual. 
	Then $\text{gr}^{\bullet} P_{(R,I)}=\Sym^{\bullet} \Omega^1_{(R,I)}$. 
	\item If $f:(R,I) \to (S,J)$ is a map 
	of PD pairs, the usual differential extends to a map $df: f^* \Omega^1_{(R,I)} \to \Omega^1_{(S,J)}$,
	and a map $f^*: f^* P_{(R,I)} \to P_{(S,J)}$ which on graded pieces is induced by $df$.
	\item For $\mathfrak{X}=\text{Spf}(\varprojlim R/I_n)$ a formal scheme with an implicit divided power ideal 
	$J \subset I_0$, let $P_{\mathfrak{X}}:=\varprojlim P_{(R/I_n,J+I_n)}$ and 
  $D_{\mathfrak{X}}\coloneqq \cup_{n} \Hom(P^n_{\mathfrak{X}},\mathcal{O}_{\mathfrak{X}})$.
  In what follows, we will only 
	use this for $\Sh^{\sharp}_x$ and $\flag^{\sharp}_{y}$.
	\end{itemize}  
	\end{defn}

	Recall that the differentials of the Grothendieck--Messing period maps induce isomorphisms 
	$T_{\Sh^{\sharp}_x }\cong \pi^*_P T_{G/P^{\wedge}_{\infty}}$ and 
	$T_{\flag^{\sharp}_y }\cong \pi^*_B T_{G/B^{\wedge}_{\infty}}$. We are again fixing a trivialization $\phi$ 
	of $M_0$ for the definition of $\pi_{B,P}$. In particular, by the property that 
  $\text{gr}^{\bullet} P_{(R,I)}=\Sym^{\bullet} \Omega^1_{(R,I)}$, the maps
   $\pi_{P}$
	and $\pi_{B}$ induce isomorphisms between their respective sheaves of crystalline differential operators 
  and their duals $P$. 
	In general let $f: X \to Y$ be a map of schemes (or PD formal schemes) such that 
$df: D_X \to f^* D_Y$ is an isomorphism, and let $V$  be a quasicoherent sheaf on $Y$. Define
the isomorphism 
\begin{equation} \label{abdc}
df_{V}: D_X \otimes f^*V \to f^*(D_Y \otimes V)
\end{equation}
as follows. 
Locally, let $f: A \to C$ be a map of rings,
$B$ an $A$-bimodule, $V$ a left $A$-module, $D$ a $C$-bimodule and $\psi: C \otimes_{A} B \cong D$
an isomorphism of left $C$ modules which respects the right $A$ module structure on both sides. Then $\psi$ 
extends to an isomorphism $\varphi: C \otimes_A (B \otimes_A V) \to D \otimes_{C} (C \otimes_A V)$ given by 
$c \otimes b \otimes v \mapsto \psi(c \otimes b) \otimes (1\otimes v)$. Define $df_{V}$ to be the inverse of $\varphi$.
On graded pieces $df_{V}$ is given by $\Sym^{\bullet} df$ on the left and the identity 
on $f^* V$.
One can check the following compatibility with composition: with $f$ and $V$ as above
the diagram 
$$ 
\begin{tikzcd} \label{composition}
D_X \otimes f^*V \arrow[rr,"df_V"] & & f^*(D_Y \otimes V)  \\
D_X \otimes D_X \otimes f^*V  \arrow[u,"\mu \otimes \text{id}"]
 \arrow[r,"\text{id} \otimes df_V"] & 
D_X \otimes f^*(D_Y \otimes V) \arrow[r,"df_{D_Y \otimes V}"] &
f^*(D_Y \otimes D_Y \otimes V) \arrow[u,"f^* (\mu \otimes \text{id})"]
\end{tikzcd}
$$
commutes. Similarly, if $f^*: P_{Y} \to f^*P_{X}$ is an isomorphism we define 
\begin{equation} \label{dual-abdc}
df^{\vee}_V:  P_X \otimes f^*V \to f^*(P_{Y}\otimes V) 
\end{equation}
in the same way as above, using that $P_{X}$ is also a 
bimodule. 
The next result defines the canonical isomorphisms of \Cref{canonical-iso-Pm} on divided power 
formal completions of points. Recall the notation $j: \Sh^{\sharp}_{x} \hookrightarrow \Sh$, 
$i: \flag^{\sharp}_{y} \hookrightarrow \flag$.
\begin{prop} \label{prop-on-PD-completions}
	Let $y \in \flag(\overline{\F}_p)$ and $x \in \Sh(\overline{\F}_p)$ its image, define $\pi_{B/P}$ for some 
	fixed choice of $\phi$ as in \Cref{GM-maps-P}. Let $R \in \{\mathcal{O},k\}$ and 
	let $V \in \Rep_{R}(P), W \in \Rep_{R}(B)$. Define the following composition of isomorphisms,
	$$
	e_{V,x}: j^*(D_{\Sh} \otimes F_P(V)) \cong D_{\Sh^{\sharp}_x} \otimes \pi_{P}^*F_{G/P}(V) \xrightarrow{d\pi_{P,V}}
	\pi^*_P(D_{G/P^{\wedge}_{\infty}} \otimes F_{G/P}(V)) \cong j^* F_P(\Ver_P(V))
	$$
	$$
	e_{W,y}: i^*(D_{\flag} \otimes F_B(W)) \cong D_{\flag^{\sharp}_{y}} 
	\otimes \pi_{B}^* F_{G/B}(W) \xrightarrow{d\pi_{B,W}}
	 \pi^*_B(D_{G/B^{\wedge}_{\infty}} \otimes F_{G/B}(W)) \cong i^* F_B(\Ver_B(W)),
	$$
	where the isomorphisms on the sides come from \Cref{GM-maps-P}(1) and \Cref{verma-flag},
  and the ones in the middle are defined 
	by \eqref{abdc}. The first two isomorphisms above induce an isomorphism
	$$
	e_{P/B,W,y}: j^*(D_{\flag/\Sh} \otimes F_B(W)) \cong D_{\flag^\sharp/\Sh^\sharp} \otimes \pi^*_B
	F_{G/B}(W) \to    j^*F_B(\Ver_{P/B}(W)).
	$$
  Similarly, for $R=\mathcal{O}/p^n$ and $V \in \mathcal{O}_{B,R}$ or $V \in \Rep_{R}(B)$ there is an isomorphism 
  $$
  e^{\vee}_{V,x}: j^*(P_{\flag} \otimes \mathcal{V}^{\vee}) \cong P_{\flag^{\sharp}_x} \otimes \pi^*_{B}
  F_{B}(V^{\vee}) \xrightarrow{df^{\vee}_{V^{\vee}}} \pi^*_{B}(P_{G/B} \otimes \mathcal{V}^{\vee})
  =j^*(F_{B} \Ver(V)^{\vee})
  $$
  over $\flag^{\sharp}_{y,R}$.
	Then $e_{V,x}$, $e_{V,y}$, $e_{P/B,W,y}$, and $e^{\vee}_{V,x}$ do not depend on the choice of $\phi$.
	If $W \in \Rep_{k}(B)$ there is also a canonical
	isomorphism 
	$$
	e^0_{W,y}: j^*F_B(\Ver^0_{P/B}(W)) \cong j^*(F^* F_* \mathcal{W}^\vee)^\vee,
	$$
	where $F: \flag \to \flag^{(p)}$ is the relative Frobenius 
	with respect to $\Shbar$.
	Furthermore, these local isomorphisms satisfy the same properties 
	$(1)$ to $(7)$ as in \Cref{canonical-iso-Pm}.
	\begin{proof}
	By the definition of $D$ on divided power formal schemes
	we get maps $d\pi_{P}: D_{\Sh^{\sharp}_x} 
	\to \pi^*_{P} D_{G/P^{\wedge}_{\infty}}$ and 
	$d\pi_{B}: D_{\flag^{\sharp}_{y}} 
	\to \pi^*_{B} D_{G/B^{\wedge}_{\infty}}$
	that respect the natural $\mathcal{O}_{G/Q^{\wedge}_{\infty}}$-bimodule 
	structure on both sides, and similarly for $P_{\Sh^{\sharp}_x}$. We 
	then define the maps $d\pi_{P,V}$ as in \eqref{abdc}.
	 Using the first diagram of \Cref{GM-KS}
	we see that on graded pieces $e_{V,x}$ is induced by the Kodaira--Spencer, so that it is an isomorphism. 
	Similarly, on graded pieces $e_{W,y}$ is induced by the symmetric powers of the isomorphism $e^1_1$
  from \Cref{KS-flag}, so it is also an isomorphism.
  Both are independent of the choice of $\phi$, since both $\pi_{Q,V}$ and the isomorphisms from \Cref{GM-maps-P}(1) 
  vary in the same way when changing $\phi$. Namely, if we use the notation $\pi_{Q,\phi}$
  to emphasize the dependence on $\phi$, then $\pi_{Q,g\phi}=g \circ \pi_{Q,\phi}$
  for $g \in G(\check{\Z}_p)$. 
	The map $e_{P/B,y}$ is induced by the cartesian diagram on \Cref{GM-maps-P}, which identifies 
  the relative differentials $D_{\flag^{\sharp}_{y}/\Sh^{\sharp}_{x}}$ and $D_{(G/B)^{\wedge}_{\infty}/(G/P)^{\wedge}_{\infty}}$
  under $\pi_{B}$, and then
	\Cref{verma-flag} identifies $D_{(G/B)^{\wedge}_{\infty}/(G/P)^{\wedge}_{\infty}}$ with $\Ver_{P/B}(1)$.  
  For $e^{\vee}_{V}$ we use \eqref{dual-abdc} to define it, it is independent of $\phi$ 
  since if we write $V=\cup V_n$ with $V_n \in \Rep(B)$, it is given by the dual of $e^{\le N}_{V_n,x}$ on all the corresponding 
  quotients. 
	
	Let $W \in \Rep_{\F_p}(B)$. There is an isomorphism $ 
	j^* F_B(\Ver^0_{P/B}(W)) \cong (G^*G_* j^*\mathcal{W}^\vee)^\vee$, where in the right-hand side 
	$G$ is the relative Frobenius with respect to $\flag^{\sharp}_{y}\to \Shbar^{\sharp}_x$.
	It follows from the corresponding isomorphism 
	on the flag variety $P/B$ in \Cref{verma-flag}(2), and the cartesian diagram in \Cref{GM-maps-P} 
	that identifies both relative Frobenius.
	To obtain $e^0_{W,y}$ we claim that
	there is a canonical isomorphism $G^* G_* j^* E \cong j^*F^* F_* E$ 
	for $E$ a vector bundle. After localizing on $\flag$ we may assume that $E=\mathcal{O}_{\flag}$.
	Write $\flag \to \Shbar$ locally as $A \to B$, $\flag^{\sharp}_{y}\to \Shbar^{\sharp}_x$
	as $\tilde{A} \to \tilde{B}$ and $j: B \to \tilde{B}$. By smoothness, we may assume that 
	$\tilde{A}=\overline{\F}_p[[\frac{x^n_i}{n!}]]$ and $\tilde{B}=\overline{\F}_p[[\frac{x^n_i}{n!},y_j]]$.
	The sheaf $G^* G_* \mathcal{O}_{\flag^{\sharp}_{y}}$ is finite locally free,
	since one can check that $G$ is flat and of finite presentation.
	There is a map $\phi: j^*F^* F_* \mathcal{O}_{\flag} \to G^* G_* \mathcal{O}_{\flag^{\sharp}_{y}}$
	as follows. Locally it is given by
	$\tilde{B} \otimes_{B} (B \otimes_{B \otimes_{A,F} A} B) \to \tilde{B} 
	\otimes_{\tilde{B} \otimes_{\tilde{A},G} \tilde{A}} \tilde{B}$ sending
	$1 \otimes x \otimes y \mapsto j(x) \otimes j(y)$. Since both the domain and target are finite locally free of 
	the same rank
	it is enough to check whether the determinant of $\phi$ is a unit in $\tilde{B}$. The map $\tilde{B}
	\to (\tilde{B}/\mathfrak{m}_{\tilde{A}}\tilde{B})$ sends non-units to non-units, so we can check it 
	after replacing $\tilde{B}$ by $\tilde{B}/\mathfrak{m}_A \tilde{B}$, and
	$B$ by $B/\mathfrak{m}_A B$. Relabelling everything,
	$\tilde{B} \to \tilde{A}=A=\kappa(x)$ becomes 
	the formal completion of $B$ at $\kappa(y)$, $F$ is the Frobenius for $B/\kappa(y)$ and 
	$G$ the Frobenius for $\tilde{B}/\kappa(y)$. We can work \'etale locally since the formal completions don't change 
	and the Frobenius pushforward commutes with \'etale base change. 
	Thus, we can reduce to affine space, where it can be checked by hand.

	Properties $(1)$ and $(2)$ follow by construction of $e_{P/B,V}$ and $e^0_V$,
	and the second diagram in \Cref{GM-KS}.
	Property $(3)$ follows directly from the associated statement on the flag variety, and $d\pi$ being functorial.
	Property $(4)$ follows from the statement on the flag variety and diagram 
	\eqref{composition}.
	Part $(5)$ follows from \Cref{GM-KS} and the statements in
	\Cref{PD-differentials}. For $(6)$
	it is sufficient to prove it after restricting to the degree at most $1$ filtered piece of the filtration,
  since the Gauss--Manin connection
	extends uniquely to a 
	HPD stratification (and $e^{\vee}_{1,y}$ is just given by the differential of $\pi_{B}$). The commutativity of the diagram is equivalent to $\nabla$ on $\Sh^{\sharp}_{x}$ being 
	the pullback connection along $\pi_{P}$ of the trivial connection on $F_{G/P}(V)\cong V \otimes
	 \mathcal{O}_{G/P^{\wedge}_{\infty}}$, since the tensor identity is precisely 
	 the one inducing the trivial connection on $F_{G/P}(V)$. We can reduce to the case where $V=\Lambda$
	 from the definition of $\nabla$,
	 where this is precisely \Cref{GM-maps-P}(2). For $(7)$ it follows from \Cref{flag-g-modules}(3) and the fact that 
   the stratification $\epsilon$ of \Cref{stratification-defn}(3) is compatible with maps that induce an isomorphism 
   on differentials, which 
   can be checked on \'etale local coordinates. 
	\end{proof}
	
	\end{prop}

\subsubsection{Definition of the canonical isomorphisms}

We construct the canonical isomorphism $e_V$ of \Cref{canonical-iso-Pm},
by reducing to the case of $\Ver_{Q}^{\le 1}(V)$ 
for  $V \in \Rep(G)$.
\begin{lemma} \label{B-rep-sub}
Let $R \in \{\mathcal{O},\mathcal{O}/p^n\}$.
Let $Q \in \{P,B\}$. Every element of $\Modfg_{R}(Q)$ is a subquotient as a $Q$-module of
the restriction to $Q$ of an element of $\Modfg_{R}(G)$. 
\begin{proof}
Let $V \in \Modfg_{R}(Q)$. Write it as $V= \oplus_{i} (R/p^i)^{n_i}$ and let 
$\pi_{i,j}: V \to R/p^i$ be the corresponding 
projections.
The map $V \to \oplus \mathcal{O}^{j}_{Q,R/p^i}$ sending $v \to (g \mapsto \pi_{i,j}(gv))$ defines an 
injective map of $Q$-representations, where $Q$ acts on $\mathcal{O}_{Q}$ by $q \cdot f(x)=f(q^{-1}x)$.
For each $i,j$ the restriction map $\phi: \mathcal{O}^j_{G,R/p^i} \to \mathcal{O}^j_{Q,R/p^i}$ is a surjection of $Q$-representations.
Let $\{v_{i,j}\}$ be the generators of $\oplus \mathcal{O}^{j}_{Q,R/p^i}$.
By \cite[\S 1.5 Prop 2]{serre-reductive} there exists $W \in \Modfg_{R}(Q)$ which 
contains some choice of elements $w_{i,j} \in \phi^{-1}(v_{i,j})$, and which is contained in $\phi^{-1}(V)$,
so that $W$ maps surjectively to $V$. Similarly, there exists $M \in \Modfg_{R}(G)$
containing $W$ and contained in $\oplus \mathcal{O}^{j}_{G,R/p^i}$. Then $W \subset M$ gives a presentation
 of $V$ as a subquotient.
\end{proof}
\end{lemma}

We prove a small lemma that will allow us to check all the properties of the canonical 
isomorphisms on PD formal completions.
To apply it in our case, we will use 
that one can cover $\Sh_{\mathcal{O}}$ by connected affine open subsets that have at least one $\overline{\F}_p$-point.
This follows from the fact that $\Sh^{\tor}$ is proper, 
so no connected component of it can live in characteristic $0$.

\begin{lemma} \label{PD-enough}
Let $R/\mathcal{O}$ be a smooth connected algebra of finite type having at least one $\overline{\F}_p$-point. 
Let $M, N$ be $R$-modules which are the sum of a finite free part and finitely many summands of the form 
$R/p^n$,
and let $\psi: M \to N$ be a map between them. Then $\psi$ is $0$ if and only if it is $0$ 
after base change to every PD formal completion of $\overline{\F}_p$-points. Alternatively, if $S/k$ is 
smooth connected of finite type, $M,N$ locally free finite $S$-modules, and $\phi: M \to N$ a map, then $\phi$ is zero if and only if 
it is zero on every $\overline{\F}_p$-point.  
\begin{proof}
The map $\psi$ decomposes into maps between free modules and maps of the form $R \to R/p^n$ or 
$R/p^m \to R/p^n$.
For the former 
the lemma follows since under the assumptions on $R$
 the composition $R\to R_{\mathfrak{m}} \to R^{\wedge}_{\mathfrak{m}} \to R^{\sharp}_{\mathfrak{m}}$ is injective for 
any $\overline{\F}_p$-point $\mathfrak{m}$. The last injection can be checked using \'etale coordinates, and will 
be of the form $\mathcal{O}[[x_i]] \hookrightarrow \mathcal{O}[[\frac{x^n_i}{n!}]]$. 
For the latter $R/p^n \to R^{\sharp}_{\mathfrak{m}}/p^n$ is not necessarily injective for a particular $\mathfrak{m}$, 
but any element in the kernel will be contained in $\mathfrak{m}$. Therefore, an element $x$
that is killed for all maximal ideals $\mathfrak{m}$ above $p$
will be contained in $p R$, since $R \otimes \F_p$ is a Jacobson ring.
Writing $x=p y$, since $R^{\sharp}_{\mathfrak{m}}$ is $p$-torsion free, we see that $n=1$ or
$y$ is also in the kernel
for every $\mathfrak{m}$, 
so that by induction $x=0$.  The second statement follows easily since $S$ is Jacobson, and $\phi$ can be locally 
written as a matrix.  
\end{proof}
\end{lemma}

We define $e^{\le 1}_V$ using its relation to the Gauss--Manin connection prescribed by 
\Cref{canonical-iso-Pm}(6). We remark that the functors $F_Q$ restricted to $\flag^{\sharp}_{y}$ are still exact over 
$\Rep_{\mathcal{O}}(Q)$, since 
divided power formal completions are $\mathcal{O}$-flat. 
\begin{prop} \label{def-less-1}
Let $R \in \{\mathcal{O},\mathcal{O}/p^n\}$.
Define $e^{\le 1}_{1}: F_{B}(\Ver^{\le 1}(1)) \cong F_{B}(1 \oplus \mathfrak{g}/\mathfrak{b}) \cong \mathcal{O} 
\oplus T_{\flag_{R}} \cong D^{\le 1}_{\flag}$
where the isomorphism in the middle is given by \Cref{basic-iso}.
\begin{enumerate}
\item 
For $V \in \Modfg_{R}(G)$ define $e^{\le 1}_V$ 
by the commutative diagram 
$$
\begin{tikzcd}
F_B(\Ver_B^{\le 1}(V)) \arrow[d,"e^{\le 1}_V"] \arrow[r,"F_B(\phi_V)"] &
F_B(V \otimes_{R} \Ver_B^{\le 1}(1)) \arrow[r,"\textnormal{id} \otimes e^{\le 1}_{1}"]
& F_B(V) \otimes D^{\le 1}_{\flag} \\ 
D^{\le 1}_{\flag} \otimes F_B(V) \arrow[rru,"\nabla"] & & 
\end{tikzcd}
$$
where $\phi_V$ is the tensor identity of \Cref{flag-g-modules}, and $\nabla$ is the dual of the 
Gauss--Manin connection on $\mathcal{V}^{\vee}/\flag_{R}$ of
\Cref{GM-definition}, seen as an isomorphism $P^{1} \otimes \mathcal{V}^{\vee} \cong \mathcal{V}^{\vee} \otimes P^1$. 
The definition makes sense since the horizontal and diagonal maps are isomorphisms. 
\item 
Let $i:V \hookrightarrow W$ be a $B$-equivariant embedding, where $W \in \Modfg_{R}(G)$ 
and $V \in \Modfg_{R}(B)$. Then 
$e^{\le 1}_W$ sends $F_B(\Ver_B^{\le 1}(V))$ to $D^{\le 1}_{\flag} \otimes F_B(V) \subseteq D^{\le 1}_{\flag} \otimes F_B(W)$.
We define $e^{\le 1}_V$ as the restriction of $e^{\le 1}_W$, by exactness of $F_B$ and $\Ver_{B}(-)$. 
It is an isomorphism independent of the embedding $i$. 
\item Let $\pi: W \twoheadrightarrow V$ be a surjection in $\Modfg_{R}(B)$ with $W$ a submodule of 
an element of $\Modfg_{R}(G)$. Then 
$e^{\le 1}_{W}$ from point $(2)$ sends $F_B(\Ver_B^{\le 1}(\Ker \pi))$ to $D^{\le 1} \otimes F_B(\Ker \pi)$, so that 
we can define $e^{\le 1}_V$ as the induced map. It is an isomorphism independent of the surjection. 
\item Let $V \in \Modfg_{R}(B)$. By \Cref{B-rep-sub} it is a subquotient of a $G$-representation, 
so we define $e^{\le 1}_V$ via the previous two points. It is an isomorphism 
independent of the presentation as a subquotient. For $V \in \Modfg_{R}(P)$ we define $e^{\le 1}_{V}$ analogously.
\end{enumerate}

\begin{proof}
We prove $(2)$ and $(3)$. The independence of embedding/quotient follows by comparing any two of them with their sum, 
since the construction of $e^{\le 1}_V$ in $(1)$ is clearly functorial. 
Let $R=\mathcal{O}$ first. By \Cref{PD-enough} we can check well-definedness 
on PD formal completions of $\overline{\F}_p$ points, e.g.
for $(2)$ we have to prove that 
$F(\Ver_B^{\le 1}(V)) \to D^{\le 1}_{\flag} \otimes F_B(W/V)$ is $0$. 
The key fact is that by \Cref{prop-on-PD-completions} the diagram in $(1)$ defining $e^{\le 1}_W$ agrees with the local 
map $e^{\le 1}_{W,x}$ at every PD formal completion. Then by naturality of the local isomorphisms
parts $(2)$ and $(3)$ hold on PD formal completions. This shows that the maps $e^{\le 1}$ 
are well-defined, and they are isomorphisms 
since one can construct an inverse by the same procedure, using 
the inverse map on $W$. For $R=k$ we use the Grothendieck--Messing map mod 
$(p,\mathfrak{m}_x^p)$, where $\mathfrak{m}_x$ is the maximal ideal of some $\overline{\F}_p$-point. Since $\mathfrak{m}_x$
has divided powers we still have that the restriction of $e^{\le 1}_V$ agrees with the isomorphisms 
in \Cref{prop-on-PD-completions}. Then we use the second part of \Cref{PD-enough} to conclude. 

\end{proof}
\end{prop}

We define $e^{\le n}_V$ by induction on $n$, by leveraging the algebra structure on the sheaf of differential 
operators. 
\begin{prop} \label{induction-dfn}
Let $R \in \{\mathcal{O},\mathcal{O}/p^n\}$, and
 $V \in \Rep_{R}(Q)$. For $n=1$ $e^{\le 1}_V$ is defined in
\Cref{def-less-1}. For $n \ge 2$ the surjection 
$\Ver_Q^{\le 1}(\Ver_Q^{\le n-1}V) \twoheadrightarrow \Ver_Q^{\le n}(V)$ induces by \Cref{def-less-1} a surjection
$$
\mu: D^{\le 1} \otimes F_Q(\Ver_{Q}^{\le n-1}(V)) \twoheadrightarrow F_Q(\Ver_{Q}^{\le n}(V)).
$$
Let $v \in F_Q(\Ver^{\le n}(V))$ be a local section on a small enough open so that it lies in the image of $\mu$,
and let $g=\sum D_i \otimes v_i$ be any element in $\mu^{-1}(v)$. 
Inductively define $e^{\le n}_V(v)\coloneqq \sum D_i \cdot e^{\le n-1}_{V}(v_i)$, where 
$\cdot: D^{\le 1} \otimes D^{\le n-1} \otimes F(V) \to D^{ \le n} \otimes F(V)$ is given by composition 
of differential operators. Then 
$e^{\le n}_V$ is a well-defined isomorphism independent of the choice of $g$. Moreover, on PD formal completions of 
$\overline{\F}_p$-points it agrees with the local isomorphisms in \Cref{prop-on-PD-completions}.

\begin{proof}
We can check well-definedness on PD formal completions of points, by \Cref{PD-enough}.
It follows from property 
$(4)$ at the end of \Cref{prop-on-PD-completions} that $e_{V}$ and $e_{V,x}$ agree on PD-formal completions
of points, so in particular $e_V$ is well-defined, and it is independent of the choice of $g$.
By its inductive construction the map $e_{V}$ is surjective. Let $N$ be its kernel.
For every $\Fpbar$-point $x$ of $\Sh$ we have that $N_{\Sh^{\sharp}_{x}}=0$ by \Cref{prop-on-PD-completions}, 
so in particular $N_{\kappa(x)}=0$. 
Since
$\Sh$ has no closed point in characteristic $0$, and $N$ is finitely generated as a $\mathcal{O}_{\Sh}$-module
being the kernel of a map between coherent sheaves, 
we get that $N=0$ by Nakayama's lemma.
Thus, the map is an isomorphism.
\end{proof}
\end{prop}

We extend the definition of the dual canonical isomorphism to $\mathcal{O}_{B,R}$.
\begin{prop} \label{dual-canonical-iso}
Let $R=\mathcal{O}/p^n$, $V \in \Modfg_{R}(B)$ or $V \in \mathcal{O}_{B,R}$. Write $V=\cup V_n$ with $V_n \in \Modfg_{R}(B)$.
Let $\mathcal{V}^{\vee}=F_{B}(V^{\vee})$.
Then there exists 
an isomorphism $e^{\vee}_{V}: F_{B}\Ver_{B}(V)^{\vee} \cong P_{\flag_{R}} \otimes \mathcal{V}^{\vee}$ fitting 
in the diagram 
$$
\begin{tikzcd}
  F_{B}\Ver_{B}(V)^{\vee} \arrow[r,"e^{\vee}_V"]  \arrow[d] & P_{\flag} \otimes \mathcal{V}^{\vee} \arrow[d] \\
  F_{B} \Ver^{\le N,\vee}_{B}(V_n) \arrow[r,"e^{\le N,\vee}_{V_n}"] & P^N_{\flag} \otimes \mathcal{V}^{\vee}_n
\end{tikzcd}
$$
for all $n,M$, 
where $e^{\le N}_{V_n}$ are defined in \Cref{induction-dfn}.
The analogous statement holds for $\mathcal{O}_{P,R}$.
\begin{proof}
Note that if $e^{\vee}_{V}$ exists and it fits in the diagram, then it is unique and it is an isomorphism. 
Using the $e^{\le N,\vee}_{V_n}$ we can define a map
 $F_{B}\Ver_{B}(V)^{\vee} \to [P_{\flag} \otimes \mathcal{V}^{\vee}]^{\wedge}$,
 where the right-hand side is the formal completion with respect to the cofiltration $P^N_{\flag} \otimes \mathcal{V}^{\vee}_n$.
 We want to prove that it factors through 
$P_{\flag} \otimes \mathcal{V}^{\vee}$. We can check this on $\Fpbar$ points. Locally we have that 
$[P_{\flag} \otimes \mathcal{V}^{\vee}]^{\wedge}/[P_{\flag} \otimes \mathcal{V}^{\vee}]$
is of the form $\prod_{I} S/\oplus_{I} S$
for $S/R$ as in 
\Cref{PD-enough}. Then we can check that an element of 
$\prod_{I} S/\oplus_{I} S$ is zero if and only if it is so on each geometric point. 
Then it follows from the properties of \Cref{prop-on-PD-completions}. 

\end{proof}
\end{prop}

We now combine the results above to prove \Cref{canonical-iso-Pm}. 
\begin{proof}[Proof of \Cref{canonical-iso-Pm}]
The isomorphisms $e_{V}$ for $\Ver_{P}(V)$ and $\Ver_{B}(V)$ are defined in \Cref{induction-dfn},
and $e^{\vee}_{V}$ is defined in \Cref{dual-canonical-iso}.
They
agree on PD formal completions with the ones defined on \Cref{prop-on-PD-completions}, by all the functoriality 
and compatibility properties that the local isomorphisms satisfy. 
For $e^{P/B}_V$, 
 let $E$ be the cokernel of 
$D_{\flag/\Sh} \otimes F(V) \to D_{\flag} \otimes F(V)$. We claim that the composition
$F_B(\Ver_{P/B}(V)) \to F_B(\Ver_B(V)) \xrightarrow{e_V} D_{\flag} \otimes F(V) \to E$ is $0$. 
By \Cref{PD-enough} it is enough to check it on PD formal completions, where it holds by
 \Cref{prop-on-PD-completions}(1).
 By exactness of $F_B$ this defines the isomorphism $e^{P/B}_V$, and it automatically 
satisfies property $(1)$. To define $e^0_V$ we use the same strategy with the square of property $(2)$: now using that we can 
check if a map of vector bundles over $\flag_{\overline{\F}_p}$ is zero on closed points,
so in particular on PD completions of points. 
This automatically proves property $(2)$. To check $(3)-(5)$ we can use \Cref{PD-enough}, by considering 
the difference of the expected map and the actual map 
(to be strict we consider each of these maps on $\Ver^{\le n}$),
and then it holds over PD formal completions by \Cref{prop-on-PD-completions}.
 Part $(6)$ for $V \in \Modfg_{R}(G)$ holds since it does on PD formal completions, for $V \in \mathcal{O}_{B,R}$
it holds since it does on every finite free projection. 
Part $(7)$ also follows from 
\Cref{prop-on-PD-completions}, since we can check it at the level of PD stratifications,
where the objects are vector bundles. 
 The isomorphisms $e^{\le 1}$ 
are Hecke equivariant since the Gauss--Manin connection and the isomorphisms of \Cref{basic-iso} are,
so the maps $e^{\le n}$
are also Hecke equivariant by 
their inductive definition.

Finally, we construct $e_V$ on toroidal compactifications. For $V \in \Modfg_{R}(G)$, $e^{\le 1}_V$
is constructed as in \Cref{def-less-1}, using the (sub)canonical extension of Gauss--Manin to a log connection. 
To construct $e^{\le 1}_V$ for general $V \in \Modfg_{R}(Q)$ we use the same procedure and the fact that one can check whether 
a map of sheaves as in \Cref{PD-enough} vanishes on the interior $\Sh \hookrightarrow \Sh^{\tor}$.
It is easy to see that they are isomorphisms, since the ones for $V \in \Modfg_{R}(G)$ are by definition. 
For $e^{\le n}$
we follow the procedure of \Cref{induction-dfn}, which is well-defined since it can be checked on the interior. 
From its construction we immediately see that $e_V$ is a surjection, and it is injective since it is so in the interior.
We define $e^{\vee}_{V}$ as in \Cref{dual-canonical-iso}, it is well-defined since we can check the lifting 
property on the interior.
One can see that all the required properties can be checked on an open dense subset of $\Sh^{\tor}$, so they 
follow from the properties on $\Sh$.
\end{proof}
	
\begin{remark}
One could have also constructed $e^0_V$ from \Cref{canonical-iso-Pm} via an appropriate base change theorem along the
cartesian diagram (see the next subsection)
$$
\begin{tikzcd}
\flag \arrow[r] \arrow[d,"\pi"] &  \GF \arrow[d]\\
\Shbar \arrow[r] & \GZip \\
\end{tikzcd}
$$
which identifies the two relative Frobenii, and then using the isomorphism of 
\Cref{verma-flag}(2) on $P/B$.
\end{remark}

\subsection{Description on the open stratum of the special fiber} \label{section2.3}
We now describe the canonical isomorphisms of \Cref{canonical-iso-Pm} more explicitly on some open dense subset 
$U \subseteq \flag_k$. As the main application, in \Cref{composition-U} 
we prove the compatibility with composition of a functor coming from
\Cref{restriction-creator-diff}, in a way that it produces differential operators over $U$. 
Consider $G/\F_p$, let $\sigma \in \Gal(\Fpbar/\F_p)$ be the $p$th power Frobenius, 
and let $P^{(p)}=P \times_{k,\sigma} k$. 
 First we introduce the notation for $G$-Zips. Recall that we have a
tuple $\mathcal{Z}=(T,P,B,G)$ of a maximal torus, parabolic, and Borel; all defined over $\F_p$
except $P$. Let $P^{-}$ be the opposite 
parabolic, and $M$ the Levi of $P$. In this subsection, fix a choice of $B_0$ fitting 
in the embedding $(G,P,B) \hookrightarrow 
(\GL(\Lambda),P_0,B_0)$ defined over $\Fpbar$. We will invoke \Cref{special-embedding} whenever a more precise 
embedding is needed,
in which case we will use such an embedding.   
\begin{defn}
A $G$-Zip of type $\mathcal{Z}:=(G,P,B)$ over some scheme $X/\F_p$ is a tuple $(I,I_{+},I_{-},\phi)$ where $I$ is a $G$-torsor
over $X$, 
$I_{+} \subseteq I$ is a $P$-torsor, $I_{-} \subseteq I$ a $P^{-(p)}$-torsor, and 
$\phi: I^{(p)}_{+}/U^{(p)}_{P} \cong I_{-}/U_{P^{-(p)}}$ is an isomorphism of $M^{(p)}$-torsors.
Here $I^{(p)}_{+}$, or on any torsor, denotes the pullback along the absolute Frobenius on $X$.
Similarly, a $G$-Zip flag of type $\mathcal{Z}$ is a tuple as before together with a $B$-subtorsor 
$I_{B} \subseteq I_{+}$. We denote by $\GZip$ and $\GF$ the stacks of $G$-Zips (flags) of type $\mathcal{Z}$.
\end{defn}

There is a cartesian square 
$$
\begin{tikzcd}
\flag \arrow[r,"\xi_{B}"] \arrow[d,"\pi"] & \GF \arrow[d] \\
\Shbar \arrow[r,"\xi_{P}"] & \GZip  
\end{tikzcd}
$$
where the maps $\xi_{P}$ and $\xi_{B}$ are smooth by \cite{Zhang-EO}. The map $\xi_{P}$ is defined by sending a point of $\Shbar$ 
to the universal $G$-zip $(G_{\dR}, P_{\dR}, P^{-(p)}_{\dR}, \phi)$ where 
$$
P^{-(p)}_\dR=\text{Isom}((\text{Ker }V \subseteq \mathcal{H}, s_{\alpha,\dR}),
 (L^{\vee} \subseteq \Lambda, s_{\alpha}) \otimes \mathcal{O}) 
$$
is a $P^{-(p)}$-torsor through the base change map $P^{-(p)} \to P^{-}$ acting on 
the right-hand side.
Then $\phi: P^{(p)}_\dR/U^{(p)}_{P} \cong P^{-(p)}_{\dR}/U_{P^{-(p)}}$ is the isomorphism induced by the Cartier 
isomorphisms $\Ker V \cong \omega^{\vee,(p)}_{A^{\vee}}$ and $\mathcal{H}/\Ker V \cong \omega^{(p)}_{A}$. 
Recall that the pushout of $B_{\dR}$ along $B \to B_0$ defines a refinement $F^{\bullet}_{\mathcal{H}}$ of the 
Hodge filtration on $\mathcal{H}$.
On $\flag$ let 
$D_{\mathcal{H},\bullet} \subseteq \mathcal{H}$ be the ascending full flag refining the conjugate filtration,
which is obtained from $F^{\bullet}_{\mathcal{H}}$ using the isomorphism of graded filtrations $\phi$. 
That is, the filtration $F^{\bullet}_{\mathcal{H}}$ on $\mathcal{H}/\omega \cong \omega^{\vee}_{A^{\vee}}$
induces a filtration on
 $\Ker V \cong \omega^{\vee,(p)}_{A^{\vee}}$, and the Hodge filtration on $\omega$ induces a filtration on 
 $\mathcal{H}/\Ker V \cong \omega^{(p)}$.
We will take $D_{\mathcal{H},0}=0$
and $D_{\mathcal{H},n}=\mathcal{H}$. 
Let 
$$
B^{z}_{\dR}=\text{Isom}_{\flag}((D_{\mathcal{H},\bullet}, s_{\alpha,\dR}), 
(F^{n-\bullet}_{\Lambda}, s_{\alpha}) \otimes \mathcal{O}_{\flag}).
$$
It is a torsor for the Borel $B^{z}:=zBz^{-1}$ generated by $(B \cap M)^{(p)}$ and $U_{P^{-(p)}}$. 
This is because $B^{z}$ is the intersection of $G$ and the Borel of $\GL_n$ generated by 
$U_{P^{-(p)}_0}$ and $B^{(p)}_0 \cap M^{(p)}_0$. 
Since $(B,T)$ are defined over $\F_p$, there exists such a $z \in W$ \cite[Lem 2.3.4]{GK-stratification}.

Let $(V,B_0,B^{z'}_0)$ be some choice of data coming from \Cref{special-embedding}.
By pushout we can define torsors $\GL(V)_{\dR}:=G_{\dR} \times^{G} \GL(V), P_{V,\dR},B_{0,\dR}, P^{-(p)}_{0,\dR}, 
B^{z'}_{0,\dR}$. 
Then $B_{0,\dR} \subseteq \GL(V)_{\dR}$ defines a descending full flag $F^{\bullet}_{\mathcal{V}}$ of $\mathcal{V}$
refining the Hodge filtration, 
and $B^{z'}_{0,\dR}$ defines a full flag $D_{\mathcal{V},\bullet}$
For a given $(V,B_0,B^{z'}_0)$ let $\LL_{V,i}:=F^{n-i}_{\mathcal{V}}/F^{n-i+1}_{\mathcal{V}}$, where $n$ is the dimension 
of $V$.
There is a map $\GF \to [B \times B^{z}
\backslash G]$, 
where $(b,b')g=bg b^{'-1}$, by only remembering the $B$-torsor, and the $B^{z}$-torsor induced 
from the $P^{-(p)}$-torsor as above. We compose with the map $[B \times B^{z}\backslash G]\to \Sbt:=[B \times 
B\backslash G]$
induced by multiplication by $z$ on the right. 
Thus, we get a map $\zeta: \flag \to \Sbt$. Similarly, let $\Sbt_{\GL}=[B_0 \times B_0\backslash \GL_n]$, 
there is a map $\Sbt \to \Sbt_{\GL}$. 
\begin{defn}
\begin{enumerate}
\item Let $B_{\Sbt}=[1\times B \backslash G \to \Sbt]$, it is a $B$-torsor such that $\zeta^{*}B_{\Sbt}=B_{\dR}$.
\item Let $\Sbt_{1}=[B \times B\backslash Bw_0B]$ be the open Schubert stratum. 
Let $(B_0,B^{z'}_0, \GL(V))$ be as \Cref{special-embedding}.
It induces a map $[B \times B^z \backslash G] \to [B_0 \times B^{z'}_0 \backslash \GL(V)]$, which 
we identify with a map $\Sbt_{G} \to \Sbt_{\GL}=[B_0 \times B_0 \backslash \GL(V)]$ after twisting appropriately.
We will make this identification from now on.
Then $\Sbt_{1}$ maps to 
$\Sbt_{\GL(V),1}=[B_0 \times B_0 \backslash B_0 w_{\GL(V),0} B_0]$.  
 Let $U=\zeta^{-1}(\Sbt_{1})$, it is 
an open dense subspace. 
\end{enumerate}
\end{defn}
By \Cref{special-embedding}
 $\Sbt_{1}$ maps to $\Sbt_{\GL(V),1}$, and we have that the extended Hodge filtration $F^{\bullet}_{\mathcal{V}}$
 and the extended conjugate filtration $D_{\mathcal{V},\bullet}$ on $\mathcal{V}$ are in general position.
  In particular there are isomorphisms 
 $H_{n-i}: F^{n-i}_{\mathcal{V}} \cap D_{\mathcal{V},n-i+1} \cong 
 {\LL'}_{\mathcal{V},n-i}:=D_{\mathcal{V},n-i+1}/D_{\mathcal{V},n-i}$ for $i=0,1, \ldots, n-1$.  
 Consequently, the line bundles $\tilde{\LL}_{\mathcal{V},n-i}:=F^{n-i}_{\mathcal{V}}
  \cap D_{\mathcal{V},n-i+1}$ provide a splitting of 
 $F^{\bullet}_{\mathcal{V}}$.

\begin{lemma}
The $B$-torsor $B_{\Sbt}$ on $\Sbt$ admits a $T$-subtorsor $T_{\Sbt} \subset B_{\Sbt}$ over $\Sbt_{1}$.
Furthermore, $B_{\Sbt}$ is the pushout of $T_{\Sbt}$ along $T \to B$. 
\begin{proof}
It is given by $T_{\Sbt}=[U_{B} \times B\backslash Bw_0B] \to [B \times B \backslash Bw_{0}B]$.
\end{proof}
\end{lemma}

\begin{defn}
Let $T_{\dR}:=\zeta^{*}T_{\Sbt}$ be the $T$-torsor over $U$ obtained by pullback from $T_{\Sbt}$, it is a $T$-reduction 
of $B_{\dR}$.
Let $F_{T}: \Rep_{k}(T) \to \text{Coh}(U)$ be the corresponding functor. 
\end{defn}

 We will use the splitting of $F^{\bullet}_{\mathcal{V}}$ to define a local section for $T_{\dR}$ in \Cref{adapted-element}. 
 Then with respect to such a section we can  explicitly describe the canonical isomorphism $e_{V}$.

\begin{lemma} \label{safety-lemma-TdR}
Choose $V$ and $B^{z'}_0$ as in \Cref{special-embedding}.
Any element $\phi \in T_{\dR}$ defines a trivialization of $\mathcal{V}$ such that $F^{\bullet}_{\mathcal{V}}$ maps  
to the split flag given by $B_0$ and its splitting corresponds to the splitting of $F^{\bullet}_{\mathcal{V}}$
by the line bundles $\tilde{\LL}_{\mathcal{V},i}$.
Further, if we modify $\phi$ by an element of $T$, it will act naturally via $T \to T_0$ on 
$\tilde{\LL}_{\mathcal{V},i}$. 
\begin{proof}
This can be checked at the level of Schubert stacks, where we can replace $G$ by $G^{\text{der}}$.
There it follows from the explicit description of $V$ in \Cref{special-embedding} 
as a sum of minuscule representations.
\end{proof}
\end{lemma}

It will be useful to consider not only the $G$-Zip structure on $\mathcal{V}$ for $V \in \Rep_{k}(G)$,
but the structure of
a \textit{de Rham $F$-gauge}. First recall the concept of $p$-curvature. 
\begin{lemma} \cite[Thm 5.1]{Katz-p-curvature}
Let $X/\overline{\F}_p$ be a scheme, and $(E,\nabla)$ a vector bundle with connection. The $p$-curvature map 
$\psi_{\nabla}: T^{(p)}_{X} \to \textnormal{End}_{X}(E)$ is defined by 
$\psi_{\nabla}(D \otimes 1)v=\nabla^p_{D}(v)-\nabla_{D^{[p]}}(v)$, for $v \in E$. 
\begin{enumerate}
\item For $V$ any vector bundle there is a canonical connection on the Frobenius pullback $V^{(p)}$
defined by $\nabla(v \otimes f)=v \otimes df$. We will denote it as the Frobenius connection. 
\item A vector bundle with connection $(E,\nabla)$ has $p$-curvature $0$ if and only if one can write 
$E=V^{(p)}$ for some 
vector bundle $V$, and $\nabla$ is the Frobenius connection. This happens if and only if Zariski locally there exists 
a basis of $E$ consisting of horizontal sections. 
\end{enumerate}
\end{lemma}

\begin{remark} \label{remark-Frobenius-connection}
  Under the equivalence of \Cref{flag-g-modules} and the isomorphism 
  $F_{G/B} \Ver^{0,\vee}_{B}(V) \cong F^*F_*\mathcal{V}^{\vee}$, we have that the right-hand side 
  is equipped with the Frobenius connection. A weight of the form $p\lambda$ has the structure of a 
  $(U\mathfrak{g},B)_{\Fpbar}$-module, by defining $x_{\gamma}$ for $\gamma \in \Phi$ to act trivially. This is 
  well-defined since the Cartan subalgebra $\mathfrak{h}$ acts trivially on $p\lambda$.
  In that case $F_{B}(p\lambda)=F^*F_{B}(\lambda)$ is also equipped with 
  the Frobenius connection. 
  \end{remark}

\begin{defn}(de Rham $F$-gauge)
  Let $X/\overline{\F}_p$ be a scheme. A de Rham $F$-gauge on $X$ consists of a tuple 
  $(E,\nabla,F^{\bullet}_{E},D_{E,\bullet},\phi)$ where 
  \begin{itemize}
  \item $(E,\nabla)$ is a vector bundle with a flat connection.
  \item  $F^{\bullet}_E$ is a descending filtration of $E$ satisfying 
  Griffiths transversality. 
  \item $D_{E,\bullet}$ is an ascending filtration on $E$ which is horizontal for $\nabla$, 
  and such that the associated connection on $\text{gr}^{\bullet}_{D} E$ has trivial $p$-curvature. 
  \item $\phi$ is a graded $\mathcal{O}_X$-linear isomorphism 
  $$
 (\text{gr}_{F}E,\theta)^{(p)} \cong (\text{gr}_{D} E, \psi),
  $$
  where  $\theta=\text{gr}_{F}\nabla: \text{gr}^{\bullet}_{F} E \to \text{gr}^{\bullet-1}_{F} E \otimes \Omega^1_X$
  is induced by Griffiths transversality, and 
  $\psi: \text{gr}^{\bullet}_{D} E \to \text{gr}^{\bullet-1}_{D} E \otimes \Omega^{1(p)}_X$
  is induced by $\psi_{\nabla}$ and the trivial $p$-curvature condition. 
  \end{itemize}
\end{defn}

One can also make sense of a de Rham $F$-gauge with a $G$-structure, by considering $G$-torsors with 
certain structures. Concretely, one gets the following. Recall \Cref{filtration-functoriality}.

\begin{prop} \label{conjugate-p-curvature}
Let $V \in \Rep_{k}(G)$, the cocharacter $\mu$ induces parabolics $P_{\GL(V)}$ and $P^{-(p)}_{\GL(V)}$ of $\GL(V)$,
which endow $\mathcal{V}=F_{P}(V)$ with filtrations $F^{\bullet}_{\mathcal{V}}$ and $D_{\mathcal{V},\bullet}$.
Our convention is that $F^0_{\mathcal{V}}=\mathcal{V},
F^{n}_{\mathcal{V}}=0$, 
$D_{\mathcal{V},0}=0,D_{\mathcal{V},n}=\mathcal{V}$, and the indices correspond to jumps in the filtration. 
 Then 
$(\mathcal{V},\nabla, F^{\bullet}_{\mathcal{V}},D_{\mathcal{V},\bullet})$ is naturally a de Rham $F$-gauge 
over $\Shbar$. Concretely, $D_{\mathcal{V},\bullet}$ is horizontal for $\nabla$, 
and the $p$-curvature $\psi$ is trivial on $\textnormal{gr}^{\bullet} D_{\mathcal{V}}$. 
Given $V \in \Rep_{\Fpbar}(G)$ let $V^{\sigma}$ be obtained by precomposing with $\sigma$ on $G_{\Fpbar}$.
There is a Cartier isomorphism 
$$
D_{\mathcal{V}^{\sigma},i}/D_{\mathcal{V}^{\sigma},i-1} \cong (F^{i-1}_{\mathcal{V}}/F^{i}_{\mathcal{V}})^{(p)}
$$
for $i \ge 1$, and the map induced from the $p$-curvature 
$$
\psi: \textnormal{gr}^{\bullet} D_{\mathcal{V}^{\sigma}}\to \textnormal{gr}^{\bullet-1}
 D_{\mathcal{V}^{\sigma}}  \otimes \Omega^{1,(p)}_{\Shbar}
$$
corresponds to the Frobenius twist of the Kodaira--Spencer map of
 $(\mathcal{V},\nabla,F^{\bullet}_{\mathcal{V}})$ 
under the Cartier isomorphism. Moreover, if we choose a pair of Borels $(B_0,B^{z'}_{0})$ of $\GL(V)$ 
such that $B \subseteq B_0 \subseteq P_{V}$, $B^{z} \subseteq  B^{z'}_{0} \subseteq P^{-(p)}_{V}$,
the 
extended conjugate filtration on $\mathcal{V}_{\flag}$ is horizontal for $\nabla$ and has trivial $p$-curvature 
on its graded pieces. 
\begin{proof}
The statement that $(\mathcal{V},\nabla, F^{\bullet}_{\mathcal{V}},D_{\mathcal{V},\bullet})$ is a de Rham F-gauge 
for $V=\Lambda$ is a classical result of Katz \cite[Thm 3.2]{Katz-p-curvature},
where he moreover proves that $\nabla$ on the graded pieces 
of the conjugate filtration $\Ker V \subseteq \mathcal{H}$ 
is identified with the Frobenius connection under the Cartier 
isomorphism. In general this is explained in \cite[\S 4.2]{de-Rham-F-gauge}.
The Cartier isomorphism comes from unpacking the isomorphism  
$P^{(p)}_\dR/U^{(p)}_{P} \cong P^{-}_{\dR}/U_{P^{-(p)}}$, since on the left-hand side it acts via 
$\sigma: M^{(p)} \to M$, 
which after pushing out to $M^{(p)}_{V}$ corresponds to considering the graded pieces of $\mathcal{\GL(V)}^{\sigma}$. 
 For the statement about the full conjugate filtration, we just observe that since
 $B^{z} \subseteq B^{z'}_{0} \subseteq P^{-(p)}_{V}$, it is constructed from 
 the full Hodge filtration via the Cartier isomorphism, so that each term in
  $\text{gr}^{\bullet} D_{\mathcal{V}}$ is 
 a Frobenius twist coming from $F^{\bullet}_{\mathcal{V}^{\sigma^{-1}}}$,
 with $\nabla$ being the Frobenius connection. This can be checked from the case $V=\Lambda$, 
 by functoriality of $\nabla$. 
\end{proof} 
\end{prop}

Using this result we can define a special kind of elements in $T_{\dR}$ over $U$, which we will use 
to describe $e_{V}$ over $U$.
\begin{lemma}(Adapted elements) \label{adapted-element}
We work over $\flag_{\overline{\F}_p}$. 
Choose $V, B_0, B^{z'}_0$ as in \Cref{special-embedding}. 
Choose some local section $\phi \in T_{\dR}$. Via the isomorphisms 
$H_i: F^{n-i}_{\mathcal{V}} \cap D_{\mathcal{V},n-i+1} \cong \LL'_{\mathcal{V},n-i}=
\textnormal{gr}^{n-i}D_{\mathcal{V},\bullet}$
 this produces
a local basis $\{v_{\bullet}\}$ of $\LL'_{\mathcal{V},\bullet}$. We can modify $\phi$ by an element of $T$ to 
ensure that
$\nabla_{\LL^{'}_{V,\bullet}}(v_{\bullet})=0$. 
We say that such a $\phi$ is an adapted element of $T_{\dR}$. 
\begin{proof}
Using \Cref{conjugate-p-curvature} and \Cref{safety-lemma-TdR} we reduce to the following claim. 
Given $x \in T \subseteq T_0$ there exists an element $t$ of $T$ 
such that all the coordinates of $tx$ in $T_0$ are $p$th powers. 
After base changing to $\overline{\F}_p$ we can assume that $T$ is split,
so that $\mathbb{G}^k_m \cong T \to T_0\cong \mathbb{G}^{n}_m$
is given by a tuple of cocharacters $\mu_i$. Then if 
$x=(x_i)$ we can take $t=(x^{p-1}_i)$, so that each $\mu_i(x^p_i)$ is a $p$th power.  

\end{proof}
\end{lemma}
Let 
$\phi \in T_{\dR}$ be an adapted element with respect to the data of \Cref{special-embedding}.
For $\gamma \in \Phi^{+}$ let $\{\omega_{\gamma} \in \Omega^1_{U}\}$ be the dual basis to
$\{e^{1}_1(\phi,x_{-\gamma}) \in T_U\}$ from \Cref{basic-iso}. 
We can describe $\nabla_{\mathcal{V}}$ succinctly in this basis. 
\begin{prop} \label{nabla-on-U}
Choose $V$ as in \Cref{special-embedding}, and choose an adapted element $\phi 
\in T_{\dR}$ with respect to it. Then for all $v \in V$ the connection $\nabla$ on $\mathcal{V}_{U}$ acts as
$$
\nabla(\phi,v)=\sum_{\gamma \in \Phi^{+}} (\phi,x_{-\gamma}v) \otimes \omega_{\gamma}. 
$$
Moreover, let $\phi \in T_{\dR}$ be an adapted element with respect to $V=\Lambda$
and some $B_0$ in \Cref{special-embedding}. Then the identity above holds for all $V \in \Rep_k(G)$ 
and $v \in V$. 
\begin{proof}
The Borel $B_0 \subseteq \GL(V)$ defines full descending 
filtration $F^{\bullet}_{V}$ of length $n$.
Let  $\{v_i: i=1,\ldots, n\}$ be a basis for $V$ such that $v_{i}$ generates 
$F^{n-i}_{V}/F^{n-i+1}_{V}$, and let $\pi_i: \mathcal{V} \to \mathcal{O}_{\flag}$ be 
the projection induced by $v_i$ and $\phi$.
 Then $(\phi,v_i)$ generates $F^{n-i}_{\mathcal{V}} \cap D_{\mathcal{V},n-i+1}$, and by construction
of $\phi$ it satisfies
$\nabla(\phi,v_i) \in \langle (\phi,v_{i+1}), \ldots, (\phi,v_n) \rangle \otimes \Omega^1_{U}$, 
where the convention is that $\nabla(\phi,v_n)=0$. Let 
$\gamma \in \Phi^{+}$, and write $x_{-\gamma}v_i=cv_{j}+y$, 
where $y$ is in the span of the rest of the basis vectors, 
for some $i <j$ and a non-zero constant $c \in \Fpbar$. 
Then the canonical map $e^1_1$ fits in the commutative diagram
$$
\begin{tikzcd}
  & F_{B}(\mathfrak{g}/\mathfrak{b}) \arrow[d,hook] \\
T_{\flag} \arrow[ur,"e^1_1"] \arrow[r] & F_{B_0}(\mathfrak{gl}(V)/\mathfrak{b}_0)
\end{tikzcd}
$$
by \Cref{KS-flag}, since $\Lie(G^{\der}) \to \Lie(\SL(V))$ is an embedding. 
Therefore, we have that 
$c\omega_{\gamma}=\pi_{j}\nabla(\phi, v_i)$. Conversely, 
if  $\pi_{j}\nabla(\phi, v_i) \neq 0 $ for some $i <j$
it must be that there exists $x \in \mathfrak{u}^{-}_{B}$
such that $v_{j}$ appears non-trivially in basis expansion of $x v_i$, this follows from the diagram above.
By writing $x$ as a sum of negative roots 
we see that there exists a $\gamma \in \Phi^{+}$ such that $v_j$ appears non-trivially 
on $x_{-\gamma}v_i$. Moreover, such a $\gamma$ is unique up to a scalar,
 since otherwise $e^1_1$ wouldn't be injective. 
 This reasoning shows that the formula holds for $V=\Lambda$.
  For a general $V \in \Rep_{k}(G)$ it follows from the previous case, 
since $\nabla_{\mathcal{V}}$ is defined by transporting $\nabla_{\mathcal{H}}$ along the action of $\mathfrak{g}$, in
\Cref{GM-definition}.
\end{proof}
\end{prop}

Furthermore, relative to an adapted element we can describe the 
canonical isomorphisms $e_V$ in a simple way. 

\begin{prop} \label{canonical-iso-U}
Let $V \in \Rep_{k}(B)$. Let $\phi \in T_{\dR}$ be an adapted element
with respect to $V=\Lambda$ in \Cref{special-embedding}. 
 For $\gamma \in \Phi^{+}$ let  $D_{\gamma}=e^{1}_{1}(\phi,x_{-\gamma}) \in T_U$. 
 Then $e_V: F_B \Ver_B(V) \cong D_{\flag} \otimes F_B(V)$ restricted to $U$ is given by 
$$
(\phi,\prod_{\gamma \in \Phi^{+}} x_{-\gamma} \otimes v) \mapsto \prod_{\gamma} D_{\gamma} \otimes (\phi,v),
$$
where in the right-hand side the product means a composition of differential operators, in the same order as the 
left-hand side. 
\begin{proof}
We can reduce to the case of $e^{\le 1}_{V}$, by the compatibility between composition of differential operators and
the maps $\Ver_B \Ver_B(V) \twoheadrightarrow \Ver_{B}(V)$. By the definition of $e^{\le 1}_{V}$ in \Cref{def-less-1}
we can reduce to the case when $V \in \Rep(G)$. 
Let $\{v_i\}$ be a basis for $V$. The dual of the tensor identity $\phi_{V}: \Ver^{\le 1,\vee}_{B}(V)
\cong V^{\vee} \otimes \Ver^{\le 1,\vee}_{B}(1)$ is given by $(x_{-\gamma} \otimes v_i)^{\vee} \mapsto v^{\vee}_i 
\otimes x^{\vee}_{-\gamma}$ and 
$$
(1 \otimes v_i)^{\vee} \mapsto v^{\vee}_i \otimes 1 + 
\sum_{\gamma \in \Phi^{+}} x_{-\gamma}v^{\vee}_i \otimes x^{\vee}_{-\gamma}. 
$$
By the definition of $e^{\le 1}_V$ we see that the proposition is then equivalent to \Cref{nabla-on-U}.
\end{proof} 
\end{prop}

Using \Cref{restriction-creator-diff} we get a functor  
$$
\Phi_{T}: \Hom_{T}(\lambda_1, \Ver_B(\lambda_2))  \to \DiffOp_{U}(\LL(-\lambda_2), \LL(-\lambda_1)).
$$
As opposed to the case of  maps of $B$-modules, it is not completely clear how to define composition on 
the left-hand side. We show that the naive way is compatible with composition on the right-hand side. 
Elements  $f \in \Hom_{T}(\lambda_1, \Ver_B(\lambda_2))$ can be described uniquely 
by an element  $x_{f}$ of $U\mathfrak{u}^{-}_{B}$,
by the PBW theorem.

\begin{theorem} \label{composition-U}The map 
$$
\Phi_{T}: \Hom_{T}(\lambda_1, \Ver_B(\lambda_2))  \to \DiffOp_{U}(\LL(-\lambda_2), \LL(-\lambda_1))
$$
is compatible with composition on both sides, where the composition on the left-hand side is defined as follows. 
Let $f \in \Hom_{T}(\lambda_1, \Ver_B(\lambda_2))$ and $g \in \Hom_{T}(\lambda_2, \Ver_B(\lambda_3))$, 
then $g \circ f$ is defined by $x_f x_g \in U \mathfrak{u}^{-}_B$.
\begin{proof}
By the proof of \Cref{creator-diff}(2) 
it is enough to prove that for a $T$-equivariant map
$f: \lambda_1 \to \Ver_{B}(\lambda_2)$
the map $h:\Ver_{B}(\lambda_1) \to \Ver_{B}(\Ver_{B}(\lambda_2))$ defined by $\prod x \otimes \lambda_1 \mapsto 
\prod x \otimes x_f \otimes \lambda_2$ is identified with $ D_{U} \otimes \LL(\lambda_1) \xrightarrow{\text{id}
\otimes F(f)} D_{U} \otimes D_{U} \otimes \LL(\lambda_2)$ after applying $F_T$ and the canonical isomorphisms from 
\Cref{canonical-iso-Pm}. Once we write both $x$ and $x_{f}$ in a PBW basis this follows immediately 
from \Cref{canonical-iso-U}.
\end{proof}
\end{theorem}

\section{The basic theta operators} \label{section3}

We construct the basic theta operators, and we prove several properties about them. We start by recalling 
the construction of Hasse invariants on $\flag$.

\subsection{The Hasse invariants}
Recall the map $\zeta: \flag \to \Sbt=[B \times B\backslash G]$ from \Cref{section2.3}. For $\alpha \in \Delta$ denote $\Sbt_{\alpha}$
the locally closed stratum of $\Sbt$ corresponding to the orbit $Bs_{\alpha} w_0B$, and let 
$\overline{\Sbt}_{\alpha}$ be its closure. Let $\Sbt_{1}$ be the open stratum corresponding to $B w_0 B$. 
Let $\overline{D}_{\alpha}=\zeta^{-1}\overline{\Sbt}_{\alpha}$, and $D_{\alpha}=\zeta^{-1} \Sbt_{\alpha}$.
The complement of $\Sbt_1$ is covered by the union of $\overline{\Sbt}_{\alpha}$ for $\alpha \in \Delta$, 
and $\Sbt \setminus \cup_{\alpha} \Sbt_{\alpha}$ has codimension at least $2$ in $\Sbt$. For a pair of weights 
$\lambda, \mu \in X^{*}(T)$ let $\LL(\lambda,\mu)$ be the line bundle of $\Sbt$ coming from $(\lambda,\mu)$
as a character of $B \times B$. If $\mu=0$ we will also write $\LL(\lambda)$.

 We define the Hasse invariants on $\flag$ as in \cite[\S 5.1]{IK-hasses}.
\begin{lemma} \label{Hasse-dfn}
Let $\alpha \in \Delta$. For $\chi_1,\chi_2 \in X^*(T)$ we have that $\H^0(\Sbt,\LL(\chi_1,\chi_2)) \neq 0$  
implies $\chi_2=-w_0\chi_1$, in which case the space is $1$-dimensional. There exists a character  
$\chi_{\alpha} \in X^{*}(T)$ such that 
$\H^0(\Sbt,\LL(-\chi_{\alpha},w_0\chi_{\alpha}))$ is $1$-dimensional, and the divisor of any non-zero section 
$\tilde{h}_{\alpha}$ is 
$$
\textnormal{div}(\tilde{h}_\alpha)=N_{\alpha}\overline{\Sbt}_{\alpha},
$$
with $N_{\alpha} \ge 1$ minimal. Moreover, the projection of $\chi_{\alpha}$ to $X^{*}(T^{\textnormal{der}})$ is unique.
More precisely, $N_{\alpha}$ is the minimal positive integer such that there exists $\chi_{\alpha} \in X^*(T)$
such that $\langle \chi_{\alpha}, \alpha^{\vee} \rangle=N_{\alpha}$ and $\langle \chi_{\alpha}, \beta^{\vee} \rangle=0$ 
for all $\beta \in \Delta \setminus \{\alpha\}$.
In particular, if $G^{\textnormal{der}}$ is simply connected then $N_{\alpha}=1$.
\begin{proof}
The first statement is \cite[Thm 2.2.1(a)]{Goldring-Koskivirta-Galois}. 
According to Chevalley's formula \cite[Thm 2.2.1(c)]{Goldring-Koskivirta-Galois}
the divisor of $\tilde{h}_{\alpha}$ is $
\sum_{\beta \in \Delta}\langle \chi_{\alpha}, \beta^\vee \rangle\overline{\Sbt}_{\beta}$.
If $G^{\text{der}}$ is simply connected we can find a $\chi_{\alpha}$ such that 
$\langle \chi_{\alpha},\beta^{\vee} \rangle$ is $1$ if $\alpha=\beta$ and zero otherwise,
and such a character is unique in $X^{*}(T^{\text{der}})$. In general, by considering the index of 
$X^*(T)$ in the weight lattice of the simply connected cover, 
there will exist a $\chi_{\alpha}$ such that $\langle \chi_{\alpha},\beta^{\vee} \rangle$
is zero for $\alpha \neq \beta$, 
and positive for $\alpha=\beta$. 
\end{proof}
\end{lemma}
We denote by $H_{\alpha}$ the pullback of some choice of $\tilde{h}_{\alpha}$ to $\flag$, then 
$N_{\alpha}$ in \Cref{Hasse-dfn} is the order of vanishing of $H_{\alpha}$. 
 The weight of 
$H_{\alpha}$ as a section of $\LL(h_{\alpha})$ is 
\begin{equation} \label{Hasse-weight}
h_{\alpha}=
-\chi_{\alpha}+pw_{0,M} \sigma^{-1}\chi_{\alpha},
\end{equation}
where 
$\sigma$ is the arithmetic Frobenius acting on $X^*(T)$. This can be computed by remembering the twist by $z \in W$
in $\flag \to \Sbt$, 
and using the formula $z=\sigma(w_{0,M})w_{0}$ from \cite[Lem 2.2.1]{IK-hasses}. 
Note that we get a different sign compared to
\cite[\S 5.1]{IK-hasses},
since 
we use different conventions for the simple roots. Also, for applications, in characteristic $0$, if $\lambda, \mu \in X^{*}(T)$
are such that their projection to $X^{*}(T^{\text{der}})$ is the same, then their associated 
line bundles $\LL(\mu)$
and $\LL(\lambda)$ on $\flag$ are isomorphic, since they are pulled back from $G/B=G^{\text{der}}/B^{\text{der}}$
via the Borel embedding. Thus, we don't lose much from $\tilde{h}_{\alpha}$ not being unique. We also remark 
that for $\alpha \in \Delta_{G_i}$ 
we can compute $h_{\alpha}$ using the weight lattice of the almost-simple quotient $G^i$.

\begin{lemma} \label{Hasse-functorial}
Consider $(V,T_0,B_0)$ as in \Cref{special-embedding}. 
Consider the map $r: \Sbt \to \Sbt_{\GL(V)}=[B_0 \times B_0 \backslash \GL(V)]$. For $\alpha_i \in \Delta_{\GL}$ one can similarly define 
Hasse invariants $\tilde{h}_{\alpha_i}$ on $\Sbt_{\GL}$. 
Then for $\alpha \in \Delta$, 
$$
\tilde{h}_{\alpha}=\prod_{\beta \in \Delta_{\GL}} r^* \tilde{h}^{n_{\beta}}_{\beta}
$$ 
for some integers $n_{\beta}$. 
\begin{proof}
Let $\LL_{\GL}(\lambda,\mu)$ be the line bundles on $\Sbt_{\GL(V)}$ given by characters of $B_0 \times B_0$. 
Then we can use $\chi_{\alpha_i}=(1,\ldots,1,0, \ldots, 0)$, where there are $i$ ones, to construct 
$\tilde{h}_{\alpha_i}$ in \Cref{Hasse-dfn}. We have 
$r^* \LL_{\GL}(\lambda,\mu)=\LL(\lambda,\mu)$ under the map $X^{*}(T_0) \to X^{*}(T)$.
The pullback map $r^*$ is injective on global 
sections, since $B w_0 B \subseteq B_0 w_{0,\GL} B_0$ by \Cref{special-embedding} (after the 
appropriate twisting). Therefore,
by the first part of \Cref{Hasse-dfn} 
$r^* \LL_{\GL}(\chi,-w_{\GL,0}\chi)=\LL(\chi,-w_{0}\chi)$ whenever the former admits a non-zero section.
Thus, we reduce to expressing $\chi_{\alpha}$ 
as a linear combination of the restrictions of $\chi_{\alpha_i}$. We conclude since the $\chi_{\alpha_i}$
span $X^*(T_0)$, and the map $X^{*}(T^{\der}_0) \to X^{*}(T^{\der})$ is surjective. This is because $V$ is a
minuscule representation, which have the property that their weights span $X^*(T^{\der})$.
\end{proof}
\end{lemma}

That is, we can construct the Hasse invariants out of Hasse invariants for general linear groups. 
Note that we are choosing $\tilde{h}_{\alpha_i}$ so that its pullback
 $H_{\alpha_i}$ on $\flag$ is $\det(F^{n-i}_{\mathcal{V}} \to \mathcal{V}/D_{\mathcal{V},n-i})$, where $\alpha_i$ 
are the simple roots of $\GL(V)$.

\subsection{Construction}
Every theta operator in $\flag$ 
will be constructed via pullback along $\zeta: \flag \to \Sbt$. For $V \in \Rep(B \times B)$ let $\underline{V}$
be the corresponding vector bundle on $\Sbt$. For $W \in \Rep(B)$ we will use $\underline{W}$ when considering $W$ as a 
representation of $B \times B$ through the first factor. Recall from \Cref{section2.3} that we have an exact functor 
$F_{T}: \Rep_{\Fpbar}(T) \to \text{Coh}(\Sbt_1) \to \text{Coh}(U)$
fitting in the diagram 
$$
\begin{tikzcd}
\Rep_{\Fpbar}(B) \arrow[r,"F_{B}"] \arrow[d] & \textnormal{Coh}(\flag) \arrow[d,"\text{res}"] \\
\Rep_{\Fpbar}(T) \arrow[r,"F_{T}"] & \text{Coh}(U).
\end{tikzcd}
$$
Therefore, by \Cref{restriction-creator-diff} we get a functor
\begin{equation} \label{diff-U}
\Hom_{T}(\lambda,\Ver_{B}(\mu)) \to \DiffOp_{U}(\LL(-\mu),\LL(-\lambda)). 
\end{equation}

\begin{defn}
Let $\gamma \in \Phi^{+}$, and $\lambda \in X^*(T)$. Consider the map of 
$T$-modules $\phi_{\gamma,\lambda}: -\lambda-\gamma \to \Ver_{B}(-\lambda)$ 
given by $-\lambda-\gamma \to 
x_{-\gamma} \otimes -\lambda$. We define $\tilde{\theta}_{\gamma}: \LL(\lambda) \to \LL(\lambda+\gamma)$
to be the differential operator on $U \subset \flag$ induced by $\phi_{\gamma,\lambda}$ and the functor 
\eqref{diff-U} above. 
\end{defn}

Thus, the $\tilde{\theta}_{\gamma}$ are  
precisely all the degree $1$ differential operators coming from
the functor above. We show how one can extend their corresponding maps on $\Sbt_{1}$ to $\Sbt$.

\begin{prop} \label{power-Hasses}
Let $\gamma \in \Phi^{+}$, $\alpha \in \Delta$ and $\lambda \in X^{*}(T)$. Let $n=n_{\gamma,\alpha,\lambda} \ge 0$
be the smallest integer such that the map of vector bundles on $\Sbt_{1}$
$$
\LL(-\lambda-\gamma) \otimes \LL(n\chi_{\alpha},-n w_0\chi_{\alpha})
\xrightarrow{\tilde{h}^n_{\alpha}} \LL(-\lambda-\gamma) \xrightarrow{\phi_{\gamma,\lambda}} 
\underline{\Ver^{\le 1}_B(-\lambda)}
$$
extends to $\Sbt_{1} \cup \Sbt_{\alpha}$. 
Let $n_{\gamma,\alpha}$ be the maximum of $n_{\gamma,\alpha,\lambda}$ along all $\lambda$.
Then 
\begin{enumerate}
\item For $\gamma \in \Delta$, $n_{\gamma, \alpha}$ is zero unless $\alpha=\gamma$, in which case it is $1$.
\item Let $p \ge 3$, and assume that $\tilde{h}_{\alpha}$ has a simple zero.
 Then $n_{\gamma,\alpha}$ is the largest integer $n$ such that 
 $\gamma-n\alpha \in \Phi^{+} \cup \{0\}$. 
\end{enumerate}
\begin{proof} 
Recall that we fix $\{x_{\gamma}\}$ to sit in a Chevalley basis. 
We follow \cite[\S 4.2]{IK-hasses}. 
We can consider $\phi_{\gamma,\lambda}$ as a section of $\underline{V}$ over $\Sbt_{1}$ for 
$V:=\Ver^{\le 1}_B(-\lambda) \otimes \LL(\lambda+\gamma) \in \Rep(B \times B)$ acting through its first factor.
Equivalently, 
as the left and right $B$-equivariant rational map $\psi: G \mapsto V$ given 
by $w_0 \mapsto v:=[x_{-\gamma} \otimes -\lambda] \otimes (\lambda+\gamma)$.
Then $n_{\gamma,\alpha,\lambda}$ can be identified 
with the pole (the highest non-trivial power of $t^{-1}$) of the element 
$F_{\alpha}(t)=u_{\alpha}(t^{-1})\alpha^{\vee}(t^{-1})v=u_{\alpha}(t^{-1}) v \in V \otimes \F_p[t,t^{-1}]$, 
since whenever $F_{\alpha}(0)$ is defined, it gives the value of $\psi$ at $s_{\alpha}w_0$. Note that it is in this identification where 
we use the assumption that $\tilde{h}_{\alpha}$ has a simple zero.
We get a slightly different formula for $F_{\alpha}(t)$ compared to \cite[Lem 4.2.1]{IK-hasses} because of our 
conventions on simple roots. 
Now, $u_{\alpha}(t^{-1})=\phi_{\alpha}\begin{pmatrix}
1 & t^{-1} \\
0 & 1 
\end{pmatrix}$ acts the formal exponential $\text{exp}(t^{-1}x_{\alpha}):=\sum_{n} \frac{1}{n!}t^{-n} x^n_{\alpha}$. 
This follows from the equivalence between $\Rep_{\Fpbar}(\SL_2)$
and $U(\SL_2)_{\Fpbar}$-modules, and using that
 $du_{\alpha}(1)=x_{\alpha}$ by definition.
Threfore, we reduce to computing 
the pole of 
$$
\sum_{n=1} \frac{1}{n!}t^{-n} x^n_{\alpha}x_{-\gamma} \otimes -\lambda \in \Ver_{B}(-\lambda)[t^{\pm1}].
$$
If $\alpha \neq \gamma$ we have that $[x_{\alpha},x_{-\gamma}]=c x_{-\gamma+\alpha}$
if $\gamma-\alpha \in \Phi^{+}$ and otherwise it is zero.
Here $c$ is the largest positive integer such that $(c-1)\alpha+\gamma \in \Phi$, 
and we can check that $c$ is at most $2$ by going through all the Dynkin diagrams.
On the other hand, by our choice of Chevalley basis we have that 
$[x_{\alpha},x_{-\alpha}]\lambda=\langle \lambda, \alpha^{\vee} \rangle \lambda$. 
It is then clear that, except in the case that $p=2$ and $c=2$, or $p \mid \langle \lambda, \alpha^{\vee} \rangle$
one has $n_{\gamma,\alpha,\lambda}=n_{\gamma,\alpha}$ and point $(2)$ holds. 
\end{proof}
\end{prop}

Since the map 
\begin{equation} \label{extension-theta}
\LL(-\lambda-\gamma) \otimes \LL(\sum_{\alpha} n_{\alpha,\gamma} \chi_{\alpha},
\sum_{\alpha} -n_{\alpha,\gamma}w_{0}\chi_{\alpha})
\xrightarrow{\prod \tilde{h}^{n_{\alpha,\gamma}}_{\alpha}} \LL(-\lambda-\gamma) \xrightarrow{\phi_{\gamma,\lambda}} 
\underline{\Ver^{\le 1}_B(-\lambda)}
\end{equation}
is defined over $\Sbt_{1} \cup_{\alpha} \Sbt_{\alpha}$, it extends to $\Sbt$,  by Hartogs lemma.  

\begin{defn} \label{defn-basic-theta}
For $\gamma \in \Phi^{+}$ via pullback by $\zeta: \flag \to \Sbt$ the dual of the map 
\eqref{extension-theta}
induces a map $F_{B}\Ver^{\le 1, \vee}(-\lambda) \to \LL(\lambda+\gamma+\sum h^{n_{\alpha,\gamma}}_{\alpha})$,
where $h_{\alpha}$ is the weight of $H_{\alpha}$.
Let $\mu_{\gamma}=\gamma+\sum n_{\alpha,\gamma}h_{\alpha}$. 
Define $\theta_{\gamma}$ to be the differential operator $\LL(\lambda) \to 
\LL(\lambda+\mu_{\gamma})$ on $\flag$ obtained from the previous map via \Cref{canonical-iso-Pm}.
\end{defn}

We make each $n_{\gamma,\alpha}$ explicit depending on the Lie algebra of $G$. We also give some 
examples in which a Hasse invariant can have a double zero.  
\begin{prop} \label{table-lie}
Let $G^i$ be the almost-simple quotients of $G^{\der}_{\Fpbar}$, 
 and let $\mathfrak{g}_i$ be
its Lie algebra. Let $\gamma \in \Phi^{+}$ lying in some $\mathfrak{g}_i$.
Assume that $\tilde{h}_{\alpha}$ has a simple zero and that $p\ge 3$. 
We give a list of all $n_{\alpha,\gamma}$,
according to the Dynkin diagram of $\mathfrak{g}_i$. We use \Cref{Hodge-embeddings} to determine what $G^i$ can be. 
\begin{enumerate}
\item ($A_{n}$) Assume that $G^i \cong \SL_{n+1}$. Its simple roots are $\alpha_1, \ldots, \alpha_n \in \mathbb{R}^{n+1}$,
given by $\alpha_{i}=e_i-e_{i+1}$, where $e_{i}$ are the basis vectors. The weight space of $G^i$
is given by the hyperplane where the sum of all coordinates is $0$.
The positive roots are of the form 
$\gamma=\alpha_i + \alpha_{i+1} + \ldots + \alpha_{j}$ for pairs $1 \le i < j \le n$. Then 
$n_{\gamma,\alpha_k}$ is $1$ if $k \in \{i,j\}$ and zero otherwise. 
\item ($C_n$) Assume that $G^i=\textnormal{Sp}_{2n}$.
The simple roots are $\alpha_{i}=e_{i}-e_{i+1}$ for $i=1,\ldots, n-1$, and $\alpha_n=2e_n$.
Then positive roots are of the form 
\begin{itemize}
\item $\gamma=e_{i}- e_j=\alpha_{i}+\ldots+\alpha_{j-1}$ for $i <j \le n$. In this case $n_{\gamma,\alpha_k}=1$ if 
$k \in \{i,j-1\}$ and otherwise it is $0$.  
\item $\gamma=2e_i=2(\alpha_{i}+\ldots+\alpha_{n-1})+\alpha_n$ for $i=1,\ldots,n$. 
In this case $n_{\gamma,\alpha_k}=2$ if $k=i<n$, it is $1$ if $k=i=n$, and otherwise it is zero. 
\item $\gamma=e_i+e_j=\alpha_{i}+\ldots+\alpha_{j-1}+2(\alpha_{j}+\ldots+\alpha_{n-1})+\alpha_n$ for $i<j$. 
In this case $n_{\gamma,\alpha_k}=1$ if $k \in \{i,j\}$, and otherwise it is zero.
\end{itemize} 
\item $(B_n)$ The positive roots are 
$\alpha_i=e_i-e_{i+1}$ for $1 \le i \le n-1$ and $\alpha_n=e_n$ inside $\mathbb{R}^n$, the coroots are 
$\alpha^{\vee}_i=e^{\vee}_i-e^{\vee}_{i+1}$ for $1 \le i \le n-1$ and $\alpha^{\vee}_n=2e^{\vee}_n$. 
If $G^{i}=\textnormal{Spin}_{2n+1}$, the weight lattice is spanned by 
$\{e_1, \ldots, e_{n-1},\frac{1}{2}(e_1 + \ldots + e_n)\}$, and the coweight lattice by the coroots. 
The positive roots are 
\begin{itemize}
\item $\gamma=\alpha_i+\alpha_{i+1}+\ldots+\alpha_n$ for $1 \le i \le n$. Here $n_{\gamma,\alpha_k}=1$ if 
$k \in \{i,n\}$ and zero otherwise. 
\item $\gamma=(\alpha_i+ \ldots +\alpha_n)+(\alpha_j+\ldots+\alpha_n)$ for $1\le i < j \le n$. If $j>i+1$ 
$n_{\gamma,\alpha_k}=1$ for $k \in \{i,j\}$ and zero otherwise. If $j=i+1$, 
$n_{\gamma,\alpha_k}=1$ for $k=i+1$ and zero otherwise. 
\item $\gamma=\alpha_i+\alpha_{i+1}+\ldots+\alpha_{j-1}$ for $1 \le i <j \le n$. Here 
$n_{\gamma,\alpha_k}=1$ if $k \in \{i,j-1\}$ and zero otherwise. 
\end{itemize}
\item ($D_n$) We can assume $n \ge 4$. The roots are $\alpha_i=e_{i}-e_{i+1}$ for $1 \le i \le n-1$ 
and $\alpha_n=e_{n-1}+e_{n}$. Suppose first that $G^i=\textnormal{Spin}_{2n}$. 
The weight lattice  is spanned by 
the $e_i$ for $1 \le i \le n-2$, $\frac{1}{2}(e_1+e_2+\ldots+e_{n-1}-e_n)$ and 
$\frac{1}{2}(e_1+e_2+\ldots+e_{n-1}+e_n)$
in  $\mathbb{R}^n$.
The cocharacter lattice is spanned by $\{\alpha^{\vee}_i\}$ under the natural pairing.
The coroots have the same coordinates under this pairing.  
If $G^i=\textnormal{SO}_{2n}$ the weight and coweight lattice are spanned
 by $\{e_i\}$ and $\{e^{\vee}_i\}$ instead.
For $i \le n-2$ $H_{\alpha_i}$ has a simple zero, with $\chi_{\alpha_i}=e_1+ \ldots +e_{i}$. For $i=n-1,n$
$H_{\alpha_i}$ has a double zero, with $\chi_{\alpha_{n-1}}=e_1+\ldots+e_{n-1}-e_n$ and 
$\chi_{\alpha_n}=e_1 + \ldots + e_n$. 
The positive roots are of the form
\begin{itemize}
\item $\gamma=e_i-e_j=\alpha_{i}+ \ldots + \alpha_{j-1}$, $1 \le i<j \le n$. Then $n_{\gamma,\alpha_k}=1$
if $k \in \{i,j-1\}$ and zero otherwise. 
\item $\gamma=e_i+e_j=\alpha_{i}+ \ldots + \alpha_{j-1}+
2(\alpha_j + \ldots +\alpha_{n-2})+\alpha_{n-1}+\alpha_n$, $1 \le i<j \le n-2$.
Then $n_{\gamma,\alpha_k}=1$ for $k=i<j-1$, $k=j$, and zero otherwise. 
\item $\gamma=e_i+e_{n-1}=\alpha_{i}+ \ldots + \alpha_{n}$, $i<n-1$. Then $n_{\gamma,\alpha_k}=1$
for $k=i< n-2$, $k=n$, $k=n-1$, and zero otherwise. 
\item $\gamma=e_i+e_n=\alpha_{i}+ \ldots + \alpha_{n-2}+\alpha_{n}$, $i<n-1$. Then $n_{\gamma,\alpha_k}=1$
for $k\in \{i,n\}$.
\item $\gamma=e_{n-1}+e_n=\alpha_n$. Then $n_{\gamma,\alpha_k}=1$ if $k=n$ and zero otherwise. 
\end{itemize}
\end{enumerate}
\end{prop}

In general if $\tilde{h}_{\alpha}$ does not have a simple zero the result above gives an upper bound 
for the $n_{\gamma,\alpha}$. We see that this upper bound might not be attained only in case that 
$G^i$ is not simply connected of $C_n$ type.

We can also relate the construction of $\tilde{\theta}_{\gamma}$ to the more classical language involving 
the Gauss--Manin connection. Recall the exact functor $F_{T}: \Rep_{\Fpbar}(T) \to \text{Coh}(U)$. 
\begin{prop} \label{easy-lemma-theta}
\begin{enumerate}
\item The map $\tilde{\theta}_{\gamma}: \LL(\lambda) \to \LL(\lambda+\gamma)$  over $U$
can also be described by the following composition. Embed $i: \lambda \hookrightarrow V$ as $T$-representations 
for some $V \in \Rep(G)$, and let $g: V \to \lambda$ be a $T$-equivariant section.
This can be done by taking $V=V(w\lambda)$
where $w\lambda$ is a dominant Weyl translate of $\lambda$. Consider the projection 
$\pi_{\gamma}: \Omega^1_{U} \to \LL(\gamma)$ given by applying $F_{T}$ to the dual of
$-\gamma \hookrightarrow \mathfrak{g}/\mathfrak{b}$. Then 
$$
\tilde{\theta}_{\gamma}: \LL(\lambda) \xrightarrow{F_{T}(i)} \mathcal{V} \xrightarrow{\nabla} \mathcal{V} \otimes \Omega^1_{U}
\xrightarrow{F_{T}(g) \otimes \pi_{\gamma}} \LL(\lambda+\gamma),
$$
in particular the composition is independent of the choice of $i$ and $g$.
\item The basic theta operators $\theta_{\gamma}$ satisfy the Leibniz rule. That is, for 
$\lambda_1, \lambda_2 \in X^*(T)$, and $f \in \LL(\lambda_1), g \in \LL(\lambda_2)$ local sections, 
$\theta_{\gamma}(fg)=f \theta_{\gamma}(g) + \theta_{\gamma}(f) g$. 
\end{enumerate}
\begin{proof}
Consider the $T$-equivariant map $\Ver^{\le 1, \vee}(-\lambda) \to \lambda+\gamma$ sending 
$(x_{-\gamma} \otimes v_{-\lambda})^{\vee}$ to $v_{\lambda+\gamma}$, and the rest of the basis elements to $0$.
After applying $F_{T}$ it defines $\tilde{\theta}_{\gamma}$.
Applying $F_{T}$ to the map 
$$
\Ver^{\le 1, \vee}(-\lambda) \xrightarrow{i^{\vee}} \Ver^{\le 1, \vee}(V^{\vee}) \xrightarrow{\phi^{\vee}_V}
V \otimes \Ver^{\le 1,\vee}(1) \xrightarrow{g \otimes \pi} \lambda \otimes (\mathfrak{g}/\mathfrak{b})^{\vee}
 \xrightarrow{\text{id} \otimes \pi_{\gamma}} \lambda+\gamma
$$
yields the (linearized) composition in the statement,
where $\phi_V$ is the tensor identity on \Cref{flag-g-modules}.
Let $v=i(v_{\lambda})$. Using that $\text{id} \otimes \pi \circ \phi^{\vee}_{V}(x_{-\beta} \otimes v^{\vee})^{\vee}=
v \otimes x^{\vee}_{-\beta}$
we obtain the equality in $(1)$. For $(2)$ it is enough to prove it for 
$\tilde{\theta}_{\gamma}$ over $U$, and then by $(1)$ it reduces to proving the statement that the Gauss--Manin connection 
satisfies the Leibniz rule for tensor products, which can be seen from its construction or from examining 
the tensor identity $\phi_V$.
\end{proof}
\end{prop}

\begin{remark}
\begin{enumerate}
\item
This also agrees with the definition given in \cite[\S 2]{paper} in the Siegel case, or in general 
in the case where the ordinary locus is non-empty. In that case after base changing to a finite \'etale cover 
$\Ig \to \Shbar^{\text{ord}}$
one can trivialize $P_{\dR}$, which can be used split the exact sequence 
\begin{equation} \label{fibration-seq}
0 \to \pi^* \Omega^1_{\Shbar^{\text{ord}}} \to \Omega^1_{\flag} \to \Omega^1_{\flag/\Shbar^{\text{ord}}} \to 0.
\end{equation}
Using our formalism one can check  
that there is a section of \eqref{fibration-seq}
on the fiber of the ordinary locus, where all the Hasse invariants 
$H_{\alpha}$ are isomorphisms for $\alpha \notin \Lie(M)$, and there exists 
a reduction $B_{M,\dR} \subseteq B_{\dR}$, for $B_{M}=M \cap B$. Both sections agree since 
the trivialization of $P_{\dR}$ induces a section of $B_{M,\dR}$, with which we can easily compare 
them over $\Ig$.
Then we can further work on $U$ to split $\Omega^1_{\flag}$ into line bundles, since each graded piece of the exact sequence \eqref{fibration-seq}
has a simple description as an automorphic vector bundle.
\item Let $\Sh$ be a PEL Shimura variety of type $A$ or $C$ such that 
the ordinary locus is non-empty. For simplicity assume that $\Sh$ is the Siegel Shimura variety.
Let $\alpha \in \Delta \backslash \Lie(M)$ so that $\pi_* H_{\alpha}$ is the classical Hasse invariant. 
In this setting \cite{Eischen-Mantovan-1} define an operator $\theta$ 
$$
\theta: \omega(\lambda) \to \omega(\lambda+w_{0,M}h_{\alpha}+w_{0,M}\alpha) 
$$
by embedding $\omega(\lambda)$ into $\mathcal{V}$ for some $V \in \Rep(G)$, then applying $\nabla$, whose image 
lands in the last non-trivial step of its Hodge filtration $F^{1}\mathcal{V}$. Lastly, they define a projection 
$F^1 \mathcal{V} \to \omega(\lambda+w_{0,M}h_{\alpha})$ using the Verschiebung and its adjugate
$\text{Ad} V: \omega^{(p)} \to \omega \otimes \text{det}^{p-1} \omega$. 
Then $\theta=\pi_*\theta_{\alpha}$. To see why we need two observations. 
First the projection $\Omega^1_{\flag} \to \LL(\alpha)$ is defined over $\flag$, 
since $-\alpha \to \mathfrak{g}/\mathfrak{b}$
is $B$-equivariant; and applying $\pi_*$ one gets the identity on $\Omega^1_{\Shbar}$. 
Second, one checks that in \Cref{easy-lemma-theta}(1) multiplying by $H_{\alpha}$ the projection 
$\mathcal{V} \to \LL(w_{0,M}\lambda)$ corresponds to $F^1 \mathcal{V} \to \omega(\lambda+w_{0,M}h_{\alpha})$
on that piece of the filtration, using the explicit description of $H_{\alpha}$ as the determinant of $V: \omega \to \omega^{(p)}$.
\end{enumerate}
\end{remark}

We can also prove some other properties typical of theta operators. 

\begin{prop} \label{theta-kills-Hasse}
	\begin{enumerate}
	\item  Let $\phi \in T_{\dR}$ be an adapted element
  with respect to the data of some $(V,B_0,B^{z'}_{0})$ from \Cref{special-embedding}. It induces a basis $\{(\phi,v_{\lambda})\}$ for each of the line bundles 
	$\LL(\lambda)$. Then with respect to this basis $\tilde{\theta}_{\gamma}$ is given by 
	$f \mapsto D_{\gamma}(f)$, where $f \in \mathcal{O}_{\flag}$, and $D_{\gamma}=e^1_1(\phi,x_{-\gamma})$.
	\item We have that $\theta_{\gamma}(H_{\alpha})=0$ for all $\gamma \in \Phi^{+}$ and $\alpha \in \Delta$. 
	Therefore, $\theta_{\gamma}$ and $H_{\alpha}$ commute as operators. 
	\item (Commutation relations) $[\tilde{\theta}_{\gamma_1},\tilde{\theta}_{\gamma_2}]=\tilde{\theta}_{[\gamma_1,\gamma_2]}$, 
	where for $m \in \overline{\F}_p$, $\tilde{\theta}_{m \gamma}:=m\tilde{\theta}_{\gamma}$. In particular, $\theta_{0}=0$.
	\end{enumerate}
	\begin{proof}
	Part $(1)$ follows directly from the way $\tilde{\theta}_{\gamma}$ is defined, and the explicit description of the canonical isomorphisms $e_{V}$
  with respect to an adapted element in 
  \Cref{canonical-iso-U}. 
	For $(2)$ it is enough to check that 
	$\tilde{\theta}_{\gamma}(H_{\alpha})=0$ over $U$. 
	By \Cref{Hasse-functorial} $H_{\alpha}$ is a combination of Hasse invariants $H_{\alpha_i}$ coming from $\Sbt_{\GL(V)}$. 
	The latter can be obtained from the isomorphisms $H_i: F^{n-i}_{\mathcal{V}}
   \cap D_{\mathcal{V},n-i+1} \cong \LL'_{\mathcal{V},n-i} \cong \LL^p_{\mathcal{V}^{\sigma^{-1}},i+1}$, i.e. 
	$H_{\alpha_i}=\text{det}(F^{n-i}_{\mathcal{V}} \to \mathcal{V}/D_{\mathcal{V},n-i})=H_{1}  \ldots H_{i}$, 
	where $\alpha_i$ are the simple roots of $\GL(V)$, and we identify
	$\LL_{\mathcal{V},i}=\tilde{\LL}_{\mathcal{V},i}:=F^{n-i}_{\mathcal{V}} \cap D_{\mathcal{V},n-i+1}$.
  Since by construction of $\phi$ the 
	$H_{i}$ are the identity (up to a $p$th power) on the basis induced by $\phi$ and $B_0$, 
  we get the conclusion from part $(1)$.
	 Part $(3)$ follows from \Cref{composition-U}.
	\end{proof}
	\end{prop}

	Part (3) also extends to $\flag$ and $\theta_{\gamma}$ by multiplying appropriately by Hasse invariants, 
	thanks to part (2). We see that for $p \ge 3$ all basic theta operators can be obtained from commutator relations from the 
	simple basic operators, but this procedure introduces superfluous powers of Hasse invariants in many cases.
	Finally, it is easy to describe which basic theta operators become linear after pushing them to $\Shbar$.

	\begin{lemma}
	Let $\gamma \in \Phi^{+} \cap \Lie(M)$. Then $\pi_* \theta_{\gamma}$ is a linear map on $\Shbar$. 
	\begin{proof}
	For such $\gamma$ we can factor $\phi_{\gamma,\lambda}: -\lambda -\gamma \to \Ver_{B}(-\lambda)$ 
  through $\Ver_{P/B}(-\lambda) \to \Ver_{B}(-\lambda)$. 
  There exists a functor $\Hom_{T}(\lambda, \Ver_{P/B}(\mu)) \to \DiffOp_{\flag/\Shbar}(\LL(-\mu), \LL(-\lambda))_{\mid U}$
  analogous to \eqref{diff-U}, by using the isomorphism $e_{P/B}$ from \Cref{canonical-iso-Pm}.
  Therefore, we can consider 
	$\tilde{\theta}_{\gamma}$ as a differential operator with respect to $\flag/\Shbar$ over $U$.
  Similarly, by mimicking the procedure of \Cref{defn-basic-theta} we see that $\theta_{\gamma}$ is a differential operator with respect to $\flag/\Shbar$.
   Now, the map $\LL(\lambda) 
	\to P_{\flag/\Shbar} \otimes \LL(\lambda)$ becomes $\mathcal{O}_{\Shbar}$-linear after applying $\pi_{*}$. 
  Therefore, $\pi_* \theta_{\gamma}: \pi_*(\LL(\lambda) 
	\to P_{\flag/\Shbar} \otimes \LL(\lambda)) \xrightarrow{\pi_*g} \LL(\lambda+\gamma)$ is linear, where 
  $g: P_{\flag} \otimes \LL(\lambda) \to \LL(\lambda+\gamma)$ is the linearized version of $\theta_{\gamma}$.
	\end{proof}
	\end{lemma}

\subsection{$p$-power relations between basic theta operators}

In characteristic $p$ we have the special behaviour that the $p$th iteration of a derivation is again a derivation. 
Hence, we might expect that for some $\gamma \in \Phi^{+}$ we can rewrite $\theta^p_{\gamma}$ in terms of another basic theta operator
and Hasse invariants. We prove that this is the case for $\gamma \in \Delta$. We do this by computing the 
$p$-curvature of certain elements on $\mathcal{V}_{U}$ for $V$ from \Cref{special-embedding}, using 
\Cref{conjugate-p-curvature} and \Cref{nabla-on-U}. We will need a small lemma first. 

\begin{lemma} \label{special-embedding-2}
Let $G_i$ be an almost-simple factor of $G^{\der}_{\Fpbar}$. There exists an irreducible minuscule representation 
$V \in \Rep_{\Fpbar}(G)$ such that $G_i \to \GL(V)$ has a finite kernel. Moreover, there exist Borel subgroups 
$B_0, B^{z'}_0 \subseteq \GL(V)$ satisfying the hypothesis of \Cref{special-embedding} together with
$B^{z'}_{0} \subseteq P^{-(p)}_{V}$.
\begin{proof}
Consider $G^{\der}$ as a quotient of $\prod H_i$ as in \Cref{Hodge-embeddings}. We take $V$ to be one of the summands 
of the Hodge embedding, so that it is non-trivial on the $H_i$ corresponding to $G_i$. We can write it as 
$V=\otimes V_i$, where each $V_i$ is a minuscule representation of $H_i$ (which can be trivial). Without loss of 
generality we may assume that they are all non-trivial.
 We follow 
the proof of \Cref{special-embedding} to choose a particular $B_0 \subseteq \GL(V)$ satisfying 
$B^{z'}_{0} \subseteq P^{-(p)}_{V}$. Since $V_i$ is minuscule 
there is an ordered basis of $T$-eigenvectors $\{v_{i,1},v_{i,2},\ldots,v_{i,n_i}\}$ of it such that its associated Borel 
(the convention is that $v_{i,1}$ generates the rank $1$ object in the flag) in $\GL(V_i)$
contains $B$. Explicitly, it satisfies that $w(v_{i,1})>w(v_{i,2})> \ldots >w(v_{i,n_i})$. 
Let $\mu$ be a dominant representative 
of the Hodge cocharacter. Since the composition of $\mu$ with $V$ must still be minuscule, we deduce that $\mu$ is trivial 
on each $H_i$ except the last one $H_k$. Then we choose an ordered basis of $V$ in a "lexicographic order" 
as follows 
\begin{align*}
\{&u_1,\ldots, u_{N}\}:=\{v_{1,1} \otimes v_{2,1} \ldots \otimes v_{k,1}, v_{1,2}\otimes v_{2,1} \ldots \otimes v_{k,1}, \ldots \\
&v_{1,n_1} \otimes v_{2,1} \ldots \otimes v_{k,1}, v_{1,2} \otimes v_{2,1} \ldots \otimes v_{k,1}, \ldots ,
v_{1,n_1} \otimes v_{2,n_2} \ldots \otimes v_{k,n_k}\}.
\end{align*}
We let $B_0$ to be the Borel associated to this flag, and $T_0$ the maximal torus 
defined by the basis vectors. Note that it is constructed so that $B \subseteq B_0 \subseteq P_{V}$.
 This choice defines a map of Weyl groups $W_{G} \to W_{\GL(V)}$. It sends 
$w_0$ to $w_{\GL(V)}$, since $w_0 v_{i,j}=v_{i,n_i-j+1}$. Thus, as in the proof of \Cref{special-embedding}
take $z'=zw_{0}w_{0,\GL(V)}=z=\sigma(w_{0,M})w_{0}$.
We check that $B^{z'}_{0} \subseteq P^{-(p)}_{V}$  is also satisfied. We have that $P^{-(p)}$ is the parabolic associated to $-\sigma \mu$, so we have 
to prove that $\langle -\sigma \mu, z' u_{r} \rangle \ge \langle -\sigma \mu,z' u_{r+1}\rangle$ for every $r$. 
Using that $z'=\sigma(w_{0,M})w_{0}$, this is equivalent to 
$\langle \mu,w_{0,M}v_{N-r} \rangle \le \langle \mu,w_{0,M}v_{N-r-1}\rangle$. This holds using three observations. 
First that $\mu$ on $\GL(V)$ is minuscule, second that $\langle \mu,v_{N-r-1} \rangle \ge \langle \mu, v_{N-r} \rangle$ and 
third that $\langle \mu, w_{0,M}\lambda \rangle=\langle \mu, \lambda \rangle$ for any $\lambda$.  
\end{proof}
\end{lemma}

\begin{theorem} \label{pth-power-relation}
  Let $G_i$ be the almost-simple factors of $G^{\der}_{\Fpbar}$. Fix some $i$ and   
 let $\alpha \in \Delta_{G_i}$. The Frobenius $\sigma$ permutes the factors $G_i$, let $G_j$
 be $\sigma(G_i)$.  
\begin{enumerate}
\item Let $\alpha \in \Phi_{M}$, then $\theta^p_{\alpha}=0$.
\item Assume that $\alpha \notin \Lie(M)$  i.e. $\alpha$ is 
the special node for $\mu$. Further assume that the almost simple quotient $G^i$ is simply connected 
or $\textnormal{SO}_{2n}$. 
Let $\delta$ be the longest root 
of $G_{j}$. Then for $G_i$ of rank at least $2$ we have
$$
\theta^p_{\alpha}=H^p_{\alpha}\theta_{\delta},
$$
and for rank $1$ ($G_i$ isogenous to $\SL_2$) we have $\theta^p_{\alpha}=H^{p+1}_{\alpha}\theta_{\delta}$.
\end{enumerate}
\begin{proof} For $(1)$ we use the compatibility of the relative Frobenius with $\Ver^0_{P/B}$ in \Cref{creator-diff}(1) and 
the compatibility with composition of the functor 
$$\Hom_{T}(\lambda,\Ver_{B}(\mu)) \to \DiffOp_{B \backslash Bw_{0}B}(\LL_{G/B}(-\mu),
\LL_{G/B}(-\lambda)),
$$ 
which can be proved via a reduction of the equivalence in \Cref{vb-on-flag} to $\Rep(T)$. Then we observe that $\tilde{\theta}^p_{\alpha}$ is given by the map of Verma modules 
associated to $x^p_{-\alpha}$, which is trivial as a map of baby Verma modules. 
Note that 
this argument also proves that $\theta^{p}_{\gamma}=0$ for $\gamma \in \Phi^{+}_{M}$.
For $(2)$, we can prove the identity over $U$. First we check that both sides have the same weight increase. 
We check it case by case using the notation of \Cref{table-lie} and \Cref{Hasse-dfn}.
Let $\theta_{\gamma}$ have weight increase
$\gamma+h_{\gamma}$, where $h_{\gamma}$ is the sum of the weights of the Hasse invariants used in $\theta_{\gamma}$.
Thus, we want to prove that $p\alpha=\delta+h_{\delta}$.
\begin{enumerate}
\item ($A_n$) Assume that $n \ge 2 $, the case of $n=1$ can be checked by hand. 
Use $\alpha_i$ for the simple roots of $G^i \cong \SL_n$
and $\beta_i$ for $G_j$. Then $\alpha=\alpha_{i}$, for some $i$, and $\delta=\beta_1+\ldots+\beta_n$, 
so that $h_{\delta}=h_{\beta_1}+h_{\beta_n}$. Let $\delta_{i}$ be the 
longest root of $G_i$. In the definition of Hasse invariants
\Cref{Hasse-dfn} we can take $\chi_{\beta_1}=e^{j}_{1}$ and $\chi_{\beta_n}=-e^j_{n}$. 
From the formula for $h_{\beta_j}$ \eqref{Hasse-weight} we reduce to proving $\alpha_{i}=w_{0,M}\sigma^{-1}\delta$. 
But $w_{0,M}\alpha_{i}=\delta_{i}$. The result follows since the Frobenius $\sigma$
always preserves the simple roots, i.e. it is an outer isomorphism 
by Lang's theorem, so that $\sigma^{-1}\delta=\delta_{i}$, by the characterization 
of the longest root.  
\item ($C_n$) Here $G^i=\textnormal{Sp}_{2n}$, and there is a unique minuscule character in this case, so that $\alpha=\alpha_n$.
We have $\delta=2e^{j}_{1}=2(\beta_1+\ldots+\beta_{n-1})+\beta_n$, so that $h_{\delta}=2h_{\beta_1}$. We can take 
$\chi_{\beta_1}=e^{j}_{1}$, so we reduce to proving $w_{0,M}\alpha_n=\sigma^{-1} \delta$. 
Since $w_{0,M}\alpha=\delta_i$
 we conclude as before. 
 \item $(B_{n})$ Let $n\ge 3$. Here $G^i=\text{Spin}_{2n+1}$, and there is a unique minuscule character, 
 so that $\alpha=\alpha_1$. We have $\delta=e^j_1+e^j_2$, 
 $h_{\delta}=h_{\beta_2}$,
  and we can take $\chi_{\beta_2}=e^{j}_1+e^j_{2}=\delta$. 
  From \eqref{Hasse-weight} we reduce to proving 
 $w_{0,M}\alpha_1=\sigma^{-1}\delta$. In this case $M^{\text{der}}=\textnormal{Spin}_{2n-1}$. 
 The longest element $w_{0} \in W_{B_n}=W_{C_n}$ satisfies $w_{0}=-\text{id}$, so we get
 $w_{0,M}\alpha_1=\delta_i$. We conclude by $\sigma^{-1}\delta=\delta_i$ as above. 
\item ($D_n$)  Let $n \ge 4$. We have 
$\delta=e^{j}_1+e^{j}_2=\beta_1+2(\beta_2+\ldots+\beta_{n-2})+\beta_{n-1}+\beta_n$, and 
$h_{\delta}=h_{\beta_2}$ with $\chi_{\beta_2}=e^{j}_1+e^{j}_2$. Note that by \Cref{table-lie}
the last equality is independent 
on whether $G^i$ is $\textnormal{Spin}_{2n}$ or $\textnormal{SO}_{2n}$.  
Again we reduce to proving 
$w_{0,M}\alpha=\sigma^{-1}\delta$. First suppose $\alpha=\alpha_n$.
In this case $M^{i} \cong \SL_{n}$ so that 
$w_{0,M}\alpha_n=e^i_1+e^i_2=\delta_i$, so it follows from $\sigma^{-1}\delta=\delta_i$. 
 Now suppose $\alpha=\alpha_1$. In this case $M^{i}=\textnormal{Spin}_{2n-2}$ or $\textnormal{SO}_{2n-2}$.
 Its Weyl group is the group of signed permutations with an even number of sign changes, and 
 $w_{0,M}(a_1,\ldots,a_n)=(a_1,-a_2,-a_3,a_4,\ldots,a_n)$.
 Thus, $w_{0,M}\alpha_1=\delta_i$, which proves the claim. The case of $\alpha=\alpha_{n-1}$ follows 
 from the case $\alpha=\alpha_{n}$, since there is an outer automorphism of $G^i$ swapping the two roots.  
\end{enumerate}
Choose some $(V,B_0,B^{z'}_0)$ as in \Cref{special-embedding-2} such that $V$ is non-trivial on $G_i$. 
 After checking that the weights agree, by \Cref{theta-kills-Hasse}(1)
we reduce to proving the identity $D^{[p]}_{\alpha}=D_{\delta}$, where we have chosen an adapted element 
$\phi \in T_{\dR}$ with respect to the data of $(V,B_0,B^{z'}_0)$. 
We will compute part of the $p$-curvature on $\mathcal{V}$ in two different ways. Let $F^{\bullet}_{\mathcal{V}}$ and 
$D_{\bullet,\mathcal{V}}$ the extended Hodge and conjugate filtration given by $B_0$ and $B^{z'}_0$. Since 
$B^{z'}_0 \subseteq P^{-(p)}$ by construction of $V$, the latter refines the one-step 
conjugate filtration on $\mathcal{V}$. Let $N$ be the dimension of $V$, and $r$ such that $D_{N-r+1,\mathcal{V}}=\mathcal{V}$
but $D_{N-r,\mathcal{V}} \neq \mathcal{V}$. 
We have the commutative diagram
$$
\begin{tikzcd} \label{p-curvature}
\mathcal{V}/D_{N-r,\mathcal{V}} \arrow[r,"\psi"] \arrow[d,"\sim"] & D_{N-r,\mathcal{V}} \otimes \Omega^{1(p)}_{\flag} 
\arrow[d,"\sim"]\\ 
F^{N-r,(p)}_{\mathcal{V}^{\sigma}} \arrow[r,"\KS^{(p)}_{\mathcal{V}}"] & 
(\mathcal{V}^{\sigma}/F^{N-r}_{\mathcal{V}^{\sigma}} \otimes \Omega^1_{\flag})^{(p)},
\end{tikzcd}
$$
by \Cref{conjugate-p-curvature}, where the vertical isomorphisms are the Cartier isomorphisms. 
Recall the ordered basis $\{v_{1,s_1} \otimes v_{2,s_2} \ldots \otimes v_{i,s_i}\}$ of $V=\otimes V_j$. 
Let $\{v_{i,1},v_{i,2},\ldots,v_{i,t}\}$ be the Hodge filtration induced on $V_i$ by $\mu$, which is
the only $V_j$ that 
has a non-trivial Hodge filtration. 
Using the diagram above we will compute $\psi_{D^{(p)}_{\alpha}}(\phi,v_{1,s_1} \otimes v_{2,s_2} \ldots \otimes v_{i,s_i})$ 
in two ways 
for certain basis vectors. We will use the following easy lemma.

\begin{lemma} \label{easy-rep-lemma}
Let $G$ be an almost-simple semisimple group of type $A,B,C$ or $D$ over an 
algebraically closed field. Let $T \subseteq B \subseteq G$ be a choice 
of maximal torus and Borel. 
 Let $V \in \Rep_{\Fpbar}(G)$ be an irreducible 
minuscule representation of $G$, and $\mu \in X_{*}(T)$ a non-trivial dominant minuscule cocharacter.
Let $\alpha \in \Delta$
be the special node for $\mu$, i.e. $\langle \mu, \alpha \rangle= 1$, and $\delta$ the longest root. 
There is a unique total order on the weights appearing in $V$: $\lambda_1> \lambda_2 > \ldots \lambda_n$, 
with $\langle \mu, \lambda_1 \rangle=1= \langle \mu, \lambda_r \rangle$ and $\langle \mu, \lambda_{r+1} \rangle=0$. 
Then there exists basis vectors $v_{\lambda_i} \in V$ with weight $\lambda_i$ such that
$$
x_{-\alpha}v_{\lambda_{r}}=v_{\lambda_{r+1}} \;\;\; x_{-\delta}v_{\lambda_1}=v_{\lambda_n}.
$$
\end{lemma}

First, we claim that under the isomorphism $\mathcal{V}/D_{N-r,\mathcal{V}} \cong F^{N-r,(p)}_{\mathcal{V}^{\sigma}}$
$(\phi,v_{1,n_1}\otimes v_{2,n_2}\ldots \otimes v_{i,1})$ is sent to $(\phi,v_{1,1} \otimes v_{2,1} \ldots \otimes 
v_{i,t})^{(p)}+\sum g_{s}(\phi,v_{1,s_1} \otimes \ldots \otimes v_{i,s_i})^{(p)}$ where in the sum $s_i<t$. 
By construction of $\phi$, we have that 
$(\phi,v_{1,n_1}\otimes v_{2,n_2}\ldots \otimes v_{i,1}) \in D_{N-\sum_{j <i} n_j+1,\mathcal{V}}$, and under the 
Cartier isomorphism 
$D_{N-\sum_{j <i} n_j+1,\mathcal{V}}/D_{N-r,\mathcal{V}} \cong F^{N-r+\sum_{j <i} n_j-1,(p)}_{\mathcal{V}^{\sigma}}$. 
The latter is generated by the vectors appearing in the expression. We can assume that 
the coefficient of $(\phi,v_{1,1} \otimes v_{2,1} \ldots \otimes 
v_{i,t})^{(p)}$ is $1$ by construction of $\phi$ (the coefficient will always be a $p$th power).

By \Cref{nabla-on-U} we have that 
$D^{(p)}_{\alpha} \circ \KS_{\mathcal{V}}(\phi,v_{1,1} \otimes v_{2,1} \ldots \otimes 
v_{i,t})=(\phi,v_{1,1} \otimes v_{2,1} \ldots \otimes 
v_{i,t+1})^{(p)}$, using that by \Cref{easy-rep-lemma} $x_{-\alpha}v_{i,t}=v_{i,t+1}$. Then under the isomorphism 
$(\mathcal{V}^{\sigma}/F^{N-r}_{\mathcal{V}^{\sigma}})^{(p)} \cong  D_{N-r,\mathcal{V}}$ we have that 
$(\phi,v_{1,1} \otimes v_{2,1} \ldots \otimes 
v_{i,t+1})^{(p)}$ corresponds to $(\phi,v_{1,n_1}\otimes v_{2,n_2} \ldots \otimes v_{i,n_i})$. This is because by 
definition of $\phi$, $(\phi,v_{1,n_1}\otimes v_{2,n_2} \ldots \otimes v_{i,n_i}) \in D_{1,\mathcal{V}}$, which 
corresponds to the rank $1$ object $(F^{N-r-1}_{\mathcal{V}^{\sigma}}/F^{N-r}_{\mathcal{V}^{\sigma}})^{(p)}$, 
which is generated by $(\phi,v_{1,1} \otimes v_{2,1} \ldots \otimes 
v_{i,t+1})$ by construction of $B_0$. Therefore, on one hand by diagram above
$$
\psi_{D^{(p)}_{\alpha}}(\phi,v_{1,n_1}\otimes v_{2,n_2}\ldots \otimes v_{i,1})=(\phi,v_{1,n_1}\otimes v_{2,n_2} \ldots \otimes v_{i,n_i}).
$$
On the other hand, by definition $\psi_{D^{(p)}_{\alpha}}=\nabla_{D^{[p]}_{\alpha}}-\nabla^{p}_{D_{\alpha}}$. From \Cref{nabla-on-U} we inmediately see that 
$\nabla^{p}_{D_{\alpha}}=0$, since $x^p_{-\alpha}$ acts trivially on $V \in \Rep(G)$. 
Locally write $D^{[p]}_{\alpha}=\sum_{\gamma \in \Phi^{+}} f_{\gamma}D_{\gamma}$ for some $f_{\gamma} \in \mathcal{O}_{U}$. 
Using \Cref{nabla-on-U} we have
\begin{align*}
&\psi_{D^{(p)}_{\alpha}}(\phi,v_{1,n_1}\otimes v_{2,n_2}\ldots \otimes v_{i,1})=
(\phi,v_{1,n_1}\otimes v_{2,n_2} \ldots \otimes v_{i,n_i})=\\
&\nabla_{D^{[p]}_{\alpha}} (\phi,v_{1,n_1}\otimes v_{2,n_2}\ldots \otimes v_{i,1}) 
=
\sum_{\gamma \in \Phi^{+}} D^{[p]}_{\alpha} (\phi, x_{-\gamma}[v_{1,n_1}\otimes v_{2,n_2}\ldots \otimes v_{i,1}]) \otimes \omega^{(p)}_{-\gamma}.
\end{align*}
By  \Cref{easy-rep-lemma} $x_{-\delta}v_{i,1}=v_{i,n_i}$ and for $\gamma \neq \delta$ $x_{-\gamma}v_{i,1}$
does not involve $v_{i,n_i}$, so we conclude that $f_{\delta}=1$.
Now let $\gamma \in \Phi^{+}$ different to $\delta$. Assume first that $\gamma \notin \Phi^{+}_{G_i}$.
Then without loss of generality there exist basis elements $v_{1,a}, v_{1,b} \in V_1$ such that $x_{-\gamma}v_{1,a}=v_{1,b}$.
Under the isomorphism $\mathcal{V}/D_{N-r,\mathcal{V}} \cong F^{N-r,(p)}_{\mathcal{V}^{\sigma}}$, 
$(\phi,v_{1,a}\otimes v_{2,1}\ldots \otimes v_{i,2})$ is sent to 
$\sum g_{s}(\phi,v_{1,s_1} \otimes \ldots \otimes v_{i,s_i})^{(p)}$ where in the sum $s_i<t$. Then for such 
$(s_j)$ we have
that $D^{(p)}_{\alpha} \KS_{\mathcal{V}^{\sigma}} (\phi,v_{1,s_1} \otimes \ldots \otimes v_{i,s_i})^{(p)}=0$. 
We conclude that $\psi_{D^{(p)}_{\alpha}}(\phi,v_{1,a}\otimes v_{2,1}\ldots \otimes v_{i,2})=0=
\nabla_{D^{[p]}_{\alpha}}(\phi,v_{1,a}\otimes v_{2,1}\ldots \otimes v_{i,2})$, which implies that 
$f_{\gamma}=0$. Now suppose that $\gamma \in \Phi^{+}_{G_i}$. There exist $v_{i,a}, v_{i,b}$ with $(a,b) \neq (1,n_i)$ 
satisfying $x_{-\gamma}v_{i,a}=v_{i,b}$. By a similar argument, unless $a=1$, 
$\psi_{D^{(p)}_{\alpha}}(\phi,v_{1,n_1}\otimes v_{2,n_2}\ldots \otimes v_{i,a})=0$. In that case we again obtain 
$f_{\gamma}=0$. If $a=1$ we use 
$\psi_{D^{(p)}_{\alpha}}(\phi,v_{1,n_1}\otimes v_{2,n_2}\ldots \otimes v_{i,1})=(\phi,
v_{1,n_1}\otimes v_{2,n_2}\ldots \otimes v_{i,n_i})$ to conclude that $f_{\gamma}=0$. Therefore, $D^{[p]}_{\alpha}=D_{\delta}$
as desired.

\end{proof}
\end{theorem}

We note that the restriction on $G^i$ in \Cref{pth-power-relation} is necessary, otherwise the weights of the Hasse invariants involved 
in $\theta_{\delta}$ might be too large for the identity to hold. In any case, the theorem proves that if the 
weight increases of $\theta^p_{\alpha}$ and $H^p_{\alpha}\theta_{\delta}$ agree, then they are equal.

\begin{remark} \label{Baby-Vermas-dont-work}
    It is not true that $F_B(\Ver^0_B(\lambda))=D^{[p]}_{\flag} \otimes \LL( \lambda)$. Otherwise, we would
     have the analogue
     of \Cref{creator-diff}$(1)$ for $\Ver^0_B(\lambda)$, which would imply that 
    $\theta^p_{\alpha}=0$ in case $(2)$ above. Another way to see this is to use Serre-Tate coordinates. 
	Take the simplest case of the modular curve. Then on the formal completion of an ordinary point $x$, theorem of 
	Katz \cite{Katz-Serre-Tate} about Serre-Tate coordinates tells us that $(\phi,x_{-}) \in \Omega^1_{\Shbar^{\wedge}_{x}}$ for an adapted element 
  is identified with $q\frac{d}{dq}$.
  Then by \Cref{canonical-iso-U}, $F_B(\Ver^0_B(1))$ is
    the quotient sheaf of $D_{\Shbar^{\wedge}_{x}}=\Fpbar[[q-1]]dq$ generated by $(q\frac{d}{dq})^{i}$ for $0 \le i \le p-1$.
	On the other hand,
    $D^{[p]}_{\Shbar^{\wedge}_{x}}$ is generated by $(\frac{d}{dq})^{i}$ for $0 \le i \le p-1$. One can explicitly 
    see that they are not isomorphic as quotients of $D_{\Shbar^{\wedge}_{x}}$.
    \end{remark}
   \begin{example}
  Let us illustrate the proposition in the case of unitary Shimura varieties.
  For simplicity take the signature $(2,1)$, $E$ to be a quadratic imaginary field with embeddings $\{\sigma,\overline{\sigma}\}$
  as in \cite[\S 2.3]{laporta}.
  Let $\LL_{\sigma} \subset \omega_{\sigma}$ be part of the Hodge filtration over $\flag$. 
  We have $\alpha=\alpha_2$. 
  Suppose that $p$ is inert in $E$. Then $G_{\overline{\F}_p}=\GL_3 \times \mathbb{G}_m$
  with Frobenius acting via the non-trivial outer automorphism. We can choose $V$ so that $\mathcal{V}=
  \mathcal{H}_{\sigma}$, and its  $\sigma$-twist gives $\mathcal{H}_{\overline{\sigma}}$. 
    The $\sigma$-isotypic parts of
  the extended Hodge and conjugate filtrations on $\flag$ are
  $$
  0 \subset \LL_{\sigma} \subset \omega_{\sigma} \subset \mathcal{H}_{\sigma} \twoheadrightarrow 
  \omega^{-1}_{\overline{\sigma}} 
  $$
  $$
  0 \subset (\omega_{\sigma}/\LL_{\sigma})^{-p} \subset \Ker(V_{\sigma})\cong \omega^{\vee(p)}_{\sigma}
   \subset \mathcal{H}_{\sigma} \twoheadrightarrow \omega^{p}_{\overline{\sigma}}, 
  $$
  where the last map is the Verschiebung. Then an adapted element $\phi$ is defined by a section
  $(\phi,f_1)$ of $\LL_{\sigma}$
  which maps under $\mathcal{H}_{\sigma} \twoheadrightarrow \omega^{p}_{\overline{\sigma}}$ to a horizontal 
  section for $\nabla$, an element $(\phi,f_2)$ of $\omega_{\sigma} \cap \Ker(V_{\sigma})$ that under
   $\omega_{\sigma} \cap \Ker(V_{\sigma}) \to \LL^{-p}_{\sigma}$ maps to a horizontal section, 
   and $(\phi,f_3)$ a horizontal section of $(\omega_{\sigma}/\LL_{\sigma})^{-p}$.
   In the $p$ split case $V$ and its $\sigma$-twist agree, and 
  we can just work with the $\sigma$-isotypic component.
  In both cases we have $\theta^p_{\alpha_2}=H^{p}_{\alpha_2}\theta_{\alpha_1+\alpha_2}$. This does not depend on 
  the splitting of $p$ on $E$, even though the structure of the conjugate filtration on $\Shbar$ changes.
   \end{example}

\subsection{Restriction to codimension $1$ strata.}
We study the restriction to codimension $1$ strata of basic theta operators, culminating in a proof 
the generic 
injectivity of basic theta operators in the $p$-restricted region in \Cref{basic-theta-injective}.
Let $\alpha \in \Delta$ and $i: \overline{D}_{\alpha} \hookrightarrow \flag$ and 
$\gamma \in \Phi^{+}$. In this subsection assume that $H_{\alpha}$ has a simple zero.
Then $\theta_{\gamma}: \LL(\lambda) \to \LL(\lambda+\mu_{\gamma})$ fits in the diagram 
$$
\begin{tikzcd}
0 \arrow[r] & \LL(\lambda-h_{\alpha}) \arrow[d,"\theta_{\gamma}"] \arrow[r,"H_{\alpha}"] & \LL(\lambda)
\arrow[d,"\theta_{\gamma}"] \arrow[r] & 
 i_{*}i^{*} \LL(\lambda) \arrow[r] \arrow[d,"\theta_{\gamma,\overline{D}_{\alpha}}"] & 0 \\
 0 \arrow[r] & \LL(\lambda-h_{\alpha}+\mu_{\gamma}) \arrow[r,"H_{\alpha}"] & \LL(\lambda+\mu_{\gamma}) \arrow[r] & 
 i_{*}i^{*} \LL(\lambda+\mu_{\gamma}) \arrow[r] & 0,
\end{tikzcd}
$$
where $\theta_{\gamma,\overline{D}_{\alpha}}$ is defined by the diagram. When taking global sections of the diagram, 
this means that for $f \in \H^0(\flag, \LL(\lambda))$, $\theta_{\gamma}(f)_{\mid \overline{D}_{\alpha}}=
\theta_{\gamma,\overline{D}_{\alpha}}(f_{\mid \overline{D}_{\alpha}})$. We use the notation 
$\Theta_{\gamma}: P^{1}_{\flag} \otimes \LL(\lambda) \to \LL(\lambda+\mu_{\gamma})$ for the linearized map
corresponding to $\theta_{\gamma}$. 

\begin{lemma} \label{linearized-theta}
Let $\alpha \in \Delta$, $\gamma, \beta \in \Phi^{+}$ such that 
$i^* \Theta_{\gamma}=B \circ i^* \Theta_{\beta}$ where $B$ is a linear map over $\overline{D}_{\alpha}$. 
Then $\theta_{\gamma,\overline{D}_{\alpha}}=i_{*}B \circ \theta_{\beta,\overline{D}_{\alpha}}$.
\begin{proof}
We give another characterization of $\theta_{\gamma,\overline{D}_{\alpha}}: i_* i^* \LL(\lambda) \to i_* i^* \LL(\mu)$. 
Take a local section $\overline{v} \in  i^* \LL(\lambda)$, and lift it to a section $v$ over $\flag$. Then 
by definition $\theta_{\gamma, \overline{D}_{\alpha}}(\overline{v})$ is the reduction of $\theta_{\gamma}(v)=\Theta_{\gamma} \circ 
g(v)$, where $g: \LL(\lambda) \to P^{1}_{\flag} \otimes \LL(\lambda)$. This makes it clear that  
$\theta_{\gamma,\overline{D}_{\alpha}}=i_{*}B \circ \theta_{\beta,\overline{D}_{\alpha}}$. 
\end{proof} 
\end{lemma}
We can describe the behaviour of simple basic theta operators $\theta_{\alpha}$ when restricted to their corresponding 
closed stratum $\overline{D}_{\alpha}$. 
\begin{theorem} \label{restriction-simple-theta}
Let $\alpha \in \Delta$, $\lambda \in X^*(T)$. Assume that $H_{\alpha}$ has a simple zero. Then
\begin{enumerate}
\item When applied to $\LL(\lambda)$ we have
 $\theta_{\alpha,\overline{D}_{\alpha}}=\langle \lambda, \alpha^{\vee} \rangle B_{\alpha}$,
where $B_{\alpha}$ is a section over $\overline{D}_{\alpha}$ that does not vanish on $D_{\alpha}$.
\item  Let $f \in \LL(\lambda)$ be a local section. Then $H_{\alpha} \mid \theta_{\alpha}(f)$
if and only if $H_{\alpha} \mid f$ or $p \mid \langle \lambda, \alpha^{\vee} \rangle$. 
\item  Assume that \Cref{assumption-vanishing}(2) holds. 
Then the map 
$\theta_{\alpha}$ is injective on $\H^0(\Shbar,\omega(\lambda))$ 
for $\lambda \in X_{1}(T)$ such that  $\langle \lambda, \sigma^{-1}\alpha^{\vee} \rangle< p-\epsilon$ where 
$\epsilon \ge 0$ depends only on $G_{\Q}$, and $p \nmid \langle w_{0,M}\lambda, \alpha^{\vee} \rangle$. 
In particular, 
it is injective 
for generic $\lambda \in X_1(T)$ as in \Cref{generic-Sh}.
\end{enumerate}
\begin{proof}
For $(1)$ we can prove the statement at the level of linearized differential operators, by \Cref{linearized-theta}.
This further reduces to a statement about vector bundles on $\overline{\Sbt}_{\alpha}$. Namely, we have to prove that 
over $\overline{\Sbt}_{\alpha}$ there is a factorization 
$$
\begin{tikzcd}
\LL(-\lambda) \otimes \LL(\chi_{\alpha}, -w_{0}\chi_{\alpha}) \arrow[rd," c B_{\alpha}"]
 \arrow[r, "H_{\alpha}"] & 
\LL(-\lambda-\alpha) \arrow[r,"\phi_{\alpha,\lambda}"] & 
\underline{\Ver_{B}(-\lambda)} \\
 & \LL(-\lambda) \arrow[ru] & 
\end{tikzcd}
$$
where $c=\langle \lambda, \alpha^{\vee} \rangle$ and $B_{\alpha}$ is non-zero.
We use the same method as in \Cref{power-Hasses}.
We see both compositions as maps
$\overline{Bs_{\alpha}w_0B} \to \Ver_{B}(-\lambda) \otimes \LL(\lambda) \otimes \LL(\chi_{\alpha},-w_{0}\chi_{\alpha})$ 
which are equivariant for the left and right $B$ action on $\overline{Bs_{\alpha}w_0B}$.
The first one extends to $G$ and it is given by
$w_0 \mapsto v=(x_{-\alpha} \otimes -\lambda) \otimes \lambda \otimes (\chi_{\alpha},- w_{0}\chi_{\alpha})$, 
and the second one by
 $s_{\alpha}w_0 \mapsto (c  \otimes -\lambda) \otimes \lambda \otimes (\chi_{\alpha},-w_{0}\chi_{\alpha})$.
 Note that a section with the same weight as $B_{\alpha}$
is unique up to a constant, since it is determined by its restriction to the open dense $Bs_{\alpha}w_0B \subseteq \overline{Bs_{\alpha}w_0B}$,
 and then by the $B$ bi-equivariance it is determined by its value at $s_{\alpha}w_0 \in Bs_{\alpha}w_0B$. 
 We compute that the value at $s_{\alpha}w_0$
 of both maps agrees. 
This would prove that such a factorization exists over $\Sbt_{\alpha}$, and since the top arrow extends to 
$\overline{\Sbt}_{\alpha}$, for any such choice of $B_{\alpha}$ it 
would automatically extend to $\overline{\Sbt}_{\alpha}$.
 Moreover, by construction $B_{\alpha}$ is non-vanishing on 
$\Sbt_{\alpha}$.
As in the proof of \Cref{power-Hasses}, the value of the top map at $s_{\alpha}w_0$ is given by the value at $0$ of $F_{\alpha}(t)=
\rho_1(u_{\alpha}(t^{-1})\alpha^{\vee}(t^{-1}))\rho_2(u_{-w_0\alpha}(-t^{-1}))v$, where $\rho_{i}$
are the actions through the respective factors of $B \times B$. Since $\rho_2$ is inflated from $T$
we get 
$F_{\alpha}(t)=t^{\langle \chi_{\alpha},\alpha^{\vee} \rangle}(v+t^{-1}x_{\alpha}v)$, so that 
$F_{\alpha}(0)=(\langle \lambda, \alpha^{\vee} \rangle \otimes -\lambda \otimes \lambda \otimes (\chi_{\alpha},- w_{0}\chi_{\alpha})$,
since $h_{\alpha}=[x_{\alpha},x_{-\alpha}]$ acts on a vector of weight $\lambda$ as the scalar 
$\langle \lambda, \alpha^{\vee} \rangle$. 

For $(2)$, the if part follows from \Cref{theta-kills-Hasse} and the observation that for 
$p \mid \langle \lambda, \alpha^{\vee} \rangle$, 
the map $\phi_{\alpha,\lambda}$ is already $B$-equivariant. If $H_{\alpha} \mid \theta_{\alpha}(f)$ and $H_{\alpha} \nmid f$, 
then since $H_{\alpha}$ has a simple zero $f_{\overline{D}_{\alpha}} \neq 0$, so that by part $(1)$
 $\theta_{\alpha}(f)_{\mid \overline{D}_{\alpha}}=
\langle \lambda, \alpha^{\vee} \rangle B_{\alpha} f_{\overline{D}_{\alpha}} \neq 0$, which implies
 $p \mid \langle \lambda, \alpha^{\vee} \rangle$.
 To prove $(3)$, we claim that any non-zero $f \in \H^0(\flag,\LL(w_{0,M}\lambda))$
with generic $\lambda \in X_{1}(T)$ cannot be divisible by $H_{\alpha}$.
Let $\beta=\sigma^{-1}\alpha \in \Delta$. 
By the compatibility of the action of Frobenius on the root system
 $\langle w_{0,M}h_{\alpha}, \beta^{\vee} \rangle=
p \langle \chi_{\alpha}, \alpha^{\vee} \rangle- \langle \sigma w_{0,M}\chi_{\alpha},\alpha^{\vee} \rangle=
p-\langle \sigma w_{0,M}\chi_{\alpha},\alpha^{\vee} \rangle$, while for $\gamma \in \Delta\setminus \{\beta\}$
$\langle w_{0,M}h_{\alpha}, \gamma^{\vee} \rangle=-\langle w_{0,M}\chi_{\alpha}, \gamma^{\vee} \rangle$. Therefore, 
$\langle w_{0,M}h_{\alpha}, \beta^{\vee} \rangle-p$ and $\langle w_{0,M}h_{\alpha}, \gamma^{\vee} \rangle$
only depends on $G_{\F_p}$ as a reductive group. 
Suppose that \Cref{assumption-vanishing}(2) holds for a 
constant $m \ge 0$.
 Then we obtain the claim for
$\lambda \in X_1(T)$ satisfying 
$\langle \lambda+w_{0,M}\chi_{\alpha}, \beta^{\vee} \rangle< p-m$ if we assume that
$\langle w_{0,M}\chi_{\alpha}, \gamma^{\vee} \rangle> -m$ for all $\gamma \in \Delta \backslash \{\alpha\}$,
which 
we can do by modifying $m$.
Moreover, by part $(2)$ if $\lambda$
also 
satisfies 
$p \nmid \langle w_{0,M}\lambda, \alpha^{\vee} \rangle$  then
the map $\theta_{\alpha}$ is injective on global sections. All these conditions on $\lambda$ only depend 
on $G_{\F_p}$, and given $G_{\Q}$ there are finitely many options for what $G_{\F_p}$ can be. 
Thus, $\theta_{\alpha}$ 
is injective for $\lambda \in X_1(T)$ generic. 
\end{proof}
\end{theorem}

\begin{remark}
 \begin{enumerate}
 \item Another way to see that $\theta_{\alpha, \overline{D}_{\alpha}}$ is linear is that since it is a degree $1$ 
differential operator satisfying the Leibniz rule,
 it is enough to check that when applied to $\lambda=0$ the restriction of $\theta_{\alpha}$ to $\overline{D}_{\alpha}$ is $0$. 
 But this is clear since in that case $\tilde{\theta}_{\alpha}$ already extends to $\flag$, 
 as the map $-\alpha \to \mathfrak{g}/\mathfrak{b}$ is $B$-equivariant.
 \item It is not necessarily the case that $B_{\alpha}$ is a partial Hasse invariant for $\overline{D}_{\alpha}$, i.e. 
 it must not vanish on $D_{\alpha}$. For instance, for $G=\GSp_4$, and $\alpha=(1,-1)$, $B_{\alpha}$ comes from a map $(\omega/\LL) \to (\omega\LL)^p$
  whose vanishing is given by the condition that $\omega$ is contained in the rank $3$ piece of the conjugate filtration. This is one of the two lower dimensional closed stratum 
 contained in $\overline{D}_{\alpha}$.
\end{enumerate}
\end{remark}

We can go further in some cases where the restriction to a stratum of a basic 
theta operator is still a differential operator. 
Even though we don't attempt to extensively understand differential operators on general
codimension $1$ strata (it is not 
even clear how to understand the tangent bundle of $\overline{D}_{\alpha}$ in general),
we can at least understand the ones that come by restriction from $\flag$. 

\begin{prop} \label{general-restriction-strata}
Let $\alpha \in \Delta$, $\gamma \in \Phi^{+}$ such that $\gamma + \alpha \in \Phi^{+}$.
Let $[x_{\alpha},x_{-\gamma-\alpha}]=c_{\alpha,\gamma}x_{-\gamma}$, and 
$c:=c_{\alpha,\gamma}\frac{(n_{\alpha,\gamma+\alpha} -1)!}{n_{\alpha,\gamma+\alpha}!}$, 
where $n_{\alpha,\gamma+\alpha}$ is as in \Cref{defn-basic-theta}.
Assume that $H_{\alpha}$ has a simple zero.
Then 
on $\LL(\lambda)$
$$
H\theta_{\gamma+\alpha, \overline{D}_{\alpha}}=c B_{\alpha} \circ \theta_{\gamma,\overline{D}_{\alpha}},
$$
where $B_{\alpha}$ is a section of an automorphic line bundle over $\overline{D}_{\alpha}$
defined in \Cref{restriction-simple-theta}, and $H$ is the product of the 
Hasse invariants that are used in the definition of $\theta_\gamma$, but not on $\theta_{\gamma+\alpha}$. 
\begin{proof}
We use \Cref{linearized-theta} to
reduce to the statement at the level of linearized operators, which reduces to the factorization of 
the following diagram of vector bundles 
$$
\begin{tikzcd}
\LL(-\lambda-\gamma-\alpha) \otimes \LL(\chi,-w_0 \chi) \arrow[d,"c B_{\alpha}"] \arrow[r, "H H_{\gamma+\alpha}"] & 
\LL(-\lambda-\gamma-\alpha) \arrow[r,"\phi_{\gamma+\alpha,\lambda}"] & 
\underline{\Ver_{B}(-\lambda)} \\
\LL(-\lambda-\gamma) \otimes \LL(\chi_{1},-w_0\chi_{1}) \arrow[r,"H_{\gamma}"]
 & \LL(-\lambda-\gamma) \arrow[ru,"\phi_{\gamma,\lambda}"] & 
\end{tikzcd}
$$
on $\overline{\Sbt}_{\alpha}$. We use $H_{\gamma}$ to denote 
the product of the Hasse invariants used in $\theta_{\gamma}$, $\chi$ for the sum of $\chi_{\alpha_i}$
used in $HH_{\gamma+\alpha}$, and similarly $\chi_1$ for $H_{\gamma}$.
 We only need to prove that the
diagram commutes on $\Sbt_{\alpha}$, since every arrow extends to $\overline{\Sbt}_{\alpha}$.
For that consider both as maps
$\overline{Bs_{\alpha}w_0B} \to \Ver_{B}(-\lambda) \otimes \LL(\lambda+\gamma+\alpha) \otimes \LL(-\chi,w_0\chi)$.
The one on top 
 extends to $G$, and it is given by
$w_0 \mapsto v_1=(x_{-\gamma-\alpha} \otimes -\lambda) \otimes (\lambda+\gamma+\alpha) \otimes (-\chi,w_0\chi)$.
As in \Cref{power-Hasses} the value at $s_{\alpha}w_0$ of the top map is given by the value at $0$ of 
$F_{\alpha}(t)=t^{\langle \chi,\alpha^{\vee} \rangle}
(v_1+t^{-1}x_{\alpha}v_1+\frac{1}{2}t^{-2}x^2_{\alpha}v_1+ \ldots)$.
For the bottom one, by the construction of $B_{\alpha}$, its value at $s_{\alpha}w_0$ is given by 
the value at $0$ of 
$t^{\langle \chi_1,\alpha^{\vee} \rangle}(u_1+t^{-1}x_{\alpha}u_1+\frac{1}{2}t^{-2}x^2_{\alpha}u_1+ \ldots)$,
where $u_{1}=c(x_{-\gamma} \otimes -\lambda) \otimes (\lambda+\gamma+\alpha) \otimes (-\chi,w_0\chi)$. 
By the assumption that $H_{\alpha}$
has a simple zero we have that
 $n_{\alpha,\gamma+\alpha}=\langle \chi,\alpha^{\vee} \rangle=\langle \chi_1,\alpha^{\vee} \rangle+1$, 
since $\theta_{\gamma+\alpha}$ uses one more copy of $H_{\alpha}$ than $\theta_{\gamma}$ by \Cref{power-Hasses}, 
and $H$ does not contain $H_{\alpha}$.
Using that $x_{\alpha}v_1=c_{\alpha,\gamma}u_1$, and multiplying by 
$c B_{\alpha}$ gives the desired equality. 
 
\end{proof}
\end{prop}

\begin{remark}
By \Cref{theta-kills-Hasse}(2)
we have $H\theta_{\gamma,\alpha}=c_{-\alpha,-\gamma}
[\theta_{\alpha},\theta_{\gamma}]$. By 
computing $\langle  h_{\alpha}+\alpha,\alpha^{\vee} \rangle$, together 
with \Cref{general-restriction-strata} and \Cref{restriction-simple-theta} it shows that 
$\theta_{\gamma,\overline{D}_{\alpha}}$ commutes with $B_{\alpha}$. 
Conversely, one could attempt to prove \Cref{general-restriction-strata} by proving that 
$\theta_{\gamma,\overline{D}_{\alpha}}$ and $B_{\alpha}$ commute. 
\end{remark}

\begin{remark} \label{compactifications-ok}
All the results of this section so far can be extended to $\flag^{\tor}$, since one can extend the Zip map $\Shbar \to \GZip$
to $\Shbar^{\tor}$ \cite[Thm I.2.5]{Goldring-Koskivirta-Galois}, so that we obtain a Zip map 
$\flag^{\tor} \to \GF$. To construct the basic theta operators on $\flag^{\tor}$ we use the extended Zip map 
and that 
\Cref{creator-diff} extends to toroidal compactifications.
The statements of \Cref{theta-kills-Hasse} and \Cref{pth-power-relation} 
can be proved on an open dense subset of $\flag^{\tor}$, 
so they automatically extend to $\flag^{\tor}$. 
The rest of the results (e.g. \Cref{restriction-simple-theta}, \Cref{general-restriction-strata} or 
\Cref{easy-lemma-theta}) are proved by working on 
$\GF \to \Sbt$ and then pulling back to $\flag$, so they 
also extend to the toroidal compactification. 
\end{remark}

Combining the previous results we get a general result about injectivity of basic theta operators 
in the $p$-restricted region. 

\begin{theorem} \label{basic-theta-injective}
Let $\gamma \in \Phi^{+}$, write it as $\gamma=\sum n_\alpha \alpha$ in terms of simple roots, with 
$n_{\alpha} \ge 1$. Suppose that every $\alpha$ appearing in the sum satisfies that $H_{\alpha}$ has a simple zero.
Further assume that \Cref{assumption-vanishing}(2) holds.
Then 
$\theta_{\gamma}$ is injective on $\H^{0}(\flag^\tor,\LL(w_{0,M}\lambda))=\H^0(\Shbar^\tor,\omega(\lambda))$ for 
 generic $\lambda \in X_1(T)$.
\begin{proof}
We use induction on the number of simple roots appearing in $\gamma$. If $\gamma \in \Delta$
this is \Cref{restriction-simple-theta}(3). Otherwise, let $\alpha \in \Delta$ appearing in $\gamma$
such that $\gamma-\alpha \in \Phi^{+}$. Suppose that 
$f \in \H^0(\flag^\tor,\LL(w_{0,M} \lambda))$ is non-zero and it is in the kernel of $\theta_{\gamma}$.
By the reasoning 
as in the proof of \Cref{restriction-simple-theta} for a generic $\lambda \in X_{1}(T)$, $f$ cannot be divisible by 
any Hasse invariant. Moreover, by \Cref{general-restriction-strata} 
$\theta_{\gamma}(f)_{\mid \overline{D}_{\alpha}}=
cB_{\alpha}\theta_{\gamma-\alpha}(f)_{\mid \overline{D}_{\alpha}}=0$, 
which implies that $H_{\alpha} \mid \theta_{\gamma-\alpha}(f)$. Here we use that
from case by case inspection $\lvert c \rvert \le 2$, 
and that we can assume that $p \ge 3$ since the statement concerns generic $\lambda$. 
We repeat the process of restricting to $\overline{D}_{\alpha}$ until $\gamma-n_{\alpha}\alpha \in \Phi^{+}$ does not contain $\alpha$,
and we obtain
$H_{\alpha} \mid \theta_{\gamma-n_{\alpha}\alpha}(f)$.
We claim that \Cref{assumption-vanishing}(2) applies to the weight of $\frac{1}{H_{\alpha}}\theta_{\gamma-n_{\alpha}\alpha}(f)$ 
for generic $\lambda \in X_1(T)$, which implies that $\theta_{\gamma-n_{\alpha}\alpha}(f)=0$.
Let $\mu=w_{0,M}(h_{\gamma-n_{\alpha}\alpha}-h_{\alpha})$ be the weight increase of $\frac{1}{H_{\alpha}}\theta_{\gamma-n_{\alpha}\alpha}$
on $\Shbar$,
where $h_{\gamma-n_{\alpha}\alpha}$ is the sum of the weights of 
the Hasse invariants appearing in $\theta_{\gamma-n_{\alpha}\alpha}$.
Then as in the proof of \Cref{restriction-simple-theta}(3)
$$
\langle \mu, \sigma^{-1}\alpha^{\vee} \rangle=-p+\langle \sigma w_{0,M}\chi_{\alpha},\alpha^{\vee} \rangle+
\sum_{\alpha'\neq \alpha} p\langle \chi_{\alpha'},\alpha^{\vee} \rangle-
\langle \sigma w_{0,M}\chi_{\alpha'},\alpha^{\vee} \rangle,
$$
where the sum runs along the simple roots of 
$\gamma-n_{\alpha}\alpha$ whose Hasse invariant appears in $h_{\gamma-n_{\alpha}\alpha}$, possibly with multiplicity.
 Since $\langle \chi_{\alpha'},\alpha^{\vee} \rangle=0$
 we see that 
$\langle \mu,\sigma^{-1}\alpha^{\vee} \rangle+p$ only depends on $G_{\F_p}$ as an reductive group.
Similarly, if $\beta^{\vee} \neq  \sigma^{-1}\alpha^{\vee}$ 
is a simple coroot 
$$
\langle \mu, \beta^{\vee} \rangle=\langle w_{0,M}\chi_{\alpha},\beta^{\vee} \rangle+
\sum_{\alpha'\neq \alpha} p\langle \chi_{\alpha'},\sigma\beta^{\vee} \rangle-
\langle w_{0,M}\chi_{\alpha'},\beta^{\vee} \rangle.
$$
We can check that $\sum_{\alpha'\neq \alpha} p\langle \chi_{\alpha'},\sigma\beta^{\vee} \rangle \le 2p$, since
$\sigma \beta \in \Delta$ and 
$H_{\sigma \beta}$ can only appear at most twice in the weight increase of a basic theta operator, by \Cref{table-lie}.
The rest of the expression only depends on $G_{\F_p}$.
 Given $G_{\Q}$, there are only finitely
many possibilities for what $G_{\F_p}$ can be, so that for generic $\lambda \in X_1(T)$, $\lambda+\mu$ satisfies 
the hypothesis of \Cref{assumption-vanishing}(2). Namely, there exists $m\ge 0$ only depending on $G_{\Q}$ such that 
$\langle \lambda+\mu, \sigma^{-1}\alpha^{\vee} \rangle< -m$ and $\langle \lambda+\mu, \beta^{\vee} \rangle< 3p+m$ for 
all other simple roots $\beta$. 
 Therefore $\theta_{\gamma-n_{\alpha}\alpha}(f)=0$,
and we conclude by the
induction hypothesis. 
\end{proof}
\end{theorem}

In fact, the injectivity result does apply to a larger set of weights than $\lambda \in X_1(T)$. For 
if $f \in \H^0(\Shbar^\tor,\omega(\lambda))$ is not divisible by any 
Hasse invariant, the result only depends on $\langle \lambda,\sigma^{-1}\alpha^{\vee}\rangle$ being small 
compared to $p$ for all but one of the simple roots $\alpha$ appearing in $\gamma$, and on the sharpness of the 
results of \Cref{assumption-vanishing}(2). It can also happen that for some of the simple roots one does not 
need to use \Cref{assumption-vanishing}(2), and one just uses that all the weights have to be $M$-dominant. 
Further, even if $f$ is divisible by a Hasse invariant $H_{\alpha}$, it might be the case that $f/H_{\alpha}$
has a weight which is generic enough for the theorem to hold.
However, the genericity conditions on $\lambda$ will keep getting worse 
as one gets further away from the $p$-restricted region, and at around the point where 
$\langle \lambda, \alpha^{\vee} \rangle$ is equal to $p^2$ for all $\alpha \in \Delta$
the theta operator will cease to be injective, since one could find sections 
in the kernel by multiplying several Hasse invariants. 
One could compute this more precisely in any particular example 
appearing in \Cref{known-vanishing}.

\section{Theta linkage maps} \label{section4}

\subsection{Construction on $\flag$.}
We construct non-zero maps between $B$-Verma modules in characteristic $p$  
whose weights are linked, hence giving the name to their associated theta linkage maps. 
We use the following notation.
\begin{defn}[Linkage relation]
Let $\lambda, \mu \in X^*(T)$.
For $\gamma \in \Phi^{+}$ and $n \in \Z$
let $s_{\gamma,n} \cdot \lambda:=\lambda+(pn-\langle \lambda+ \rho,\gamma^{\vee} \rangle) \gamma$.
We say that $\lambda \uparrow_{\gamma} \mu$ if there exists $n \in \Z$ such that $\mu=s_{\gamma,n} \cdot \lambda$,
$\lambda \le \mu$, and $n$ is the smallest integer $m$ such that $\lambda \le s_{\gamma,m} \cdot \lambda$. 
Then $\lambda \uparrow \mu$ if there exists a chain 
$\lambda \uparrow_{\gamma_1} \lambda_1 \ldots \uparrow_{\gamma_n} \mu$. 
\end{defn}
Geometrically $\lambda \uparrow_{\gamma} \mu$
if and only if $\mu$ is the reflection of $\lambda$ in the positive root direction across the closest wall
 defined by $\gamma^{\vee}$.
We start by recalling the classical result for maps of Verma modules in characteristic $0$.

\begin{prop}[\cite{Verma}, \cite{BGG}] \label{Hom-Verma-char0}
Let $\lambda_{1}, \lambda_2 \in X^*(T)$, then 
$\Hom_{(U\mathfrak{g},B)_{\mathbb{C}}}(\Ver_{B}(\lambda_1), \Ver_{B}(\lambda_2))$ 
	is non-zero precisely when there is a sequence $\lambda_1 \le \mu_1 \ldots \le \mu_n=\lambda_2$ 
	such that $\mu_{i+1}=s_{\gamma_i} \cdot \mu_{i}$ for $\gamma_i \in \Phi^{+}$, 
	and in that case the Hom space is $1$-dimensional. 
\end{prop}

In particular 
	 there is no non-zero map between Verma modules of different weights, both of which are $G$-dominant, 
   since $W$ acts transitively on the Weyl chambers. However, in characteristic $p$
	 we get many more maps. 

\begin{prop} \label{Verma-linkage}
Let $\lambda, \mu \in X^{*}(T)$. 
\begin{enumerate}
\item The Hom space of Verma modules in characteristic $p$ is 
 invariant under $p$-translation, that is, for any $\nu \in X^{*}(T)$
$$
\Hom_{(U\mathfrak{g},B)_{\Fpbar}}(\Ver_B(\lambda),\Ver_B(\mu))=
\Hom_{(U\mathfrak{g},B)_{\Fpbar}}(\Ver_B(\lambda+p\nu),\Ver_B(\mu+p\nu)).
$$
\item Assume that there exists $\gamma \in \Phi^{+}$ such that 
$\lambda \uparrow_{\gamma,n} \mu$, and that there exists some $\nu \in X^{*}(T)$
such that $\langle \nu, \gamma^{\vee} \rangle \mid n$. In particular, such a $\nu$ always exists 
if $G^{\textnormal{der}}$ is simply connected. 
 Then we can construct an explicit $1$-dimensional subspace of 
$$
\langle \phi(\lambda \uparrow \mu) \rangle \subseteq \Hom_{(U\mathfrak{g},B)_{\Fpbar}}(\Ver_B(\lambda),\Ver_B(\mu)) \neq 0.
$$
In general if $G^{\text{der}}$ is simply connected and $\lambda \uparrow \mu$, then 
$\Hom_{(U\mathfrak{g},B)_{\Fpbar}}(\Ver_B(\lambda),\Ver_B(\mu)) \neq 0$.
\end{enumerate}

\begin{proof}
For $(1)$, let $f \in \Hom_{(U\mathfrak{g},B)_{\Fpbar}}(\Ver_B(\lambda),\Ver_B(\mu))$ be given by 
$x_{f} \in U\mathfrak{u}^{-}_{B}$. We claim that $x_f$ also defines a map in 
$\Hom_{(U\mathfrak{g},B)_{\Fpbar}}(\Ver_B(\lambda+p\nu),\Ver_B(\mu+p\nu))$, which proves the claim by reversing the process. 
We have to prove that $g: \lambda + p\nu \to \Ver_{B}(\mu +p\nu)$ is $B$-equivariant,
 which is the same as being equivariant for the algebra of distributions $U(B)$. 
 We can reduce it to being equivariant for elements of the form 
$\frac{y^n}{n!}$ for positive roots $y$, and $H \choose n$ for $H \in \mathfrak{h}$. For the latter 
it follows since the weight increase of both maps is the same. For the former, we have to prove that 
$\frac{y^n}{n!} x_{f} \otimes (\mu+p\nu)=0$. Write $\frac{y^n}{n!}x_f=x_g x_h$ for 
$x_g \in U\mathfrak{u}^{-}_B$ and $x_h \in U(B)$ by the following procedure. 
By induction on the weight of $y^n$ and the degree of $x_f$ (with respect to a choice of PBW basis) it suffices to prove that $\frac{y^n}{n!}x \in 
x\frac{y^n}{n!}+U(B)U\mathfrak{u}^{-}_{B}$
for $x \in \mathfrak{u}^{-}_{B}$, where the terms in $U(B)$ have smaller weight than $y^n$.
This follows from the PBW theorem on $U(G)$, and by
observing that in the expression for $[y^n,x]$ in a PBW basis
the degree of any element of $\mathfrak{u}^{-}_{B}$ appearing can be at most $1$. Since $f$ is equivariant 
we have that $x_{g} \otimes x_{h} \mu=0$. Then 
$\frac{y^n}{n!} x_{f} \otimes (\mu+p\nu)=x_{g} \otimes x_{h} (\mu+p\nu)=0$, since 
$x_{h} p\nu=0$. 

For $(2)$
choose $\nu \in X^{*}(T)$ such that $n=\langle \nu, \gamma^{\vee} \rangle m$
for some $m \ge 1$. Consider $\tilde{\lambda}\coloneqq \lambda-pm\nu$, then there exists a unique 
up to scalar non-zero map $\phi: \tilde{\lambda} \to \Ver_B(s_{\gamma,0} \cdot \tilde{\lambda})$ over 
$\overline{\Q}_p$
by \Cref{Hom-Verma-char0},
given by an element $\tilde{x}_f \in U \mathfrak{u}^{-}_{B,\overline{\Q}_p}$. 
After rescaling $\tilde{x}_f$ in some PBW basis, choose $x_f$ over the ring of integers $\mathcal{O}$
such that $x_f$ is 
not $0$ mod $p$, this is unique up to a unit. Then $x_f$ defines a map $\tilde{\lambda} \to \Ver_B(s_{\gamma,0} \cdot \tilde{\lambda})$
of $B$-modules over $\mathcal{O}$ since the Verma modules are $p$-torsion free.
 Since $\tilde{\lambda}-\lambda=
s_{\gamma,0} \cdot \tilde{\lambda}-s_{\gamma,n} \cdot \lambda=-pm\nu$, part $(1)$ ensures that the map 
$\lambda \to \Ver_B(\mu)_{\Fpbar}$ defined by $\lambda \mapsto x_f \otimes \mu$ is $B$-equivariant.
We claim that the one-dimensional subspace generated by this map is independent  
of the choice of $\nu$. This amounts to the fact that in characteristic $0$ the maps 
$\Ver_{B}(\lambda) \to \Ver_{B}(s_{\gamma} \cdot \lambda)$ are given by an element of $U \mathfrak{b}^{-}$ 
for all $\lambda$ lying in a hyperplane $\langle \lambda, \gamma^{\vee} \rangle=r$ \cite[4.12]{Humphreys}. 
For two valid $\nu_{1,2}$
the maps in characteristic zero
$\Ver_{B}(\lambda-pn\nu_{i})\to \Ver_{B}(s_{\gamma}\cdot(\lambda-pn\nu_{i}))$ are then given by 
the same element $h=h_1 \cdot h_2 \in U \mathfrak{b}^{-}=U\mathfrak{u}^{-}_{B} \cdot U\mathfrak{h}$.
After rescaling to get a non-zero map over $\Fpbar$, we have that $h_2$ acts by the same scalar on both 
$s_{\gamma}(\lambda-pn\nu_{i})$,
since the difference of their weights is a multiple of $p$.
Thus, we end up with the same element in $U \mathfrak{u}^{-}_{B,\Fpbar}$. Finally, the last statement in $(2)$ follows from the first 
since $U\mathfrak{u}^{-}_{B,\Fpbar}$ has no non-trivial zero-divisors,
so we can compose the previously defined maps. 
\end{proof}
\end{prop}

\begin{remark}
\begin{enumerate}
\item If there is a non-zero map $\Ver_B(\lambda) \to \Ver_B(\mu)$ over $\overline{\F}_p$ then 
$\mu \in W_{\text{aff}} \cdot \lambda$ by considering action of the Harish-Chandra center of $U\mathfrak{g}$.
Moreover, if the projection to $\Ver^0_B(\lambda) \to \Ver^0_B(\mu)$ is non-zero and 
$\lambda \in X_{1}(T)$, then $\lambda \uparrow \mu$ by 
the linkage principle on $G_1B$-representations \cite[Cor 9.12]{Janzten-book},
since the simple $L(\lambda)$ is the head of $\Ver^0_B(\lambda)$.
\item \Cref{Verma-linkage} 
doesn't prove that for two linked weights $\lambda \uparrow \mu$ the map constructed 
above is independent of a choice of sequence 
$\lambda \uparrow_{\gamma_1} \lambda_1 \uparrow_{\gamma_2} \ldots \uparrow_{\gamma_n} \mu$.
In general, we don't know if the Hom space is at most
one-dimensional as in the characteristic $0$ case. A key difference in characteristic $p$ is that 
Verma modules are not of finite length as $(U\mathfrak{g},B)$-modules. E.g. for $G=\GL_2$ 
for each $n \ge 0$ the submodule generated by $\{x^i_{-} \otimes \lambda : i \ge np \} \subseteq \Ver_{B}(\lambda)$
is a $(U\mathfrak{g},B)$-module, forming
an infinite descending chain of submodules. 
\end{enumerate}
\end{remark}

We define theta linkage map by plugging the maps of \Cref{Verma-linkage} into \Cref{creator-diff}. 
\begin{defn}(Theta linkage maps)
Let $\lambda \uparrow_{\gamma} \mu$ as in part $(2)$ of \Cref{Verma-linkage}. Recall the functor 
$\Phi_{B}: \Hom_{(U\mathfrak{g},B)_{\Fpbar}}(\Ver_{B}(\lambda),\Ver_{B}(\mu)) \to 
\DiffOp_{\flag_{\Fpbar}}(\LL(-\mu),\LL(-\lambda))$ from \Cref{creator-diff}.
We define
$$
\theta^{\lambda \uparrow \mu}:=\Phi_{B}(\phi(\lambda \uparrow \mu)): \LL(-\mu) \to \LL(-\lambda)
$$
for some non-zero choice of $\phi(\lambda \uparrow \mu)$ in \Cref{Verma-linkage}(2). 
 We will refer to these as theta linkage maps. They extend to a toroidal compactification $\flag^{\tor}$.
\end{defn}

\begin{cor} \label{linkage-map-composition-basic}
After post-composing with some powers of Hasse invariants, the theta linkage maps $\theta^{\lambda \uparrow \mu}$ are combinations of the basic theta operators.
Namely, given $\theta^{\lambda \uparrow \mu}$
there exists some $f \in U\mathfrak{u}_{B}$ which we see as a non-commutative polynomial on the $\{x_{\gamma}: \gamma \in \Phi^{+}\}$, and some 
integers $n_{\alpha} \ge 0$ such that 
$$
\prod_{\alpha \in \Delta} H^{n_{\alpha}}_{\alpha} \circ \theta^{\lambda \uparrow \mu}=f([\theta_{\gamma}]_{\gamma \in \Phi^{+}}).
$$
\begin{proof}
This follows from \Cref{composition-U}, which identifies $\theta^{\lambda \uparrow \mu}$ over $U$ with a combination 
of basic theta operators $\tilde{\theta}_{\gamma}$. Since we know that this automatically extends to $\flag$ 
it must be the case that the same combination of $\theta_{\gamma}$ becomes divisible by all the Hasse invariants 
used in the combination. 
\end{proof}
\end{cor}

\begin{notation}[Reparametrization of theta linkage maps according to Serre weights] \label{linkage-notation}
  As alluded to in the introduction the BGG decomposition over $\overline{\C}_p$ implies that for 
$\lambda$ a dominant weight  
$\H^0(\Sh^\tor_{\overline{\C}_p},\omega(\lambda+\eta))\neq 0 \implies 
\H^{d}_{\et}(\Sh_{\overline{\Q}_p},V(\lambda)_{\overline{\C}_p}) \neq 0$,
where $\eta$ is the sum of roots appearing in $\Lie(U_{P})$.   
The same will hold with mod $p$ coefficients provided that we can lift coherent cohomology to characteristic $0$. 
This motivates the following reparametrization of the theta linkage maps. 
Given two weights $\lambda \uparrow_{\gamma,n} \mu$ for some $\gamma \in \Phi^{+}\backslash \Phi^{+}_{M}$ we have 
that $-w_{0,M}\mu-\eta \uparrow_{w_{0,M}\gamma,-n} -w_{0,M}\lambda-\eta$.
 This comes down to the fact 
that $\rho+w_{0,M}\rho=\eta$ and $w_{0,M}\gamma \in \Phi^{+}$. This follow from the identities
 $\rho=\frac{1}{2}\eta+\rho_{M}$, $w_{0,M}\rho_{M}=-\rho_{M}$ and that $w_{0,M}$ preserves 
 the set $\Phi^{+}\backslash \Phi^{+}_{M}$. 
Then we define $\theta_{\lambda \uparrow \mu}$ to be the pushforward along $\pi$ of the operator induced by 
$\phi(- w_{0,M}\mu-\eta \uparrow  -w_{0,M}\lambda-\eta)$
$$
\theta_{\lambda \uparrow \mu}:=\pi_*\Phi_{B}[\phi(- w_{0,M}\mu-\eta \uparrow  -w_{0,M}\lambda-\eta)]: 
\omega(\lambda+\eta) \to \omega(\mu+\eta). 
$$
If $\gamma \in \Phi^{+}_{M}$ then $w_{0,M}\gamma$ is a negative simple root  and 
$-w_{0,M}\lambda-\eta \uparrow_{-w_{0,M}\gamma,n} -w_{0,M}\mu-\eta$. We define $\theta_{\lambda \uparrow \mu}$ 
to be the pushforward along $\pi$ of the operator induced by 
$\phi(- w_{0,M}\lambda-\eta \uparrow  -w_{0,M}\mu-\eta)$
$$
\theta_{\lambda \uparrow \mu}:=\pi_*\Phi_{B}[\phi(- w_{0,M}\lambda-\eta \uparrow  -w_{0,M}\mu-\eta)]
: \omega(\mu+\eta) \to \omega(\lambda+\eta),
$$
so we see that in this case the theta linkage map goes in the opposite direction. 
  \end{notation}

\begin{remark}($p$-translation property of theta linkage maps)
Let $\lambda,\mu, \nu \in X^*(T)$ such that $\lambda \uparrow \mu$. Let $F$ be the absolute Frobenius on $\flag$, 
then by the projection formula $F_{*}\LL(-\lambda-p\nu)=F_*\LL(-\lambda) \otimes \LL(-\nu)$, 
we have
$$
F_* \theta^{\lambda+p\nu \uparrow \mu+p\nu}: F_*\LL(-\mu) \otimes \LL(-\nu) 
\xrightarrow{F_* \theta^{\lambda \uparrow \mu} \otimes \text{id}} F_*\LL(-\lambda) \otimes \LL(-\nu)=F_* \LL(-\mu-p\nu).
$$
This is because by construction and \Cref{linkage-map-composition-basic} 
both $\theta^{\lambda+p\nu \uparrow \mu+p\nu}$ and $\theta^{\lambda \uparrow \mu}$ are constructed using the same 
combination of 
basic theta operators. Thus, to prove the above property we can replace the linkage map by a basic theta operator, 
in which case it follows from the Leibniz rule. Therefore, in order to compute kernels of linkage maps as maps of 
sheaves we can assume that the domain is in $X_1(T)$.
\end{remark}

\begin{example}(Simple theta linkage maps)
In the case that $\lambda \uparrow_{\alpha,n} \mu$ with $\alpha \in \Delta$, the associated map of 
Verma modules is easy to describe. Write $N_{\alpha}=pn-\langle \lambda+\rho,\alpha^{\vee} \rangle$ for some integer $n$ 
with 
the property that $1 \le N_{\alpha} \le p$.
Then the map $-\mu \to \Ver_{B}(-\lambda)$ is given by $-\mu \mapsto x^{N_{\alpha}}_{-\alpha} \otimes -\lambda$. 
One can explicitly check that this is $B$-equivariant using the commutation relation 
$x_{\alpha}x^{N}_{-\alpha}=x^{N}_{-\alpha}x_{\alpha}+Nx^{N-1}_{-\alpha}h_{\alpha}-N(N-1)x^{N-1}_{-\alpha}$
and that $\langle \rho, \alpha^{\vee} \rangle=1$. Therefore,
by \Cref{composition-U} 
$\theta^{\lambda \uparrow \mu}=\frac{1}{H^{N_{\alpha}}_{\alpha}} \theta^{N_{\alpha}}_{\alpha}$. 
We will call these \textit{simple theta linkage maps}, they are the ones that can be understood explicitly 
in terms of the basic theta operators. By composing two such simple affine Weyl reflections 
we also get that $\lambda \uparrow_{\alpha} \lambda+p\alpha$, so that $H^p_{\alpha} \mid \theta^{p}_{\alpha}$. 
Then \Cref{pth-power-relation} says that this is $0$ in some cases, or another basic theta operator in others.  
\end{example}

We also note, that even though that $x^{p}_{-\gamma}$ for $\gamma \in \Phi^{+}$ is always a central element in $U\mathfrak{g}$, it does 
not define a map of Verma modules unless $\gamma$ is simple. This is because the resulting map won't be $B$-equivariant, one can see that this already 
fails by considering equivariance for $\frac{x^p_{\alpha}}{p!}$ for $\alpha \in \Delta$ appearing in $\gamma$.

\begin{remark} \label{remark-injective}
We could also prove that in the example above
$\frac{1}{H^{N_{\alpha}}_{\alpha}}\theta^{N_{\alpha}}_{\alpha}$ is well-defined without appealing to 
\Cref{composition-U} and its relation to 
$\theta^{w_{0,M}\lambda \uparrow w_{0,M}\mu}$. We can prove this by combining
\Cref{restriction-simple-theta}(2) and \Cref{pth-power-relation} as in \cite[Thm 3.10]{paper}, 
in a way closely resembling the proof of the classical theta cycle map. 
Let $f \in \LL(\lambda)$ be a local section, and suppose first that it is not divisible by $H_{\alpha}$.
 By the first one 
we know that for $n<N_{\alpha}$ we have
$H_{\alpha} \nmid \theta^n_{\alpha}(f)$, and some power $k$ of $H_{\alpha}$ divides 
$\theta^{N_{\alpha}}_\alpha(f)$. 
We also know by \Cref{pth-power-relation} 
that $H^p_{\alpha} \mid \theta^p_{\alpha}(f)$, so by repeating the argument we get that 
$k \ge N_{\alpha}$. The case of $H_{\alpha} \mid f$ 
follows by induction on $\langle \lambda,\alpha^{\vee} \rangle
\bmod{p}$, with the case where $p \mid \langle \lambda,\alpha^{\vee} \rangle$ being the base case. 
\end{remark}

In the case of simple linkage maps we can understand their behaviour on global sections. 

\begin{theorem} \label{injective-linkage}
Let $p \ge 3$.
Assume that $G^{\der}$ is simply connected, and that \Cref{assumption-vanishing}(2) holds.
Let $G_i$ be an almost-simple factor of $G^{\text{der}}_{\Fpbar}$, 
$\alpha \in \Delta_{G_i}$ which is not on $\Lie(M)$. 
 Let $\lambda \in X_1(T)$ and $\mu \in X^*(T)$
such that $\lambda \uparrow_{w_{0,M}\alpha} \mu$. 
Then 
$$
\theta_{\lambda \uparrow \mu} : \H^0(\Shbar^\tor,\omega(\lambda+\eta)) \to \H^0(\Shbar^\tor,\omega(\mu+\eta))
$$
is injective for generic $\lambda \in X_1(T)$. 
\begin{proof}
From \Cref{linkage-notation} we see that $\theta_{\lambda \uparrow \mu}$ is a simple theta linkage map. 
Using the identities $\frac{1}{H^{p}_\alpha}\theta^p_{\alpha}=\theta_{\mu \uparrow \lambda+pw_{0,M}\alpha} \circ 
\theta_{\lambda \uparrow \mu}$ from \Cref{linkage-map-composition-basic} and 
$\frac{1}{H^{p}_\alpha}\theta^p_{\alpha}=\theta_{\delta}$ from 
\Cref{pth-power-relation} (the case of semisimple rank $1$ is analogous) we can reduce to proving injectivity of $\theta_{\delta}$. 
But then it follows from 
\Cref{basic-theta-injective}.
\end{proof}
\end{theorem}

Note that in the notation from above $w_{0,M}\alpha$ is always 
the longest root of $G_j$, see the beginning of the proof of \Cref{pth-power-relation}. Therefore,
when considering unitary Shimura varieties of different signatures, one always finds
the same $\gamma \in \Phi$ such that \Cref{injective-linkage} applies to $\theta_{\lambda \uparrow_{\gamma} \mu}$.

\subsection{Relation between Borel and parabolic linkage maps} \label{parabolic-thetas}
  We can also define a kind of theta linkage maps on $\Shbar$ by constructing maps between parabolic Verma modules. 
  These are harder to understand, but they occur naturally when considering the de Rham cohomology of 
  the Shimura variety. We prove that we can construct maps of parabolic Verma modules 
  such that the associated operator over $\Shbar$ comes from pushforward 
  from a theta linkage map on $\flag$, which we can understand more explicitly
  in terms of basic theta operators. 

  \begin{lemma} \label{pushforward-Verma}
  Let $\lambda, \mu \in X^*(T)$, and let $f: \lambda \to \Ver_{B}(\mu)$ a $B$-equivariant map. 
  Then there is a unique $P$-equivariant map $g: W(-w_{0,M}\lambda)^{\vee} \to \Ver_{P}W(-w_{0,M}\mu)^{\vee}$
  fitting in the diagram of $B$-equivariant maps
  $$
  \begin{tikzcd}
  \lambda \arrow[r,"f"] \arrow[d,"\phi"] &\Ver_{B}(\mu) \arrow[d] \\
  W(-w_{0,M}\lambda)^{\vee} \arrow[r,"g"] & \Ver_{P}W(-w_{0,M}\mu)^{\vee}.
  \end{tikzcd}
  $$
  The map $\phi$ is given by a choice of a highest weight vector of $W(-w_{0,M}\lambda)^{\vee}$, 
  using that $-\lambda=w_{0,M}(-w_{0,M}\lambda)$ is a lowest weight for $W(-w_{0,M}\lambda)$.
  On the other hand, $\Ver_{B}(\mu) \to \Ver_{P}W(-w_{0,M}\mu)^{\vee}$ is defined by $yx \otimes \mu 
  \mapsto y \otimes x v^{\vee}_{-\mu}$ for $y \in U\mathfrak{u}^{-}_{P}$, $x \in U\mathfrak{u}^{-}_{M}$, 
  and $v_{-\mu}$ a vector of $W(-w_{0,M}\mu)$ of weight $-\mu$.
  \begin{proof}
  The map $f$ induces a differential operator $P_{G/B} \otimes \LL(-\mu) \to \LL(-\lambda)$ by 
  \Cref{verma-flag}. Let $\pi: G/B \to G/P$.
   Its pushforward 
  precomposed with $P_{G/P} \otimes W(-w_{0,M}\mu) \to \pi_{*}[P_{G/B} \otimes \LL(-\mu)]$ (using that $-\mu$ is a lowest 
  weight for $W(-w_{0,M}\mu)$)
  defines a $G$-equivariant differential operator 
  $h: P_{G/P} \otimes W(-w_{0,M}\mu) \to W(-w_{0,M}\lambda)$. 
   We define $g$ as the associated map of Verma modules. 
  Let $\epsilon: \pi^* \pi_{*} \to \text{id}$ be the unit of the adjunction. 
  $$
  \begin{tikzcd}
	\pi^*[P_{G/P} \otimes W(-w_{0,M}\mu)] \arrow[d,"\phi"] \arrow[r,"\pi^*\pi_* \phi"] & \pi^* \pi_*[P_{G/B} \otimes \LL(-\mu)] 
	\arrow[r,"\pi^*\pi_* f^{\vee}"] \arrow[ld,"\epsilon"] & W(-w_{0,M}\lambda) \arrow[d,"\epsilon"]\\
	P_{G/B} \otimes \LL(-\mu) \arrow[rr,"f^{\vee}"] & & \LL(-\lambda)
  \end{tikzcd}
  $$
  To prove $g$ makes the diagram commute we must prove that the big square above commutes, where the composition 
  of the top horizontal arrows corresponds to $\pi^* g$. 
  By naturality of $\epsilon$ we reduce to proving the commutativity of the small triangle, which is clear.
   Conversely, if $g$ fits 
  in the diagram it corresponds to the composition 
  $\pi^*(P_{G/P} \otimes W(-w_{0,M}\mu)) \to \pi^*\pi_{*}(P_{G/B} \otimes \LL(-\mu)) 
  \xrightarrow{\pi^* \pi_* f} \pi^{*} \pi_* \LL(-\lambda)$. Applying $\pi_*$ and the projection formula recovers
   the previous construction, so this proves the uniqueness. 
  \end{proof}
  \end{lemma}

  We will denote such a map $g$ as $\pi_* f$. The maps of parabolic Verma modules that are pushforwards from 
  Borel Verma modules are sometimes denoted as \textit{standard} in the literature,
   e.g. \cite[\S 3]{maps-parabolic-verma}. It is known that there can be
  more maps than the standard one, see \cite[Prop 10.1]{thesis-parabolic-verma} for some examples in characteristic $0$
  for the group  
  $G=\GL_4$.
  \begin{prop} \label{compatibility-pushforward}
	Let $\lambda, \mu \in X^*(T)$, and let $f: \lambda \to \Ver_{B}(\mu)$ a $B$-equivariant map. It induces 
	a differential operator $\theta(f): \LL(-\mu) \to \LL(-\lambda)$. Consider $\pi_* f$ from 
	\Cref{pushforward-Verma}, it induces a differential operator $\theta(\pi_* f):
	 \omega(-w_{0,M}\mu) \to \omega(-w_{0,M} \lambda)$. Then 
	 $$
	 \pi_* \theta(f)=\theta(\pi_*f). 
	 $$

	\begin{proof}
	 We can prove it on PD formal completions of $\Fpbar$-points, where we can use the Grothendieck--Messing 
   period map to reduce to the case of the flag variety, where it follows by definition of $\pi_*f$. 
	\end{proof}
  \end{prop}

  \begin{remark}
  The above shows that given $\lambda \uparrow_{\gamma} \mu$ for $\gamma \in \Phi^{+}_{M}$, 
  $\theta_{\lambda \uparrow \mu}: \omega(\mu+\eta) \to \omega(\lambda+\eta)$ is linear, and it is induced 
  by the map of $P$ representations induced by $\pi_* \phi(-w_{0,M}\lambda-\eta \uparrow 
  -w_{0,M}\mu-\eta)$. 
  \end{remark}

  As a corollary of the above we get a relation between theta linkage maps on $\flag$ and $\Shbar$. 

  \begin{cor} \label{uniqueness-parabolic-linkage}
  Let $\lambda, \mu \in X^*(T)$ such that $\lambda \uparrow_{\gamma} \mu$ for some $\gamma \in \Phi^{+}$.
  \begin{enumerate}
  \item There is a $P$-equivariant map $g:W(-w_{0,M}\mu)^{\vee} \to \Ver_{P}W(-w_{0,M}\lambda)^{\vee}$, 
  and it satisfies 
  $\pi_*\theta^{\lambda \uparrow \mu}=\Phi_{P}(g)$. 
  \item  Moreover, if $\gamma=\alpha \in \Delta$, $\alpha \notin \Lie(M)$, 
  and $W(-w_{0,M}\mu)$ is irreducible, then there is only one $P$-equivariant map 
  $W(-w_{0,M}\mu)^{\vee} \to \Ver_{P}W(-w_{0,M}\lambda)^{\vee}$.
  \end{enumerate}
  \begin{proof}
	For $(1)$ we take $g$ the pushforward of 
	the map of Borel Verma modules $\phi(-\mu \uparrow -\lambda)$ in \Cref{Verma-linkage}.
	 To prove uniqueness in 
	the second case, note that
   $\mu-\lambda$ is a multiple of $\alpha$.  Since $W(-w_{0,M}\mu)^{\vee}$ is generated by the 
  highest weight vector $v_{\mu}$ \cite[2.13(b)]{Janzten-book},
  such a map $g$ must send $v_{\mu}$ to $x^{N_{\alpha}}_{-\alpha} \otimes v_{\lambda}$, 
  where $v_{\lambda}$ is the unique vector of weight $\lambda$. 
  \end{proof}
  \end{cor}

\subsection{Kernel of flag-like linkage maps}
By \Cref{Frobenius-differentials}
a crystalline differential operator $\theta: \LL(\lambda) \to \LL(\mu)$ over $\flag$ is always Frobenius linear, 
so it is natural to try to understand the kernel of $F_{*} \LL(\lambda) \to F_{*}\LL(\mu)$, where $F$ 
is the absolute Frobenius. If one understood this one could reduce the study of whether a theta operator is injective on 
global sections to understanding the sections of the kernel. Sadly we currently don't understand much about this,
for instance, we don't even know if $F_{*}\LL(\lambda)$ is an automorphic vector bundle on $\flag$. 
We can only understand this systematically 
in the case that $\theta$ is a linkage map such its pushforward to $\Shbar$ is linear, and replacing $F$ by the relative Frobenius. 

\begin{defn}
We say that a $B$-equivariant map $\lambda \to \Ver_{B}(\mu)$ is flag-like if it factors through $\Ver_{P/B}(\mu)$.
We say that a theta linkage map is flag-like if its associated map of Verma modules is flag-like.
\end{defn}

In fact, flag-like theta linkage maps are precisely the ones that become linear on $\Shbar$. 

\begin{lemma}
A theta linkage map $\theta$ is flag-like if and only if $\pi_* \theta$ is a linear map. Moreover, in that case 
$\pi_* \theta$ comes from a map of $M$-representations. 
\begin{proof}
Let $\theta$ be given by a map $f: \lambda \to \Ver_{B}(\mu)$. Suppose first that $\phi$ is flag-like 
so that $f$ factors through $\Ver_{P/B}(\mu)$.
Then by weight considerations the map
$\pi_*f: W(-w_{0,M}\lambda)^{\vee} \to \Ver_{P}W(-w_{0,M}\mu)^{\vee}$ must come from a map of $P$-representations 
$g: W(-w_{0,M}\lambda)^{\vee} \to W(-w_{0,M}\mu)^{\vee}$. Then \Cref{compatibility-pushforward} implies that 
$\pi_*\theta=\theta(\pi_*f)$ 
is linear and comes from the map $g$. Conversely, suppose that $\pi_*\theta=\theta(\pi_*f)$ is linear. 
Then $\pi_* f$ comes from a map $g$ as before. By the diagram in \Cref{pushforward-Verma} we see that 
$\mu-\lambda$ must be in the $\Z^{\ge 0}\Phi^{+}_{M}$, which implies that $f$ factors through $\Ver_{P/B}(\mu)$. 
\end{proof}
\end{lemma}

\begin{prop} \label{kernel-flag-linkages}
Let $\phi: \Ver_{P/B}(\lambda) \to \Ver_{P/B}(\mu)$ be a $B$-equivariant map, with associated differential operator 
$\theta: \LL(-\mu) \to \LL(-\lambda)$. Let $F: \flag \to \flag^{(p)}$ be the relative Frobenius with respect to $\Shbar$. 
Let also $\phi^0: \Ver^0_{P/B}(\lambda) \to \Ver^0_{P/B}(\mu)$ be the map induced by $\phi$. Then 
$$
\Ker F^*F_{*} \theta=F_{B}\textnormal{Ker} (\phi^0)^{\vee},
$$
where we identify $F^*F_* \LL(\lambda)$ with $\Ver^{0,\vee}_{P/B}(-\lambda)$ using $e^{0}_{P/B}$ in 
\Cref{canonical-iso-Pm}. 
Moreover, let $\tilde{F}: P/B \to P/B$ be the absolute Frobenius, 
then $\phi$ induces a $G$-equivariant map $\tilde{\theta}: \tilde{F}_{*}\LL_{P/B}(-\mu) \to \tilde{F}_{*}\LL_{P/B}(-\lambda)$. 
Then we can write $\Ker \tilde{\theta}=F_{P/B}(V_{\phi})$ for some $V_{\phi} \in \Rep(B)$. Then 
$$
\Ker F_{*} \theta=\sigma^{*} F_{B}(V_{\phi}),
$$
where $\sigma: \flag^{(p)} \to \flag$ is the base change morphism.
\begin{proof}
The first part follows immediately from \Cref{creator-diff}(1). The second part follows in the same way, by transporting 
it from the flag variety via the Grothendieck--Messing map. 
\end{proof}
\end{prop}

\begin{example} \label{example-kernel-gsp4}
Take the example of $G=\GSp_4$ with $P$ the Siegel parabolic. In the notation of 
\Cref{table-lie} let $\alpha=\alpha_1$. 
 Let $\phi: (-l',-k') \to \Ver_{P/B}(-l,-k)$ be the map inducing the linkage map 
for $\alpha \in \Lie(M)$. Then if $k-l=ap+b$ with $0 \le b \le p-1$ we have 
$$
\Ker \phi^{0,\vee}=(ap,0) \otimes W(l+b,l),
$$
with $\mathfrak{p}$ acting naturally on $W(l+b,l)$ and trivially on $(ap,0)$. By the $P/B$ 
version of \Cref{remark-Frobenius-connection} this ensures that 
$V_{\phi}=(a,0)\otimes W(l+b,l)$. 
\end{example}

\Cref{kernel-flag-linkages} together with the divisibility criteria for simple basic theta operators 
gives a description of the kernel of a flag-like simple basic theta operator. 
\begin{cor} \label{kernel-simple-flag-like}
Let $\alpha \in \Delta \cap \Lie(M)$ and $\lambda \in X^*(T)$.
Let $F: \flag \to \flag^{(p)}$ be the relative Frobenius with respect to $\Shbar$. 
Assume that $H_{\alpha}$ has a simple zero.  Write 
$\langle \lambda, \alpha^{\vee} \rangle=pa+b$ for $1 \le b \le p$.
Then 
$$
F_* \Ker(\theta_{\alpha}: \LL(\lambda) \to \LL(\lambda+h_{\alpha}+\alpha))=
F_* \Ker(\theta_{\alpha}/H_{\alpha}: \LL(\lambda-(p-b)h_{\alpha})\to \LL(\lambda-(p-b)h_{\alpha}+\alpha)),
$$
where the inclusion to $F_{*}\LL(\lambda)$ is given by multiplication by $H^{p-b}_{\alpha}$, and the right-hand side 
can be computed by \Cref{kernel-flag-linkages}.
\begin{proof}
Given a local section $f \in \LL(\lambda)$ in the kernel of $\theta_{\alpha}$ we have to show that 
$H^{p-b}_{\alpha} \mid f$. We proceed by induction on $b$, with 
$b=p$ being the base case. Then we use \Cref{restriction-simple-theta}, \Cref{theta-kills-Hasse}(2)
and the observation that 
$\langle h_{\alpha},\alpha^{\vee}\rangle=\langle pw_{0,M}\sigma^{-1}\chi_{\alpha}-\chi_{\alpha},\alpha^{\vee} \rangle=-1 \bmod {p}$.
\end{proof}
\end{cor}

We currently don't have a good understanding of the kernel of a "non flag-like"
simple basic theta operator, even just on global sections. Essentially we can only understand the 
following case, which is 
a direct generalization of the $\GL_2$ case. 
\begin{example} \label{kernel-parallel-wt}
Let $\Shbar$ be the unitary Shimura variety of signature $(n-1,1)$ as in
\Cref{unitary-compact} with $p$ split 
in the quadratic imaginary field. Let $\alpha=\alpha_{n-1} \in \Delta \setminus \Lie(M)$. 
Let $k_1,k_2 \in \mathbb{Z}$, and write $\omega(\underline{k_1},k_2):=\omega(k_1,k_1,\ldots,k_1,k_2)$. 
Let $F$ be the absolute Frobenius on $\Shbar$. Assume for simplicity that $1 \le k_1 \le p$, then 
$$
\Ker(F_*\pi_*\theta_{\alpha}:F_* \omega(\underline{k_1},pk_2)\to F_* \omega(\underline{k_1}+(p,p-1, \ldots, p-1),pk_2-1))
=\omega(\underline{k_1-p+1},k_2)
$$
via $\omega(\underline{k_1-p+1},k_2) \to F_* \omega(p\underline{(k_1-p+1)},pk_2)
\xrightarrow{H^{p-k_1}_{\alpha}} F_* \omega(\underline{k_1},pk_2)$, where the first map induced by the adjunction 
$F^{*} F_*$. This follows by \Cref{restriction-simple-theta}(2), and the observation that by the Leibniz rule,
on a local basis 
$\theta_{\alpha}$ acting on $\LL(p\underline{(k_1-p+1)},pk_2)$ can be identified with $\theta_{\alpha}$ acting 
on $\mathcal{O}_{\flag}$. Therefore, we can check that
 $\pi_*\theta_{\alpha}$ can be identified with $d:\mathcal{O}_{\Shbar}  \to \Omega^1_{\Shbar}$,
 e.g. using \Cref{compatibility-pushforward},
  whose kernel is precisely the functions that are $p$th powers. 
The same kind of example will work for a unitary Shimura variety of any signature with $p$ totally split
in the CM field,
or in the Siegel case. 
\end{example}

\section{Duality} \label{section5}
We prove a compatibility result between the functors in \Cref{creator-diff} and Serre duality.

\subsection{Duality of differential operators}
Fix $X/S$ a smooth scheme whose fibers are equidimensional of dimension $d$. Let $\omega_X$ be its 
relative canonical bundle. 
\begin{defn}
For a vector bundle $E$ over $X$ let $E^{*}=E^{\vee} \otimes \omega_X$. Let $f: P_{X} \otimes E \to V$
be a (HPD) crystalline differential operator between two vector bundles.
By tensor-hom adjunction, and since $E,V$ are reflexive, we get a map $V^{*} \otimes P_{X} \otimes E \to \omega_X$.
Again by tensor-hom adjunction we get a map 
$f^*: V^{*} \otimes P_X \to E^{*}$ which is $\mathcal{O}$-linear for the 
right $\mathcal{O}$-structure on $P_{X}$, so it defines 
an element of $\DiffOp(V^*, E^*)$. We say that $f^*$ is the dual differential operator of $f$.
The same definition extends to log differential operators. We note that $(f^*)^*=f$ canonically.
\end{defn}

The goal of this subsection is to state the following compatibility with Serre duality. 
\begin{theorem} \label{Serre-duality-diff-op}
Let $S$ an affine scheme and
let $\pi: X \to S$ be a smooth proper scheme whose fibers are equidimensional of dimension $d$. Let $E,V$ be vector bundles 
on $X$, and $f: P \otimes E \to V$ a PD differential operator. Then the following diagram commutes 
$$
\begin{tikzcd}
\textnormal{RHom}(R\pi_* E, \mathcal{O}_S)[-d] \arrow[r,"\sim"] & R\pi_{*}E^{*} \\ 
\textnormal{RHom}(R\pi_* V, \mathcal{O}_S)[-d] \arrow[r,"\sim"] \arrow[u, "\textnormal{RHom}(f)"] 
& R\pi_{*}V^{*} \arrow[u,"f^*"],
\end{tikzcd}
$$
where the horizontal isomorphisms are given by Serre's duality. 
If $S$ is the spectrum of a field, talking cohomology for $0 \le i \le d$ we get that the 
following diagram commutes
$$
\begin{tikzcd}
	\H^i(X,E)^{\vee} \arrow[r,"\sim"] & \H^{d-i}(X,E^{*}) \\ 
	\H^i(X,V)^{\vee} \arrow[r,"\sim"] \arrow[u, "\H^i(f)^{\vee}"] & \H^{d-i}(X,E^{*}) \arrow[u,"f^*"].
	\end{tikzcd}
	$$
The same statements hold for log differential operators for a pair $(X,D)$ over $S$. 
\begin{proof}
The commutativity of the first diagram in the case of Grothendieck's differential operators 
is  proved in  
\cite[\S 7.1]{duality-diff-op}. 
 There is a map 
$\DiffOp(E,V) \to \Hom_{\mathcal{O}_X}(\tilde{P} \otimes E,V) \to \Hom_{\mathcal{O}_{S}}(E,V)$, 
where the first map to Grothendieck's differential operators is compatible with duality.
Since the commutativity of the diagram 
here only concerns 
the image of this composition, it follows from that case. 
To obtain the result on cohomology we take $\H^i$ of the diagram above, and since $S$ is a field 
every complex in $D(\mathcal{O}_{S})$ is a perfect complex with trivial differentials, so we get 
the correct cohomology groups. The statement for log differential operators 
is immediate since they map to usual differential operators. 
\end{proof}
\end{theorem}

\subsection{Dual linkage maps}
Let $Q \in \{P,B\}$ with Levi $M$, $R=\mathcal{O}/p^n$ and $V_{i} \in \Rep_{R}(Q)$. 
Given $\phi: \Ver_{Q}(V_1) \to \Ver_{Q}(V_2)$ a $Q$-equivariant map, it induces a $G$-equivariant 
differential operator $\psi: P_{G/Q} \otimes \mathcal{V}^{\vee}_2 \to \mathcal{V}^{\vee}_1$ by \Cref{verma-flag}.
Let $K_{G/Q}$ denote the canonical bundle. 
 Then $\phi^*: P_{G/Q} 
\otimes \mathcal{V}_1 \otimes K_{G/Q} 
\to \mathcal{V}_2 \otimes K_{G/Q}$ is still $G$-equivariant since the tensor-hom adjunction defining it 
preserves $G$-equivariance. Therefore, it induces a map 
$\phi^*: \Ver_{Q}(V^{\vee}_2 \otimes W(2\rho_{M}-2\rho)) \to \Ver_{Q}(V^{\vee}_1 \otimes W(2\rho_{M}-2\rho))$, 
where $\rho_{M}$ is the half sum of positive roots of $M$, so that $K_{G/Q}=\mathcal{W}(2\rho-2\rho_{M})$.

\begin{defn} \label{Serre-dual-verma}
Given a $Q$-equivariant map $\phi: \Ver_{Q}(V_1) \to \Ver_{Q}(V_2)$ we define its dual as $\phi^*$ constructed above. 
Denote by $V^{*}:=V^{\vee} \otimes W(2\rho-2\rho_{M})$.
\end{defn}

\begin{prop}
Let $R=\mathcal{O}/p^n$, $Q \in \{P,B\}$, $V_{i} \in \Rep_{R}(Q)$,
 and $\phi: \Ver_{Q}(V_1) \to \Ver_{Q}(V_2)$ a $Q$-equivariant map over $R$. 
 Then 
 $$
 \Phi_{Q}(\phi)^{*}=\Phi_{Q}(\phi^*),
 $$
where $\Phi_{Q}$ is the functor of \Cref{creator-diff}. On toroidal compactifications 
we have 
$$
\Phi^{\can}_{Q}(\phi)^{*}=\Phi^{\sub}_{Q}(\phi^*)
$$
and the other way around. 
\begin{proof}
We can check everything at the level of $\overline{\F}_p$-points since it is a map of 
vector bundles. Thus, for $x \in \flag_{Q}(\overline{\F}_p)$, we use the description of 
the canonical isomorphisms defining $\Phi_{Q}$ in terms the Grothendieck--Messing maps 
over $\flag^{\sharp}_{Q,x,R}$ of \Cref{induction-dfn}.  
Since $\phi^*$ by definition produces the dual differential operator 
on $G/Q$ we reduce to the following lemma. We say that a PD formal scheme $\text{Spf}(A)$ (whose underlying 
reduced is a point) over $R$ 
is good if $A \cong R[[\frac{x^{n_i}_i}{n_i !}, y_j]]$, where $(x_i)$ is the implicit ideal with 
divided powers. Even if $A$ is not smooth can define crystalline differential operators on $A$ by 
using $P_{A,\delta}$, which is a $A$-bimodule. 
We can also define $\omega_{A,\delta}$ and dual differential operators. 
\begin{lemma}
Let $\pi: \mathfrak{X}_1=\textnormal{Spf}(R_1) \to \mathfrak{X}_2=\textnormal{Spf}(R_2)$ be a map of affinoid good $PD$ formal schemes over $R$, 
such that it induces an isomorphism $\pi^* \Omega^1_{\mathfrak{X}_2,\delta} \cong \Omega^1_{\mathfrak{X}_1,\delta}$. 
  Let $E,V/R_2$ be locally free modules, and let $f: P_{\mathfrak{X}_2,\delta} \otimes E \to V$ be 
  a differential operator. Consider the isomorphism $d\pi^{\vee}_{V}: \pi^*(P_{\mathfrak{X}_2,\delta} \otimes V)
  \cong P_{\mathfrak{X}_1,\delta} \otimes \pi^*V$ of \eqref{abdc}. Together $d\pi^{\vee}_{E}$ and $\pi^*(f^*)$ 
  define a differential operator $f_1: P_{\mathfrak{X}_1} \otimes \pi^*E \to \pi^*V$. Similarly, 
  $d\pi^{\vee}_{V^{*}}$ and $\pi^* f$ define a differential operator $f_2: 
  P_{\mathfrak{X}_1,\delta} \otimes (\pi^*V )^* \cong P_{\mathfrak{X}_1,\delta} \otimes \pi^*V^{*} 
  \to \pi^*E^{*} \cong (\pi^*E)^*$. Then 
  $$
  f^*_1=f_2.
  $$
  \begin{proof}
  This follows formally from the compatibility of tensor-hom adjunction with $\pi^*$, and the observation 
  that we can reduce to the case $E=V=\mathcal{O}$, by the way
  $d\pi^{\vee}_{V}$ is defined. 
  \end{proof}
\end{lemma}
Finally, the statement on toroidal compactifications can be proved on the interior, being a map between vector bundles. 

\end{proof}
\end{prop}

Duality for $B$ and $\flag$ is not very interesting, as automorphic line bundles are self-dual, it does not change
the weight shift of 
a theta linkage map, and the basic theta operators are self-dual.
\begin{remark}
In the case of $\GSp_4$ Yamauchi \cite[Prop 3.13]{Yamauchi-1} defined $3$ theta operators on $\Shbar$, one of which corresponds to 
$\pi_* \theta_{\beta}$. By \Cref{easy-lemma-theta}(1)
and the fact that the dual of $\nabla$ on $\mathcal{V}$ is 
$\mathcal{V} \otimes \Omega^{d-1}_{\Shbar}  \to \mathcal{V} \otimes \Omega^{d}_{\Shbar}$,
we can check that the one denoted as $\theta^{\underline{k}}_1$  is the dual differential operator of 
$\pi_* \theta_{\beta}$. In that paper they are only defined 
for certain weights, since in general these dual maps will be only maps between the duals of automorphic 
vector bundles. 
\end{remark}

\subsection{Serre duality on $\mathcal{O}_{P,R}$}
For the applications in \Cref{section7} it will be convenenient to extend the duality $\Ver_{P}W \mapsto \Ver_{P}W^{*}$ of \Cref{Serre-dual-verma}
to $D^{b}(\mathcal{O}_{P,R})$. It will be immediately clear that one needs to work over the derived category. 
Let $R=W/p^n$, for $W=W(\Fpbar)$. 

\begin{lemma} \label{lemma-complex-verma}
Let $C \in D^{b}(\mathcal{O}_{P,R})$, then $C$ is quasi-isomorphic to a complex whose terms are of the form $\Ver_{P} W$ for $W \in \Modfg_{R}(P)$.
\begin{proof}
By virtue of being finitely generated as a $U\mathfrak{g}$-module every object in $\mathcal{O}_{P,R}$ admits a surjection from a Verma module as in the statement. Then the lemma follows from 
\cite[Thm 12.7]{buehler2009exactcategories}. 
\end{proof}
\end{lemma}

\begin{defn}
Let $V \in D^{b}(\mathcal{O}_{P,R})$. Let $V$ be quasi-isomorphic to $C=[\ldots \to \Ver_{P} W_i \xrightarrow{f_i} \Ver_{P} W_{i+1}\to \ldots]$ as in \Cref{lemma-complex-verma}.
Define $V^{*}$ as $[\ldots \to \Ver_{P}W^{*}_{i+1} \xrightarrow{f^*_i} \Ver_{P}W^{*}_i \to \ldots]$. 
Then $V^*$ is independent of the choice of $C$.
\begin{proof}
Suppose that $C_1=[\Ver_{P}W_i],C_2=[\Ver_{P} V_i]$ are two quasi-isomorphic complexes of Verma modules. Then using \Cref{verma-flag} and the exactness of the admissible dual we get two quasi-isomorphic complexes 
$[P_{G/P}\otimes \mathcal{W}^{\vee}_i]\cong  [P_{G/P}\otimes \mathcal{V}^{\vee}_i]$ of $G$-equivariant sheaves on $G/P$. Then we take the Serre dual of both complexes to get a quasi-isomorphism 
$[P_{G/P}\otimes \mathcal{W}_i \otimes K_{G/P}]\cong  [P_{G/P}\otimes \mathcal{V}_i \otimes K_{G/P}]$.
Using \Cref{vb-on-flag} we get a quasi-isomorphism of complexes 
$[\Ver^{\vee}(W^{*}_i)] \cong [\Ver^{\vee}(V^{*}_i)]$. Applying duality on $\mathcal{O}^{+}_{P,R}$ and using that the double dual is isomorphic to the identity \Cref{lemma-category-O}(2),
 we get that $C^{*}_1 \cong C^*_2$, 
since dual maps of Verma modules are defined via this process. 
\end{proof}
\end{defn}

As an example, given $V \in \Rep_{\Fpbar}(G)$, it has a standard resolution by 
the standard complex $C=[\Ver_{P}(\wedge^{\bullet} \mathfrak{g}/\mathfrak{p} \otimes V)]$. One can explicitly see that for $p$ large enough with respect to $G$, $V^{*}=V^{\vee}[-d]$ where $d$ is the dimension of $G/P$.
However, Serre duality preserves a certain class of objects in $\mathcal{O}_{P,R}$. 
We endow $D^{b}(\mathcal{O}_{P,R})$ with the natural $t$-structure, which enables to see $\mathcal{O}_{P,R}$ as the subcategory of objects whose cohomology is concentrated in degree $0$. 
\begin{lemma} \label{lemma-serre-dual-acyclic}
Let $V \in \mathcal{O}_{P,R}$ have a finite filtration by Verma modules. Then $[V]^*$ lies in $\mathcal{O}_{P,R}$, and we simply denote it by $V^*$. Similarly, 
if $C=[V_1 \to V_2 \to \ldots V_n] \in D^{b}(\mathcal{O}_{P,R})$ is a complex with each $V_i$ having a finite filtration by Verma modules, then 
$$
C^*=[V^{*}_n \to V^{*}_{n-1}\to \ldots \to V^*_1].
$$
\begin{proof}
By induction it is enough to prove it when $0 \to \Ver_{P} W_1 \to V \to \Ver_{P} W_2 \to 0$ is an extension of two Verma modules. We get an exact triangle 
$[\Ver_{P}W^{*}_2] \to V^* \to [\Ver_{P}W^{*}_1]$, which shows that $V^*$ has cohomology concentrated in degree $0$, so it lies in the heart of the $t$-structure. 
The second statement follows from a general fact about complexes whose terms are acyclic for a derived functor. 
\end{proof}
\end{lemma}

\section{The weight part of Serre's conjecture: entailments} \label{section6}
\subsection{Generalities about the weight part of Serre's conjecture} \label{generalities-SW}
We explain the setup of the weight part of Serre's conjecture for Shimura varieties in some generality. 
Let $\Sh/\mathcal{O}$ be the integral model of a Hodge type Shimura variety at hyperspecial level, 
with reductive model $G/\Z_p$, and reflex field $E_0/\Q$. There is an action of a spherical 
Hecke algebra $\mathbb{T}/\Z_p$ on coherent or \'etale cohomology of the Shimura variety.
Let $F/\Q$ be a finite extension such that $G_{\Q}$ splits over $F$, so that $E_0 \subseteq F$. Let $\check{G}$ be the dual group 
of the split group $G_{F}$, as a group over $\Z_p$.
We will work with the following assumption. 
\begin{itemize}
\item Given a mod $p$ Hecke eigensystem $\mathfrak{m} \subseteq \mathbb{T}$
appearing in the coherent or \'etale cohomology of the Shimura variety, 
we can attach a semisimple mod $p$ Galois representation $\overline{r}_{\mathfrak{m}} : 
\Gal_{F} \to \check{G}(\overline{\F}_p)$
as in 
\cite[Conj 2.2.1]{Eischen-Mantovan-1}.
\end{itemize}

Recall that $\mathcal{O} \otimes \Q_p=E$ is the completion of the reflex field $E_0$ at some prime $\mathfrak{p}$ above $p$.
For simplicity let us assume that we can take $F=E_0$ and that $G^\der_{\overline{\Q}}$ is almost-simple.
This holds, for instance, for the Siegel threefold and unitary Shimura varieties where the CM field 
is quadratic and the signature is not $(n,n)$. From now on we will assume that $\Sh$ is one of these, 
and in the case of unitary Shimura varieties we will also assume that it is compact \footnote{This is mainly so that 
\Cref{assumption-vanishing}(1) holds by the work of \cite{Deding-unitary}.}.  
Let $\overline{\rho}_{\mathfrak{m}}$ be the restriction of $\overline{r}_{\mathfrak{m}}$ to the decomposition group 
of $\mathfrak{p}$. Then we define the set of modular Serre weights as 
$$
W(\overline{\rho}_{\mathfrak{m}}):=\{F(\lambda): \H^{\bullet}_{\et}(\Sh_{\overline{\Q}_p},F(\lambda))_{\mathfrak{m}} \neq 0 \}. 
$$
We will often drop $\mathfrak{m}$ from the notation and write $W(\overline{\rho})$.
In the case of $\GSp_4$ let $\mathfrak{m} \subset \mathbb{T}$ be a maximal ideal occurring in the coherent cohomology of
the toroidal compactification $\Shbar^{\tor}$, e.g. occuring in $\H^0(\Shbar,\omega(k,l))$
for $k \ge l \ge 0$.
Then one can attach a semisimple Galois representation $\overline{r}_{\mathfrak{m}}: \Gal_{\Q} \to \gsp_4(\overline{\F}_p)$
 by the construction of Galois representations for automorphic representations of regular weight \cite{Galois-reps},
 and the use of generalized Hasse invariants in \cite{boxer-thesis} or \cite{Goldring-Koskivirta-Galois}. 
We say that $\mathfrak{m}$ is non-Eisenstein if $\overline{r}_{\mathfrak{m}}$ is irreducible 
as a $\GL_4$-valued representation. For some compact unitary Shimura varieties, 
as in Harris--Taylor Shimura varieties from \Cref{unitary-non-compact}, Galois representations for characteristic $0$
automorphic forms were constructed 
in \cite{Harris-Taylor} and \cite{Shin-galois}. For mod $p$ eigenforms, they are constructed via congruences as in the $\GSp_4$ case.

We briefly explain Herzig's 
recipe \cite{Herzig1} \cite{Herzig-Tilouine} for $W(\overline{\rho})$ in the case that $\overline{\rho}$ is semisimple
(in the $\GSp_4$ case this means semisimple as a $4$-dimensional representation), 
as it will be relevant for the weight shifts of the operators. 
For the two cases above we may assume that $E=\Q_{p^n}$ is unramified.
For $t \ge 1$ let $\omega_{t}: \Gal_{\Q_{p^{nt}}} \to \Fpbar^{\times}$
be a niveau $tn$ fundamental character. Note that in the two cases that we are considering $G$ is self-dual. 
Given $\overline{\rho}: \Gal_{E} \to \check{G}(\overline{\F}_p)$ we can
write the semisimplification of 
the restriction to inertia $\tau\coloneqq \overline{\rho}^{\text{ss}}_{\mid I_{\mathfrak{p}}}$ as
$$
\tau=\tau(\mu, w)\coloneqq(\overline{\mu}+p^nw^{-1}\overline{\mu}+ \ldots +
p^{n(t-1)}w^{1-t}\overline{\mu})\omega_{t}
$$
where  $w \in W$ is such that 
$\check{w} \tau \sim \tau^{p^n}$, 
 $t\ge 1$ is the order of $w$, $\mu \in X^{*}(T)$, and $\overline{\mu} \in X_{*}(\check{T})$. We are fixing
the isomorphism $G=\check{G}$ swapping $T$ and $\check{T}$,
which on Weyl groups we denote by $w \to \check{w}$, and sending $\mu$ to $\overline{\mu}$. 
Then any semisimple $\tau: I_{\mathfrak{p}} \to G(\overline{\F}_p)$ that lifts 
to $\Gal_{E}$ (a tame inertial parameter)
can always be written as $\tau(\mu,w)$, and for generic $\tau$ we can assume that $\mu \in X_1(T)$. 
 Herzig proposed a conjectural set of Serre weights $W^{?}(\overline{\rho})$ when $\overline{\rho}$ is semisimple
(restricting to generic weights),
based 
on the compatibility with inertial mod $p$ Langlands. 
Further restricting to generic 
weights, the set proposed by Herzig is \cite[Prop 6.28]{Herzig1} \cite[Prop 10.1.13]{Gee-Herzig-Savitt}
\begin{align} \label{Herzig-recipe}
  W^{?}(\overline{\rho})=\{ F(\mu) \text{ such that  } & \exists \; \mu' \uparrow \mu \text{ with } \mu'+\rho 
  \text{ dominant and} \\
  & w \in W,
   \text{ such that }
  \overline{\rho}_{\mid I}
  =\tau(\mu'+\rho,w)\}.
  \end{align}
  We say that $F(\mu) \in W^{?}(\overline{\rho})$ is an \textit{obvious weight} if we can take $\mu'=\mu$ above, 
  and a \textit{shadow 
  weight} otherwise. The obvious weights take their name since for $\overline{\rho}$ semisimple, given 
  any obvious weight in $W^{?}(\overline{\rho})$,  
  it admits an "obvious" crystalline lift (by lifting the fundamental characters)
   whose Hodge-Tate weights are associated to the obvious weight. 
  For general $\overline{\rho}$ one expects that $W(\overline{\rho}) \subseteq W^{?}(\overline{\rho}^{\text{ss}})$. 
   We immediately see that the combinatorics of the shadow weights are related to the combinatorics 
   of the theta linkage maps.
  In fact, we will use theta linkage maps to prove the existence of certain shadow weights 
  in this section and the next one. 
  We recall our definition of generic weights in \Cref{generic-Sh}, since we will make heavy use of it 
  in this section.

\subsection{A generic entailment for $\GL_4$}
By studying Breuil--Mezard cycles \cite{Le-Hung-Lin} can prove that there should not be
generic entailments for the group $\GL_3/\Q_p$, so the next 
natural candidate for a generic entailment beyond $\GSp_4/\Q_p$ is $\GL_4/\Q_p$. 
For this, we can work with a unitary group $G/\Q$ with respect 
to a quadratic imaginary $E/\Q$ with $p$ split in $E$. We can either work with signature $(3,1)$ or $(2,2)$.
\begin{notation} \label{notation-unitary}
Consider $\Sh$, the unitary Shimura variety of signature $(a,b)$ for a quadratic imaginary $E/\Q$. 
Then $G_{\overline{\Q}}=\GL_n \times \mathbb{G}_m$. Weights are of the form $(\lambda,c)$ with $c$ a character of 
$\mathbb{G}_m$. We will suppress the similitude character $c$ from the notation from now on, since the representation 
theory of $G_{\Fpbar}$ is equivalent to that of $\GL_n$.
For $\GL_n$ we will use the notation for roots as in \Cref{table-lie}, and take 
$\rho=(n-1,n-2,\ldots,0)$. 
\end{notation}
For $\GL_4/\Fpbar$ we have $6$ $p$-restricted alcoves consisting of $(a,b,c,d)-\rho \in X^*(T)$ satisfying 
	\begin{align*}
		C_0:\quad & 0 < a-b,\;b-c,\;c-d;\ a-d < p, \\
		C_1:\quad & 0 < b-c;\ p < a-d;\ a-c,\;b-d < p, \\
		C_2:\quad & 0 < c-d;\ p < a-c;\ a-b,\;b-d < p, \\
		C_3:\quad & 0 < a-b;\ p < b-d;\ c-d,\;a-c < p, \\
		C_4:\quad & p < a-c,\;b-d;\ b-c < p;\ a-d < 2p, \\
		C_5:\quad & 2p < a-d;\ a-b,\;b-c,\;c-d < p, \\
		C_{0'}:\quad & 0 < b-c,\;c-d;\ p < a-b;\ a-d < 2p, \\
		C_{0''}:\quad & 0 < a-b,\;b-c;\ p < c-d;\ a-d < 2p.
	  \end{align*}
  Here $C_{0'}$ and $C_{0''}$ are not $p$-restricted but appear as constituents of $p$-restricted Weyl modules. 
	We can also describe the decomposition of dual Weyl modules into irreducible constituents in this case. We  
	only state the ones that we will use. 
	\begin{prop} \label{Weyl-decomposition-GL4} \cite{Jantzen74}
	Let $\lambda_0=(a,b,c,d) \in C_0$, and let $\lambda_i \in C_i$ be the unique element of 
	$W_{\textnormal{aff}} \cdot \lambda_0$ in $C_i$. Then in the Grothendieck group of representations of 
	$\GL_4/\overline{\F}_p$ we have the following relations. 
	\begin{itemize}
	\item $V(\lambda_4)_{\Fpbar}=L(\lambda_4)+L(\lambda_3)+L(\lambda_2)+L(\lambda_1)+L(\lambda_0)$
	\item $V(\lambda_5)_{\Fpbar}=L(\lambda_5)+L(\lambda_4)+L(\lambda_3)+L(\lambda_2)+L(\lambda_1)+L(\lambda_{0'})
	 +L(\lambda_{0''})$
	\end{itemize}
	where $\lambda_{0'}=(p+b-1,c-1,d-1,a-p+3)$ and $\lambda_{0''}=(p+d-3,a+1,b+1,c-p+1)$.
	In the Grothendieck group of $\Rep_{\Fpbar}\GL_4(\F_p)$ we have that for generic $\lambda_0 \in C_0$
	\begin{align*}
	F(\lambda_{0'})&=F(b,c-1,d-1,a-p+3)+F(b-1,c,d-1,a-p+3)\\
	&+F(b-1,c-1,d,a-p+3)+F(b-1,c-1,d-1,a-p+4)
	\end{align*}
	and 
	\begin{align*}
	F(\lambda_{0''})&=F(p+d-2,a-p+2,b-p+2,c-p+1)+F(p+d-2,a-p+1,b-p+2,c-p+2)\\
	&+F(p+d-2,a-p+2,b-p+1,c-p+2)+F(p+d-3,a-p+2,b-p+2,c-p+2)
	\end{align*}
  Similarly, in the Grothendieck group of $\Rep_{\Fpbar}U_4(\F_p)$ we have that for generic $\lambda_0 \in C_0$
  \begin{align*}
    F(\lambda_{0'})&=F(b-1,c-1,d-1,a-p+2)+F(b-1,c-1,d-2,a-p+3)\\
    &+F(b-1,c-2,d-1,a-p+3)+F(b-2,c-1,d-1,a-p+3)
    \end{align*}
    and 
    \begin{align*}
    F(\lambda_{0''})&=F(p+d-4,a-p,b-p,c-p+1)+F(p+d-4,a-p+1,b-p,c-p)\\
    &+F(p+d-4,a-p,b-p+1,c-p)+F(p+d-3,a-p,b-p,c-p)
    \end{align*}
  
	\begin{proof}
	The decomposition of dual Weyl modules is from \cite{Jantzen74}. The decompositions for $F(\lambda_{0'})$
	follow from Steinberg's tensor product, the fact that $L(\lambda)=V(\lambda)$ for $\lambda \in C_0$,
	and the Littlewood-Richardson rule for tensor products of highest weight representations. 
  For $\GL_4(\F_p)$ we use that $L(1,0,0,0)^{(p)}(\F_p)=F(1,0,0,0)$ and $L(1,1,1,0)^{(p)}(\F_p)
  =F(1,1,1,0)$.
  For $U_4(\F_p)$ we use that 
  $L(1,0,0,0)^{(p)}(\F_p)=F(0,0,0,-1)$, and  $L(1,1,1,0)^{(p)}(\F_p)=F(0,-1,-1,-1)$
    since in that case Frobenius acts as 
  $\sigma(a_1,a_2,a_3,a_4)=(-a_4,-a_3,-a_2,-a_1)$ swapping $\alpha_1$ and $\alpha_3$. 
	\end{proof}
	\end{prop}

	We have included only these two dual Weyl modules because
  they are the only pair of $p$-restricted $\lambda_{i}, \lambda_{j}$
	with $\lambda_i \uparrow \lambda_j$ such that $V(\lambda_i)_{\Fpbar}$ contains an irreducible factor which is not contained 
	in $V(\lambda_j)_{\Fpbar}$, which is the key phenomenon we look for to produce entailments. The full list can be checked 
	in \cite[Prop 9.3]{Herzig1} (in the thesis version). 
	Also, crucially for us $\lambda_4 \uparrow_{\delta} \lambda_5$, for $\delta$ the longest root of $\GL_4$, 
	which is the relevant one in \Cref{injective-linkage}. 
  We make \Cref{known-vanishing} explicit in the $\GL_4$ case. See \Cref{unitary-compact} for 
  the explicit description of the automorphic vector bundles. 
  \begin{prop} \label{vanishing-GL4}
  Let $\Shbar$ be a compact unitary Shimura variety with $p$ split in the quadratic imaginary field,
  and $\lambda=(a_1,a_2,a_3,a_4) \in X^{*}(T)$.
   Then 
  \begin{enumerate}
  \item If $G_{\Q}$ has signature $(2,2)$ then if $\H^{0}(\Shbar,\omega(\lambda)) \neq 0$ we have 
  $p(a_2-a_3)+(a_1-a_4) \ge 0$. If $G_{\Q}$ has signature $(3,1)$ the condition is that 
  $\lambda$ is in the cone generated by the weights of the Hasse invariants and the cone 
  $\{(b_1,b_2,b_3,b_4) : (b_1-b_4)+p(b_2-b_4)+p^2(b_3-b_4)\}$
 \cite[Thm 4.2.7/8]{general-cone-conjecture1}.
  \item Let $p \ge  7$. 
  If $G_{\Q}$ has signature $(3,1)$ then $\H^1(\Shbar,\omega(\lambda))=0$ whenever $a_1\ge a_2 \ge a_3 \ge 7$ and 
  $(a_1-a_4-2)+p(a_4-a_3+6)>0$. If $G_{\Q}$ has signature $(2,2)$ then we replace the last condition by 
  $(a_1-a_4-2)+p(a_3-a_2+2)>0$ \cite{Deding-unitary}.
  \end{enumerate}
  \end{prop}

    To complete the proof of the entailment we need a way to relate coherent cohomology of $\Shbar^\tor$
    and mod $p$ \'etale cohomology of $\Sh_{\overline{\Q}_p}$.
	We will work in the Harris--Taylor setting of \Cref{unitary-compact} with signature $(n-1,1)$, for simplicity. 
  Also from now on we assume that these satisfy the assumptions of 
  \cite[\S I.7]{Harris-Taylor} for some auxiliary prime $l$. We will provide a different proof
  from the one for $\GSp_4$ in \cite[Prop 5.1]{paper},
  using that the Galois representation is contained in the \'etale cohomology. 
	\begin{prop} \label{coherent-to-betti-GL4}
		Let $\lambda \in X^*(T)$ be a dominant weight. Let $E$ be some imaginary quadratic number field, and 
    $\Sh$ the associated compact Harris--Taylor Shimura variety of signature $(n,1)$ of \Cref{unitary-compact}.
    Let $\mathfrak{m} \subset \mathbb{T}$ be a non-Eisenstein generic 
		maximal ideal. Let $w \in W^{M}$, then
    $$
     \H^n_{\text{\'et}}(\Sh_{\Q_p},V(\lambda)^{\vee}_{\overline{\F}_p})_{\mathfrak{m}} \neq 0  \implies 
      \H^{n-l(w)}(\Shbar,\omega(w \cdot \lambda)^{\vee})_{\mathfrak{m}} \neq 0.
    $$
    Similarly,
		\begin{equation} \label{why}
		\H^n_{\text{\'et}}(\Sh_{\Q_p},V(\lambda)_{\overline{\F}_p})_{\mathfrak{m}} \neq 0  \implies 
		\H^0(\Shbar,\omega(\lambda+\eta))_{\mathfrak{m}} \neq 0,
    \end{equation}
    where $\eta$ is the weight of $\Omega^{n}_{\Shbar}$.
		Moreover, $\H^{*}_{\text{\'et}}(\Sh_{\Q_p},V(\lambda) \otimes \overline{\F}_p)_{\mathfrak{m}}$ is concentrated 
		in degree $n$. Further, assuming \Cref{assumption-vanishing}(1) 
    the reverse implications
    \begin{align} \label{reverse}
     &\H^{0}(\Shbar,\omega(\lambda+\eta))_{\mathfrak{m}} \neq 0 \implies 
     \H^n_{\text{\'et}}(\Sh_{\overline{\Q}_p},V(\lambda)_{\overline{\F}_p})_{\mathfrak{m}} \neq 0 \\ 
     & \H^{n}(\Shbar,\omega(-w_0\lambda)^{\vee})_{\mathfrak{m}} \neq 0 \implies 
     \H^n_{\text{\'et}}(\Sh_{\overline{\Q}_p},V(\lambda)_{\overline{\F}_p})_{\mathfrak{m}} \neq 0
    \end{align}
    holds for generic $\lambda \in X_1(T)$. 
    \Cref{assumption-vanishing}(1) holds when $p$ is split in the CM field. For $n=2$ \eqref{reverse}
    holds for $\lambda=(a_1,a_2,a_3) \in X_1(T)$ satisfying $a_1>a_2+2$.
    \begin{proof}
    The concentration of 
    \'etale cohomology in middle degree was proved in \cite{caraiani-scholze-compact},
    in fact under a weaker genericity condition 
    than the one in \cite{Hamann-Lee}. First we prove the forward implication. By torsion freeness of \'etale cohomology 
    we can choose 
     $\tilde{\mathfrak{m}} \subseteq \mathbb{T}$ corresponding to an integral system of Hecke eigenvalues lifting 
    $\mathfrak{m}$. Then under the restriction that $\mathfrak{m}$ is non-Eisenstein, 
    $r_{\tilde{\mathfrak{m}}}$ is a summand of 
    $\H^n_{\et}(\Sh_{\overline{\Q}},V(\lambda)^{\vee}_{\overline{\Q}_p})_{\tilde{\mathfrak{m}}}$.
    This follows from \cite[Prop VII.1.8]{Harris-Taylor}. Let $\pi^{\infty}$ be a representation of $G(\mathbb{A}^{\infty})$
    appearing in the \'etale cohomology, whose eigensystem is
    $\tilde{\mathfrak{m}}$, and 
    $\pi=\pi^{\infty} \otimes \pi_{\infty}$ a choice of automorphic representation of $G(\mathbb{A})$.
    For the cited result to apply we need to check that 
    the Jacquet-Langlands transform of the base change of $\pi$ is cuspidal as a representation 
    of $\GL_{n+1}(\mathbb{A}_{E_0})$,
    this follows since the Galois representations of Eisenstein forms are known to be reducible. 
    Moreover, by the \'etale-crystalline comparison theorem
    $\rho_{\tilde{\mathfrak{m}}}$ is crystalline with regular Hodge-Tate weights, 
    so each graded piece of $D_{\crys}(\rho_{\tilde{\mathfrak{m}}})$ is one-dimensional. 
    By the rational \'etale-crystalline comparison theorem and the BGG decomposition the graded pieces of 
    $D_{\crys}(\rho_{\tilde{\mathfrak{m}}})$  are labelled by $W^{M}$ and they are a summand of 
    $\H^{n-l(w)}(\Sh_{\overline{\Q}_p},\omega(w \cdot \lambda)^{\vee})_{\tilde{\mathfrak{m}}}$ for each $w \in W^{M}$. 
    Therefore, $\H^{n-l(w)}(\Sh_{\overline{\Q}_p},\omega(w \cdot \lambda)^{\vee})_{\tilde{\mathfrak{m}}} \neq 0$
    for each $w \in W^M$.
    By rescaling an eigenform in $\H^{n-l(w)}(\Sh_{\overline{\Q}_p},\omega(w \cdot \lambda)^{\vee})_{\tilde{\mathfrak{m}}}$
    we can assume that it lies in 
    $\H^{n-l(w)}(\Sh_{\check{\Z}_p},\omega(w \cdot \lambda)^{\vee})_{\tilde{\mathfrak{m}}}$, 
    and then we reduce mod $p$ to obtain 
     $\H^{n-l(w)}(\Shbar,\omega(w \cdot \lambda)^{\vee})_{\mathfrak{m}} \neq 0$.
     To prove \eqref{why}, let $\lambda^{\vee}=-w_{0}\lambda$, then in characteristic $0$ we have that $\omega(w \cdot \lambda^{\vee})^{\vee}=\omega(-w_{0,M}(w \cdot \lambda^{\vee}))$
     and $V(\lambda^{\vee})^{\vee}=V(\lambda)$. 
     For $w=w_{0,M}w_{0} \in W^M$ of length $n$ we can compute that the former is equal to $\omega(\lambda+\eta)$. Therefore, we obtain that 
     $\H^{0}(\Sh_{\check{\Z}_p},\omega(\lambda+\eta))_{\tilde{\mathfrak{m}}} \neq 0$. 

    Now we prove the first implication in \eqref{reverse}. 
    The key is that by \Cref{known-vanishing}(1), $\H^{1}(\Shbar,\omega(\lambda+\eta))_{\mathfrak{m}}=0$, so that 
    we can lift any eigenclass to $\H^{0}(\Sh_{\check{\Z}_p},\omega(\lambda+\eta))_{\mathfrak{m}}$, and this space is torsion-free.
    Thus, $\H^{0}(\Sh_{\overline{\Q}_p},\omega(\lambda+\eta))_{\mathfrak{m}}=
    \H^{0}(\Sh_{\overline{\Q}_p},\omega(w_{0,M}w_{0} \cdot \lambda^{\vee})^{\vee})_{\mathfrak{m}} \neq 0$. 
    Then we use the BGG in characteristic $0$ and torsion-freeness of \'etale cohomology 
    to get that $\H^n_{\text{\'et}}(\Sh_{\Q_p},V(\lambda)_{\check{\Z}_p})_{\mathfrak{m}} \neq 0$. We get the statement 
    reducing modulo $p$. The genericity conditions for $n=2,3$ are given in 
     \Cref{vanishing-GL4}(2).  For the second implication in \eqref{reverse}, the lifting of 
      $\H^{n}(\Shbar,\omega(-w_0\lambda)^{\vee})_{\mathfrak{m}}$ to characteristic $0$ is immediate, and the torsion-freeness 
      is given by the vanishing of $\H^{n-1}(\Shbar,\omega(-w_0\lambda)^{\vee})_{\mathfrak{m}} 
      \cong \H^1(\Shbar,\omega(-w_{0}\lambda+\eta))^{\vee}_{\mathfrak{m'}}$. Here we are using that Serre duality is 
      Hecke equivariant up to a twist \cite[Prop 4.2.9]{higher-coleman}, and $\mathfrak{m}'$ is that twisted eigensystem,
       which is still non-Eisenstein generic.

    \end{proof}
	\end{prop}

	With that we can exhibit a general method to construct generic weak entailments, which depends on the structure of the irreducible constituents of certain dual Weyl modules.
  Then using \Cref{Weyl-decomposition-GL4} we will get a generic weak
   entailment for $\GL_4/\Q_p$. 
    
	\begin{theorem} \label{general-entailment}
  Let $\Sh/\mathcal{O}$ be a compact unitary Shimura variety of signature $(n-1,1)$ with $p$ split
  in the quadratic imaginary field. Assume that \Cref{assumption-vanishing}(2) holds.
  Let $\mathfrak{m} \subseteq \mathbb{T}$ be a generic non-Eisenstein maximal ideal
	appearing in the cohomology of $\Sh$. Let $\lambda, \mu \in X_1(T)$ such that 
  $\lambda \uparrow_{\delta} \mu$, where $\delta$ is the longest root. Assume that there exists some
  $F(\chi) \in \textnormal{JH}[V(\lambda)_{\Fpbar}]$ such that $F(\chi) \notin \textnormal{JH}[V(\mu)_{\Fpbar}]$.
  Let $\{F(\mu_i)\}$ be the Jordan--Holder factors of $V(\mu)_{\Fpbar}$.
  Then the statement 
  $$
  F(\chi) \in W(\overline{\rho}_{\mathfrak{m}}) \implies F(\mu_i) \in W(\overline{\rho}_{\mathfrak{m}}) \textnormal{ for some } i
  $$
  holds for generic $\lambda \in X_1(T)$. Note that $F(\mu_i) \neq F(\chi)$ for all $i$.
  The same result holds for $p$ inert if \Cref{assumption-vanishing}(1) holds.
  \begin{proof}
  The method is the same as in \cite[Thm 5.2]{paper}. By definition $F(\chi) \in W(\overline{\rho}_{\mathfrak{m}})$ implies that 
  $\H^{n-1}_{\et}(\Sh_{\overline{\Q}_p},F(\chi))_{\mathfrak{m}} \neq 0$. By the concentration of \'etale cohomology in 
  \Cref{coherent-to-betti-GL4} and the assumption that $F(\chi) \in \textnormal{JH}[V(\lambda)_{\Fpbar}]$ we obtain 
  $\H^{n-1}_{\et}(\Sh_{\overline{\Q}_p},V(\lambda)_{\Fpbar})_{\mathfrak{m}} \neq 0$. By 
  \Cref{coherent-to-betti-GL4} this implies that 
  $\H^0(\Shbar^\tor,\omega(\lambda+\eta))_{\mathfrak{m}} \neq 0$. By \Cref{injective-linkage}
  and \Cref{known-vanishing}(1) $\theta_{\lambda \uparrow \mu}$ is injective 
  on global sections for generic $\lambda$, so that 
  $\H^0(\Shbar^\tor,\omega(\mu+\eta))_{\mathfrak{m}} \neq 0$. Applying \Cref{coherent-to-betti-GL4} in the other direction 
  yields $\H^{n-1}_{\et}(\Sh_{\overline{\Q}_p},V(\mu)_{\Fpbar})_{\mathfrak{m}} \neq 0$. Finally, we conclude 
  from the the decomposition of $V(\mu)_{\Fpbar}$ into irreducible components. 
  \end{proof}
	\end{theorem}
  Note that the theorem above is quite robust, and it should extend to more general settings, such as the Siegel case, 
  where \Cref{assumption-vanishing}(2) is known. The only ingredient that we have not written in this generality is 
  \Cref{coherent-to-betti-GL4}.
  In the $\GL_{4,\Q_p}$ case we get an unconditional result. Moreover, we also get (assuming the mild \Cref{assumption-vanishing}(1)) 
  a generic entailment for the non-split group $GU(4)$ with respect to $\Q_{p^2}$, 
  where much less is known in general.
  \begin{cor} \label{entailment-GL4}
  Let $\Sh/\mathcal{O}$ be a compact unitary Shimura variety of signature $(3,1)$ 
  with $p$ split in the quadratic imaginary field.
  Let $\mathfrak{m} \subseteq \mathbb{T}$ be a generic non-Eisenstein maximal ideal.
  Then the statement
  $$
	F(\lambda_0) \in W(\overline{\rho}_{\mathfrak{m}}) \implies F(\mu) \in W(\overline{\rho}_{\mathfrak{m}})
  $$
  for some $F(\mu) \in \textnormal{JH}[V(\lambda_5)_{\overline{\F}_p}]$ holds for $\lambda_0 \in C_0$ generic.
  Moreover,
	$F(\mu) \neq F(\lambda_0)$.
  The same result holds for $p$ inert if \Cref{assumption-vanishing}(1) holds.
  \begin{proof}
  It follows from \Cref{general-entailment}. We use as input that \Cref{assumption-vanishing}(2) is known in this case by 
  \Cref{known-vanishing}(2), and \Cref{Weyl-decomposition-GL4} for the Jordan--Holder factors of 
  $V(\lambda_4)_{\Fpbar}$ and $V(\lambda_5)_{\Fpbar}$. We also use the observation that $\lambda_4 \uparrow_{\delta} \lambda_5$.
  \end{proof}
  \end{cor}

  \begin{remark} \label{heuristic-BM}
  Using local arguments, \cite{Le-Hung-Lin} prove that if the Breuil--Mezard 
  conjecture is true, then $\mathcal{Z}_{F(\lambda_0)}
  \subseteq \mathcal{Z}_{F(\lambda_5)}$ for generic $\lambda_0 \in C_0$.  
  We give a heuristic argument that suggests that \Cref{entailment-GL4} implies that. 
  The idea is to use Kisin--Taylor--Wiles patching \cite{CEGGPS}
  on the injection $\H^0(\Shbar, \omega(\lambda_4+\eta)) \hookrightarrow \H^0(\Shbar,\omega(\lambda_5+\eta))$
  from \Cref{injective-linkage}. One would  
  still obtain an injection of patched modules $M_{\infty}(\lambda_4) \hookrightarrow M_{\infty}(\lambda_5)$
  over a thickening $R_{\infty}$ of the local Galois deformation $R_{\overline{\rho}}$ ring for $\overline{\rho}$.
  One key ingredient in order to patch coherent cohomology (although this is not written in the literature, 
  some of the tools can be found in \cite{Harris} \cite{Atanasov-Harris})
  is that both $\omega(\lambda_i+\eta)$ have coherent cohomology concentrated in degree $0$ by \Cref{known-vanishing}.
  The heuristic is that by \Cref{coherent-to-betti-GL4}, and assuming a form of Fontaine-Mazur's conjecture, 
   $\text{supp}_{R_{\overline{\rho}}}M_{\infty}(\lambda_i)$ should be (set-theoretically) 
   $\mathcal{Z}_{V(\lambda_i)_{\Fpbar}}(\overline{\rho})$.
   Then the injection of patched modules would imply that 
   $\mathcal{Z}_{V(\lambda_4)_{\Fpbar}}(\overline{\rho}) \subseteq \mathcal{Z}_{V(\lambda_5)_{\Fpbar}}(\overline{\rho})$
   for every $\overline{\rho}$, so that using \Cref{Weyl-decomposition-GL4} we get that
   $\mathcal{Z}_{F(\lambda_0)} \subseteq \mathcal{Z}_{F(\lambda_5)}+\sum_{\mu}\mathcal{Z}_{F(\mu)}$
   where each $\mu$ is in the lowest alcove, coming from \Cref{Weyl-decomposition-GL4}.
   One can easily discard the possibility that $\mathcal{Z}_{F(\lambda_0)}$ is contained in any cycle of a different
	$C_0$ weight 
	by Fontaine--Laffaille theory, so the above would suggest that 
	$$
	\mathcal{Z}_{F(\lambda_0)} \subseteq \mathcal{Z}_{F(\lambda_5)},
	$$
	which agrees with the local predictions/results of \cite{Le-Hung-Lin} for $\GL_{4,\Q_p}$.
	Moreover, this is the only predicted generic entailment for $\GL_{4,\Q_p}$. 
  Assuming that \Cref{assumption-vanishing} holds this heuristic generalizes. 
  Therefore, given a weak entailment proved as in \Cref{entailment-GL4}, it can help us guess 
  if there exists an entailment refining it. Indeed, in the $\GSp_4$ case this heuristic also correctly predicts the 
  entailment predicted by \cite{Le-Hung-Lin}. In that case, see \cite[\S 4]{le2025serreweightconjecturesmathrmgsp4} for a similar reasoning.
  \end{remark}

	\subsection{Non-generic entailments for $\GL_3$} 
	Even though there are no generic entailments for $\GL_{3,\Q_p}$, 
	we can still obtain some non-generic entailments. We will say that an entailment is non-generic if it holds 
  for a family of weights which is not generic. Nevertheless, we will give explicit conditions on the weights, which make the 
  non-genericity clear.
  There are two $p$-restricted alcoves: the lowest 
	alcove $C_0=\{(a,b,c): 0 \le a-b,b-c ; a-c<p-2\}$ and the upper alcove 
	$C_1=\{(a,b,c): a-b,b-c <p-1 ; a-c >p-2\}$. For $\lambda_0 \in C_0$ let $\lambda_1 \in C_1$ be its corresponding 
	reflection. 

	\begin{lemma} \label{decomposition-GL3}
		For $\lambda_0 \in C_0$, $V(\lambda_0)=L(\lambda_0)$. For $\lambda_1 \in C_1$ we have an exact sequence
		$$
		0 \to L(\lambda_1) \to V(\lambda_1) \to L(\lambda_0)
		$$
		If $\lambda \in X_{1}(T)$ and is in the boundary of $C_1$ then $V(\lambda)=L(\lambda)$.
		\end{lemma}

    We specify \Cref{known-vanishing} and \Cref{injective-linkage} to the $U(2,1)$ $p$ split case. 
    \begin{prop} \label{vanishing-GL3}
    Let $\Shbar$ be a compact unitary Shimura variety of signature $(2,1)$ with $p$ split in the quadratic 
    imaginary field. Then 
    \begin{enumerate}
    \item If $\H^0(\Shbar,\omega(a,b,c))\neq 0$ then $a \ge b$ and $(a-c)+p(b-c)\ge 0$. 
    \item For $p\ge 5$ we have $\H^1(\Shbar,\omega(a,b,c))=0$ for $a \ge b \ge c+5$ satisfying 
    $p(b-c-4)+(a-c-2)>0$. 
    \item For $\lambda_0=(a,b,c)\in C_0$ the linkage map 
$\theta_{\lambda_0 \uparrow \lambda_1} : \H^0(\Shbar,\omega(a+1,b+1,c-2))_{\mathfrak{m}} 
\to \H^0(\Shbar,\omega(p+c-1,b+1,a-p))_{\mathfrak{m}}$
is injective for $\mathfrak{m}$ non-Eisenstein. 
    \end{enumerate}
    \begin{proof}
    Part $(1)$ follows from \cite[Thm 5.1.1, Fig 1]{general-cone-conjecture0}. Part $(2)$ 
    is \cite[Thm 1.5/6]{Deding-unitary}. For $(3)$, 
    let $f$ be an element of the kernel, by \Cref{pth-power-relation} $\theta_{\alpha_1+\alpha_2}(f)=0$.
    By \Cref{general-restriction-strata} we know $H_{\alpha_1} \mid \theta_{\alpha_2}(f)$. 
    The resulting weight of $\theta_{\alpha_2}(f)/H_{\alpha_1}$ would not be $M$-dominant, so $\theta_{\alpha_2}(f)=0$.
     Restricting to $\overline{D}_{\alpha_2}$,
    and using 
    \Cref{restriction-simple-theta} we get $H_{\alpha_2} \mid f$, since $p \nmid a-c+3$. Unless $a=b=c+p-4$, part $(1)$ 
    implies that $f=0$. Otherwise, $f/H_{\alpha_2} \in \H^0(\Shbar, \omega(c-2,c-2,c-2))$, which only contains 
    the trivial eigensystem. 
    \end{proof}
    \end{prop}

	\begin{theorem} \label{non-generic-entailment-GL3}
  Let $\Sh$ be a compact unitary Shimura variety of signature $(2,1)$ with 
  $p$ split in the quadratic imaginary field.
  Let $\mathfrak{m}$ be a generic non-Eisenstein eigensystem 
	with associated local Galois representation $\overline{\rho}$.

	\begin{enumerate}
	\item Let $\lambda=(a,b,b) \in X_1(T)$ satisfying $b \ge 5$.
	 Then $F(\lambda) \in W(\overline{\rho}_{\mathfrak{m}}) \implies F(a+p-1,b+p-1,b) \in W(\overline{\rho}_{\mathfrak{m}})$. 
	 \item Let $\lambda=(a,a,b) \in X_1(T)$. Then $F(\lambda) \in W(\overline{\rho}_{\mathfrak{m}}) \implies F(a+p-1,a,b) 
	 \in W(\overline{\rho}_{\mathfrak{m}})$.
   \item Let $p\ge 5$. Let $(a,b)$ such that $0 \le a-b \le p-3,  b \ge 0$ 
    and  $\H^0(\Shbar,\omega(a+1,b+1,b-1))_{\mathfrak{m}} \neq 0$. Then 
    $F(p+b-1,b,a-p+2) \in W(\overline{\rho}_{\mathfrak{m}})$ and $F(p+a,p+b-1,b) \in W(\overline{\rho}_{\mathfrak{m}})$. 
	\end{enumerate}
	\begin{proof}
  We use that $L(\mu)=V(\mu)$ for all $\mu$ appearing in the proposition, by \Cref{decomposition-GL3}.
	For $(1)$, by 
	\Cref{coherent-to-betti-GL4} $\H^0(\Shbar,\omega(a+1,b+1,b-2))_{\mathfrak{m}} \neq 0$. 
	Applying $H_{\alpha_2}$ 
	we get $\H^0(\Shbar,\omega(a+p,b+p,b-2))_{\mathfrak{m}} \neq 0$, we conclude by 
 \Cref{vanishing-GL3} and 
	\Cref{coherent-to-betti-GL4} again. For $(2)$  
	by the proof of \Cref{coherent-to-betti-GL4} $\H^2(\Shbar,\omega(a,b,a))_{\mathfrak{m}} \neq 0$.
  There is a map 
	$H_{\alpha_2}: \H^2(\Shbar,\omega(a,b,a+p-1)) \to  \H^2(\Shbar,\omega(a,b,a))$ which is 
	Serre dual to an injection of global sections, so by the trivial case of \Cref{Serre-duality-diff-op}
	it is surjective, so that 
	$\H^2(\Shbar,\omega(a,b,a+p-1))_{\mathfrak{m}} \neq 0$. We are using that we can take $H_{\alpha_2}$
  to be of weight $(0,0,1-p)$.
   From \eqref{reverse} in \Cref{coherent-to-betti-GL4} we
	deduce that $\H^2_{\text{\'et}}(\Sh,V(a+p-1,a,b)_{\overline{\F}_p})_{\mathfrak{m}} \neq 0$.
  For part $(3)$, we use that the linkage map $\H^0(\Shbar,\omega(a+1,b+1,b-1)) \to 
  \H^0(\Shbar,\omega(p+b,b+1,a-p))$  and $\theta_{\alpha_2}$ are injective .
  These follow as in \Cref{vanishing-GL3}(3),
  using that $p \nmid a-b+2$. Then the result follows from 
  \Cref{coherent-to-betti-GL4} and that $L(p+b-1,b,a-p+2)=V(p+b-1,b,a-p+2)$. 
	\end{proof}
	\end{theorem}

  One could see part $(3)$ as suggesting congruences between a form of non-regular Hodge-Tate weights 
  $\{a+2,b+1,b+1\}$ and one of regular Hodge-Tate weights $\{p-b+1,b+1,a-p+2\}$ or 
  $\{p+a+2,b+1,a-p+2\}$. However, in general it is not known if the Galois representation 
  of a form in $\H^0(\Sh_{\overline{\Q}_p},\omega(a+1,b+1,b-1))$ is crystalline of those Hodge-Tate weights.

	\begin{remark}
	The local analogue of parts $(1)$ and $(2)$ at the level of Breuil--Mezard cycles is proved
  (without the extra genericity assumption in $(1)$) in upcoming work of Levin, Le, Le Hung and Morra. 
  They prove that if we define candidate Breuil--Mezard cycles by the property that 
  $\mathcal{Z}_{F(\lambda)}$ is its corresponding irreducible component except in the cases
  $\mathcal{Z}_{F(a+p-1,a,b)}=C_{F(a+p-1,a,b)}+C_{F(a,a,b)}$ and
  $\mathcal{Z}_{F(a+p-1,b+p-1,b)}=C_{F(a+p-1,b+p-1,b)}+C_{F(a,b,b)}$;
  then these cycles 
   satisfy the Breuil--Mezard
  conjecture for Hodge-Tate weights $\{2,1,0\}$ and tame types.
	 Therefore, one can reach almost all possible entailments geometrically. 
	\end{remark}

	In general, for $\GL_n/\Q_p$ one expects the following "obvious non-generic entailments": for $1 \le k \le n-1$
	$F(\lambda) \in W(\overline{\rho}_{\mathfrak{m}}) \implies F(\lambda+(p-1,p-1,\ldots,p-1,0,\ldots,0)) \in W(\overline{\rho}_{\mathfrak{m}})$
	where there are $k$ $(p-1)$s and $\lambda$ is of the form $
	(\lambda_1,\ldots, \lambda_{k},\lambda_k, \lambda_{k+2},\ldots, \lambda_n)$. Using Hasse invariants
  and duality for unitary Shimura varieties of different signatures 
	we can reach those weight increases at the level of coherent cohomology. However, in general 
	we won't have $F(\lambda+(p-1,p-1,\ldots,p-1,0,\ldots,0))=V(\lambda+(p-1,p-1,\ldots,p-1,0,\ldots,0))$, 
	so we can't deduce the result at the level of Serre weights. 
\subsection{A remark on the injectivity of theta linkage maps}
The reader might wonder in what generality are theta linkage maps (generically) injective. Let us briefly explain the conjectural picture 
in the case of $\GSp_4$, which is the simplest case where we don't understand this for theta linkage maps within the $p$-restricted region (and coherent cohomology in degree $0$).
There are $4$ $p$-restricted $\rho$-shifted alcoves $C_i$ for $i=0,1,2,3$, with $C_0$ being the lowest alcove. Let $\lambda_0 \in C_0$ and let $\lambda_i \in C_i$ be its affine Weyl reflections. 
Then, we have $3$ theta linkage maps in this region:  $\theta_{\lambda_0 \uparrow \lambda_1}, \theta_{\lambda_1 \uparrow \lambda_2}$, and $\theta_{\lambda_2 \uparrow \lambda_3}$. 
For $\theta_{\lambda_1 \uparrow \lambda_2}$, \Cref{injective-linkage} says that it is injective on global sections for generic $\lambda_0$. 
\begin{enumerate}
\item We expect, although we cannot prove, $\theta_{\lambda_0 \uparrow \lambda_1}$ to be injective on global sections for generic $\lambda_0$. This is based on the conjectural generic middle degree concentration for certain spaces 
of de Rham cohomology, which we study in forthcoming work \cite{de-Rham-paper}. If $\theta_{\lambda_0 \uparrow \lambda_1}$ were not injective its kernel would contribute to degree $2$ de Rham cohomology with coefficients 
in $L(\lambda_1)$. 
\item We know that $\theta_{\lambda_2 \uparrow \lambda_3}$ is not injective on global sections, even for generic $\lambda_0$. We give two arguments. One is that if that was true, then using the general machinery to produce entailments in 
\Cref{general-entailment} we would get a weak entailment between $F(\lambda_1)$ and $F(\lambda_2)$ or $F(\lambda_3)$, using that in the Grothendieck group
$[V(\lambda_{i})_{\Fpbar}]=[F(\lambda_i)]+[F(\lambda_{i-1})]$ for $i=2,3$. However, from the local results of \cite{Le-Hung-Lin} on Breuil--Mezard cycles we know this shouldn't happen. 
Less conjecturally, one can prove that the image of $\theta_{\lambda_1 \uparrow \lambda_2}$ on global sections is contained in the kernel of $\theta_{\lambda_2 \uparrow \lambda_3}$
after possibly 
choosing another non-zero map of Verma modules inducing an operator of the same weight increase
as $\theta_{\lambda_2 \uparrow \lambda_3}$. This is the case because over $\flag$
the composition of $\LL(w_{0,M}\lambda_1+\eta) \to \LL(w_{0,M}\lambda_2+\eta) \to \LL(w_{0,M}\lambda_3+\eta)$ 
and $\LL(w_{0,M}\lambda_1+\eta) \to \LL(w_{0,M}s_{\alpha}\cdot \lambda_1+\eta) \to \LL(w_{0,M}\lambda_3+\eta)$ agree, 
but $\pi_*\LL(w_{0,M}s_{\alpha}\cdot \lambda_1+\eta)=0$. This commutativity can be seen by $p$-translating a piece of 
the BGG complex over $\flag$ for a certain $\lambda \in C_0$, see \cite[Rmk 7.2.0.5]{thesis}. 
Alternatively, it would follow if we knew that the Hom space of Verma modules is $1$-dimensional. 
\item One could also ask the question about the behaviour of theta linkage maps in higher cohomology. We currently don't understand much about these. The two main problems is that higher coherent cohomology is 
not understood at all in general, and that restricting to strata becomes a less viable strategy for higher coherent cohomology. 
At least in some cases the
conjectures in \cite{de-Rham-paper} can help guide us in this regard. 
\end{enumerate}

\section{A de Rham realization functor} \label{section7}
For the purposes of our companion paper \cite{de-Rham-paper} on computing de Rham cohomology of $\Shbar^\tor$, we describe a functor that computes de Rham cohomology
and whose domain is $D^{b}(\mathcal{O}_{P,\Fpbar})$. This will allow us to compute generalized BGG decompositions in  
\cite{de-Rham-paper}.
\subsection{The realization functor}
Let $Q \in \{P,B\}$, and let $\flag_{P}=\Shbar$, $\flag_{B}=\flag$. Even though we have mostly focused on Verma modules 
it is natural to work in the larger category $\mathcal{O}_{Q,R}$ of $(U\mathfrak{g},Q)$-modules. 
As a corollary of \Cref{stratification-defn}(4) and \Cref{nilpotent-connection} we get the following. 
\begin{prop} \label{O-to-D-modules}
Let $R=\mathcal{O}/p^n$ and
$V \in \mathcal{O}_{Q,R}$. Then $\mathcal{V}^{\vee}=F_{Q}(V^{\vee})$ is naturally equipped with a 
HPD stratification over $\Sh_{R}$. It extends to a log HPD stratification on $\Sh^{\tor}_{R}$.
\end{prop}

Let $X/\mathcal{O}$ be a smooth proper scheme. Let $R=\mathcal{O}/p^n$ and let $X_n:=X_{R}$. We will need some important facts about crystals, 
see \cite{Berthelot-Ogus}. We recall that $\mathcal{O}=W(k)$. 
\begin{prop} \label{crystalline-theory}
\begin{enumerate}
\item The abelian category of $\mathcal{O}_{X_n}$-modules equipped with a HPD stratification with respect to $R$ and 
$\mathcal{O}_{X_n}$-linear maps that respect the HPD stratifications as morphisms
is equivalent 
to the category of crystals on $(X_n/R)_{\textnormal{crys}}$, the crystalline topos of $X_n$ with respect to $R$.
\item There is an exact functor $L$ from the category of $\mathcal{O}_{X_n}$-modules (with HPD differential operators with respect 
to $R$
as morphisms) to the category of crystals on $(X_n/R)_{\crys}$. Under the equivalence of $(1)$ it is given by $E \mapsto P_{X_n} \otimes E$
with the HPD stratification of \Cref{stratification-defn}(3).
\item Consider the derived pushforward $Ru_{*}: (X_n/R)_{\text{crys}} \to \Sh(X_{n,\text{zar}})$
to the topos of Zariski sheaves on $X_n$. It satisfies 
$Ru_*L(E)=E$ for $E$ a $\mathcal{O}_{X_{n}}$-module. Further, if $f: E_1 \to E_2$ is a 
HPD differential operator with respect to $R$,
then $Ru_*[L(E_1) \xrightarrow{L(f)} L(E_2)]=[E_1 \xrightarrow{f} E_2]$.
\item All the statements above extend to the setup of log-HPD stratifications and log-crystals as in 
\cite[\S 4]{Mokrane-Tilouine}.
\end{enumerate}

\end{prop}

Putting the two propositions together we get the realization functor which we will use to compute de Rham 
cohomology of $\Shbar$. 
\begin{theorem} \label{realization-functor}
Let $R=\mathcal{O}/p^n$ and $V \in \Rep_{R}(G)$. For a scheme $X/S$ let $\mathcal{C}_{X/S}$
be the category of $\mathcal{O}_X$-modules with $\mathcal{O}_S$-linear maps. We can define an exact functor
$$
\Psi_{Q,R}:=Ru_{*} \circ F_{Q} \circ (-)^{\vee} : D^{b}(\mathcal{O}_{Q,R}) \to D^{b}(\mathcal{C}_{\flag_{Q}/R}),
$$
where we are using \Cref{crystalline-theory}(1) and \Cref{O-to-D-modules} to make sense of it. 
If $C^{\bullet}$ is a complex of
$(U\mathfrak{g},Q)$-modules which is quasi-isomorphic to the complex $[V]$ concentrated in degree $0$,
$\Psi_{Q,R}(C^{\bullet})$ is quasi-isomorphic to the de Rham complex $\mathcal{V} \otimes \Omega^{\bullet}_{\flag_Q}$.
Moreover, if each 
term $C^{-k}$ is an extension of Verma modules $\Ver_{Q}(V_{i,k})$, then $\Psi_{Q}(C^{\bullet})^{k}$ is an 
extension of $\mathcal{V}^{\vee}_{i,k}$, and whenever a map between the graded pieces is defined it is given 
by the functor in \Cref{creator-diff}. The functor $\Psi_{Q,R}$ extends to functors $\Psi^{\can}_{Q,R}$ and 
$\Psi^{\sub}_{Q,R}$ on toroidal compactifications, satisfying the same properties. 
\begin{proof}
By \Cref{O-to-D-modules} $F_{Q} \circ (-)^{\vee}$  sends $C^{\bullet}$ to a complex of 
crystals on $(\Sh_{R}/R)_{\crys}$. Since it is an exact functor it descends to the respective bounded derived categories, so 
$F_{Q}(C^{\bullet,\vee})$ is quasi-isomorphic to $\mathcal{V^{\vee}}$ considered as a crystal, 
thus it computes 
its de Rham cohomology. Also, $\mathcal{V^{\vee}}$ is quasi-isomorphic to the linearized de Rham complex 
$L(\mathcal{V^{\vee}} \otimes \Omega^{\bullet}_{\flag_Q})$ \cite[Thm 6.12]{Berthelot-Ogus}. After applying $Ru_{*}$
and \Cref{crystalline-theory}(3) and we get the quasi-isomorphism 
in $D(\mathcal{C}_{\flag/R})$. For the last statement, we have that 
$\Psi_{Q,R}$ sends $\Ver_{Q}(V)$ to the crystal $L(\mathcal{V}^{\vee})$, by \Cref{canonical-iso-Pm}(7).
Thus, if $M$ is an extension of $\Ver_{Q}(V_1)$ and $\Ver_{Q}(V_2)$, then $F_{Q}(M^{\vee})$ is an extension of 
$L(V^{\vee}_2)$ and $L(V^{\vee}_1)$. After applying 
$Ru_{*}$ and using by \Cref{crystalline-theory}(3) again, it becomes an extension of $\mathcal{V}^{\vee}_2$
and $\mathcal{V}^{\vee}_1$. Finally, for a map $\phi: \Ver_Q V_1 \to \Ver_Q V_2$, under the isomorphism from
\Cref{canonical-iso-Pm}
$F_{Q}(\phi^{\vee}): L(\mathcal{V}^{\vee}_2) \to L(\mathcal{V}^{\vee}_1)$ 
is also identified with $L(\Phi_{P}(\phi))$, by construction of $\Phi_{P}$. Therefore, we conclude by \Cref{crystalline-theory} 
again. 
\end{proof}
\end{theorem}

\subsection{The lowest alcove BGG complex}
As an illustration of the realization functor, we reprove the lowest alcove integral BGG decomposition of 
\cite{Lan-Polo} and \cite{Tilouine-Polo}. Some of our technical inputs are different from theirs, and as opposed to  \cite{Lan-Polo}
we exhibit a quasi-isomorphism of the whole de Rham complex, not just of the graded pieces of its Hodge filtration. 
Let $\lambda \in C_0$. Then the dual Weyl module $V(\lambda)$ over $\Fpbar$ is irreducible.
Consider the standard complex $\text{Std}_{P}(V(\lambda))$ for $V(\lambda) \in \Rep_{\check{\Z}_p}(P)$
$$
0 \to \Ver_{P}(\bigwedge^{d}(\mathfrak{g}/\mathfrak{p}) \otimes_{\Z_p} V(\lambda)) \to
\ldots \to  \Ver_{P}((\mathfrak{g}/\mathfrak{p}) \otimes V(\lambda)) \to \Ver_{P}V(\lambda) \to 0,
$$
defined as in \cite[2.2]{Tilouine-Polo}, where $d$ is the dimension of
$\mathfrak{g}/\mathfrak{p}$. It is also
the dual of the de Rham complex of $F_{G/P}(V(\lambda)^{\vee})$ over $G/P$, under the equivalence of
\Cref{vb-on-flag}.
They are complexes of $(U\mathfrak{g},P)$-modules,
exact except at the term $\Ver_{P}V(\lambda)$, which has homology $V(\lambda)$.
Let $W(a) \subset W$ be the elements of length $a$.  We will need a small lemma, let 
$W_{\text{ext}}=pX^*(T) \rtimes W$ be the extended affine Weyl group, containing 
$W_{\text{aff}}=p\mathbb{\Z}\Phi^{+} \rtimes W$. They both act by the dot action on $X^*(T)$
as in \Cref{rep-theory}. Let $(U\mathfrak{g}_{\Fpbar})^{G}$ be the
$G$-invariant elements under the adjoint action. It is a commutative algebra, and there is 
a natural map of algebras $(U\mathfrak{g}_{\check{\Z}_p})^{G} \to (U\mathfrak{g}_{\Fpbar})^{G}$. 
\begin{lemma} \label{harish-chandra}
Assume that $p$ is greater than the Coxeter number of $G$. Let $\lambda,\mu \in X^*(T_{\check{\Z}_p})$, and let 
$\chi_{\lambda}$ be the character by which $(U\mathfrak{g}_{\Fpbar})^{G}$ acts on 
$\Ver_{B}(\lambda)_{\check{\Z}_p}$. Then 
$\Ver_{B}(\mu)_{\chi_{\lambda}}$ is non-zero, and equal to $\Ver_{B}(\mu)$, if and only if 
$\mu \in W_{\textnormal{ext}}
\cdot \lambda$. The same statement holds for 
$\Ver_{P}(L_{M}(\mu))_{\chi_{\lambda}}$. 
\begin{proof}
Under the assumption on $p$ we have that  $(U\mathfrak{g}_{\Fpbar})^{G}$ is isomorphic to $S(\mathfrak{h})^{W}$,
where $S(\mathfrak{h})$
 is the symmetric algebra of the Cartan subalgebra, and the action of $W$ is given by the dot action
 \cite{mod-p-Harish-Chandra}.
Note that any element of $(U\mathfrak{g}_{\Fpbar})^{G}$ also acts as $\chi_{\mu}$ on $\Ver_{P}(L_{M}(\mu))_{\Fpbar}$,
since $L_{M}(\mu)$ is cyclic as a $P$-module.
Therefore,
the localization of the Verma modules by $\chi_{\lambda}$ are non-zero if
 and only if $\chi_{\lambda}=\chi_{\mu}$. 
By the Harish-Chandra isomorphism above we have that $\chi_{\lambda}=\chi_{\mu}$ if and only if $\lambda \in W \cdot \mu + pX^{*}(T)$, i.e. 
$\lambda \in W_{\text{ext}} \cdot \mu$. 
\end{proof}
\end{lemma}
In all applications we will have $\lambda-\mu \in \Z \Phi^{+}$, so we can replace 
the extended affine Weyl group by the affine Weyl group. Also, note that the assumption on $p$ 
is needed for the lowest alcove to be non-empty, so we will often not mention it. 
With that we prove the BGG resolution for $\lambda \in C_0$ in the category of $(U\mathfrak{g},Q)$-modules. 
Recall that a weight $\mu \in X^*(T)$ is $p$-small 
if $|\langle \mu+\rho, \gamma^{\vee} \rangle| <p$ for all roots $\gamma$.
Let $\tilde{\mathcal{O}}=W(k')$ such that $G$ splits over $k'$.
\begin{theorem} \label{BGG-Verma}
Let $\lambda \in X^*(T)$ in $C_0$.
There exists a complex of $(U\mathfrak{g},P)_{\tilde{\mathcal{O}}}$-modules $\BGG_{P,\lambda}$, quasi-isomorphically 
embedded as a summand $\BGG_{P,\lambda} \subset \Std_{P}(V(\lambda))_{\tilde{\mathcal{O}}}$, such that  
$$
\BGG^{k}_{P,\lambda}=\oplus_{w \in W(k) \cap W^{M}} \Ver_{P}W(w \cdot \lambda)
$$
and the differentials go in the direction $\BGG^{d}_{P,\lambda} \to \BGG^{d-1}_{P,\lambda}$.

\begin{proof}
Let $\chi_{\lambda}$ be the character by which $(U\mathfrak{g}_{\Fpbar})^{G}$
acts on $V(\lambda)_{\Fpbar}$, since it is an irreducible representation. Then define 
$$
\BGG_{P,\lambda}\coloneqq \Std_{Q}(V(\lambda))_{\chi_{\lambda}}
$$
as the isotypic component of $\chi_{\lambda}$. It is a summand of $\Std_{P}(V(\lambda))$
since the latter is a direct sum over all the characters of $(U\mathfrak{g}_{\Fpbar})^{G}$. The homology 
of $\Std_{P}(V(\lambda))$ is $V(\lambda)$ concentrated in one degree, which is killed by every other
character $\chi\neq \chi_{\lambda}$
of $(U\mathfrak{g}_{\Fpbar})^{G}$. Therefore,
 $\BGG_{P,\lambda}$ is quasi-isomorphic to $\Std_{P}(V(\lambda))$.
The rest of 
the proposition for $\BGG_{P,\lambda}$ is proved in \cite[Thm D]{Tilouine-Polo} and \cite[Thm 5.2]{Lan-Polo}. 
We work with slightly more general groups, but the key input that we need is \Cref{harish-chandra}.
\end{proof}
\end{theorem}

For $\lambda \in C_0$ we have the identity $V(\lambda)^\vee=V(-w_0 \lambda)$, and 
$W(\lambda)^\vee=W(-w_{0,M} \lambda)$.  
Denote $\lambda^\vee=-w_0 \lambda$.
Let $\dR_{\flag}(\lambda)$ be the de Rham complex of 
$\mathcal{V}(\lambda)$. Recall that $\tilde{\mathcal{O}}=W(k')$ for some $k'/\F_p$ such that $G$ splits over $k'$.
Define $\BGG_{\Sh}(\lambda)$ as a complex on $\Sh_{\tilde{\mathcal{O}}}$ obtained 
by applying the functor $\Phi_{P}$ of 
\Cref{creator-diff} to $\BGG_{P,\lambda^\vee}$. Note that this is well-defined since $\BGG_{P,\lambda^\vee}$ is a complex whose terms are direct sums of Verma modules, not 
just an element of the derived category.
By \Cref{BGG-Verma}, it comes equipped with maps $\BGG_{\Sh}(\lambda) \to 
\dR_{\Sh}(\lambda)$. This is because $\Psi_{P,W(k')/p^n}$ of the standard complex for $V(\lambda)^{\vee}$ is 
$\dR_{\Sh}(\lambda)_{W(k')/p^n}$.
Similarly, let $\BGG^{?}_{\Sh}(\lambda)$ and $\dR^{?}_{\Sh^\tor}(\lambda)$ for $? \in \{\can,\sub\}$
be their (sub)canonical extensions to $\Sh^{\tor}$. Both the de Rham and BGG complexes are equipped with 
a natural Hodge filtration, as in \cite[Def 3.10]{Lan-Polo}.

\begin{theorem} \label{quasi-iso-BGG}
	Let $\lambda \in C_0$ and $? \in \{\can,\sub\}$. Then
  the map 
  $$
  \BGG^{?}_{\Sh^\tor}(\lambda) \to 
  \dR^{?}_{\Sh^\tor}(\lambda)
  $$ 
  is a quasi-isomorphism of filtered complexes over $\Sh_{\tilde{\mathcal{O}}}$. Moreover, 
  it admits a section.
  \begin{proof}
  The fact that the base change of the map above to each $\tilde{\mathcal{O}}/p^n$ is a quasi-isomorphism follows from \Cref{realization-functor} and \Cref{BGG-Verma}, 
  by the compatibility of the functors $\Phi_{P}$ and $\Psi_{P}$. Since both are complexes of vector bundles 
  it implies that the map over $\tilde{\mathcal{O}}$ is a quasi-isomorphism. 
  The fact that they are summands follows from functoriality of \Cref{realization-functor}.
  The compatibility with the Hodge filtrations follows in the same way, since as defined on each 
  $D(\mathcal{O}_{Q,\tilde{\mathcal{O}}/p^n})$
  both Hodge filtrations are quasi-isomorphic. Explicitly, for a Verma module $\Ver_{P}W(\chi)$ with $W(\chi) \in \Rep(P)$ irreducible, 
  its $n$-th Hodge filtration is $\Ver_{P}W(\chi)$ itself if $\langle \chi, \mu \rangle \le n$, and otherwise it is zero. 
   This defines a filtration on $\BGG_{P,\lambda,R}$. On the standard complex we define its Hodge filtration 
  as the one induced 
  by the Hodge filtration on $V(\lambda) \otimes \wedge^{\bullet} \mathfrak{g}/\mathfrak{p}$. 
  From its definition the quasi-isomorphism $\BGG_{P,\lambda,\tilde{\mathcal{O}}/p^n} \cong \Std_{P}V(\lambda)$ preserves these filtrations. 
  Finally, by the fact that $\Psi_{P,R} \Ver_{P}V=\mathcal{V}^{\vee}[0]$ we see that the associated filtration 
  on $\dR_{\Sh^\tor}(\lambda)_{\tilde{\mathcal{O}}/p^n}=\Psi_{P,\tilde{\mathcal{O}}/p^n}\Std_{P}V(\lambda)^{\vee}$
  is the Hodge filtration induced by Griffiths transversality. 
  \end{proof}
  \end{theorem}

\bibliographystyle{alpha}
\bibliography{bib}

\end{document}